\documentclass[11pt]{amsart}
\usepackage{amscd,amssymb}
\usepackage{amsthm,amsmath,amssymb}
\usepackage[matrix,arrow]{xy}
\usepackage{fancyhdr}

\makeatletter\@addtoreset{equation}{section} \makeatother

\newtheorem{theorem}[equation]{Theorem}

\newtheorem{proposition}[equation]{Proposition}
\newtheorem{lemma}[equation]{Lemma}
\newtheorem{corollary}[equation]{Corollary}

\theoremstyle{definition}
\newtheorem{example}[equation]{Example}

\theoremstyle{remark}
\newtheorem{remark}[equation]{Remark}

\author{Ivan Cheltsov}

\title{Elliptic structures on weighted\break\hfill three-dimensional Fano hypersurfaces}

\begin{document}

\begin{abstract}
We classify birational maps into elliptic fibrations of a general
quasismooth hypersurface in
$\mathbb{P}(1,a_{1},a_{2},a_{3},a_{4})$ of degree
$\sum_{i=1}^{4}a_{i}$ that has terminal singularities.
\end{abstract}

\maketitle\thispagestyle{empty}%

\section {Introduction.}
\label{section:introduction}

Let $X$ be a quasismooth hypersurface in
$\mathbb{P}(1,a_{1},a_{2},a_{3},a_{4})$ of degree
$d=\sum_{i=1}^{4}a_{i}$ that has terminal singularities, where
$a_{1}\leqslant a_{2}\leqslant a_{3}\leqslant a_{4}$. Then the
hypersurface $X$ is a Fano threefold, and there are exactly $95$
possibilities for the four-tuple $(a_{1},a_{2},a_{3},a_{4})$. We
use the notation~$n$~for the  entry numbers of these families,
which are ordered in the same way as in \cite{IF00}.

Suppose that the hypersurface $X$ is general. The following result
is proved in \cite{CPR}.

\begin{theorem}
\label{theorem:CPR} The hypersurface $X$ is birationally
rigid\footnote{Let $V$ be a Fano variety of Picard rank $1$ having
terminal $\mathbb{Q}$-factorial singularities. Then
$V$~is~said~to~be~bira\-ti\-onally rigid if it is not birational
to any other Mori fiber space (see \cite{Ch05umn}). The variety
$V$ is said to be birationally superrigid if it is birationally
rigid and $\mathrm{Bir}(V)=\mathrm{Aut}(V)$.}.
\end{theorem}

For every $n$ there are involutions
$\tau_{1},\ldots,\tau_{k_{n}}\in\mathrm{Bir}(X)$ that generates
the group $\mathrm{Bir}(X)$ up to bi\-re\-gu\-lar automorphisms
(see \cite{CPR}). In the case when $n\not\in \{7,20,60\}$ and
$k_{n}>0$, the hypersurface~$X$~can be bi\-ra\-ti\-onally
transformed into an elliptic fibration that is invariant under the
induced action of the group $\mathrm{Bir}(X)$, which is used in
\cite{ChPa05} to find relations between
$\tau_{1},\ldots,\tau_{k_{n}}$.

It is natural to try to classify all  birational transformations
of the hypersurface $X$ into elliptic fibrations, which is
equivalent to the following problem: find all rational maps
$X\dasharrow\mathbb{P}^{2}$ whose general fiber is birational to
an elliptic curve. Let us consider few examples.

\begin{example}
\label{example:elliptic-1} Let $n=1$. Then $X$ is a quartic
threefold. Let $\xi\colon X\dasharrow\mathbb{P}^{2}$ be a
projection from a line that is contained in $X$. Then a general
fiber of the map $\xi$ is an elliptic curve.
\end{example}

\begin{example}
\label{example:elliptic-2} Let $n=2$. Then $X$ is a hypersurface
in $\mathbb{P}(1,1,1,1,2)$ of degree $5$, which has one singular
point of type $\frac{1}{2}(1,1,1)$. There is a commutative diagram
$$
\xymatrix{
Y\ar@{->}[d]_{\pi}\ar@{->}[rr]^{\gamma}&&Z\ar@{->}[d]^{\eta}\\
X\ar@{-->}[rr]_{\psi}&&\mathbb{P}^{3}},
$$
where $\psi$ is the natural projection, $\pi$ is a weighted blow
up of the singular point of the hypersurface $X$ with weights
$(1,1,1)$, $\gamma$ is a birational morphism that contracts $15$
irreducible smooth rational curves $C_{1},\ldots, C_{15}$, and
$\eta$ is a double cover. Put $\xi=\chi\circ\psi$, where
$\chi\colon\mathbb{P}^{3}\dasharrow\mathbb{P}^{2}$ is a projection
from the point $\eta\circ\gamma(C_{i})$. Then a general fiber of
the map $\xi$ is an elliptic curve.
\end{example}

\begin{example}
\label{example:elliptic-17} Let $n=17$. Then $X$ is a hypersurface
in $\mathbb{P}(1,1,3,4,4)$ of degree $12$, whose singularities
consist of three singular points of type $\frac{1}{4}(1,1,3)$.
There is a commutative diagram
$$
\xymatrix{
&U\ar@{->}[d]_{\alpha}&&W\ar@{->}[ll]_{\beta}\ar@{->}[d]^{\omega}&\\%
&X\ar@{-->}[rr]_{\xi}&&\mathbb{P}(1,1,4),&}
$$
where $\xi$ is a projection, $\alpha$ is a weighted blow up of a
singular point of type $\frac{1}{4}(1,1,3)$ with
weights~$(1,1,3)$, $\beta$ is a weighted blow up with weights
$(1,1,2)$ of the singular point that is contained in the
exceptional divisor of the morphism $\alpha$, and $\omega$ is an
elliptic fibration.
\end{example}

\begin{example}
\label{example:elliptic-26} Let $n=26$. Then $X$ is a hypersurface
in $\mathbb{P}(1,1,3,5,6)$ of degree $15$ that has two singular
points of type $\frac{1}{3}(1,1,2)$. There is a commutative
diagram
$$
\xymatrix{
&&U\ar@{->}[dl]_{\sigma}\ar@{->}[dr]^{\omega}&&\\
&X\ar@{-->}[rr]_{\xi}&&\mathbb{P}(1,1,6),&}
$$
where $\xi$ is a projection, $\sigma$ is a weighted blow up of a
singular point of type $\frac{1}{3}(1,1,2)$ with
weights~$(1,1,2)$, and $\omega$ is a morphism given by the linear
system $|-6K_{X}|$. Then the normalization of a general fiber of
the rational map $\xi$ is an elliptic curve.
\end{example}

\begin{example}
\label{example:elliptic-31} Let $n=31$. Then $X$ is a hypersurface
in $\mathbb{P}(1,1,4,5,6)$ of degree $16$ that has a singular
point of type $\frac{1}{5}(1,1,4)$, and there is a commutative
diagram
$$
\xymatrix{
&U\ar@{->}[d]_{\alpha}&&W\ar@{->}[ll]_{\beta}\ar@{->}[d]^{\omega}&\\%
&X\ar@{-->}[rr]_{\xi}&&\mathbb{P}(1,1,6),&}
$$
where $\xi$ is a projection, $\alpha$ is a weighted blow up of the
singular point of type $\frac{1}{5}(1,1,4)$ with weights
$(1,1,4)$, $\beta$ is a weighted blow up with weights $(1,1,3)$ of
the singular point that is contained in the exceptional divisor of
$\alpha$, and $\omega$ is an elliptic fibration.
\end{example}

\begin{example}
\label{example:elliptic-7-11-19} Let $n\in \{7, 11, 19\}$. Then
$a_{2}=a_{3}$, and the hypersurface $X$ has  $d/a_{2}$ singular
points of type $\frac{1}{a_{2}}(1,1,a_{2}-1)$. Let  $\xi\colon
X\dasharrow\mathbb{P}(1,a_{1},a_{2})$ be a rational map induced by
a linear subsystem in the linear system $|-a_{2}K_{X}|$ consisting
of surfaces that pass through a given singular points of type
$\frac{1}{a_{2}}(1,1,a_{2}-1)$. Then the normalization of a
general fiber of $\xi$ is an elliptic curve.
\end{example}

\begin{example} \label{example:elliptic-7-9-20-30-36-44-49-51-64}
Let $n\in \{7, 9, 20, 30, 36, 44, 49, 51, 64\}$. Then $X$ can be
given by
$$
tz^{k}+\sum_{i=0}^{k-1}g_{i}(x,y,t,w)z^{i}=0\ \text{or}\
w^{2}t+wg(x,y,z,t)+f(x,y,z,t)=0\subset\mathrm{Proj}\Big(\mathbb{C}[x,y,z,t,w]\Big),
$$
where $\mathrm{wt}(x)=1$, $\mathrm{wt}(y)=a_{1}$,
$\mathrm{wt}(z)=a_{2}$, $\mathrm{wt}(t)=a_{3}$,
$\mathrm{wt}(w)=a_{4}$, and $g_{i}$ is a quasihomogeneous
polynomial. Let $\xi\colon X\dasharrow\mathbb{P}(1,a_{1},a_{3})$
be a rational map given by a linear system consisting of surfaces
that are cut out by  $f(x,y,t)=0$, where $f(x,y,t)$ is a
quasihomogeneous polynomial of degree $a_{3}$. Then the
normalization of a general fiber of the map $\xi$ is an elliptic
curve.
\end{example}

\begin{example}
\label{example:elliptic-main} Let $n\not\in \{1, 2, 3, 7, 11, 19,
60, 75, 84, 87, 93\}$, and $\xi\colon
X\dasharrow\mathbb{P}(1,a_{1},a_{2})$ be the natural projection.
Then the normalization of a general fiber of the map  $\xi$ is an
elliptic curve.
\end{example}

The purpose of this paper is to prove the following
result\footnote{In the case when $n\in\{1, 3, 14, 22, 28, 34, 37,
39, 52, 53, 57, 59, 60, 66, 70, 72, 73, 75, 78, 81, 84, 86, 87,
88, 89, 90, 92,\allowbreak 93, 94, 95\}$, the claim of
Theorem~\ref{theorem:main} is proved in \cite{Ch00a},
\cite{Ch03a}, \cite{ChPa05}.}.

\begin{theorem}
\label{theorem:main} Let $\rho\colon X\dasharrow\mathbb{P}^{2}$ be
a rational map whose general fiber is birational to an elliptic
curve. Then there is a commutative diagram
$$
\xymatrix{
&X\ar@{-->}[d]_{\xi}\ar@{-->}[rr]^{\sigma}&&X\ar@{-->}[d]^{\rho}&\\%
&\mathbb{P}(1,a_{1},a_{i})\ar@{-->}[rr]_{\phi}&&\mathbb{P}^{2},&}
$$
where $\phi$ is a birational map, $\sigma\in\mathrm{Bir}(X)$, and
$\xi$ is one of the dominant rational maps constructed in
Examples~\ref{example:elliptic-1}, \ref{example:elliptic-2},
\ref{example:elliptic-17}, \ref{example:elliptic-26},
\ref{example:elliptic-31}, \ref{example:elliptic-7-11-19},
\ref{example:elliptic-7-9-20-30-36-44-49-51-64} and
\ref{example:elliptic-main}.
\end{theorem}

\begin{corollary}
\label{corollary:main} Let $\rho\colon X\dasharrow\mathbb{P}^{2}$
be a rational map whose general fiber is birational to an elliptic
curve. Suppose that $n\not\in\{1, 2, 7, 9, 11, 17, 19, 20, 26, 30,
31, 36, 44, 49, 51, 64\}$. Then there is a birational map
$\psi\colon\mathbb{P}(1,a_{1},a_{2})\dasharrow\mathbb{P}^{2}$ such
that the diagram
$$
\xymatrix{
&&X\ar@{-->}[dl]_{\psi}\ar@{-->}[dr]^{\rho}&&&\\%
&\mathbb{P}(1,a_{1},a_{2})\ar@{-->}[rr]_{\phi}&&\mathbb{P}^{2}&}
$$
commutes, where $\psi$ is the natural projection.
\end{corollary}

\begin{corollary}
\label{corollary:3-60-75-84-87-93} The hypersurface $X$ can be
birationally transformed into an elliptic fibration if and only if
$n\not\in\{3, 60, 75, 84, 87, 93\}$.
\end{corollary}

To illustrate our technique let us prove the following result.

\begin{proposition}
\label{proposition:n-14}
The claim of Theorem~\ref{theorem:main} holds for $n=14$.%
\end{proposition}

\begin{proof}
Let $n=14$. Then $X$ is a hypersurface in $\mathbb{P}(1,1,1,4,6)$
of degree $12$, the singularities of the hypersurface $X$ consist
of a singular point of type $\frac{1}{2}(1,1,1)$. Let $\psi\colon
X\dasharrow\mathbb{P}^{2}$ be the natural projection, and
$\pi\colon U\to X$ be a weighted blow up with weights $(1,1,1)$ of
the singular point of the hypersurface $X$. Then $\psi\circ\pi$ is
an elliptic fibration.

Let $\rho\colon X\dasharrow\mathbb{P}^{2}$ be a rational map such
that the normalization of a general fiber of the rational map
$\rho$ is an irreducible elliptic curve. Let us consider
commutative diagram
$$
\xymatrix{
&&V\ar@{->}[ld]_{\alpha}\ar@{->}[rd]^{\beta}&&\\%
&X\ar@{-->}[rr]_{\rho}&&\mathbb{P}^{2},&}
$$ %
where $V$ is smooth, $\alpha$ is a birational morphism, and
$\beta$ is a morphism. Let $\mathcal{M}$ be the proper transform
of $|\beta^{*}(\mathcal{O}_{\mathbb{P}^{2}}(1))|$ on $X$. There is
a natural number $k>0$ such that $\mathcal{M}\sim-kK_{X}$. Then
$$
K_{V}+\frac{1}{k}\Big|\beta^{*}\Big(\mathcal{O}_{\mathbb{P}^{2}}(1)\Big)\Big|\sim_{\mathbb{Q}}\alpha^{*}\Big(K_{X}+\frac{1}{k}\mathcal{M}\Big)+\sum_{i=1}^{\delta}a_{i}E_{i}\sim_{\mathbb{Q}}\sum_{i=1}^{\delta}a_{i}E_{i},
$$
where $E_{i}$ is an $\alpha$-exceptional divisor, $a_{i}$ is a
rational number, and $\delta$ is the number of exceptional
divisors of the birational morphism $\alpha$. It follows from
\cite{CPR} that $a_{i}\geqslant 0$ for every $i$, but that there
is an index $j$ such that $a_{j}\leqslant 0$ by
Lemma~\ref{lemma:Noether-Fano}. Put $Z_{j}=\alpha(E_{j})$.

Suppose that $Z_{j}$ is a smooth point of $X$. Then
$\mathrm{mult}_{Z_{j}}(S_{1}\cdot S_{2})\geqslant 4k^{2}$ by
Lemma~1.10~in~\cite{Co00}, but the linear system $|-4K_{X}|$
induces a double cover $X\to\mathbb{P}(1,1,1,4)$. Thus, we have
$$
2k^{2}=H\cdot S_{1}\cdot S_{2}\geqslant\mathrm{mult}_{Z_{j}}\Big(S_{1}\cdot S_{2}\Big)\geqslant 4k^{2},%
$$
where $S_{1}$ and $S_{2}$ are general surfaces of the linear
system $\mathcal{M}$, and $H$ is a sufficiently general divisor in
the linear system $|-4K_{X}|$ that passes through the point
$Z_{j}$, which is a contradiction.

It follows from Corollary~\ref{corollary:curves} that $Z_{j}$ is
not a curve, which implies that $Z_{j}$ is the unique singular
point of the hypersurface $X$. Let $\mathcal{D}$ be the proper
transform of the linear system $\mathcal{M}$ on the variety $U$.
Then $\mathcal{D}\sim_{\mathbb{Q}}-kK_{U}$ by
Theorem~\ref{theorem:Kawamata}, which implies that $\mathcal{D}$
lies in the fibers of the elliptic fibration $\psi\circ\pi$, which
concludes the proof.
\end{proof}

Let us describe the structure of the paper. We consider auxiliary
results in Section~\ref{section:preliminaries}, the first steps of
the proof of Theorem~\ref{theorem:main} is done in
Section~\ref{section:start}, where we also prove
Theorem~\ref{theorem:main} in the case when $n\in\{1, 3, 5, 11,
14, 22, 28, 34, 37, 39, 52, 53, 57, 59, 60, 66, 70, 72, 73, 75,
78, 81, 84, 86, 87, \allowbreak 88, 89, 90, 92, 93, 94, 95\}$.
Then we prove Theorem~\ref{theorem:main} in all other cases.

\medskip

The author is very grateful to Max-Plank-Institute f\"ur
Mathematik at Bonn for the hospitality and excellent work
conditions. The author would like to thank A.\,Co\-rti,
M.\,Gri\-nen\-ko, V.\,Is\-kov\-skikh, Yu.\,Ma\-nin, J.\,Park,
Yu.\,Pro\-kho\-rov, A.\,Pukh\-li\-kov, V.\,Sho\-ku\-rov and
D.\,Stepanov for useful conversations.

\section{Preliminaries.}
\label{section:preliminaries}

Let $X$ be a threefold having terminal $\mathbb{Q}$-factorial
singularities, $\mathcal{M}$ be a linear system on the threefold
such that $\mathcal{M}$ does not have fixed components, and
$\lambda$ be an arbitrary non-negative rational number. In this
section we consider technical results describing properties of the
mobile\footnote{Elementary properties of mobile log pairs can be
found in \cite{Ch05umn}.} log pair $(X, \lambda\mathcal{M})$,
which are used in the proof of Theorem~\ref{theorem:main}. As
usual, the set of centers of canonical singularities of the mobile
log pair $(X, \lambda\mathcal{M})$ is denoted as $\mathbb{CS}(X,
\lambda\mathcal{M})$.

\begin{lemma}
\label{lemma:Noether-Fano}  Let $\rho\colon
X\dasharrow\mathbb{P}^{2}$ be a rational map whose general fiber
is birational to an elliptic curve, and $\pi\colon V\to X$ be a
resolution of the indeterminacies of $\rho$. Suppose that the
$\mathcal{M}$~is~a~proper transform of the linear system
$|\rho\circ\pi^{*}(\mathcal{O}_{\mathbb{P}^{2}}(1))|$, the divisor
$-K_{X}$ is nef and big, and the  equi\-va\-lence
$K_{X}+\lambda\mathcal{M}\sim_\mathbb{Q} 0$ holds. Then the
singularities of the log pair $(X,
\lambda\mathcal{M})$~are~not~terminal.
\end{lemma}

\begin{proof}
See the proof of Theorem~1.4.4 in \cite{Ch05umn}.
\end{proof}

The following well-known result is proved in \cite{Ka96}.

\begin{theorem}
\label{theorem:Kawamata} Let $O$ be a singular point of the
threefold $X$ of type $\frac{1}{r}(1,a,r-a)$, where~$a$~and~$r$
are coprime natural numbers such that $r>a$, and
$\mathrm{mult}_{O}(\mathcal{M})$ be a rational number such that
$$
\mathcal{D}\sim_{\mathbb{Q}} \pi^{*}(\mathcal{M})-\mathrm{mult}_{O}(\mathcal{M})G,%
$$
where $\pi\colon U\to X$ is a weighted blow up of $O$ with weights
$(1,a,r-a)$, $G$ is the ex\-cep\-ti\-onal divisor of $\pi$,
$\mathcal{D}$ is a proper transform of $\mathcal{M}$ on the
variety $U$. Suppose that $\mathbb{CS}(X, \lambda\mathcal{M})$
contains either the point $O$, or a curve passing through the
point $O$. Then $\mathrm{mult}_{O}(\mathcal{M})\geqslant
1/(r\lambda)$.
\end{theorem}

In the proof of Theorem~\ref{theorem:Kawamata} implies the
following result.

\begin{lemma}
\label{lemma:Cheltsov-Kawamata} Under the assumptions and
notations of Theorem~\ref{theorem:Kawamata}, suppose that the
singularities of $(X, \lambda\mathcal{M})$ are canonical and
$\mathbb{CS}(X, \lambda\mathcal{M})=\{O\}$, but $\mathbb{CS}(U,
\lambda\mathcal{D})\ne\varnothing$. Then
\begin{itemize}
\item the set $\mathbb{CS}(U, \lambda\mathcal{D})$ does not contain smooth points of the surface $G\cong\mathbb{P}(1,a,r-a)$,%
\item if the set $\mathbb{CS}(U, \lambda\mathcal{D})$ contains a
curve $L$, then $L\in|\mathcal{O}_{\mathbb{P}(1,\,a,\,r-a)}(1)|$,
and every singular point of the surface $G$ is contained in the
set $\mathbb{CS}(U, \lambda\mathcal{D})$.
\end{itemize}
\end{lemma}

\begin{proof}
We consider only the case when $r=5$ and $a=2$. Thus, we have
$G\cong\mathbb{P}(1,2,3)$.

Let $P$ and $Q$ be singular points of $G$, and $L$ be the curve in
$|\mathcal{O}_{\mathbb{P}(1,\,2,\,3)}(1)|$. Then $L$ passes
through the points $P$ and $Q$, but
$\mathrm{mult}_{O}(\mathcal{M})=1/(5\lambda)$ by
Theorem~\ref{theorem:Kawamata}, which implies
$\mathcal{D}\vert_{G}\sim_{\mathbb{Q}}\lambda L$.

Suppose that the set $\mathbb{CS}(U, \lambda\mathcal{D})$ contains
a subvariety $Z$ of subvariety $U$ that is different from the
curve $L$ and the points $P$ and $Q$. Then $Z\subset G$.

Suppose that $Z$ is a point. Then the point $Z$ is smooth on the
variety $U$, which implies the inequality
$\mathrm{mult}_{Z}(\mathcal{D})>1/\lambda$. Let $C$ be a general
curve in the linear system
$|\mathcal{O}_{\mathbb{P}(1,\,2,\,3)}(6)|$ that passes through the
point $Z$. Then the curve $C$ is not contained in the base locus
of the linear system $\mathcal{D}$. Hence, we have
$1/\lambda=C\cdot\mathcal{D}\geqslant\mathrm{mult}_{Z}(C)\mathrm{mult}_{Z}(\mathcal{D})>1/\lambda$,
which is a contradiction.

Therefore, the subvariety $Z$ is a curve. Then
$\mathrm{mult}_{Z}(\mathcal{D})\geqslant 1/\lambda$. Let $C$ be a
sufficiently general curve in the linear system
$|\mathcal{O}_{\mathbb{P}(1,\,2,\,3)}(6)|$. Then
$$
\frac{1}{\lambda}=C\cdot\mathcal{D}\geqslant\mathrm{mult}_{Z}(\mathcal{D})C\cdot Z\geqslant \frac{C\cdot Z}{\lambda},%
$$
which implies that $C\cdot Z=1$. Hence, the curve $Z$ is contained
in $|\mathcal{O}_{\mathbb{P}(1,\,2,\,3)}(1)|$.
\end{proof}

\begin{lemma}
\label{lemma:curves} Let $C$ be a curve on $X$ such that
$C\in\mathbb{CS}(X, \lambda\mathcal{M})$. Suppose that the
complete linear system $|-mK_{X}|$ is base-point-free for some
natural number $m>0$. Then $-K_{X}\cdot C\leqslant -K_{X}^{3}$.
\end{lemma}

\begin{proof}
Let $M_{1}$ and $M_{2}$ be general surfaces in $\mathcal{M}$. Then
the inequalities
$$
\mathrm{mult}_{C}\Big(M_{1}\cdot M_{2}\Big)\geqslant\mathrm{mult}_{C}\big(M_{1}\big)\mathrm{mult}_{C}\big(M_{2}\big)\geqslant \frac{1}{\lambda^{2}}%
$$
hold. Let $H$ be a general surface in $|-mK_{X}|$. Then
$$
\frac{-mK_X^3}{\lambda^{2}}=H\cdot M_{1}\cdot M_{2}\geqslant
\big(-mK_{X}\cdot C\big)\mathrm{mult}_{C}\Big(M_{1}\cdot
M_{2}\Big)\geqslant\frac{-mK_{X}\cdot C}{\lambda^{2}},
$$
which implies that $-K_{X}\cdot C\leqslant -K_{X}^{3}$.
\end{proof}

\begin{lemma}
\label{lemma:Cheltsov} Suppose that the linear system
$\mathcal{M}$ is not composed from a pencil. Then there~is~no
proper Zariski closed subset $\Sigma\subsetneq X$ such that
$$
\mathrm{Supp}(S_{1})\cap\mathrm{Supp}(S_{2})\subset\Sigma\subsetneq X,%
$$
where $S_{1}$ and $S_{2}$ are general divisors in the linear
system $\mathcal{M}$.
\end{lemma}

\begin{proof}
Suppose that there is a proper Zariski closed subset
$\Sigma\subset X$ such that the set-theoretic intersection of the
sufficiently general divisors $S_{1}$ and $S_{2}$  of the linear
system $\mathcal{M}$ is contained in the set $\Sigma$. Let
$\rho\colon X\dasharrow\mathbb{P}^{n}$ be a rational map induced
by the linear system $\mathcal{M}$, where $n$ is the dimension of
the linear system $\mathcal{M}$. Then there is a commutative
diagram
$$
\xymatrix{
&&W\ar@{->}[ld]_{\alpha}\ar@{->}[rd]^{\beta}&&\\%
&X\ar@{-->}[rr]_{\rho}&&\mathbb{P}^{n},&}
$$ %
where $W$ is a smooth variety, $\alpha$ is a birational morphism,
and $\beta$ is a morphism. Let $Y$ be the image of the morphism
$\beta$. Then $\mathrm{dim}(Y)\geqslant 2$, because $\mathcal{M}$
is not composed from a pencil.

Let $\Lambda$ be a Zariski closed subset of the variety $W$ such
that the morphism
$$
\alpha\vert_{W\setminus\Lambda}:W\setminus\Lambda\longrightarrow X\setminus\alpha(\Lambda)%
$$
is an isomorphism, and $\Delta$ be a union of the subset
$\Lambda\subset W$ and the closure of the proper transform of the
set $\Sigma\setminus\alpha(\Lambda)$ on the variety $W$. Then
$\Delta$ is a Zariski closed proper subset of  $W$.

Let $B_{1}$ and $B_{2}$ be general hyperplane sections of the
variety $Y$, and $D_{1}$ and $D_{2}$ be proper transforms of the
divisors $B_{1}$ and $B_{2}$  on the variety $W$ respectively.
Then $\alpha(D_{1})$ and $\alpha(D_{2})$ are general divisors of
the linear system $\mathcal{M}$. Hence, in the set-theoretic sense
we have
\begin{equation}
\label{equation:mobile-log-pairs-set-theoretic-intersection}
\varnothing\ne\beta^{-1}\Big(\mathrm{Supp}(B_{1})\cap\mathrm{Supp}(B_{2})\Big)=\mathrm{Supp}(D_{1})\cap\mathrm{Supp}(D_{2})\subset\Delta\subsetneq W,%
\end{equation}
because $\mathrm{dim}(Y)\geqslant 2$. However, the set-theoretic
identity~\ref{equation:mobile-log-pairs-set-theoretic-intersection}
is an absurd.
\end{proof}

The following result is implied by Lemma~0.3.3 in \cite{KMM} and
Lemma~\ref{lemma:Cheltsov}.

\begin{corollary}
\label{corollary:Cheltsov} Let $X$ be a three-dimensional variety
with canonical singularities, $D$ be divisor on the variety $X$
that is big and nef, $\mathcal{M}$ be a linear system on the
variety $X$ that does not have fixed components and is not
composed from a pencil, and $S_{1}$ and $S_{2}$ be sufficiently
general surfaces of the linear system $\mathcal{D}$. Then the
inequality $D\cdot S_{1}\cdot S_{2}>0$ holds.
\end{corollary}

The proof of Lemma~\ref{lemma:Cheltsov} implies the following
result.

\begin{lemma}
\label{lemma:Cheltsov-II} Let $X$ be a variety, $\mathcal{M}$ be a
linear system on the variety $X$ that does not have fixed
components and is not composed from a pencil, and $\mathcal{D}$ be
a linear system on $X$ that does not have fixed components. Then
there is no Zariski closed subset $\Sigma\subsetneq X$ such that
$$
\mathrm{Supp}(S)\cap\mathrm{Supp}(D)\subset\Sigma\subsetneq X,
$$
where $S$ and $D$ are sufficiently general divisors of the linear
system $\mathcal{M}$ and $\mathcal{D}$ respectively.
\end{lemma}

The claim of Lemma~\ref{lemma:Cheltsov} and the proof of
Lemma~\ref{lemma:curves}  imply the following result.

\begin{corollary}
\label{corollary:curves} Under the assumptions and notations of
Lemma~\ref{lemma:curves}, suppose that the linear system
$\mathcal{M}$ is not composed from a pencil, and the divisor
$-K_{X}$ is big. Then $-K_{X}\cdot C<-K_{X}^{3}$.
\end{corollary}

Many applications of Lemma~\ref{lemma:Cheltsov-II} use the
following simple result.

\begin{lemma}
\label{lemma:normal-surface} Let $S$ be a surface,  $D$ be an
effective divisor on $S$ such that
$D\equiv\sum_{i=1}^{r}a_{i}C_{i}$, where $a_{i}\in\mathbb{Q}$, and
$C_{1},\ldots, C_{r}$ are irreducible curves on $S$ whose
intersection form is negatively defined. Then
$D=\sum_{i=1}^{r}a_{i}C_{i}$.
\end{lemma}

\begin{proof}
Let $D=\sum_{i=1}^{k}c_{i}B_{i}$, where $B_{i}$ is an irreducible
curve on $S$, and $c_{i}$ is a nonnegative rational number.
Suppose that
$$
\sum_{i=1}^{k}c_{i}B_{i}\ne \sum_{i=1}^{r}a_{i}C_{i},
$$
and the curve $B_{i}$ is not one of the curves among
$C_{1},\ldots, C_{r}$ for every  $i$. We have
$$
0\geqslant\Big(\sum_{a_{i}>0}a_{i}C_{i}\Big)\cdot\Big(\sum_{a_{i}>0}a_{i}C_{i}\Big)=\Big(\sum_{i=1}^{k}c_{i}B_{i}\Big)\cdot\Big(\sum_{a_{i}>0}a_{i}C_{i}\Big)-\Big(\sum_{a_{i}\leqslant 0}a_{i}C_{i}\Big)\cdot\Big(\sum_{a_{i}>0}a_{i}C_{i}\Big)\geqslant 0,%
$$
which gives $\sum_{c_{i}\geqslant
0}c_{i}B_{i}\equiv\sum_{a_{i}\leqslant 0}a_{i}C_{i}$. Hence, we
have $c_{i}=0$ and $a_{i}=0$ for every $i$.
\end{proof}

\section{Beginning of classification.}
\label{section:start}

Let us use the notations and assumptions of
Section~\ref{section:introduction}. In this section we begin to
prove the claim of Theorem~\ref{theorem:main}. Suppose that there
is a birational map $\rho\colon X\dasharrow V$ and an elliptic
fibration $\nu\colon V\to\mathbb{P}^{2}$ such that $V$ is smooth,
and fibers of $\nu$ are connected. We must show that there is a
commutative diagram
\begin{equation}
\label{equation:diagram-main} \xymatrix{
&X\ar@{-->}[d]_{\xi}\ar@{-->}[rr]^{\sigma}&&X\ar@{-->}[rr]^{\rho}&&V\ar@{->}[d]^{\nu}&\\%
&\mathbb{P}(1,a_{1},a_{i})\ar@{-->}[rrrr]_{\zeta}&&&&\mathbb{P}^{2},&}
\end{equation}
where $\phi$ is a birational map, $\sigma$ is a birational
automorphism of $X$, and $\xi$ is one of the dominant rational
maps constructed in Examples~\ref{example:elliptic-1},
\ref{example:elliptic-2}, \ref{example:elliptic-17},
\ref{example:elliptic-26}, \ref{example:elliptic-31},
\ref{example:elliptic-7-11-19},
\ref{example:elliptic-7-9-20-30-36-44-49-51-64} and
\ref{example:elliptic-main}.

The commutative diagram~\ref{equation:diagram-main} implies the
commutative diagram
\begin{equation}
\label{equation:diagram-auxiliary} \xymatrix{
&X\ar@{-->}[d]_{\xi}\ar@{-->}[rr]^{\rho}&&V\ar@{->}[d]^{\nu}&\\%
&\mathbb{P}(1,a_{1},a_{i})\ar@{-->}[rr]_{\zeta}&&\mathbb{P}^{2}&}
\end{equation}
in the case when $\xi\circ\sigma=\chi\circ\xi$ for every
$\sigma\in\mathrm{Bir}(X)$, where
$\chi\in\mathrm{Bir}(\mathbb{P}(1,a_{1},a_{i}))$.

\begin{example}
\label{example:automorphisms-and-projection} Let $\psi\colon
X\dasharrow\mathbb{P}(1,a_{1},a_{2})$ be a projection, and
$\sigma$ be any  birational automorphism of the threefold $X$.
Suppose that $n\not\in \{1, 2, 3, 7, 11, 19, 20, 36, 60, 75, 84,
87, 93\}$. Then it follows from \cite{CPR} that there is a
birational automorphism $\chi$ of  $\mathbb{P}(1,a_{1},a_{2})$
such that $\psi\circ\sigma=\chi\circ\psi$.
\end{example}

Let $\mathcal{M}$ be a proper transform of the linear system
$|\nu^{*}(\mathcal{O}_{\mathbb{P}^{2}}(1))|$ on $X$. Then
$\mathcal{M}\sim_{\mathbb{Q}}-kK_{X}$ for some natural $k$, but
the singularities of the log pair $(X, \frac{1}{k}\mathcal{M})$
are not terminal by Lemma~\ref{lemma:Noether-Fano}.

\begin{remark}
\label{remark:canonical-singularities} It follows from \cite{CPR}
that there is a birational automorphism $\sigma\in\mathrm{Bir}(X)$
such that the singularities of the log pair $(X,
\frac{1}{k^{\prime}}\sigma(\mathcal{M}))$ are canonical, where
$k^{\prime}\in\mathbb{N}$ such that
$\mathcal{M}\sim_{\mathbb{Q}}-k^{\prime}K_{X}$.
\end{remark}

We may assume that the singularities of the log pair $(X,
\frac{1}{k}\mathcal{M})$ are canonical.

\begin{theorem}
\label{theorem:smooth-points} The set $\mathbb{CS}(X,
\frac{1}{k}\mathcal{M})$ does not contain smooth points of $X$ if
$n\ne 1$ and $n\ne 2$.
\end{theorem}

\begin{proof}
The  claim follows from the proof of Theorem~5.1.2 in \cite{CPR}.
\end{proof}

The following corollary is implied by Lemma~\ref{lemma:curves}.

\begin{corollary}
\label{corollary:cheap-curves} The set of centers of canonical
singularities $\mathbb{CS}(X, \frac{1}{k}\mathcal{M})$ does not
contain curves that do not contain singular points of the
hypersurface $X$ in the case when $n\geqslant 6$.
\end{corollary}

The set $\mathbb{CS}(X, \frac{1}{k}\mathcal{M})$ contains a
singular point of $X$ in the case $n\geqslant 6$ by
Theorem~\ref{theorem:Kawamata}.

\begin{proposition}
\label{proposition:singular-points} Suppose that the set
$\mathbb{CS}(X, \frac{1}{k}\mathcal{M})$ contains a singular point
$O$ of the hypersurface $X$ that is a singularity of type
$\frac{1}{r}(1,r-a,a)$, where $a$ and $r$ are coprime natural
numbers and $r>a$. Let $\pi\colon Y\to X$ be a weighted blow up of
$O$ with weights $(1,a,r-a)$. Then $-K_{Y}^{3}\geqslant 0$.
\end{proposition}

\begin{proof}
Suppose that $-K_{Y}^{3}<0$. Let $E$ be the $\pi$-exceptional
divisor, and $\mathcal{B}$ be a proper transform of $\mathcal{M}$
on the variety $Y$. Then
$-K_{Y}^{3}=-K_{X}^{3}-1/\big(ra(r-a)\big)$, but
$\mathcal{B}\sim_{\mathbb{Q}}-kK_{Y}$ by
Theorem~\ref{theorem:Kawamata}.

Let $\overline{\mathbb{NE}}(Y)$ be a closure in $\mathbb{R}^{2}$
of the cone generated by effective one-dimensional cycles of the
variety $Y$. Then $-E\cdot E$ generates an extremal ray of
$\overline{\mathbb{NE}}(Y)$, but Corollary~5.4.6 in \cite{CPR}
implies that there are integers $b>0$ and $c\geqslant 0$ such that
the cycle
$$
-K_{Y}\cdot \Big(-bK_{Y}+cE\Big)
$$
is numerically equivalent to an
effective, irreducible and reduced curve $\Gamma$ on the variety
$Y$ that generates the extremal ray of the cone
$\overline{\mathbb{NE}}(Y)$ different from the ray generated by
$-E\cdot E$.

Let $S_{1}$ and $S_{2}$ be general surfaces in $\mathcal{B}$. Then
$S_{1}\cdot S_{2}\in \overline{\mathbb{NE}}(Y)$, but $S_{1}\cdot
S_{2}\equiv k^{2}K^{2}_{Y}$, which implies that the cycle
$S_{1}\cdot S_{2}$ generates an extremal ray of the cone
$\overline{\mathbb{NE}}(Y)$ that contains the curve $\Gamma$.
Moreover, for every effective cycle $C\in \mathbb{R}^{+}\Gamma$ we
have
$$
\mathrm{Supp}(C)=\mathrm{Supp}\Big(S_{1}\cdot S_{2}\Big),
$$
because $S_{1}\cdot \Gamma<0$ and $S_{2}\cdot\Gamma<0$, which
contradicts Lemma~\ref{lemma:Cheltsov}.
\end{proof}

The following result is implied by
Proposition~\ref{proposition:singular-points}.

\begin{proposition}
\label{proposition:n-14-22-28-92-94-95} The claim of
Theorem~\ref{theorem:main} holds for $n\in\{14, 22, 28, 34, 37,
39, 52, 53, 57, 59,\allowbreak 66, 70, 72, 73, 78, 81, 86, 88, 89,
90, 92, 94, 95\}$.
\end{proposition}

\begin{proof}
We must show the existence of the commutative diagram
\begin{equation}
\label{equation:n-14-22-28-92-94-95-commutative-diagram}
\xymatrix{
&X\ar@{-->}[d]_{\psi}\ar@{-->}[rr]^{\rho}&&V\ar@{->}[d]^{\nu}&\\%
&\mathbb{P}(1,a_{1},a_{2})\ar@{-->}[rr]_{\phi}&&\mathbb{P}^{2},&}
\end{equation}
where $\psi$ is the natural projection, and $\phi$ is a birational
map.

It follows from Theorems~\ref{theorem:Kawamata} and
\ref{theorem:smooth-points} and Lemma~\ref{lemma:curves} that
$\mathbb{CS}(X, \frac{1}{k}\mathcal{M})$ contains a singular point
$P\in X$ of type $\frac{1}{r}(1,a,r-a)$, where $a$ and $r$ are
coprime natural numbers and $r>a$.

Let $\pi\colon Y\to X$ be a weighted blow up of  $P$ with weights
$(1,a,r-a)$, and $\mathcal{B}$ be the proper transform of the
linear system $\mathcal{M}$ on variety $Y$. Then $-K_{Y}^{3}=0$ by
Proposition~\ref{proposition:singular-points}.

It is easy to see that the linear system $|-rK_{Y}|$ does not have
base points for $r\gg 0$ and induces the morphism $\eta\colon
Y\to\mathbb{P}(1,a_{1},a_{2})$ such that the diagram
$$
\xymatrix{
&&&Y\ar@{->}[lld]_{\pi}\ar@{->}[rrd]^{\eta}&&&\\%
&X\ar@{-->}[rrrr]_{\psi}&&&&\mathbb{P}(1,a_{1},a_{2})&}
$$
is commutative. The equivalence
$\mathcal{B}\sim_{\mathbb{Q}}-kK_{Y}$ holds by
Theorem~\ref{theorem:Kawamata}, which implies the existence of the
commutative
diagram~\ref{equation:n-14-22-28-92-94-95-commutative-diagram}.
\end{proof}

The following result implies
Corollary~\ref{corollary:3-60-75-84-87-93}.

\begin{lemma}
\label{lemma:n-3} The claim of Theorem~\ref{theorem:main} holds
for $n\not\in\{3,60,75,84,87,93\}$.
\end{lemma}

\begin{proof}
It follows from Proposition~\ref{proposition:singular-points} that
$n\not\in\{75, 84, 87, 93\}$.

Suppose that $n=3$. Then the hypersurface $X$ is smooth, and the
set $\mathbb{CS}(X, \frac{1}{k}\mathcal{M})$ contains an
irreducible  curve $\Gamma$ such that $-K_{X}\cdot\Gamma=1$ by
Lemma~\ref{lemma:curves}. In particular, the curve $\Gamma$ is
smooth.

Let $\gamma\colon \bar{X}\to X$ be a blow up of $\Gamma$, $G$ be
the exceptional divisor of $\gamma$, and $\bar{S}_{1}$ and
$\bar{S}_{2}$ are proper transforms on $\bar{X}$ of general
surfaces in $\mathcal{M}$. Then the divisor
$\gamma^{*}(-3K_{X})-G$ is nef and big, but
$$
\Big(\eta^{*}(-3K_{X})-G\Big)\cdot\bar{S}_{1}\cdot\bar{S}_{1}=0,
$$
which contradicts Corollary~\ref{corollary:Cheltsov}.

We have $n=60$. The threefold $X$ is a hypersurface in
$\mathbb{P}(1,4,5,6,9)$ of degree $24$.

It is easy to check that $\mathbb{CS}(X, \frac{1}{k}\mathcal{M})$
does not contain curves by Corollary~\ref{corollary:curves}, which
implies that the set $\mathbb{CS}(X, \frac{1}{k}\mathcal{M})$
consists of the singular point $O$ of the hypersurface $X$ that is
a quotient singularity of type $\frac{1}{9}(1,4,5)$ by
Proposition~\ref{proposition:singular-points}. Let $\pi\colon Y\to
X$ be a weighted blow up of the singular point $O$ with weights
$(1,4,5)$, and $\mathcal{D}$ be a proper transform of the linear
system $\mathcal{M}$ on the threefold $Y$. Then
$\mathcal{D}\sim_{\mathbb{Q}}-kK_{Y}$ by
Theorem~\ref{theorem:Kawamata}.

Let $P$ and $Q$ be the points of $Y$ contained in the
$\pi$-exceptional divisor that are singularities of types
$\frac{1}{4}(1,1,3)$ and $\frac{1}{5}(1,1,4)$ respectively. Then
$\mathbb{CS}(Y, \frac{1}{k}\mathcal{D})\subseteq\{P, Q\}$ by
Lemmas~\ref{lemma:Cheltsov-Kawamata} and \ref{lemma:curves}.

Suppose that the set $\mathbb{CS}(Y, \frac{1}{k}\mathcal{D})$
contains the point $Q$. Let $\alpha\colon U\to Y$ be a weighted
blow up of the point $Q$ with weights $(1,1,4)$, and $\mathcal{B}$
be a proper transform of the linear system $\mathcal{M}$ on the
variety $U$. Then $\mathcal{B}\sim_{\mathbb{Q}}-kK_{U}$ by
Theorem~\ref{theorem:Kawamata}, the linear system $|-4K_{U}|$ is a
proper transform of the pencil $|-4K_{X}|$, and the base locus of
the pencil $|-4K_{U}|$ consists of an irreducible reduced curve
$Z$ on the variety $U$ such that the curve $\pi\circ\alpha(Z)$ is
a base curve of the pencil $|-4K_{X}|$.

Let $H$ be a general surface in $|-4K_{U}|$. Then $Z^{2}=-1/30$ on
$H$, but $\mathcal{B}\vert_{H}\sim_{\mathbb{Q}} kZ$, which
contradicts  Lemmas~\ref{lemma:Cheltsov-II} and
\ref{lemma:normal-surface}. We see that $\mathbb{CS}(Y,
\frac{1}{k}\mathcal{D})=\{P\}$ by Lemmas~\ref{lemma:Noether-Fano}
and \ref{lemma:Cheltsov-Kawamata}.

Let $\beta\colon W\to Y$ be a weighted blow up of the point $P$
with weights $(1,1,3)$, and $D$ be the proper transform of the
general surface in $|-5K_{X}|$ on $W$. Then $D$ is nef and big,
but the equality $D\cdot H_{1}\cdot H_{2}=0$ holds, where $H_{1}$
and $H_{2}$ are general surfaces of the proper transform of the
linear system $\mathcal{M}$ on the variety $W$, which contradicts
Corollary~\ref{corollary:Cheltsov}.
\end{proof}

Now we prove the following very simple result.

\begin{proposition}
\label{proposition:n-11} Suppose that $n=11$. Then the
diagram~\ref{equation:diagram-auxiliary} exists, where $\xi$ is
one~of~the~five rational maps constructed in
Example~\ref{example:elliptic-7-11-19}.
\end{proposition}

\begin{proof}
The threefold $X$ is a hypersurface in $\mathbb{P}(1,1,2,2,5)$ of
degree $10$, whose singularities consist of points $P_{1}$,
$P_{2}$, $P_{3}$, $P_{4}$ and $P_{5}$ that are singularities of
types $\frac{1}{2}(1,1,1)$. The diagram
$$
\xymatrix{
&&U_{i}\ar@{->}[dl]_{\pi_{i}}\ar@{->}[dr]^{\eta_{i}}&\\
&X\ar@{-->}[rr]_{\xi_{i}}&&\mathbb{P}(1,1,2)&&}
$$
commutes, where $\xi_{i}$ is a projection, $\pi_{i}$ is the
weighted blow up of $P_{i}$ with weights $(1,1,1)$, and $\eta_{i}$
is an elliptic fibration. It follows from
Theorems~\ref{theorem:Kawamata} and \ref{theorem:smooth-points}
and Lemma~\ref{lemma:curves} that $P_{i}\in\mathbb{CS}(X,
\frac{1}{k}\mathcal{M})$ for some $i\in\{1,2,3,4,5\}$.The
diagram~\ref{equation:diagram-auxiliary} exists for $\xi=\xi_{i}$
by Theorem~\ref{theorem:Kawamata}.
\end{proof}

The following result is due to \cite{Ry02}.

\begin{theorem}
\label{theorem:Ryder} Suppose that
$\mathbb{CS}(X,\frac{1}{k}\mathcal{M})$ contains an irreducible
curve $C$, and $n\geqslant 3$. Then
$$
\mathrm{Supp}(C)\subset\mathrm{Supp}\Big(S_{1}\cdot S_{2}\Big),
$$
where
$S_{1}$ and $S_{2}$ are different surfaces of the linear system
$|-K_{X}|$.
\end{theorem}

Let us prove the following result of Daniel Ryder.

\begin{proposition}
\label{proposition:n-5} The claim of Theorem~\ref{theorem:main}
holds for $n=5$.
\end{proposition}

\begin{proof}
Let $n=5$. Then $X$ is a hypersurface in $\mathbb{P}(1,1,1,2,3)$
of degree $7$, whose singularities consist of the points $P$ and
$Q$ such that $P$ is a singular point of type
$\frac{1}{2}(1,1,1)$, and $Q$~is~a~singular point of type
$\frac{1}{3}(1,1,2)$. The hypersurface $X$ can be given by the
equation
$$
w^{2}f_{1}(x,y,z)+f_{4}(x,y,z,t)w+f_{7}(x,y,z,t)=0\subset\mathbb{P}(1,1,1,2,3)\cong\mathrm{Proj}\Big(\mathbb{C}[x,y,z,t,w]\Big),
$$
where $\mathrm{wt}(x)=\mathrm{wt}(y)=\mathrm{wt}(z)=1$,
$\mathrm{wt}(t)=2$, $\mathrm{wt}(w)=3$, and $f_{i}$ is a
quasihomogeneous~poly\-no\-mi\-al of degree $i$. There is a
commutative diagram
$$
\xymatrix{
&&&&W\ar@{->}[dll]_{\gamma}\ar@{->}[d]^{\alpha}&&&\\%
&&Y\ar@{->}[dr]_{\eta}&&X\ar@{-->}[dl]_{\xi}\ar@{-->}[dr]^{\psi}&&U\ar@{->}[ld]^{\omega}\ar@{->}[llu]_{\beta}&\\
&&&\mathbb{P}(1,1,1,2)\ar@{-->}[rr]_{\chi}&&\mathbb{P}^{2},&&}
$$
where $\chi$ and $\xi$ are the natural projections, the morphism
$\omega$ is an elliptic fibration, the morphism $\alpha$ is a
weighted blow up of the point $Q$ with weights $(1,1,2)$, the
morphism $\gamma$ is the birational morphism that contracts $14$
smooth irreducible rational curves $C_{1},\ldots, C_{14}$ into
$14$ isolated ordinary double points $P_{1},\ldots, P_{14}$ of the
variety $Y$ respectively, the morphism $\eta$ is a double cover
branched over the surface $R\subset\mathbb{P}(1,1,1,2)$ that is
given by the equation
$$
f_{4}(x,y,z,t)^{2}-4f_{1}(x,y,z)f_{7}(x,y,z,t)=0\subset\mathbb{P}(1,1,1,2)\cong\mathrm{Proj}\Big(\mathbb{C}[x,y,z,t]\Big)
$$
and has $14$ isolated ordinary double points
$\eta(P_{1}),\ldots,\eta(P_{14})$, and $\beta$ is the composition
of the weighted blow ups with the weights $(1,1,1)$ of two
singular points of the variety $W$ that are singularities of types
$\frac{1}{2}(1,1,1)$. It follows from Theorem~\ref{theorem:Ryder}
that the set $\mathbb{CS}(X, \frac{1}{k}\mathcal{M})$ does not
contain curves (see the proof of Lemma~\ref{lemma:n-8-curves}).

Suppose that the set $\mathbb{CS}(X, \frac{1}{k}\mathcal{M})$
consists of the point $Q$. Let $O$ be the singular point of the
threefold $W$ that dominates the singular  point $Q$, and
$\mathcal{D}$ be the proper transform of the linear system
$\mathcal{M}$ on the threefold $W$. Then $O\in\mathbb{CS}(W,
\frac{1}{k}\mathcal{D})$ by Theorem~\ref{theorem:Kawamata} and
Lemmas~\ref{lemma:Noether-Fano} and \ref{lemma:Cheltsov-Kawamata}.

Let $\mathcal{B}$ be the proper transform of $\mathcal{D}$ on $Y$.
Then it follows from Theorem~\ref{theorem:Kawamata} and
Lemmas~\ref{lemma:Cheltsov-Kawamata} and \ref{lemma:curves} that
$\mathbb{CS}(Y, \frac{1}{k}\mathcal{B})$ contains a curve $C$ such
that $-K_{Y}\cdot C=1/2$, and $\chi\circ\eta(C)$ is a point.

There is an irreducible curve $Z$ on the variety $Y$ such that the
curve $Z$ is different from the curve $C$, but $\eta(Z)=\eta(C)$.
Let $S$ be a general surface of the linear system $|-K_{Y}|$ that
contains the curve $C$. Then $Z^{2}<0$ on $S$, but
$\mathcal{B}\vert_{S}\sim_{\mathbb{Q}} kC+kZ$, which contradicts
to Lemma~\ref{lemma:Cheltsov-II}.

It follows from Theorem~\ref{theorem:Kawamata} and
Lemmas~\ref{lemma:Cheltsov-Kawamata} and \ref{lemma:curves} that
$$
\mathbb{CS}\Big(X, \frac{1}{k}\mathcal{M}\Big)=\big\{P, Q\big\},
$$
but $O\in\mathbb{CS}(W, \frac{1}{k}\mathcal{D})$ by
Corollary~\ref{corollary:curves} and
Lemmas~\ref{lemma:Noether-Fano} and \ref{lemma:Cheltsov-Kawamata}.
Hence, the proper transform of the linear system $\mathcal{M}$ on
the variety $U$ is contained in the fibers of $\omega$ by
Theorem~\ref{theorem:Kawamata}.
\end{proof}

The claim of Theorem~\ref{theorem:Ryder} implies the following
result.

\begin{lemma}
\label{lemma:Ryder} Suppose that $a_{2}\ne 1$. Then the set
$\mathbb{CS}(X,\frac{1}{k}\mathcal{M})$ does not contains curves.
\end{lemma}

\begin{proof}
Suppose that $\mathbb{CS}(X,\frac{1}{k}\mathcal{M})$ contains a
curve $C$. Then the claim of Theorem~\ref{theorem:Ryder} implies
that there are surfaces $S_{1}$ and $S_{2}$ of the pencil
$|-K_{X}|$ such that the curve $C$ is an irreducible component of
the irreducible and reduced curve $S_{1}\cap S_{2}$, which
contradicts Lemma~\ref{lemma:curves}.
\end{proof}

Let us illustrate Lemma~\ref{lemma:Ryder} by proving the following
result.

\begin{proposition}
\label{proposition:n-18} The claim of Theorem~\ref{theorem:main}
holds for $n=18$.
\end{proposition}

\begin{proof}
Let $n=18$. Then $X$ is a hypersurface in $\mathbb{P}(1,2,2,3,5)$
of degree $12$, whose singularities consist of  the points
$O_{1}$, $O_{2}$, $O_{3}$, $O_{4}$, $O_{5}$ and $O_{6}$ that are
quotient singularities of type $\frac{1}{2}(1,1,1)$, and the point
$P$ that is a quotient singularity of type $\frac{1}{5}(1,2,3)$.
There is a commutative diagram
$$
\xymatrix{
&U\ar@{->}[d]_{\alpha}&&W\ar@{->}[ll]_{\beta}\ar@{->}[d]^{\eta}&\\%
&X\ar@{-->}[rr]_{\psi}&&\mathbb{P}(1,2,2),&}
$$
where $\psi$ is a projection, $\alpha$ is the weighted blow up of
$P$ with weights $(1,2,3)$, $\beta$ is the weighted blow up with
weights $(1,1,2)$ of the singular point of the variety $U$ that is
a quotient singularity of type $\frac{1}{3}(1,1,2)$, and $\eta$ is
an elliptic fibration.

It follows from Theorem~\ref{theorem:smooth-points},
Lemma~\ref{lemma:Ryder} and
Proposition~\ref{proposition:singular-points} that $\mathbb{CS}(X,
\frac{1}{k}\mathcal{M})=\{P\}$.

Let $\mathcal{D}$ be the proper transform of $\mathcal{M}$ on the
variety $U$, and $Q$ and $O$ be the singular points of the variety
$U$ contained in the exceptional divisor of the birational
morphism $\alpha$ that are singularities of types
$\frac{1}{3}(1,1,2)$ and $\frac{1}{2}(1,1,1)$ respectively. Then
$\mathcal{D}\sim_{\mathbb{Q}}-kK_{U}$ and
$$
\varnothing\ne\mathbb{CS}\Big(U,
\frac{1}{k}\mathcal{D}\Big)\subseteq\big\{Q, O\big\}
$$
by Theorem~\ref{theorem:Kawamata} and
Lemma~\ref{lemma:Noether-Fano}, \ref{lemma:Cheltsov-Kawamata} and
\ref{lemma:curves}.

Suppose that the set $\mathbb{CS}(U, \frac{1}{k}\mathcal{D})$
contains the point $O$. Let $\pi\colon Y\to U$ be the weighted
blow up of $O$ with weights $(1,1,1)$, $F$ be the
$\pi$-exceptional divisor, and $\mathcal{H}$ and $\mathcal{P}$ be
the proper transforms of $\mathcal{M}$ and $|-3K_{U}|$ on $Y$
respectively. Then $\mathcal{H}\sim_{\mathbb{Q}}-kK_{Y}$ by
Theorem~\ref{theorem:Kawamata}, but
$$
\mathcal{P}\sim_{\mathbb{Q}}\pi^{*}\big(-3K_{U}\big)-\frac{1}{2}F,
$$
and the base locus of the linear system $\mathcal{P}$ consists of
the irreducible curve $Z$ such that $\alpha\circ\pi(Z)$ is the
base curve of the linear system $|-3K_{X}|$. Moreover, for a
general surface $S$ of the linear system $\mathcal{P}$, the
inequality $S\cdot Z>0$ holds, which implies that the divisor
$\pi^{*}(-6K_{U})-F$ is nef and big. On the other hand, for
 general surfaces $D_{1}$ and $D_{2}$ of the linear
system $\mathcal{H}$, we have
$$
\Big(\pi^{*}\big(-6K_{U}\big)-F\Big)\cdot D_{1}\cdot D_{2}=\Big(\pi^{*}\big(-6K_{U}\big)-F\Big)\cdot\Big(\pi^{*}\big(-kK_{U}\big)-\frac{k}{2}F\Big)^{2}=0,%
$$
which contradicts Corollary~\ref{corollary:Cheltsov}.

Therefore, the set $\mathbb{CS}(U, \frac{1}{k}\mathcal{D})$
contains the point $Q$. Let $\mathcal{B}$ be the proper transform
of the linear system $\mathcal{M}$ on the variety $W$. Then the
equivalence $\mathcal{B}\sim_{\mathbb{Q}}-kK_{W}$ holds by
Theorem~\ref{theorem:Kawamata}, which easily implies that the
claim of Theorem~\ref{theorem:main} holds for the hypersurface
$X$.
\end{proof}

We conclude this section by proving the following result, which is
proved in \cite{Ch00a} and \cite{Ch03a}.

\begin{theorem}
\label{theorem:quartic-threefold-elliptic-fibrations} The claim of
Theorem~\ref{theorem:main} holds for $n=1$.
\end{theorem}

\begin{proof}
Let $X$ be a general hypersurface in $\mathbb{P}^{4}$ of degree
$4$. Then we must show that there is a line $L\subset X$ such that
there is a commutative diagram
\begin{equation}
\label{equation:quartic-threefold-commutative-diagram} \xymatrix{
&X\ar@{-->}[d]_{\psi}\ar@{-->}[rr]^{\rho}&&V\ar@{->}[d]^{\nu}&\\%
&\mathbb{P}^{2}\ar@{-->}[rr]_{\sigma}&&\mathbb{P}^{2},&}
\end{equation}
where $\psi$ is the projection from the $L$, and $\sigma$ is a
birational map.

Suppose that $\mathbb{CS}(X, \frac{1}{k}\mathcal{M})$ contains a
point $P$ of the quartic $X$. Let $H$ be a general hyperplane
section of $X$ passing through the point $P$. Then it follows from
Lemma~1.10 in \cite{Co00} that
$$
4k^{2}\leqslant\mathrm{mult}_{P}\Big(D_{1}\cdot D_{2}\Big)\leqslant D_{1}\cdot D_{2}\cdot H=4k^{2},%
$$
where $D_{1}$ and $D_{2}$ are general surfaces of the linear
system $\mathcal{M}$. Therefore, the support of the effective
one-dimensional cycle $D_{1}\cdot D_{2}$ is contained in the union
of a finite number of lines on the quartic $X$ that pass through
the point $P$, which contradicts Lemma~\ref{lemma:Cheltsov}.

Therefore, the set $\mathbb{CS}(X, \frac{1}{k}\mathcal{M})$
contains a curve $C$. Thus, the inequality
$\mathrm{mult}_{C}(\mathcal{M})\geqslant k$ holds, but it follows
from Lemma~\ref{lemma:curves} that $\mathrm{deg}(C)\leqslant 3$.

Suppose that $C$ is not contained in any plane in
$\mathbb{P}^{4}$. Then the curve $C$ is either a smooth curve of
degree $3$ or $4$, or a rational curve of degree $4$ having one
double point.

Suppose that $C$ is smooth. Let $\alpha\colon U\to X$ be the blow
up of the curve $C$, $F$ be the the exceptional divisor of
$\alpha$, and $\mathcal{D}$ be the proper transform of
$\mathcal{M}$ on the variety $U$. Then the base locus of the
linear system $|\alpha^{*}(-\mathrm{deg}(C)K_{X})-F|$ does not
contain curves, but
$$
\Big(\alpha^{*}\big(-\mathrm{deg}(C)K_{X}\big)-F\Big)\cdot D_{1}\cdot D_{2}<0,%
$$
where $D_{1}$ and $D_{2}$ are general surfaces of the linear
system $\mathcal{D}$, which is a contradiction.

Thus, the curve $C$ is a quartic curve with a double point $P$.
Let $\beta\colon W\to X$ be a composition of the blow up of $P$
with the blow up of  he proper transform of $C$. Let $G$ and $E$
be the exceptional divisors of $\beta$ such that $\beta(E)=C$ and
$\beta(G)=P$. Then the base locus of the linear system
$|\beta^{*}(-4K_{X})-E-2G|$ does not contain curves, but
$$
\Big(\beta^{*}\big(-4K_{X}\big)-E-2G\Big)\cdot D_{1}\cdot D_{2}<0
$$
where $D_{1}$ and $D_{2}$ are general surfaces of the linear
system $\mathcal{D}$, which is a contradiction.

Hence, we see that the curve $C$ is contained in a plane in
$\mathbb{P}^{4}$.

Suppose that $\mathrm{deg}(C)\ne 1$. Then we have the following
possibilities:
\begin{itemize}
\item the curve $C$ is a smooth conic;%
\item the curve $C$ is a smooth plane cubic;%
\item the curve $C$ is a singular plane cubic.%
\end{itemize}

Suppose that $C$ is smooth. Let $\alpha\colon U\to X$ be a blow up
of the curve $C$, $F$ be the exceptional divisor of the morphism
$\alpha$, and $\mathcal{D}$ be the proper transform of
$\mathcal{M}$ on the variety $U$. Then one can easily check that
the base locus of the linear system
$|\alpha^{*}(-\mathrm{deg}(C)K_{X})-F|$ does not contain curves.
Therefore, the divisor $\alpha^{*}(-\mathrm{deg}(C)K_{X})-F$ is
nef and big, but
$$
\Big(\alpha^{*}\big(-\mathrm{deg}(C)K_{X}\big)-F\Big)\cdot D_{1}\cdot D_{2}=0,%
$$
where $D_{1}$ and $D_{2}$ are general surfaces of  $\mathcal{D}$,
which contradicts Corollary~\ref{corollary:Cheltsov}.

Hence, the curve $C$ is a plane cubic with a double point $P$. Let
$\beta\colon W\to X$ be a composition of the blow up of $P$ with
the blow up of the proper transform of $C$. Let $G$ and $E$ be the
exceptional divisors of the morphism $\beta$ such that
$\beta(E)=C$ and $\beta(G)=P$. Then the base locus of the linear
system $|\beta^{*}(-3K_{X})-E-2G|$ does not contain curves, which
implies that the divisor $\beta^{*}(-3K_{X})-E-2G$ is nef and big.
On the other hand, the inequality
$$
\Big(\beta^{*}\big(-3K_{X}\big)-E-2G\Big)\cdot D_{1}\cdot
D_{2}\leqslant 0
$$
holds, where $D_{1}$ and $D_{2}$ are general surfaces of
$\mathcal{D}$, which is impossible by
Corollary~\ref{corollary:Cheltsov}.

Thus, we see that the curve $C$ is a line. The equality
$\mathrm{mult}_{C}(D)=k$ implies the existence of the commutative
diagram~\ref{equation:quartic-threefold-commutative-diagram},
where $L=C$.
\end{proof}

\section{Case $n=2$, hypersurface of degree $5$ in $\mathbb{P}(1,1,1,1,2)$.}
\label{section:n-2}

We use the notations and assumptions of
Section~\ref{section:start}. Let $n=2$. Then $X$ is a sufficiently
general hypersurface in $\mathbb{P}(1,1,1,1,2)$ of degree $5$, the
equality $-K_{X}^{3}=5/2$ holds, and the singularities of the
hypersurface $X$ consist of a point $O$ that is a quotient
singularity of type $\frac{1}{2}(1,1,1)$.

The hypersurface $X$ can be given by the equation
$$
w^{2}f_{1}(x,y,z,t)+f_{3}(x,y,z,t)w+f_{5}(x,y,z,t)=0\subset\mathbb{P}(1,1,1,1,2)\cong\mathrm{Proj}\Big(\mathbb{C}[x,y,z,t,w]\Big),
$$
where
$\mathrm{wt}(x)=\mathrm{wt}(y)=\mathrm{wt}(z)=\mathrm{wt}(t)=1$
and $\mathrm{wt}(w)=2$, $f_{i}(x,y,z,t)$ is a homogeneous
polynomial of degree $i$, and the point $O$ is given by the
equations $x=y=z=t=0$.

Let $\psi\colon X\dasharrow\mathbb{P}^{3}$ be the natural
projection. Then there is a commutative diagram
$$
\xymatrix{
&&&Y\ar@{->}[dll]_{\pi}\ar@{->}[dr]^{\gamma}&&W_{i}\ar@{->}[ll]_{\alpha_{i}}\ar@{->}[dr]^{\omega_{i}}&&\\%
&X\ar@{-->}[drr]_{\psi}&&&Z\ar@{->}[dl]^{\eta}&&U_{i}\ar@{->}[ll]_{\beta_{i}}\ar@{->}[ld]^{\xi_{i}}&\\
&&&\mathbb{P}^{3}\ar@{-->}[rr]_{\chi_{i}}&&\mathbb{P}^{2}&&}
$$
where $\pi$ is the weighted blow up of the point $O$ with weights
$(1,1,1)$, the morphism $\gamma$ is the birational morphism that
contracts $15$ smooth rational curves $C_{1},\ldots, C_{15}$ to
$15$ isolated ordinary double points $P_{1},\ldots, P_{15}$ of the
variety $Z$ respectively, $\eta$ is a double cover branched over
the surface $R\subset\mathbb{P}^{3}$ of degree $6$ that is given
by the equation
$$
f_{3}(x,y,z,t)^{2}-4f_{1}(x,y,z,t)f_{5}(x,y,z,t)=0\subset\mathbb{P}^{3}\cong\mathrm{Proj}\Big(\mathbb{C}[x,y,z,t]\Big)
$$
and has $15$ isolated ordinary double points $\eta(P_{1}),\ldots,
\eta(P_{15})$, $\alpha_{i}$ is a blow up of $C_{i}$, $\beta_{i}$
is a blow up of the point $P_{i}$, $w_{i}$ is a birational
morphism, $\chi_{i}$ is a projection from  $\eta(P_{i})$, and
$\xi_{i}$ is an elliptic fibration. Moreover, the points
$\eta(P_{1}),\ldots, \eta(P_{15})$ are given by the equations
$f_{3}=f_{1}=f_{5}=0$.

\begin{remark}
\label{remark:n-2-singularities-are-canonical} It follows from
\cite{CPR} that a general fiber of $\chi_{i}\circ\psi$ is
$\mathrm{Bir}(X)$-invariant.
\end{remark}

In the rest of the section we prove the following result.

\begin{proposition}
\label{proposition:n-2} There is a commutative diagram
\begin{equation}
\label{equation:n-2-commutative-diagram} \xymatrix{
&U_{i}\ar@{->}[d]_{\xi_{i}}\ar@{-->}[rr]^{\rho}&&V\ar@{->}[d]^{\nu}&\\%
&\mathbb{P}^{2}\ar@{-->}[rr]_{\phi}&&\mathbb{P}^{2}&}
\end{equation}
for some $i\in\{1,\ldots,15\}$, where $\phi$ is a birational map.
\end{proposition}

Let $\mathcal{B}_{i}$ be the proper transform of the linear system
$\mathcal{M}$ on the variety $W_{i}$. Then in order to prove the
existence of the commutative
diagram~\ref{equation:n-2-commutative-diagram} it is enough to
show that $\mathcal{B}_{i}$ lies in the fibers of the elliptic
fibration $\xi_{i}\circ\omega_{i}$, which is implied by the
equivalence $\mathcal{B}_{i}\sim_{\mathbb{Q}}-kK_{W_{i}}$.

\begin{lemma}
\label{lemma:n-2-smooth-points} Suppose that $\mathbb{CS}(X,
\frac{1}{k}\mathcal{M})$ contains a smooth point of $X$. Then the
commutative diagram~\ref{equation:n-2-commutative-diagram} exists
for some $i\in\{1,\ldots,15\}$.
\end{lemma}

\begin{proof}
Suppose that $\mathbb{CS}(X, \frac{1}{k}\mathcal{M})$ contains a
smooth point $P$ of $X$. Let $S$ be a sufficiently general surface
in $|-K_{X}|$ that passes through $P$. In the case when
$P\not\in\cup_{i=1}^{15}\pi(C_{i})$, the surface $S$ does not
contain irreducible components of the effective cycle $D_{1}\cdot
D_{2}$, where $D_{1}$ and $D_{2}$ are general surfaces of the
linear system $\mathcal{M}$. Therefore, in the latter case we have
$$
\mathrm{mult}_{P}\Big(D_{1}\cdot D_{2}\Big)\leqslant D_{1}\cdot D_{2}\cdot S=-k^{2}K_{X}^{3}=\frac{5}{2}k^{2},%
$$
which contradicts Lemma~1.10 in \cite{Co00}. Thus, the point $P$
is contained in $\pi(C_{i})$ for some $i$.

We use arguments of \cite{CPR}. Put $C=\pi(C_{i})$ and
$\mathcal{M}\vert_{S}=\mathcal{L}+\mathrm{mult}_{C}(\mathcal{M})C$,
where $\mathcal{L}$ is a linear system on the surface $S$ without
fixed components. Then the log pair
$$
\Big(S,
\frac{1}{k}\mathcal{L}+\frac{\mathrm{mult}_{C}\big(\mathcal{M}\big)}{k}C\Big)
$$
is not log terminal in $P$ by Theorem~7.5 in \cite{Ko97}.
Therefore, we have
$$
\mathrm{mult}_{P}\Big(L_{1}\cdot L_{2}\Big)\geqslant
4\Big(1-\frac{\mathrm{mult}_{C}\big(\mathcal{M}\big)}{k}\Big)k^{2}
$$
by Theorem~3.1 in \cite{Co00}, where $L_{1}$ and $L_{2}$ are
general curves in $\mathcal{L}$. On the other hand, the equality
$$
L_{1}\cdot L_{2}=\frac{5}{2}k^{2}-\mathrm{mult}_{C}\big(\mathcal{M}\big)k-\frac{3}{2}\mathrm{mult}^{2}_{C}\big(\mathcal{M}\big)%
$$
holds on the surface $S$, because $C^{2}=-3/2$. Hence, we have
$$
\frac{5}{2}k^{2}-\mathrm{mult}_{C}\big(\mathcal{M}\big)k-\frac{3}{2}\mathrm{mult}^{2}_{C}\big(\mathcal{M}\big)\geqslant
4\Big(1-\frac{\mathrm{mult}_{C}\big(\mathcal{M}\big)}{k}\Big)k^{2},
$$
which implies $\mathrm{mult}_{C}(\mathcal{M})=k$. Thus, the curve
$\pi(C_{i})$ is contained in $\mathbb{CS}(X,
\frac{1}{k}\mathcal{M})$, and the equivalence
$\mathcal{B}_{i}\sim_{\mathbb{Q}}-kK_{W_{i}}$ follows from
Theorem~\ref{theorem:Kawamata}, which concludes the proof.
\end{proof}

We may assume that $\mathbb{CS}(X, \frac{1}{k}\mathcal{M})$ does
not contain smooth points of $X$.

\begin{lemma}
\label{lemma:n-2-curves-in-nonsingular-locus} Let $C$ be a curve
on $X$ such that $C\cap\mathrm{Sing}(X)=\varnothing$. Then
$C\not\in\mathbb{CS}(X, \frac{1}{k}\mathcal{M})$.
\end{lemma}

\begin{proof}
Suppose that the set $\mathbb{CS}(X, \frac{1}{k}\mathcal{M})$
contains the curve $C$. Then $\mathrm{mult}_{C}(\mathcal{M})=k$.
Let $H$ be a very ample divisor on $X$. Then
$H\sim_{\mathbb{Q}}-\lambda K_{X}$ holds for a natural number
$\lambda$. Thus, we have
$$
\frac{5\lambda k^{2}}{2}=-\lambda k^{2}K_{X}^{3}=H\cdot S_{1}\cdot S_{2}\geqslant\mathrm{mult}^{2}_{C}\big(\mathcal{M}\big)H\cdot C\geqslant -\lambda k^{2}K_{X}\cdot C,%
$$
where $S_{1}$ and $S_{2}$ are general surfaces in $\mathcal{M}$.
Therefore, we have the following possibilities:
\begin{itemize}
\item the equality $-K_{X}\cdot C=1$ holds, and the curve $C$ is smooth and rational;%
\item the equality $-K_{X}\cdot C=2$ holds, and the curve $C$ is smooth and rational;%
\item the equality $-K_{X}\cdot C=2$ holds, and the arithmetic genus of the curve $C$ is $1$.%
\end{itemize}

Let $\sigma\colon \check{X}\to X$ be the blow up of the ideal
sheaf of the curve $C$, and $G$ be the exceptional divisor of the
birational morphism $\sigma$. Then the variety $\check{X}$ is
smooth  in the neighborhood of the divisor $G$ whenever the curve
$C$ is smooth. Moreover, in the case when the curve $C$ has an
ordinary double point the singularities of $\check{X}$ in the
neighborhood of the divisor $G$ consist of a single isolated
ordinary double point. In the case when the curve $C$ has a
cuspidal singularity the singularities of the variety $\check{X}$
in the neighborhood of  $G$ consist of an isolated double point
such that $\check{X}$ in the neighborhood of this point is locally
isomorphic to the hypersurface
$$
x_{1}^{2}+x_{2}^{2}+x_{3}^{2}+x_{4}^{3}=0\subset\mathbb{C}\cong\mathrm{Spec}\Big(\mathbb{C}[x_{0},x_{1},x_{2},x_{3}]\Big).
$$

Let $\check{S}_{1}$ and $\check{S}_{2}$ be the  proper transforms
of the surfaces $S_{1}$ and $S_{2}$ on the variety $\check{X}$
respectively.

Suppose that $-K_{X}\cdot C=1$. Then the curve $C$ is cut in the
set-theoretic sense by the surfaces of the linear system
$|-2K_{X}|$ that pass through  $C$. Moreover, the scheme-theoretic
intersection of two general surfaces of the linear system
$|-2K_{X}|$ passing through $C$ is reduced in a general point of
$C$. Thus, the divisor $\sigma^{*}(-2K_{X})-G$ is nef and big (see
Lemma~5.2.5 in \cite{CPR}), but
$$
\Big(\sigma^{*}(-2K_{X})-G\Big)\cdot\check{S}_{1}\cdot\check{S}_{2}=0,%
$$
which is impossible by Corollary~\ref{corollary:Cheltsov}.

Suppose that $-K_{X}\cdot C=2$, and $C$ is smooth and rational.
Then $\sigma^{*}(-2K_{X})-G$ is nef, because the curve $C$ is cut
in the set-theoretic sense by the surfaces of the linear system
$|-2K_{X}|$ that pass through the curve $C$, but the
scheme-theoretic intersection of two general surfaces of the
linear system $|-2K_{X}|$ passing through$C$ is reduced in a
general point of $C$. We have
$$
0>-3k^{2}=\Big(\sigma^{*}(-2K_{X})-G\Big)\cdot\check{S}_{1}\cdot\check{S}_{2}\geqslant 0.%
$$

Hence,  the arithmetic genus of the curve $C$ is $1$ and
$-K_{X}\cdot C=2$. The curve $C$ is a set-theoretic intersection
of the surfaces in $|-4K_{X}|$ that pass through $C$. Moreover,
the scheme-theoretic intersection of two general surfaces of the
linear system $|-4K_{X}|$ passing through $C$ is reduced in a
general point $C$. Hence, the divisor $\sigma^{*}(-4K_{X})-G$ is
nef and big, but
$$
\Big(\sigma^{*}(-4K_{X})-G\Big)\cdot\check{S}_{1}\cdot\check{S}_{2}=0,
$$
which contradicts Corollary~\ref{corollary:Cheltsov}.
\end{proof}

It follows from Theorem~\ref{theorem:Kawamata} that
$O\in\mathbb{CS}(X, \frac{1}{k}\mathcal{M})$. Let $\mathcal{D}$ be
the proper transform of $\mathcal{M}$ on the variety $Y$. Then
Theorem~\ref{theorem:Kawamata} implies that
$\mathcal{D}\sim_{\mathbb{Q}}-kK_{Y}$. Thus, the commutative
diagram~\ref{equation:n-2-commutative-diagram} exists in the case
when $C_{i}\in\mathbb{CS}(Y, \frac{1}{k}\mathcal{D})$. Therefore,
we may assume that
$$
\mathbb{CS}\Big(Y,
\frac{1}{k}\mathcal{D}\Big)\cap\big\{C_{1},\ldots,C_{15}\big\}=\varnothing.
$$

\begin{lemma}
\label{lemma:n-2-smooth-points-again} The set $\mathbb{CS}(Y,
\frac{1}{k}\mathcal{D})$ does not contain smooth points of the
variety $Y$.
\end{lemma}

\begin{proof}
The set $\mathbb{CS}(X, \frac{1}{k}\mathcal{M})$ does not contain
smooth points of $X$. Therefore, to conclude the proof it is
enough to show that the set $\mathbb{CS}(Y,
\frac{1}{k}\mathcal{D})$ does not contain points of the
exceptional divisor of the morphism $\pi$, which is implied by
Lemma~\ref{lemma:Cheltsov-Kawamata}.
\end{proof}

Put $\mathcal{H}=\gamma(\mathcal{D})$. Then
$\mathcal{H}\sim_{\mathbb{Q}}-kK_{Z}$, the singularities of the
log pair $(Z, \frac{1}{k}\mathcal{H})$ are canonical, and the
claim of Lemma~\ref{lemma:Noether-Fano} implies that the set
$\mathbb{CS}(Z, \frac{1}{k}\mathcal{H})$ is not empty.

\begin{lemma}
\label{lemma:n-2-points-downstairs} The set $\mathbb{CS}(Z,
\frac{1}{k}\mathcal{H})$ does not contain points of the variety
$Z$.
\end{lemma}

\begin{proof}
It follows from Lemma~\ref{lemma:n-2-smooth-points-again} that
smooth points of the variety $Z$ are not contained in the set
$\mathbb{CS}(Z, \frac{1}{k}\mathcal{H})$. The condition
$P_{i}\in\mathbb{CS}(Z, \frac{1}{k}\mathcal{H})$ implies that the
set $\mathbb{CS}(Y, \frac{1}{k}\mathcal{D})$ contains either the
curve $C_{i}$, or a point on the curve $C_{i}$, which is
impossible.
\end{proof}

Thus, there is a curve $\Gamma$ on $Z$ that is contained in the
set $\mathbb{CS}(Z, \frac{1}{k}\mathcal{H})$, and
$\mathrm{mult}_{\Gamma}(\mathcal{H})=k$.

\begin{lemma}
\label{lemma:n-2-double-space-with-15-double-points-again} The
equality $-K_{Z}\cdot\Gamma=1$ holds.
\end{lemma}

\begin{proof}
Let $H$ be a general divisor of the linear system $|-K_{Z}|$. Then
$$
2k^{2}=H\cdot D_{1}\cdot D_{2}\geqslant
\mathrm{mult}_{\Gamma}\Big(D_{1}\cdot D_{2}\Big)H\cdot
\Gamma\geqslant -k^{2}K_{Z}\cdot\Gamma,
$$
where $D_{1}$ and $D_{2}$ are sufficiently general surfaces of the
linear system $\mathcal{H}$. Therefore,  the inequality
$-K_{Z}\cdot\Gamma\leqslant 2$ holds. Moreover, the equality
$-K_{Z}\cdot\Gamma=2$ implies that the support of the effective
cycle $D_{1}\cdot D_{2}$ coincides with the curve $\Gamma$, which
contradicts Lemma~\ref{lemma:Cheltsov}.
\end{proof}

The curve $\eta(\Gamma)$ is a line in $\mathbb{P}^{3}$, and
$\eta\vert_{\Gamma}\colon\Gamma\to\eta(\Gamma)$ is an isomorphism.
However, the arguments used in the proof of
Lemma~\ref{lemma:n-2-double-space-with-15-double-points-again}
easily imply that $\mathbb{CS}(Z,
\frac{1}{k}\mathcal{H})=\{\Gamma\}$ and
$\mathrm{mult}_{\Gamma}(D_{1}\cdot D_{2})<2k^{2}$, where $D_{1}$
and $D_{2}$ are general surfaces in $\mathcal{H}$.

\begin{lemma}
\label{lemma:n-2-double-space-with-15-double-points-again-and-again}
The line $\eta(\Gamma)$ is contained in the ramification surface
$R$ of the double cover $\eta$.
\end{lemma}

\begin{proof}
Suppose that the line $\eta(\Gamma)$ is not contained in the
surface $R$. Let $S$ be a general surface of the linear system
$|-K_{Z}|$ that passes through the curve $\Gamma$. Then
$$
\mathcal{H}\vert_{S}=\mathrm{mult}_{\Gamma}(\mathcal{H})\Gamma+\mathrm{mult}_{\Omega}(\mathcal{H})\Omega+\mathcal{L},
$$
where $\mathcal{L}$ is a linear system on the surface $S$ that
does not have fixed components, and $\Omega$ is a smooth rational
curve on the variety $Z$ such that $\eta(\Omega)=\eta(\Gamma)$,
but $\Omega\ne\Gamma$. We have
$$
\mathrm{Sing}(Z)\cap\Gamma=\big\{P_{i_{1}},\ldots,P_{i_{r}}\big\}\subsetneq\big\{P_{1},\ldots,P_{15}\big\}=\mathrm{Sing}(Z),
$$
but $P_{i_{j}}$ is an ordinary double point of the surface $S$.
The equalities $\Gamma^{2}=\Omega^{2}=-2+r/2$ hold on the surface
$S$, but $r\leqslant 3$. Hence, the inequality $\Omega^{2}<0$
holds on the surface $S$, and
$$
\Big(k-\mathrm{mult}_{\Omega}(\mathcal{H})\Big)\Omega^{2}=\Big(\mathrm{mult}_{\Gamma}(\mathcal{H})-k\Big)\Gamma\cdot \Omega+L\cdot\Omega=L\cdot\Omega\geqslant 0,%
$$
where $L$ is a general curve in the linear system $\mathcal{L}$.
Therefore, the inequality
$\mathrm{mult}_{\Omega}(\mathcal{H})\geqslant k$ holds, which
implies that $\Omega\in\mathbb{CS}(Z, \frac{1}{k}\mathcal{H})$,
which is impossible.
\end{proof}

Let $H$ be a general hyperplane in $\mathbb{P}^3$ passing through
the line $\eta(\Gamma)$. Then the curve
$$
\Delta=H\cap R=\eta(\Gamma)\cup\Upsilon
$$
is reduced and $\eta(\Gamma)\not\subset\mathrm{Supp}(\Upsilon)$,
where $\Upsilon$ is a plane curve of degree $5$. Moreover, the
reducible curve $\Delta$ is singular in every singular point
$\eta(P_{i})$ of the surface $R$ that lies on the line
$\eta(\Gamma)$, but the set $\eta(\Gamma)\cap\Upsilon$ contains at
most $5$ points. On the other hand, we have
$$
\mathrm{Sing}(\Delta)\cap\eta(\Gamma)=\Upsilon\cap\eta(\Gamma),
$$
which implies that $|\mathrm{Sing}(Z)\cap\Gamma|\leqslant 5$.
Moreover, the surface $R$ is given by the equation
$$
f_{3}(x,y,z,t)^{2}=4f_{1}(x,y,z,t)f_{5}(x,y,z,t)\subset\mathbb{P}^{3}\cong\mathrm{Proj}\Big(\mathbb{C}[x,y,z,t]\Big),
$$
and singular points of $R$ are given by the equations
$f_{1}=f_{3}=f_{5}=0$. We may assume that the equations
$f_{1}=f_{3}=0$ and $f_{1}=f_{5}=0$ defines irreducible curves in
$\mathbb{P}^{3}$, which implies that at most $3$ points of the
subset $\mathrm{Sing}(R)\subset\mathbb{P}^{3}$ can lie on a single
line. Therefore, the Bertini theorem implies that
$\eta(\Gamma)\cap\Upsilon$ contains different points $O_{1}$ and
$O_{2}$ that are not contained in $\mathrm{Sing}(R)$.

\begin{remark}
\label{remark:n-2-tangent-hyperplane-to-sextic-surface} The
hyperplane $H$ tangents the surface $R$ in the points $O_{1}$ and
$O_{2}$.
\end{remark}

Let $L_{j}$ be a  general line on the plane $H$ that passes
through the point $O_{j}$, $\tilde{O}_{j}$ be a proper transform
of the point $O_{j}$ on the variety $Z$, and $\tilde{L}_{j}$ be
the proper transform of the line $L_{j}$ on the variety $Z$. Then
$L_{j}$ tangents the surface $R$ in the point $O_{j}$, and the
curve $\tilde{L}_{j}$ is irreducible and singular in the point
$\tilde{O}_{j}$, but $-K_{Z}\cdot \tilde{L}_{j}=2$. Let
$\tilde{H}$ be the proper transform of $H$ on the variety $Z$, and
$D$ be a general surface of the linear system $\mathcal{H}$. Then
$$
D\vert_{\tilde{H}}=\mathrm{mult}_{\Gamma}(\mathcal{H})\Gamma+\Psi,
$$
where $\Psi$ is an effective divisor on $\tilde{H}$ such that
$\Gamma\not\subset\mathrm{Supp}(\Psi)$. Let
$\Lambda_{j}=(D\backslash\Gamma)\cap \tilde{L}_{j}$. Then
$$
2k=D\cdot\tilde{L}_{j}\geqslant
\mathrm{mult}_{\tilde{O}_{j}}(\tilde{L}_{j})\mathrm{mult}_{\Gamma}(D)+\sum_{P\in
\Lambda_{j}} \mathrm{mult}_{P}(D)\cdot
\mathrm{mult}_{P}(\tilde{L}_{j})\geqslant 2k+\sum_{P\in
\Lambda_{j}} \mathrm{mult}_{P}(D)\cdot
\mathrm{mult}_{P}(\tilde{L}_{j}),
$$
which implies that $D\cap \tilde{L}_{j}\subset\Gamma$. On the
other hand, when we vary the lines $L_{1}$ and $L_{2}$ on the
plane $H$ the curves $\tilde{L}_{1}$ and $\tilde{L}_{2}$ span two
different pencils on the surface $\tilde{H}$, whose base locus
consist of the points $\tilde{O}_{1}$ and $\tilde{O}_{1}$
respectively. Hence, we have
$$
\mathrm{Supp}(D)\cap\mathrm{Supp}(\tilde{H})=\mathrm{Supp}(\Gamma),
$$
where $\tilde{H}$ is a general divisor in $|-K_{Z}|$ that passes
through the curve $\Gamma$, and $D$ is a general divisor in
$\mathcal{H}$, which contradicts Lemma~\ref{lemma:Cheltsov-II}.
The claim of Proposition~\ref{proposition:n-2} is
proved\footnote{It is easy to see that the given proof implies the
claim of Proposition~\ref{proposition:n-2} under the weaker
assumption that the hypersurface $X$ is quasismooth, and the
projection $\psi\colon X\dasharrow\mathbb{P}^{3}$ contracts $15$
different curves.}.

\section{Case   $n=4$, hypersurface of degree $6$ in $\mathbb{P}(1,1,1,2,2)$.}
\label{section:n-4}

We use the notations and assumptions of
Section~\ref{section:start}. Let $n=4$. Then $X$ is a sufficiently
general hypersurface in $\mathbb{P}(1,1,1,2,2)$ of degree $6$, the
equality $-K_{X}^{3}=3/2$ holds, and the singularities of the
hypersurface $X$ consist of points $P_{1}$, $P_{2}$, $P_{3}$ that
are quotient singularities of types $\frac{1}{2}(1,1,1)$.

Let $\psi\colon X\dasharrow\mathbb{P}^{2}$ be the natural
projection. Then a general fiber of the rational map $\psi$ is an
elliptic curve, and the composition $\psi\circ\eta$ is a morphism,
where $\eta\colon Y\to X$ is a composition of the weighted blow
ups of the singular points $P_{1}$, $P_{2}$ and $P_{3}$ with
weights $(1,1,1)$.

\begin{proposition}
\label{proposition:n-4} The claim of Theorem~\ref{theorem:main}
holds for $n=4$.
\end{proposition}

Let us prove Proposition~\ref{proposition:n-4}. We must show the
existence of the commutative diagram
\begin{equation}
\label{equation:n-4-commutative-diagram} \xymatrix{
&X\ar@{-->}[d]_{\psi}\ar@{-->}[rr]^{\rho}&&V\ar@{->}[d]^{\nu}&\\%
&\mathbb{P}^{2}\ar@{-->}[rr]_{\phi}&&\mathbb{P}^{2},&}
\end{equation}
where $\phi\in\mathrm{Bir}(\mathbb{P}^{2})$. Let $Z$ be an element
of $\mathbb{CS}(X, \frac{1}{k}\mathcal{M})$. We have the following
possibilities:
\begin{itemize}
\item the subvariety $Z$ is a curve that is contained in $X\setminus\mathrm{Sing}(X)$;%
\item the subvariety $Z$ is a curve that contains a singular point of the hypersurface $X$;%
\item the subvariety $Z$ is a singular point of the hypersurface $X$.%
\end{itemize}

Suppose that $Z$ is an irreducible  curve such that $Z$ does not
contain singular points of the hypersurface $X$. Then the equality
$-K_{X}\cdot Z=1$ holds by Lemma~\ref{lemma:curves}, which implies
that the curve $Z$ is smooth. Let $\gamma\colon W\to X$ be the
blow up of the curve $Z$, and $G$ be the exceptional divisor of
the morphism $\gamma$. Then the divisor $\gamma^{*}(-4K_{X})-G$ is
nef, which implies that
$$
\Big(\gamma^{*}(-4K_{X})-G\Big)\cdot\bar{S}_{1}\cdot\bar{S}_{2}\geqslant 0,%
$$
where $\bar{S}_{1}$ and $\bar{S}_{2}$ are proper transforms on the
variety $W$ of sufficiently general surfaces of the linear system
$\mathcal{M}$. We have
$0\leqslant(\gamma^{*}(-4K_{X})-G)\cdot\bar{S}_{1}\cdot\bar{S}_{2}=-k^{2}$,
which is a contradiction.

Suppose that $Z$ is an irreducible curve that passes through some
singular point of the hypersurface $X$. Then the inequality
$-K_{X}\cdot Z\leqslant 1$ holds by Lemma~\ref{lemma:curves}. The
curve $C$ is contracted by the rational map $\psi$ to a point, and
either $-K_{X}\cdot Z=1/2$, or $-K_{X}\cdot Z=1$.

Let $F$ be a sufficiently general surface of the linear system
$|-K_{X}|$ that passes through the curve $Z$. Then the surface $F$
is smooth outside the points $P_{1}$, $P_{2}$ and $P_{3}$, which
are isolated ordinary double points of $F$. Let $\tilde{Z}$ be a
fiber of $\psi$ over the point $\psi(Z)$. Then the generality of
the hypersurface $X$ implies that the curve $Z$ is an irreducible
component of $\tilde{Z}$.

Suppose that  $\tilde{Z}$ consists of irreducible curves $Z$ and
$\bar{Z}$. Then the inequality $\bar{Z}^{2}<0$ holds on the
surface $F$, but $\mathcal{M}\vert_{F}\sim_{\mathbb{Q}}
kZ+k\bar{Z}$. On the other hand, we have
$$
\mathcal{M}\vert_{F}=\mathrm{mult}_{Z}(\mathcal{M})Z+\mathrm{mult}_{\bar{Z}}(\mathcal{M})\bar{Z}+\mathcal{F},
$$
where $\mathcal{F}$ is a linear system on the surface $F$ that
does not have fixed components, but
$$
\Big(k-\mathrm{mult}_{\bar{Z}}\big(\mathcal{M}\big)\Big)\bar{Z}\sim_{\mathbb{Q}}\mathcal{F},
$$
which implies that $\mathrm{mult}_{\bar{Z}}(\mathcal{M})=k$ and
$\mathrm{Supp}(S_{1}\cdot S_{2})\subset Z\cup\bar{Z}$, which
contradicts Lemma~\ref{lemma:Cheltsov}.

Suppose that the fiber $\tilde{Z}$ consists of irreducible curves
$Z$, $\hat{Z}$ and $\check{Z}$. Then
$$
-K_{X}\cdot\hat{Z}=-K_{X}\cdot\check{Z}=-K_{X}\cdot Z=\frac{1}{2},
$$
but the intersection form of the curves $\hat{Z}$ and $\check{Z}$
on the surface $F$ is negatively defined, which implies that the
support of  $S_{1}\cdot S_{2}$ is contained in
$Z\cup\hat{Z}\cup\check{Z}$, where $S_{1}$ and $S_{2}$ are general
surfaces of the linear system $\mathcal{M}$, which is impossible
by Lemma~\ref{lemma:Cheltsov}.

Hence, we have $\mathbb{CS}(X,
\frac{1}{k}\mathcal{M})\subseteq\{P_{1},P_{2},P_{3}\}$.

Let $\pi\colon U\to X$ be a composition of the weighted blow ups
with weights $(1,1,1)$ of the singular points of the hypersurface
$X$ that are contained in the set $\mathbb{CS}(X,
\frac{1}{k}\mathcal{M})$, and $\mathcal{D}$ be the proper
transform of the linear system $\mathcal{M}$ on the variety $U$.
Then it follows from Theorem~\ref{theorem:Kawamata} that the
equivalence $\mathcal{D}\sim_{\mathbb{Q}}-kK_{U}$ holds, but the
divisor $-K_{U}$ is nef.

Suppose that $\mathbb{CS}(U, \frac{1}{k}\mathcal{D})$ contains a
subvariety $\Delta\subset U$. Then $\Delta$ is contained in some
exceptional divisor of $\pi$. Let $G$ be an exceptional divisor of
$\pi$ that contains $\Delta$. Then $\Delta$ is a line on the
surface $G\cong\mathbb{P}^{2}$ by
Lemma~\ref{lemma:Cheltsov-Kawamata}. The linear system
$|\pi^{*}(-2K_{X})-G|$ does not have base points, and the divisor
$\pi^{*}(-2K_{X})-G$ is nef and big. It follows from Lemma~0.3.3
in \cite{KMM} that there is a proper Zariski closed subset
$\Lambda\subset U$ such that $\Lambda$ contains all curves on $U$
having trivial intersections with the divisor
$\pi^{*}(-2K_{X})-G$. We have
$$
2k^{2}=\Big(\pi^{*}\big(-2K_{X}\big)-G\Big)\cdot
\tilde{S}_{1}\cdot\tilde{S}_{2}\geqslant\mathrm{mult}^{2}_{\Delta}\big(\mathcal{D}\big)\Big(\pi^{*}\big(-2K_{X}\big)-G\Big)\cdot\Delta=2k^{2},
$$
where $\tilde{S}_{1}$ and $\tilde{S}_{2}$ are the proper
transforms of $S_{1}$ and $S_{2}$ on $U$ respectively. Thus, the
support of the effective cycle $\tilde{S}_{1}\cdot\tilde{S}_{2}$
is contained in $\Lambda\cup\Delta$, which is impossible by
Lemma~\ref{lemma:Cheltsov}.

Hence, the set $\mathbb{CS}(U, \frac{1}{k}\mathcal{D})$ is empty,
and $-K_{U}^{3}=0$ by Lemma~\ref{lemma:Noether-Fano}. Then
$\pi=\eta$, which implies that $\mathcal{D}$ lies in the fibers of
$\psi\circ\pi$. Therefore, the commutative
diagram~\ref{equation:n-4-commutative-diagram} exists.

\section{Case $n=6$, hypersurface of degree $8$ in $\mathbb{P}(1,1,1,2,4)$.}
\label{section:n-6}

We use the notations and assumptions of
Section~\ref{section:start}. Let $n=6$. Then $X$ is a sufficiently
general hypersurface in $\mathbb{P}(1,1,1,2,4)$ of degree $8$, the
equality $-K_{X}^{3}=1$ holds, and the singularities of the
hypersurface $X$ consists of points $P_{1}$ and $P_{2}$ that are
singularities of type $\frac{1}{2}(1,1,1)$.

There is a commutative diagram
\begin{equation}
\label{equation:n-6-commutative-diagram-and-projection} \xymatrix{
&&U\ar@{->}[dl]_{\pi}\ar@{->}[dr]^{\eta}&&\\%
&X\ar@{-->}[rr]_{\psi}&&\mathbb{P}^{2},&}
\end{equation}
where $\psi$ is the natural projection, $\pi$ is a composition of
the weighted blow ups of the singular points $P_{1}$ and $P_{2}$
with weights $(1,1,1)$, and $\eta$ is an elliptic fibration.

\begin{proposition}
\label{proposition:n-6} The claim of Theorem~\ref{theorem:main}
holds for $n=6$.
\end{proposition}

Let us prove Proposition~\ref{proposition:n-6}. It follows from
Theorem~\ref{theorem:smooth-points} that the set $\mathbb{CS}(X,
\frac{1}{k}\mathcal{M})$ does not contain smooth points of the
hypersurface $X$. Therefore, the set $\mathbb{CS}(X,
\frac{1}{k}\mathcal{M})$ contains a singular point of the
hypersurface $X$ by Corollary~\ref{corollary:cheap-curves} and
Theorem~\ref{theorem:Kawamata}.

\begin{remark}
\label{remark:n-6-two-singular-points} Suppose that the set
$\mathbb{CS}(X, \frac{1}{k}\mathcal{M})$ contains both points
$P_{1}$ and $P_{2}$. Then it easily follows from
Theorem~\ref{theorem:Kawamata} that the claim of the
Theorem~\ref{theorem:main} holds for the hypersurface $X$.
\end{remark}

We may assume that $P_{1}\in\mathbb{CS}(X,
\frac{1}{k}\mathcal{M})$ and $P_{2}\not\in\mathbb{CS}(X,
\frac{1}{k}\mathcal{M})$.

\begin{lemma}
\label{lemma:n-6-curves} The set $\mathbb{CS}(X,
\frac{1}{k}\mathcal{M})$ does not contain curves.
\end{lemma}

\begin{proof}
Suppose that $\mathbb{CS}(X, \frac{1}{k}\mathcal{M})$ contains a
curve $C$. Then $-K_{X}\cdot C=1/2$ by Lemma~\ref{lemma:curves}.

Let $\check{C}$ be the proper transform of the curve $C$ on the
variety $U$. Then $-K_{U}\cdot\check{C}=0$, which implies that the
curve $\check{C}$ is a component of a fiber of $\eta$. Therefore,
the curve $C$ is contracted by the rational map $\psi\colon
X\dasharrow\mathbb{P}^{2}$ to a point. In particular, the curve
$C$ is smooth and rational.

Let $S$ be a general surface of the linear system $|-K_{X}|$ that
contains the curve $C$. Then the surface $S$ is smooth outside the
points $P_{1}$ and $P_{2}$, which are isolated ordinary double
points on the surface $S$. Let $F$ be a fiber of the rational map
$\psi$ over the point $\psi(C)$. Then $F$ consists of two
irreducible components such that the curve $C$ is one of them. Let
$Z$ be the component of the curve $F$ that is different from the
curve $C$. Then $Z^{2}<0$ on the surface $S$, but
$$
\mathcal{M}\vert_{S}=\mathrm{mult}_{C}(\mathcal{M})C+\mathrm{mult}_{Z}(\mathcal{M})Z+\mathcal{L},
$$
where $\mathcal{L}$ is a linear system having no fixed components.
We have
$(k-\mathrm{mult}_{Z}(\mathcal{M}))Z\sim_{\mathbb{Q}}\mathcal{L}$,
which implies that $\mathrm{mult}_{Z}(\mathcal{M})=k$. The latter
is impossible by  Lemmas~\ref{lemma:Cheltsov-II} and
\ref{lemma:normal-surface}.
\end{proof}

Thus, the set $\mathbb{CS}(X, \frac{1}{k}\mathcal{M})$ consists of
the point $P_{1}$. Let $\zeta\colon Y\to X$ be the weighted blow
up of the point $P_{1}$ with weights $(1,1,1)$, and $\mathcal{D}$
be the proper transform of the linear system $\mathcal{M}$ on the
variety $Y$. Then $\mathcal{D}\sim_{\mathbb{Q}}-kK_{Y}$ by
Theorem~\ref{theorem:Kawamata}, and $\mathbb{CS}(Y,
\frac{1}{k}\mathcal{D})\ne\varnothing$ by
Lemma~\ref{lemma:Noether-Fano}.

Let $T$ be a subvariety of $Y$ that is contained in
$\mathbb{CS}(Y, \frac{1}{k}\mathcal{D})$. Then
$\zeta(T)\in\mathbb{CS}(X, \frac{1}{k}\mathcal{M})$, which is
impossible by Lemmas~\ref{lemma:Cheltsov-Kawamata} and
\ref{lemma:curves}.

\section{Case $n=7$, hypersurface of degree $8$ in $\mathbb{P}(1,1,2,2,3)$.}
\label{section:n-7}

We use the notations and assumptions of
Section~\ref{section:start}. Let $n=7$. Then $X$ is a
sufficiently~general hypersurface in $\mathbb{P}(1,1,2,2,3)$ of
degree $8$, whose singularities consist of the point $Q$ that is a
quotient singularity of type $\frac{1}{3}(1,1,2)$, and the points
$P_{1}$, $P_{2}$, $P_{3}$ and $P_{4}$ that are quotient
singularities of type $\frac{1}{2}(1,1,1)$.  There is a
commutative diagram
$$
\xymatrix{
&U_{i}\ar@{->}[d]_{\alpha_{i}}&&Y_{i}\ar@{->}[ll]_{\beta_{i}}\ar@{->}[drr]^{\eta_{i}}&&&\\
&X\ar@{-->}[rrrr]_{\xi_{i}}&&&&\mathbb{P}(1,1,2),&}
$$
where $\xi_{i}$ is a projection, $\alpha_{i}$ is the weighted blow
up of $P_{i}$ with weights $(1,1,1)$, $\beta_{i}$ is the weighted
blow up of the proper transform of  $Q$ on the variety $U_{i}$
with weights $(1,1,2)$, and $\eta_{i}$ is an elliptic fibration.
There is a commutative diagram
$$
\xymatrix{
&U_{0}\ar@{->}[d]_{\alpha_{0}}&&Y_{0}\ar@{->}[ll]_{\beta_{0}}\ar@{->}[drr]^{\eta_{0}}&&&\\
&X\ar@{-->}[rrrr]_{\xi_{0}}&&&&\mathbb{P}(1,1,2),&}
$$
where $\xi_{0}$ is a projection, $\alpha_{0}$ is the weighted blow
up of $Q$ with weights $(1,1,2)$, $\beta_{0}$ is the blow up with
weights $(1,1,1)$ of the singular point of $U_{0}$ that dominates
$Q$, and $\eta_{0}$ is an elliptic fibration.

\begin{proposition}
\label{proposition:n-7} There is a commutative diagram
\begin{equation}
\label{equation:n-7-commutative-diagram} \xymatrix{
&X\ar@{-->}[d]_{\xi_{i}}\ar@{-->}[rr]^{\sigma}&&X\ar@{-->}[rr]^{\rho}&&V\ar@{->}[d]^{\nu}&\\%
&\mathbb{P}(1,1,2)\ar@{-->}[rrrr]_{\phi}&&&&\mathbb{P}^{2}&}
\end{equation}
for some $i\in\{0,1,2,3,4\}$, where $\sigma$ and $\phi$ are
birational maps.
\end{proposition}

Let us prove Proposition~\ref{proposition:n-7}. The set
$\mathbb{CS}(X, \frac{1}{k}\mathcal{M})$ does not contains smooth
points of the hypersurface $X$, and it follows from
Lemma~\ref{lemma:Ryder} that $\mathbb{CS}(X,
\frac{1}{k}\mathcal{M})\subseteq\{P_{1},P_{2},P_{3},P_{4},Q\}$.

\begin{lemma}
\label{lemma:n-7-two-singular-points-of-index-2} The set
$\mathbb{CS}(X, \frac{1}{k}\mathcal{M})$ does not contain two
points of the set $\{P_{1},P_{2},P_{3},P_{4}\}$.
\end{lemma}

\begin{proof}
Suppose that the set $\mathbb{CS}(X, \frac{1}{k}\mathcal{M})$
contains the points $P_{1}$ and $P_{2}$. Let $\pi\colon W\to X$ be
the composition of the weighted blow ups of $P_{1}$ and $P_{2}$
with weights $(1,1,1)$, and $\mathcal{B}$ be the proper transform
of the linear system $\mathcal{M}$ on the variety $W$. Then
$\mathcal{B}\sim_{\mathbb{Q}}-kK_{W}$ by
Theorem~\ref{theorem:Kawamata}, but the base locus of $|-K_{W}|$
consists of a curve $C$ such that $C^{2}=-1/3$ on general surface
of the pencil $|-K_{W}|$, which contradicts
Lemma~\ref{lemma:Cheltsov-II}.
\end{proof}

Let $\mathcal{D}_{i}$ be the proper transform of the linear system
$\mathcal{M}$ on the variety $U_{i}$.

\begin{remark}
\label{remark:n-7-nef-and-big} The set $\mathbb{CS}(U_{0},
\frac{1}{k}\mathcal{D}_{0})$ is not empty by
Lemma~\ref{lemma:Noether-Fano}, if the set $\mathbb{CS}(X,
\frac{1}{k}\mathcal{M})$ contains the point $Q$. Similarly, the
set $\mathbb{CS}(U_{i}, \frac{1}{k}\mathcal{D}_{i})$ is not empty,
if $\mathbb{CS}(X, \frac{1}{k}\mathcal{M})$ contains the point
$P_{i}$.
\end{remark}

\begin{lemma}
\label{lemma:n-7-one-singular-points-of-index-two} The set
$\mathbb{CS}(X, \frac{1}{k}\mathcal{M})$ does not consists of the
point $P_{i}$.
\end{lemma}

\begin{proof}
Suppose that $\mathbb{CS}(X, \frac{1}{k}\mathcal{M})=\{P_{i}\}$.
Then the set $\mathbb{CS}(U_{i}, \frac{1}{k}\mathcal{D}_{i})$
contains an irreducible subvariety $Z\subset U_{i}$ by
Lemma~\ref{lemma:Noether-Fano}. Let $E_{i}$ be the exceptional
divisor of $\alpha_{i}$. Then $Z$ is a line on the surface
$E_{i}\cong\mathbb{P}^{2}$ by Lemma~\ref{lemma:Cheltsov-Kawamata},
which contradicts Lemma~\ref{lemma:curves}.
\end{proof}

Therefore, the condition $P_{i}\in\mathbb{CS}(X,
\frac{1}{k}\mathcal{M})$ implies that $\mathbb{CS}(X,
\frac{1}{k}\mathcal{M})=\{P_{i},Q\}$, which easily implies that
the proper transform of the linear system $\mathcal{M}$ on the
variety $Y_{i}$ lies in the fibers of the elliptic fibration
$\eta_{i}$. The latter implies the existence of the commutative
diagram~\ref{equation:n-7-commutative-diagram}.

We may assume that the set $\mathbb{CS}(X,
\frac{1}{k}\mathcal{M})$ consists of the point $Q$.

Let $O$ be a singular point of the variety $U_{0}$ that is
contained in the exceptional divisor of the morphism $\alpha_{0}$.
Then $O$ is contained in $\mathbb{CS}(U_{0},
\frac{1}{k}\mathcal{D}_{0})$ by
Lemma~\ref{lemma:Cheltsov-Kawamata}, which implies the existence
of the commutative diagram~\ref{equation:n-7-commutative-diagram}
by Theorem~\ref{theorem:Kawamata}. The claim of
Proposition~\ref{proposition:n-7} is proved.

\section{Case $n=8$, hypersurface of degree $9$ in $\mathbb{P}(1,1,1,3,4)$.}
\label{section:n-8}

We use the notations and assumptions of
Section~\ref{section:start}. Let $n=8$. Then $X$ is a sufficiently
general hypersurface in $\mathbb{P}(1,1,1,3,4)$ of degree $9$, the
equality $-K_{X}^{3}=3/4$ holds, and the singularities of the
hypersurface $X$ consist of the singular point $O$ that is a
quotient singularity of type $\frac{1}{4}(1,1,3)$.

\begin{proposition}
\label{proposition:n-8} The claim of Theorem~\ref{theorem:main}
holds for $n=8$.
\end{proposition}

Let us prove Proposition~\ref{proposition:n-8}. The hypersurface
$X$ can be given by the equation
$$
w^{2}z+f_{5}(x,y,z,t)w+f_{9}(x,y,z,t)=0\subset\mathbb{P}(1,1,1,3,5)\cong\mathrm{Proj}\Big(\mathbb{C}[x,y,z,t,w]\Big),
$$
where $\mathrm{wt}(x)=\mathrm{wt}(y)=\mathrm{wt}(z)=1$,
$\mathrm{wt}(t)=3$, $\mathrm{wt}(w)=4$, and $f_{i}$ is a
quasihomogeneous po\-ly\-no\-mi\-al of degree $i$. The point $O$
is given by $x=y=z=t=0$. There is a commutative~diagram
$$
\xymatrix{
&&&W\ar@{->}[dll]_{\sigma}\ar@{->}[d]^{\alpha}&&U\ar@{->}[ll]_{\beta}&&&\\%
&Y\ar@{->}[dr]_{\omega}&&X\ar@{-->}[dl]_{\xi}\ar@{-->}[drrr]^{\psi}&&&&Z\ar@{->}[ld]^{\eta}\ar@{->}[llu]_{\gamma}&\\
&&\mathbb{P}(1,1,1,3)\ar@{-->}[rrrr]_{\chi}&&&&\mathbb{P}^{2},&&}
$$
where rational maps $\xi$, $\psi$ and $\chi$ are the natural
projections, $\alpha$ is the weighted blow up of the singular
point $O$ with weights $(1,1,3)$, $\beta$ is the weighted blow up
with weights $(1,1,2)$ of the singular point of the variety $W$
that is a quotient singularity of type $\frac{1}{3}(1,1,2)$,
$\gamma$ is the weighted blow up with weights $(1,1,1)$ of the
singular point of the variety $U$ that is a quotient singularity
of type $\frac{1}{2}(1,1,1)$, $\eta$ is an elliptic fibration,
$\sigma$ is a birational morphism that contracts $15$ smooth
rational curves to $15$ isolated ordinary double points
$P_{1},\ldots, P_{15}$ of $Y$ respectively, $\omega$ is a double
cover branched over the surface $R\subset\mathbb{P}(1,1,1,3)$ that
is given by the equation
$$
f_{5}(x,y,z,t)^{2}-4zf_{9}(x,y,z,t)=0\subset\mathbb{P}(1,1,1,3)\cong\mathrm{Proj}\Big(\mathbb{C}[x,y,z,t]\Big)
$$
and has $15$ isolated ordinary double points
$\omega(P_{1}),\ldots,\omega(P_{15})$ given by $z=f_{5}=f_{9}=0$.

Let $G$ be the exceptional divisor of the morphism $\alpha$, $F$
be the exceptional divisor of the morphism $\beta$, $P$ be the
singular point of $W$, $Q$ be the singular point of $U$,
$\mathcal{D}$ be the proper transform of $\mathcal{M}$ on $W$,
$\mathcal{H}$ be the proper transform of $\mathcal{M}$ on the
variety $U$, $\mathcal{P}$ be the proper transform of
$\mathcal{M}$ on the variety $Z$, and $\mathcal{B}$ be the proper
transform of $\mathcal{M}$ on the variety $Y$.

\begin{remark}
\label{remark:n-8-generic} The divisors $-K_{W}$ and $-K_{U}$ are
nef and big, and $\omega\circ\sigma(G)$ is given by $z=0$.
\end{remark}

\begin{lemma}
\label{lemma:n-8-curves} The set $\mathbb{CS}(X,
\frac{1}{k}\mathcal{M})$ does not contain curves.
\end{lemma}

\begin{proof}
Let $L$ be a curve in $\mathbb{CS}(X, \frac{1}{k}\mathcal{M})$. It
follows from Theorem~\ref{theorem:Ryder} that there are two
different surfaces $D$ and $D^{\prime}$ in the linear system
$|-K_{X}|$ such that the irreducible curve $L$ is a component of
the cycle $D\cdot D^{\prime}$. The cycle $D\cdot D^{\prime}$ is
reduced and contains at most two components.

Let $\mathcal{P}$ be the pencil in $|-K_{X}|$ that is generated by
the surfaces $D$ and $D^{\prime}$. Then we may assume that $D$ is
a sufficiently general surface of the pencil $\mathcal{P}$.
Applying Lemma~\ref{lemma:Cheltsov-II} together with the proof of
Lemma~\ref{lemma:curves} to the linear system $\mathcal{M}$ and
the pencil $\mathcal{P}$, we immediately obtain a contradiction in
the case when the equality $-K_{X}\cdot L=3/4$ holds. Therefore,
we may assume that either $-K_{X}\cdot L=1/4$ or $-K_{X}\cdot
L=1/2$. We consider only the case $-K_{X}\cdot L=1/4$, because the
case $-K_{X}\cdot L=1/2$ is simpler and very similar.

Let $S_{W}$ be he proper transform on $W$ of the surface given by
$z=0$, and $L_{W}$ be the proper transform on $W$ of the curve
$L$. Then $S_{W}$ must contain $L_{W}$, because
$$
S_{W}\sim_{\mathbb{Q}} \alpha^*\big(-K_{X}\big)-\frac{5}{4}G,
$$
but $G\cdot L_{W}\geqslant 1/3$. Thus, either the curve $L_{W}$ is
contracted by $\sigma$ or the curve $\omega(L_{Y})$ is a ruling of
the cone $\mathbb{P}(1,1,1,3)$ contained in the surface
$\omega\circ\sigma(G)$, where $L_Y$ is the image of $L_W$ by
$\sigma$.

Suppose that the curve $L_{W}$ is not contracted by $\sigma$.
Then, the curve $\omega(L_{Y})$ is not contained in the surface
$R$, which implies that $\omega(L_{Y})$ contains at most one
singular point of $R$ that is different from the point
$\omega\circ\sigma(P)$. Moreover, the curve $\omega(L_{Y})$ must
contain a singular point of the surface $R$ different from
$\omega\circ\sigma(P)$ because $-K_{X}\cdot L=3/4$ otherwise.
Thus, we may assume that the curve  $\omega(L_{Y})$ contains the
point $\omega(P_{1})$.

Let $D_Y$ and $D^{\prime}_Y$ be the proper transforms on $Y$ of
the surfaces $D$ and $D^{\prime}$ respectively. Then the point
$P_{1}$ is an isolated ordinary double point of the surface
$D_{Y}$. Thus, we see that the proper transform of the curve
$L^{\prime}$ on the threefold $W$ is contracted to the point
$P_{1}$ by $\sigma$ and
$$
D_{Y}\cdot D^{\prime}_{Y}=L_{Y}+\bar{L}_{Y},
$$
where $\bar{L}_{Y}$ is a ruling of the cone
$G\cong\mathbb{P}(1,1,3)$. In particular, we have $-K_{X}\cdot
L^{\prime}=1/4$, which contradicts the equality $-K_{X}\cdot
L=1/4$. Hence, the curve $L_{W}$ is contracted by $\sigma$.

Let $L^{\prime}_{W}$ be the proper transform on $W$ of the curve
$L^{\prime}$, and $L^{\prime}_{Y}=\sigma(L^{\prime}_{W})$. Then
$\omega(L^{\prime}_{Y})$ is a ruling of the cone
$\mathbb{P}(1,1,1,3)$ that is contained in the surface
$\omega\circ\sigma(G)$. The curve $\omega(L^{\prime}_{Y})$ is not
contained in the surface $R$, which implies that
$\omega(L^{\prime}_{Y})$ contains at most one singular point of
the surface $R$ different from the point $\omega\circ\sigma(P)$.
The curve $\omega(L^{\prime}_{Y})$ must contain a singular point
of the surface $R$ different from $\omega\circ\sigma(P)$ because
$-K_{X}\cdot L^{\prime}\ne 3/4$. Thus, we may assume that the
curve $\omega(L^{\prime}_{Y})$ contains the point $\omega(P_{1})$.

The point $P_{1}$ is an isolated ordinary double point of $D_{Y}$
and $\sigma(L_{W})=P_{1}$. Hence, we have
$$
D_{Y}\cdot D^{\prime}_{Y}=L^{\prime}_{Y}+\bar{L}^{\prime}_{Y},
$$
where $\bar{L}^{\prime}_{Y}$ is a ruling of
$G\cong\mathbb{P}(1,1,3)$. The intersection $L_{W}\cap
L_{W}^{\prime}$ consists of a point $O^{\prime}$ such that
$O^{\prime}\not\in G$. Hence, the intersection $L\cap L^{\prime}$
contains the point $\alpha(O^{\prime})$ that is different from
$O$.

The surface $D$ is smooth at the point $\alpha(O^{\prime})$, but
$(L+L^{\prime})\cdot L^{\prime}=1/2$ on the surface $D$, which
implies that $L^{\prime}\cdot L^{\prime}<0$ on the surface $D$.
Therefore, we have
$$
\mathcal{M}\vert_{D}=m_{1}L+m_{2}L^{\prime}+\mathcal{L}\equiv  kL+kL^{\prime},%
$$
where $\mathcal{L}$ is a linear system on $D$ that does not have
fixed components, and $m_1$ and $m_{2}$ are natural numbers such
that $m_{1}\geqslant k$. In particular, we have
$$
0\leqslant(m_{1}-k)L^{\prime}\cdot L+\mathcal{L}\cdot L^{\prime}=(k-m_{2})L^{\prime}\cdot L^{\prime},%
$$
which implies that $m_{2}=m_{1}=k$ and
$\mathcal{M}\vert_{D}=kL+kL^{\prime}$, which is impossible by
Lemma~\ref{lemma:Cheltsov-II}.
\end{proof}

It follows from Theorem~\ref{theorem:smooth-points} that
$\mathbb{CS}(X, \frac{1}{k}\mathcal{M})=\{O\}$. Hence, the claim
of Theorem~\ref{theorem:Kawamata} implies that
$\mathcal{D}\sim_{\mathbb{Q}}-kK_{W}$, but $P\in\mathbb{CS}(W,
\frac{1}{k}\mathcal{D})$ by Lemmas~\ref{lemma:Noether-Fano} and
\ref{lemma:Cheltsov-Kawamata}.

\begin{lemma}
\label{remark:n-8-second-floor} Suppose that the set
$\mathbb{CS}(W, \frac{1}{k}\mathcal{D})$ does not contain
subvarieties of $W$ that are different from the point $P$. Then
the claim of Theorem~\ref{theorem:main} holds for the hypersurface
$X$.
\end{lemma}

\begin{proof}
The equivalence $\mathcal{H}\sim_{\mathbb{Q}}-kK_{U}$ holds by
Theorem~\ref{theorem:Kawamata}. Hence, the set $\mathbb{CS}(U,
\frac{1}{k}\mathcal{H})$ contains the point $Q$ by
Lemmas~\ref{lemma:Noether-Fano} and \ref{lemma:Cheltsov-Kawamata}.
We see that $\mathcal{P}$ lies in the fibers of $\eta$ by
Theorem~\ref{theorem:Kawamata}.
\end{proof}

We may assume that the set $\mathbb{CS}(W,
\frac{1}{k}\mathcal{D})$ contains a subvariety $Z$ that is
different from the point $P$. Then $Z$ is a curve in $G$ by
Lemma~\ref{lemma:Cheltsov-Kawamata}, which is impossible by
Corollary~\ref{corollary:curves}.

\section{Case $n=9$, hypersurface of degree $9$ in $\mathbb{P}(1,1,2,3,3)$.}
\label{section:n-9}

We use the notations and assumptions of
Section~\ref{section:start}. Let $n=9$. Then $X$ is a sufficiently
general hypersurface in $\mathbb{P}(1,1,2,3,3)$ of degree $9$, the
equality $-K_{X}^{3}=1/2$ holds, and the singularities of the
hypersurface $X$ consist of the point $O$ that is a quotient
singularity of type $\frac{1}{2}(1,1,1)$, and the points $P_{1}$,
$P_{2}$ and $P_{3}$ that are quotient singularities of type
$\frac{1}{3}(1,1,2)$.

There is a commutative diagram
$$
\xymatrix{
&&&U_{ij}\ar@{->}[dl]_{\beta_{ij}}&&Y\ar@{->}[ll]_{\gamma_{ij}}\ar@{->}[ddr]^{\eta}&&\\%
&&U_{i}\ar@{->}[dr]_{\alpha_{i}}&&&&&\\
&&&X\ar@{-->}[rrr]_{\xi}&&&\mathbb{P}(1,1,2),&}
$$
where $\xi_{i}$ is a projection, $\alpha_{i}$ is the blow up of
$P_{i}$ with weights $(1,1,2)$, $\beta_{ij}$ is the weighted blow
up with weights $(1,1,2)$ of the proper transform of  $P_{j}$ on
the variety $U_{i}$, $\gamma_{ij}$ is the weighted blow up with
weights $(1,1,2)$ of the proper transform of $P_{k}$ on the
$U_{ij}$, and $\eta$ is an elliptic fibration, where $i\ne j$ and
$k\not\in\{i,j\}$.

\begin{remark}
\label{remark:n-9-nef-and-big} The divisors $-K_{U_{i}}$ and
$-K_{U_{ij}}$ are nef and big.
\end{remark}

There is a commutative diagram
$$
\xymatrix{
&&Z\ar@{->}[dl]_{\pi}\ar@{->}[dr]^{\sigma}&&\\
&X\ar@{-->}[rr]_{\chi}&&\mathbb{P}(1,1,3),&}
$$
where $\chi$ is a projection, $\pi$ is the blow up of $Q$ with
weights $(1,1,1)$, and $\sigma$ is an elliptic fibration.

\begin{proposition}
\label{proposition:n-9} Either there is a commutative diagram
\begin{equation}
\label{equation:n-9-commutative-diagram-I} \xymatrix{
&X\ar@{-->}[d]_{\xi}\ar@{-->}[rr]^{\rho}&&V\ar@{->}[d]^{\nu}&\\%
&\mathbb{P}(1,1,2)\ar@{-->}[rr]_{\phi}&&\mathbb{P}^{2},&}
\end{equation}
or there is a commutative diagram
\begin{equation}
\label{equation:n-9-commutative-diagram-II} \xymatrix{
&X\ar@{-->}[d]_{\chi}\ar@{-->}[rr]^{\omega}&&X\ar@{-->}[rr]^{\rho}&&V\ar@{->}[d]^{\nu}&\\%
&\mathbb{P}(1,1,3)\ar@{-->}[rrrr]_{\upsilon}&&&&\mathbb{P}^{2},&}
\end{equation}
where $\phi$, $\omega$ and $\upsilon$ are birational maps.
\end{proposition}

Let us prove Proposition~\ref{proposition:n-9}, which implies the
claim of Theorem~\ref{theorem:main} for $n=9$.

\begin{lemma}
\label{lemma:n-9-singular-point-of-index-2} Suppose that
$Q\in\mathbb{CS}(X, \frac{1}{k}\mathcal{M})$. Then there is a
commutative diagram~\ref{equation:n-9-commutative-diagram-II}.
\end{lemma}

\begin{proof}
Let $\mathcal{B}$ be the proper transform of $\mathcal{M}$ on $Z$.
Then $\mathcal{B}\sim_{\mathbb{Q}}-kK_{X}$ by
Theorem~\ref{theorem:Kawamata}, which implies the equality $S\cdot
C=0$, where $S$ is a general surface of the linear system
$\mathcal{M}$, and $C$ is a general fiber of the morphism
$\sigma$. Thus, the commutative
diagram~\ref{equation:n-9-commutative-diagram-II} exists.
\end{proof}

\begin{lemma}
\label{lemma:n-9-three-singular-points-of-index-3} The commutative
diagram~\ref{equation:n-9-commutative-diagram-I} exists whenever
$\mathbb{CS}(X, \frac{1}{k}\mathcal{M})=\{P_{1},P_{2},P_{3}\}$.
\end{lemma}

\begin{proof}
See the proof of Lemma~\ref{lemma:n-9-singular-point-of-index-2}.
\end{proof}

We may assume that $\mathbb{CS}(X,
\frac{1}{k}\mathcal{M})\subseteq\{P_{2},P_{3}\}$ by
Theorem~\ref{theorem:smooth-points} and Lemma~\ref{lemma:Ryder}.

\begin{lemma}
\label{lemma:n-9-one-singular-points-of-index} The set
$\mathbb{CS}(X, \frac{1}{k}\mathcal{M})$ contains the point
$P_{2}$.
\end{lemma}

\begin{proof}
Suppose that the set $\mathbb{CS}(X, \frac{1}{k}\mathcal{M})$
consists of the point $P_{1}$. Let $\mathcal{B}$ be the proper
transform of the linear system $\mathcal{M}$ on the variety
$U_{1}$. Then $\mathcal{B}\sim_{\mathbb{Q}}-kK_{U_{1}}$ by
Theorem~\ref{theorem:Kawamata}, and Lemma~\ref{lemma:Noether-Fano}
implies that the set $\mathbb{CS}(U_{1}, \frac{1}{k}\mathcal{B})$
is not empty.

Let $Z$ be an element of the set $\mathbb{CS}(U_{1},
\frac{1}{k}\mathcal{B})$, and $G$ be the exceptional  divisor of
the birational morphism $\alpha_{1}$. Then it follows from
Lemmas~\ref{lemma:Cheltsov-Kawamata} and \ref{lemma:curves} that
$Z$ is a singular point of $G$, which is a quotient singularity of
type $\frac{1}{2}(1,1,1)$ on the variety $U_{1}$.

Let $\delta\colon W\to U_{1}$ be the weighted blow up of the point
$Z$ with weights $(1,1,1)$, $\mathcal{D}$ be the proper transform
of the linear system $\mathcal{M}$ on the variety $W$, and $F$ be
a general surface of the linear system $|-K_{W}|$. Then  the
inequality $\Delta^{2}=-1/2$ holds on $F$, but
$\mathcal{D}\vert_{F}\sim_{\mathbb{Q}}k\Delta$ by
Theorem~\ref{theorem:Kawamata}, which is impossible by
Lemmas~\ref{lemma:normal-surface} and \ref{lemma:Cheltsov-II}.
\end{proof}

Now we can apply the arguments of the proof of
Lemma~\ref{lemma:n-9-one-singular-points-of-index} to get a
contradiction.

\section{Case $n=10$, hypersurface of degree $10$ in $\mathbb{P}(1,1,1,3,5)$.}
\label{section:n-10}

We use the notations and assumptions of
Section~\ref{section:start}. Let $n=10$. Then $X$ is a
sufficiently general hypersurface in $\mathbb{P}(1,1,1,3,5)$ of
degree $10$, the equality $-K_{X}^{3}=2/3$ holds, the
singularities of the hypersurface $X$ consist of the point $O$
that is a quotient singularity of type $\frac{1}{3}(1,1,2)$.

It should be pointed out that $X$ is birationally superrigid.

\begin{proposition}
\label{proposition:n-10} The claim of Theorem~\ref{theorem:main}
holds for $n=10$.
\end{proposition}

There is a commutative diagram
$$
\xymatrix{
&&&&W\ar@{->}[dll]_{\alpha}&&&\\%
&&X\ar@{->}[dl]_{\xi}\ar@{-->}[drrr]^{\psi}&&&&Y\ar@{->}[ld]^{\eta}\ar@{->}[llu]_{\beta}&\\
&\mathbb{P}(1,1,1,3)\ar@{-->}[rrrr]_{\chi}&&&&\mathbb{P}^{2},&&}
$$
where $\xi$, $\psi$ and $\chi$ are projections, $\alpha$ is the
weighted blow up of  $O$ with weights $(1,1,1)$, $\beta$ is the
weighted blow up with weights $(1,1,1)$ of the singular point of
$W$, and $\eta$ is an elliptic fibration.

In the rest of the section we prove
Proposition~\ref{proposition:n-10}. Let $Q$ be the unique singular
point of the variety $W$, $\mathcal{D}$ be the proper transforms
of the linear system $\mathcal{M}$ on the variety $W$, and
$\mathcal{H}$ be the proper transforms of the linear system
$\mathcal{M}$ on the variety $Y$.

\begin{lemma}
\label{lemma:n-10-curves-exclusion} The set $\mathbb{CS}(X,
\frac{1}{k}\mathcal{M})$ does not contains curves.
\end{lemma}

\begin{proof}
Suppose that $\mathbb{CS}(X, \frac{1}{k}\mathcal{M})$ contains a
curve $C$. Then $-K_{X}\cdot C=1/3$ by Lemma~\ref{lemma:curves},
which implies that the curve $C$ is contracted by the rational map
$\psi$ to a point.

There is a curve $Z$ on the variety $X$ such that $Z\ne C$, but
$\xi(Z)=\xi(C)$. Let $S$ be a general surface of the linear system
$|-K_{X}|$ that passes through $C$ and $Z$. Then $Z^{2}<0$ on the
normal surface $S$, but $\mathcal{M}\vert_{S}\sim_{\mathbb{Q}}
kC+kZ$, which is impossible by Lemmas~\ref{lemma:Cheltsov-II} and
\ref{lemma:normal-surface}.
\end{proof}

The set $\mathbb{CS}(X, \frac{1}{k}\mathcal{M})$ consists of the
point $O$ by Theorem~\ref{theorem:smooth-points}, the set
$\mathbb{CS}(W, \frac{1}{k}\mathcal{D})$ contains the point $Q$ by
Lemmas~\ref{lemma:Noether-Fano} and \ref{lemma:Cheltsov-Kawamata},
but $\mathcal{H}\sim_{\mathbb{Q}}-kK_{Y}$  by
Theorem~\ref{theorem:Kawamata}, which implies that $\mathcal{H}$
is contained in the fibers of $\eta$. The claim of
Proposition~\ref{proposition:n-10} is proved.

\section{Case $n=12$, hypersurface of degree $10$ in $\mathbb{P}(1,1,2,3,4)$.}
\label{section:n-12}

We use the notations and assumptions of
Section~\ref{section:start}. Let $n=12$. Then
$X$~is~a~general~hyper\-sur\-face in $\mathbb{P}(1,1,2,3,4)$ of
degree $10$, whose singularities consist of the points $P_{1}$ and
$P_{2}$ that are singularities of type $\frac{1}{2}(1,1,1)$, the
point $P_{3}$ that is a  singularity of type $\frac{1}{3}(1,1,2)$,
and the point $P_{4}$ that is a  singularity of type
$\frac{1}{4}(1,1,3)$. There is a commutative diagram
$$
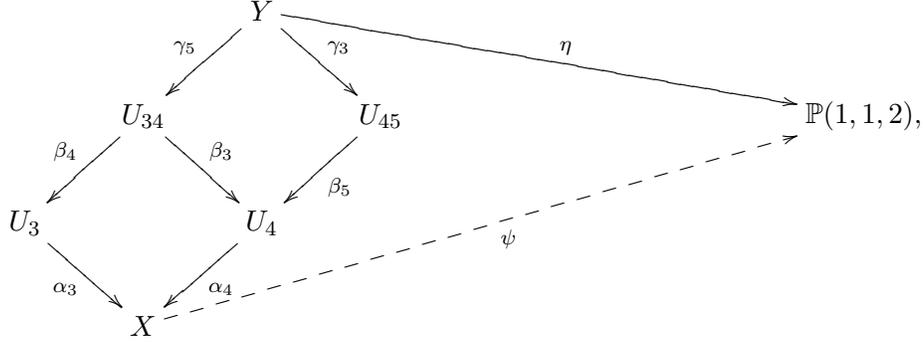
\xymatrix{
&&&Y\ar@{->}[dl]_{\gamma_{5}}\ar@{->}[dr]^{\gamma_{3}}\ar@{->}[drrrrrr]^{\eta}&&&&&&\\
&&U_{34}\ar@{->}[dl]_{\beta_{4}}\ar@{->}[dr]^{\beta_{3}}&&U_{45}\ar@{->}[dl]^{\beta_{5}}&&&&&\mathbb{P}(1,1,2),\\
&U_{3}\ar@{->}[dr]_{\alpha_{3}}&&U_{4}\ar@{->}[dl]^{\alpha_{4}}&&&&&&\\
&&X\ar@{-->}[rrrrrrruu]_{\psi}&&&&&&&}
$$
where $\psi$ is a projection, $\alpha_{3}$ is the weighted blow up
of $P_{3}$ with weights $(1,1,2)$, $\alpha_{4}$ is the weighted
blow up of $P_{4}$ with weights $(1,1,3)$, $\beta_{4}$ is the
weighted blow up with weights $(1,1,3)$ of the proper transform of
the point $P_{4}$ on  $U_{3}$, $\beta_{3}$ is the weighted blow up
with weights $(1,1,2)$ of the proper transform of the point
$P_{3}$ on $U_{4}$, $\beta_{5}$ is the weighted blow up with
weights $(1,1,2)$ of the singular point of $U_{4}$ that is
contained in the exceptional divisor of $\alpha_{4}$, $\gamma_{3}$
is the weighted blow up with weights $(1,1,2)$ of the proper
transform of the point $P_{3}$ on the variety $U_{45}$,
$\gamma_{5}$ is the weighted blow up with weights $(1,1,2)$ of the
singular point of $U_{34}$ that is contained in the exceptional
divisor of the morphism $\beta_{4}$, and $\eta$ is an elliptic
fibration.

\begin{proposition}
\label{proposition:n-12} The claim of Theorem~\ref{theorem:main}
holds for $n=12$.
\end{proposition}

In the rest of this sectiop we prove
Proposition~\ref{proposition:n-12}.

\begin{remark}
\label{remark:n-12-big-nef} The divisors $-K_{U_{3}}$,
$-K_{U_{4}}$, $-K_{U_{34}}$ and $-K_{U_{45}}$ are nef and big.
\end{remark}

It follows from Theorem~\ref{theorem:smooth-points} and
Lemma~\ref{lemma:Ryder} that $\mathbb{CS}(X,
\frac{1}{k}\mathcal{M})\subseteq\{P_{1},P_{2},P_{3},P_{4}\}$.

\begin{lemma}
\label{lemma:n-12-points-of-index-2} The set $\mathbb{CS}(X,
\frac{1}{k}\mathcal{M})$ does not contain $P_{1}$ and $P_{2}$.
\end{lemma}

\begin{proof}
Suppose that $P_{1}\in \mathbb{CS}(X, \frac{1}{k}\mathcal{M})$.
Let $\pi\colon W\to X$ be the weighted blow up of the singular
point $P_{1}$ with weights $(1,1,1)$, $\mathcal{B}$ be the proper
transform of $\mathcal{M}$ on  $W$, $S$ be a general surface of
the pencil $|-K_{W}|$, and $Z$ be the base curve of the pencil
$|-K_{W}|$. Then $Z^{2}=-1/24$ on the surface $S$, but
$\mathcal{B}\vert_{S}\sim_{\mathbb{Q}} kZ$, which contradicts
Lemmas~\ref{lemma:normal-surface} and \ref{lemma:Cheltsov-II}.
\end{proof}

Let $\mathcal{D}_{3}$, $\mathcal{D}_{4}$, $\mathcal{D}_{34}$ and
$\mathcal{D}_{45}$ be the proper transforms of $\mathcal{M}$ on
$U_{3}$, $U_{4}$, $U_{34}$ and $U_{45}$ respectively, then it
follows from Lemma~\ref{lemma:Noether-Fano} that
$\mathbb{CS}(U_{\mu}, \frac{1}{k}\mathcal{D}_{\mu})\ne\varnothing$
in the case when $\mathcal{D}_{\mu}\sim_{\mathbb{Q}}
-kK_{U_{\mu}}$.

\begin{lemma}
\label{lemma:n-12-centers-on-U-4} Suppose that the set
$\mathbb{CS}(X, \frac{1}{k}\mathcal{M})$ contains the point
$P_{4}$. Let $\bar{P}_{3}$ be the proper transform of the point
$P_{3}$ on the variety $U_{4}$, and $\bar{P}_{5}$ be the singular
point of the variety $U_{4}$ that dominates the point $P_{4}$.
Then $\mathcal{D}_{4}\sim_{\mathbb{Q}} -kK_{U_{4}}$ and
$\mathbb{CS}(U_{4},
\frac{1}{k}\mathcal{D}_{4})\subseteq\{\bar{P}_{3},\bar{P}_{5}\}$.
\end{lemma}

\begin{proof}
Suppose that $\mathbb{CS}(U_{4},
\frac{1}{k}\mathcal{D}_{4})\not\subseteq\{\bar{P}_{3},\bar{P}_{5}\}$.
Let $C$ be an element of $\mathbb{CS}(U_{4},
\frac{1}{k}\mathcal{D}_{4})$ different from the points
$\bar{P}_{3}$ and $\bar{P}_{5}$. Then $\alpha_{4}(C)=P_{4}$, which
contradicts Lemma~\ref{lemma:Cheltsov-Kawamata} and
Corollary~\ref{corollary:curves}.
\end{proof}

\begin{lemma}
\label{lemma:n-12-centers-on-U-3} Suppose that the set
$\mathbb{CS}(X, \frac{1}{k}\mathcal{M})$ contains the point
$P_{3}$. Let $\bar{P}_{4}$ be the proper transform of the point
$P_{4}$ on the variety $U_{3}$. Then
$\mathcal{D}_{3}\sim_{\mathbb{Q}}-kK_{U_{3}}$ and
$\mathbb{CS}(U_{3}, \frac{1}{k}\mathcal{D}_{3})=\{\bar{P}_{4}\}$.
\end{lemma}

\begin{proof}
See the proof of
Lemma~\ref{lemma:n-9-one-singular-points-of-index}.
\end{proof}

It follows from Theorem~\ref{theorem:Kawamata} that either
$\mathcal{D}_{34}\sim_{\mathbb{Q}}-kK_{U_{34}}$ or
$\mathcal{D}_{45}\sim_{\mathbb{Q}}-kK_{U_{45}}$. Let $\mathcal{D}$
be the proper transform of the linear system $\mathcal{M}$ on the
variety $Y$.

\begin{lemma}
\label{lemma:n-12-centers-on-U-34} Suppose that
$\mathcal{D}_{34}\sim_{\mathbb{Q}}-kK_{U_{34}}$. Then
$\mathcal{D}\sim_{\mathbb{Q}} -kK_{Y}$.
\end{lemma}

\begin{proof}
Let $F$ be the exceptional  divisor of the morphism $\beta_{3}$,
$G$ be the exceptional  divisor of the morphism $\beta_{4}$,
$\check{P}_{5}$ be the singular point of the surface $G$, and
$\check{P}_{6}$ be the singular point of the surface $F$. Then
$\mathcal{D}\sim_{\mathbb{Q}}-kK_{Y}$ by
Theorem~\ref{theorem:Kawamata} if the set $\mathbb{CS}(U_{34},
\frac{1}{k}\mathcal{D}_{34})$ contains the point $\check{P}_{5}$.

Suppose that $\check{P}_{5}\not\in\mathbb{CS}(U_{34},
\frac{1}{k}\mathcal{D}_{34})$. Then
$\check{P}_{6}\in\mathbb{CS}(U_{34}, \frac{1}{k}\mathcal{D}_{34})$
by Lemmas~\ref{lemma:Noether-Fano} and
\ref{lemma:Cheltsov-Kawamata}.

Let $\pi\colon W\to U_{34}$ be the weighted blow up of the point
$\check{P}_{6}$ with weights $(1,1,1)$, $\mathcal{B}$ be the
proper transform of the linear system $\mathcal{M}$ on the variety
$W$, and $S$ be a general surface of the linear system $|-K_{W}|$.
Then the base locus of $|-K_{W}|$ consists of the irreducible
curve $\Delta$ such that
$\mathcal{B}\vert_{S}\sim_{\mathbb{Q}}k\Delta$, but $\Delta^{2}<0$
holds on $S$, which contradicts Lemmas~\ref{lemma:normal-surface}
and \ref{lemma:Cheltsov-II}.
\end{proof}

The proof of Lemma~\ref{lemma:n-12-centers-on-U-34} implies that
$\mathcal{D}\sim_{\mathbb{Q}}-kK_{Y}$ whenever
$\mathcal{D}_{45}\sim_{\mathbb{Q}}-kK_{U_{45}}$, but it follows
from the equivalence $\mathcal{D}\sim_{\mathbb{Q}}-kK_{Y}$ that
the linear system $\mathcal{D}$ lies in the fibers of the elliptic
fibrations $\eta$. Hence, the claim of
Proposition~\ref{proposition:n-12} is proved.

\section{Case  $n=13$, hypersurface of degree $11$ in $\mathbb{P}(1,1,2,3,5)$.}
\label{section:n-13}

We use the notations and assumptions of
Section~\ref{section:start}. Let $n=13$. Then $X$ is a
sufficiently general hypersurface in $\mathbb{P}(1,1,2,3,5)$ of
degree $11$, the singularities of the hypersurface $X$ consist of
the point $P_{1}$ that is a quotient singularity of type
$\frac{1}{2}(1,1,1)$, the point $P_{2}$ that is a singularity of
type $\frac{1}{3}(1,1,2)$, the point $P_{3}$ that is a singularity
of type $\frac{1}{5}(1,2,3)$, and $-K_{X}^{3}=11/30$.

\begin{proposition}
\label{proposition:n-13} The claim of Theorem~\ref{theorem:main}
holds for  $n=13$.
\end{proposition}

Let us prove Proposition~\ref{proposition:n-13}. There is a
commutative diagram
$$
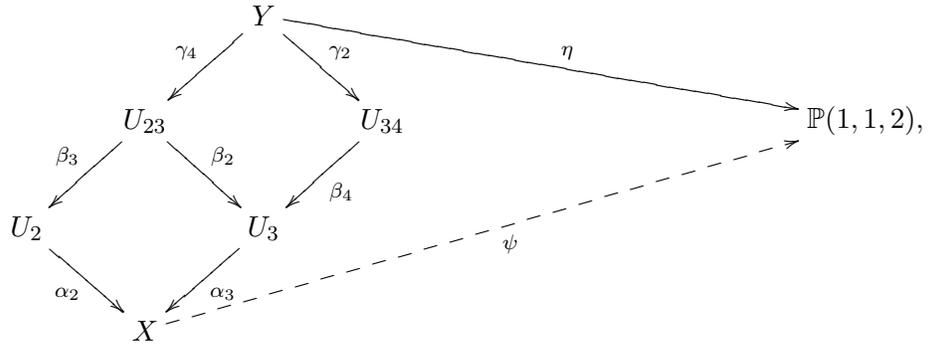
\xymatrix{
&&&Y\ar@{->}[dl]_{\gamma_{4}}\ar@{->}[dr]^{\gamma_{2}}\ar@{->}[drrrrrr]^{\eta}&&&&&&\\
&&U_{23}\ar@{->}[dl]_{\beta_{3}}\ar@{->}[dr]^{\beta_{2}}&&U_{34}\ar@{->}[dl]^{\beta_{4}}&&&&&\mathbb{P}(1,1,2),\\
&U_{2}\ar@{->}[dr]_{\alpha_{2}}&&U_{3}\ar@{->}[dl]^{\alpha_{3}}&&&&&&\\
&&X\ar@{-->}[rrrrrrruu]_{\psi}&&&&&&&}
$$
where $\psi$ is a projection, $\alpha_{2}$ is the weighted blow up
of  $P_{2}$ with weights $(1,1,2)$, $\alpha_{3}$ is the weighted
blow up of $P_{3}$ with weights $(1,2,3)$, $\beta_{3}$ is the
weighted blow up with weights $(1,2,3)$ of the proper transform of
the point $P_{3}$ on the variety $U_{2}$, $\beta_{2}$ is the
weighted blow up with weights $(1,1,2)$ of the proper transform of
$P_{2}$ on the variety $U_{3}$, $\beta_{4}$ is the weighted blow
up with weights $(1,1,2)$ of the singular point of the variety
$U_{3}$ that is the quotient singularity of type
$\frac{1}{3}(1,1,2)$ contained in the exceptional divisor of
$\alpha_{3}$, $\gamma_{2}$ is the weighted blow up with weights
$(1,1,2)$ of the proper transform of $P_{2}$ on the variety
$U_{34}$, $\gamma_{4}$ is the weighted blow up with weights
$(1,1,2)$ of the singular point of the variety $U_{23}$ that is
the quotient singularity of type $\frac{1}{3}(1,1,2)$ contained in
the exceptional divisor of the morphism $\beta_{3}$, and $\eta$ is
an elliptic fibration.

It follows from Theorem~\ref{theorem:smooth-points},
Lemma~\ref{lemma:Ryder} and
Proposition~\ref{proposition:singular-points} that $\mathbb{CS}(X,
\frac{1}{k}\mathcal{M})\subseteq\{P_{2},P_{3}\}$.

\begin{lemma}
\label{lemma:n-13-P3-is-a-center} The set $\mathbb{CS}(X,
\frac{1}{k}\mathcal{M})$ contains the point $P_{3}$.
\end{lemma}

\begin{proof}
Suppose that $P_{3}\not\in\mathbb{CS}(X, \frac{1}{k}\mathcal{M})$.
Then $\mathbb{CS}(X, \frac{1}{k}\mathcal{M})=\{P_{2}\}$. Let
$\mathcal{D}_{2}$ be the proper transform of $\mathcal{M}$ on the
variety $U_{2}$. Then it follows from
Theorem~\ref{theorem:Kawamata} that
$\mathcal{D}_{2}\sim_{\mathbb{Q}}-kK_{U_{2}}$, but the set
$\mathbb{CS}(U_{2}, \frac{1}{k}\mathcal{D}_{2})$ is not empty by
Lemma~\ref{lemma:Noether-Fano}.

Let $E$ be the exceptional  divisor of the morphism $\alpha_{2}$.
Then $E$ can be identified with a cone over the smooth rational
curve in $\mathbb{P}^{3}$ of degree $3$. Let $Z$ be a subvariety
of the variety $U_{2}$ that is contained in the set
$\mathbb{CS}(U_{2}, \frac{1}{k}\mathcal{D}_{2})$. Then it follows
from Lemmas~\ref{lemma:Cheltsov-Kawamata} and \ref{lemma:curves}
that $Z$ is the vertex of the cone $E$, which is a quotient
singularity of type $\frac{1}{2}(1,1,1)$ on the variety $U_{2}$.

Let $\pi\colon W\to U_{2}$ be the  blow up of the point $Z$ with
weights $(1,1,1)$, and $S$ be a sufficiently general surface of
the pencil $|-K_{W}|$. Then the base locus of the pencil
$|-K_{W}|$ consists of an irreducible curve $\Delta$ such that
$\Delta^{2}=-3/10$ on the surface $S$, but
$\mathcal{B}\vert_{S}\sim_{\mathbb{Q}} k\Delta$, where
$\mathcal{B}$ is the proper transform of $\mathcal{M}$ on $W$. We
have $\mathcal{B}\vert_{S}=k\Delta$, which is impossible by
Lemma~\ref{lemma:Cheltsov-II}.
\end{proof}

\begin{lemma}
\label{lemma:n-13-P2-is-a-center} The set $\mathbb{CS}(X,
\frac{1}{k}\mathcal{M})$ contains the point $P_{2}$.
\end{lemma}

\begin{proof}
Suppose that $\mathbb{CS}(X, \frac{1}{k}\mathcal{M})$ does not
contain the point $P_{2}$. Then $\mathbb{CS}(X,
\frac{1}{k}\mathcal{M})$ consists of the point $P_{3}$. Let
$\mathcal{D}_{3}$ be the proper transform of $\mathcal{M}$ on
$U_{3}$. Then $\mathbb{CS}(U_{3},
\frac{1}{k}\mathcal{D}_{3})\ne\varnothing$~by~Lemma~\ref{lemma:Noether-Fano}.

Let $G$ be the exceptional  divisor of the morphism $\alpha_{3}$,
$P_{4}$ be the singular point of  $G$ that is a quotient
singularity of type $\frac{1}{3}(1,1,2)$ on the variety $U_{3}$,
and $P_{5}$ be the singular point of the surface $G$ that is a
quotient singularity of type $\frac{1}{2}(1,1,1)$ on $U_{3}$. Then
$G\cong\mathbb{P}(1,2,3)$, and it follows from
Lemma~\ref{lemma:Cheltsov-Kawamata} that either
$\mathbb{CS}(U_{3}, \frac{1}{k}\mathcal{D}_{3})=\{P_{4}\}$, or
$P_{5}\in\mathbb{CS}(U_{3}, \frac{1}{k}\mathcal{D}_{3})$.

Suppose that the set $\mathbb{CS}(U_{3},
\frac{1}{k}\mathcal{D}_{3})$ contains the point $P_{5}$. Let
$\pi\colon W\to U_{3}$ be the weighted blow up of $P_{5}$ with
weights $(1,1,1)$, $\mathcal{B}$ be the proper transform of
$\mathcal{M}$ on $W$, $L$ be the curve on the surface $G$ that is
contained in  $|\mathcal{O}_{\mathbb{P}(1,\,2,\,3)}(1)|$,
$\bar{L}$ be the proper transform of the curve $L$ on the variety
$W$, and $S$ be a general surface in $|-K_{W}|$. Then
$\mathcal{B}\vert_{S}\sim_{\mathbb{Q}}k\Delta+k\bar{L}$ by
Theorem~\ref{theorem:Kawamata}, the base locus of $|-K_{W}|$
consists of the curves $\bar{L}$ and $\Delta$ such that
$\alpha\circ\pi(\Delta)$ is the base curve of the pencil
$|-K_{X}|$. The equalities
$$
\Delta^{2}=-5/6,\ \bar{L}^{2}=-4/3,\ \Delta\cdot\bar{L}=1
$$
hold on the surface $S$, which imply that the intersection form of
the curves $\Delta$ and $\bar{L}$ on the surface $S$ is negatively
defined. The latter is impossible by
Lemmas~\ref{lemma:Cheltsov-II} and \ref{lemma:normal-surface}.

Hence, the set $\mathbb{CS}(U_{3}, \lambda\mathcal{D}_{3})$
consists of the point $P_{4}$. Let $D_{34}$ be the proper
transform of the linear system $\mathcal{M}$ on the variety
$U_{34}$. Then $\mathbb{CS}(U_{34},
\frac{1}{k}\mathcal{D}_{34})\ne\varnothing$ by
Lemma~\ref{lemma:Noether-Fano}, because the equivalence
$D_{34}\sim_{\mathbb{Q}}-kK_{U_{34}}$ holds by
Theorem~\ref{theorem:Kawamata}.

Let $E$ be the exceptional divisor of the morphism $\beta_{4}$,
and $P_{6}$ be the singular point of the surface $E$. Then
$E\cong\mathbb{P}(1,1,3)$, the point $P_{6}$ is a quotient
singularity of type $\frac{1}{2}(1,1,1)$ on the variety $U_{34}$,
and $P_{6}\in\mathbb{CS}(U_{34}, \frac{1}{k}\mathcal{D}_{34})$ by
Lemma~\ref{lemma:Cheltsov-Kawamata}.

Let $\zeta\colon Z\to U_{34}$ be the weighted blow up of $P_{6}$
with weights $(1,1,1)$, $\mathcal{H}$ be the proper transform of
$\mathcal{M}$ on $Z$, and $F$ be a general surface of the pencil
$|-K_{Z}|$. Then the base locus of the pencil $|-K_{Z}|$ consists
of irreducible curves $\check{L}$ and $\check{\Delta}$ such that
the equalities
$$
\check{\Delta}^{2}=-5/6,\ \check{L}^{2}=-3/2,\ \check{\Delta}\cdot\check{L}=1%
$$
hold on $F$. Thus, the intersection form of the curves
$\check{\Delta}$ and $\check{L}$ on the surface $F$ is negatively
defined, but
$\mathcal{H}\vert_{F}\sim_{\mathbb{Q}}k\check{\Delta}+k\check{L}$
by Theorem~\ref{theorem:Kawamata}, which contradicts
Lemmas~\ref{lemma:Cheltsov-II} and \ref{lemma:normal-surface}.
\end{proof}

Let $\mathcal{D}_{23}$ be the proper transform of  $\mathcal{M}$
on $U_{23}$. Then $\mathbb{CS}(U_{23},
\frac{1}{k}\mathcal{D}_{23})\ne\varnothing$ by
Lemma~\ref{lemma:Noether-Fano}.

\begin{remark}
\label{remark:n-13-final} The proof of
Lemma~\ref{lemma:n-13-P3-is-a-center} implies that the set
$\mathbb{CS}(U_{23}, \frac{1}{k}\mathcal{D}_{23})$ does no contain
the singular point of the variety $U_{23}$ that is contained in
the exceptional divisor of $\beta_{2}$, but the proof of
Lemma~\ref{lemma:n-13-P2-is-a-center} implies that
$\mathbb{CS}(U_{23}, \frac{1}{k}\mathcal{D}_{23})$ does not
contain the singular point of $U_{23}$ that is a singularity of
type $\frac{1}{2}(1,1,1)$ contained in the exceptional divisor of
$\beta_{3}$.
\end{remark}

Therefore, the set $\mathbb{CS}(U_{23},
\frac{1}{k}\mathcal{D}_{23})$ contains the singular point the
variety $U_{23}$ that is a singularity of type
$\frac{1}{3}(1,1,2)$ contained in the $\beta_{3}$-exceptional
divisor, which implies that the proper transform of $\mathcal{M}$
on the variety $Y$ is contained in the fibers of $\eta$ by
Theorem~\ref{theorem:Kawamata}.

\section{Case   $n=15$, hypersurface of degree $12$ in $\mathbb{P}(1,1,2,3,6)$.}
\label{section:n-15}

We use the notations and assumptions of
Section~\ref{section:start}. Let $n=15$. Then $X$ is a
sufficiently general hypersurface in $\mathbb{P}(1,1,2,3,6)$ of
degree $12$, the equality $-K_{X}^{3}=1/3$ holds, and the
singularities of the hypersurface  $X$ consist of the points
$P_{1}$ and $P_{2}$ that are quotient singularities of type
$\frac{1}{2}(1,1,1)$, and the points $P_{3}$ and $P_{4}$ that are
quotient singularities of type $\frac{1}{3}(1,1,2)$.

There is a commutative diagram
$$
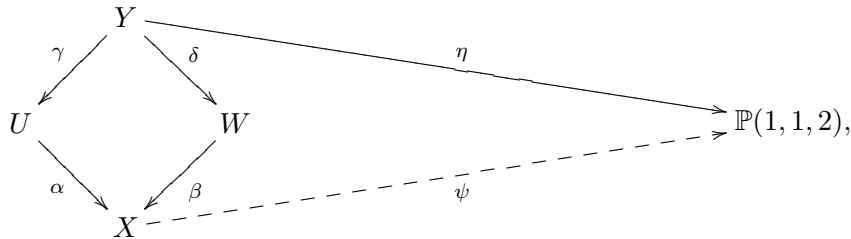
\xymatrix{
&&Y\ar@{->}[dl]_{\gamma}\ar@{->}[dr]^{\delta}\ar@{->}[drrrrrrr]^{\eta}&&&&&&&\\
&U\ar@{->}[dr]_{\alpha}&&W\ar@{->}[dl]^{\beta}&&&&&&\mathbb{P}(1,1,2),\\
&&X\ar@{-->}[rrrrrrru]_{\psi}&&&&&&&}
$$
where $\psi$ is a projection, $\alpha$ is the weighted blow up of
the point $P_{3}$ with weights $(1,1,2)$, $\beta$ is the weighted
blow up of $P_{4}$ with weights $(1,1,2)$, $\gamma$ is the
weighted blow up with weights $(1,1,2)$ of the proper transform of
$P_{4}$ on the variety $U$, $\delta$ is the weighted blow up with
weights $(1,1,2)$ of the proper transform of the point $P_{3}$ on
the variety $W$, and $\eta$ is an elliptic fibration.

\begin{proposition}
\label{proposition:n-15} The claim of Theorem~\ref{theorem:main}
holds for $n=15$.
\end{proposition}

\begin{proof}
We have $\mathbb{CS}(X,
\frac{1}{k}\mathcal{M})\subseteq\{P_{3},P_{4}\}$ by
Theorem~\ref{theorem:smooth-points},
Proposition~\ref{proposition:singular-points} and
Lemma~\ref{lemma:Ryder}.

Let $\mathcal{H}$ be the proper transform of the linear system
$\mathcal{M}$ on the variety $Y$. To conclude the proof we must
show that $\mathcal{H}$ lies in the fibers of the morphism $\eta$,
which easily follows from the condition $\mathbb{CS}(X,
\frac{1}{k}\mathcal{M})=\{P_{3},P_{4}\}$ by
Theorem~\ref{theorem:Kawamata}. We may assume that
$P_{4}\not\in\mathbb{CS}(X, \frac{1}{k}\mathcal{M})$.

Let $\mathcal{B}$ be the proper transform of $\mathcal{M}$ on $U$,
and $Q$ be the singular point of $U$ that is contained in the
exceptional divisor of the morphism $\alpha$. Then
$\mathcal{B}\sim_{\mathbb{Q}}-kK_{U}$ by
Theorem~\ref{theorem:Kawamata}, and~it~follows from
Lemmas~\ref{lemma:Noether-Fano} and \ref{lemma:Cheltsov-Kawamata}
that the set $\mathbb{CS}(U, \frac{1}{k}\mathcal{B})$ contains the
singular point $Q$.

Let $\zeta\colon Z\to U$ be the weighted blow up of the point $Q$
with weights $(1,1,1)$, $\mathcal{D}$ be the proper transform of
$\mathcal{M}$ on the threefold $Z$, and $S$ be a general surface
in $|-K_{Z}|$. Then the base locus of the pencil $|-K_{Z}|$
consists of an irreducible curve $\Delta$ such that $\Delta^{2}<0$
on the surface $S$, which is impossible by
Lemmas~\ref{lemma:normal-surface} and \ref{lemma:Cheltsov-II},
because $\mathcal{D}\vert_{S}\sim_{\mathbb{Q}}-k\Delta$ by
Theorem~\ref{theorem:Kawamata}.
\end{proof}

\section{Case $n=16$, hypersurface of degree $12$ in $\mathbb{P}(1,1,2,4,5)$.}
\label{section:n-16}

We use the notations and assumptions of
Section~\ref{section:start}. Let $n=16$. Then $X$ is a
sufficiently general hypersurface in $\mathbb{P}(1,1,2,4,5)$ of
degree $12$, the equality $-K_{X}^{3}=3/10$ holds, and the
singularities of the hypersurface $X$ consist of the points
$P_{1}$, $P_{2}$ and $P_{3}$ that are quotient singularities of
type $\frac{1}{2}(1,1,1)$, and the point $P_{4}$ that is a
quotient singularity of type $\frac{1}{5}(1,1,4)$.

\begin{proposition}
\label{proposition:n-16} The claim of Theorem~\ref{theorem:main}
holds for $n=16$.
\end{proposition}

Let us prove Proposition~\ref{proposition:n-16}. There is a
commutative diagram
$$
\xymatrix{
&U\ar@{->}[d]_{\alpha}&&W\ar@{->}[ll]_{\beta}&&Y\ar@{->}[ll]_{\gamma}\ar@{->}[d]^{\eta}&\\
&X\ar@{-->}[rrrr]_{\psi}&&&&\mathbb{P}(1,1,2)&}
$$
where $\psi$ is the natural projection, $\alpha$ is the weighted
blow up of  $P_{4}$ with weights $(1,1,4)$, $\beta$ is the
weighted blow up with weights $(1,1,3)$ of the singular point of
the variety $U$ that is contained in the exceptional divisor of
$\alpha$, $\gamma$ is the weighted blow up with weights $(1,1,2)$
of the singular point of $W$ that is contained in the exceptional
divisor of $\beta$, and $\eta$ is an elliptic fibration.

\begin{remark}
\label{remark:n-16-big-nef} The divisors $-K_{U}$ and $-K_{W}$ are
nef and big.
\end{remark}

It follows from Theorem~\ref{theorem:smooth-points},
Lemma~\ref{lemma:Ryder} and
Proposition~\ref{proposition:singular-points} that $\mathbb{CS}(X,
\frac{1}{k}\mathcal{M})$ consists of the point $P_{4}$. Let $G$ be
the exceptional divisor of the morphism $\alpha$, $\bar{P}_{5}$ be
the singular point of the surface $G$, and $\mathcal{D}$ be the
proper transform of $\mathcal{M}$ on the variety $U$. Then $G$ is
a cone over a smooth rational curve of degree $4$, and
$\bar{P}_{5}$ is a  singularity of type $\frac{1}{4}(1,1,3)$ on
the variety $U$.

We see that $\mathcal{D}\sim_{\mathbb{Q}}-kK_{U}$ by
Theorem~\ref{theorem:Kawamata} and $\mathbb{CS}(U,
\frac{1}{k}\mathcal{D})\ne\varnothing$ by
Lemma~\ref{lemma:Noether-Fano}.

\begin{lemma}
\label{lemma:n-16-point-P5} The set $\mathbb{CS}(U,
\frac{1}{k}\mathcal{D})$ consists of the point $\bar{P}_{5}$.
\end{lemma}

\begin{proof}
Suppose that the set $\mathbb{CS}(U, \frac{1}{k}\mathcal{D})$
contains a subvariety $C$ of the variety $U$ that is different
from the point $\bar{P}_{5}$. Then it follows from
Lemma~\ref{lemma:Cheltsov-Kawamata} that $C$ is a ruling of the
cone $G$, which is impossible by Corollary~\ref{corollary:curves}.
\end{proof}

Let $\mathcal{H}$ be the proper transform of $\mathcal{M}$ on $W$.
Then it follows from Theorem~\ref{theorem:Kawamata} and
Lemmas~\ref{lemma:Noether-Fano} and \ref{lemma:Cheltsov-Kawamata}
that the set $\mathbb{CS}(W, \frac{1}{k}\mathcal{H})$ contains the
singular point of the variety $W$ that is contained in the
exceptional divisor of $\beta$. Therefore, the proper transform of
the linear system $\mathcal{M}$ on the variety $Y$ lies in the
fibers of the morphism $\eta$ by Theorem~\ref{theorem:Kawamata}.

\section{Case $n=17$, hypersurface of degree $12$ in $\mathbb{P}(1,1,3,4,4)$.}
\label{section:n-17}

We use the notations and assumptions of
Section~\ref{section:start}. Let $n=17$. Then
$X$~is~a~general~hyper\-surface in $\mathbb{P}(1,1,3,4,4)$ of
degree $12$, whose singularities consist of the points $P_{1}$,
$P_{2}$ and $P_{3}$ that are singularities of type
$\frac{1}{4}(1,1,3)$. There is a commutative diagram
$$
\xymatrix{
&&Z\ar@{->}[dl]_{\pi}\ar@{->}[dr]^{\omega}&&\\
&X\ar@{-->}[rr]_{\psi}&&\mathbb{P}(1,1,3),&}
$$
where $\psi$ is a projection, $\pi$ is a composition of the
weighted blow ups of $P_{1}$, $P_{2}$ and $P_{3}$ with weights
$(1,1,3)$, and $\omega$ is and elliptic fibration. There is a
commutative diagram
$$
\xymatrix{
&&&Y_{i}\ar@{->}[dl]_{\beta_{i}}\ar@{->}[ddrr]^{\eta_{i}}&&\\%
&&U_{i}\ar@{->}[dr]_{\alpha_{i}}&&&&\\
&&&X\ar@{-->}[rr]_{\xi_{i}}&&\mathbb{P}(1,1,4),&}
$$
where $\xi_{i}$ is a projection, $\alpha_{i}$ is the blow up of
$P_{i}$ with weights $(1,1,3)$, $\beta_{i}$ is the weighted blow
up with weights $(1,1,2)$ of the singular point the variety
$U_{i}$ that is contained in the exceptional divisor of the
morphism $\alpha_{i}$, and $\eta_{i}$ is an elliptic fibration.

\begin{proposition}
\label{proposition:n-17} Either there is a commutative diagram
\begin{equation}
\label{equation:n-17-commutative-diagram-I} \xymatrix{
&X\ar@{-->}[d]_{\psi}\ar@{-->}[rr]^{\rho}&&V\ar@{->}[d]^{\nu}&\\%
&\mathbb{P}(1,1,3)\ar@{-->}[rr]_{\phi}&&\mathbb{P}^{2},&}
\end{equation}
or there is a commutative diagram
\begin{equation}
\label{equation:n-17-commutative-diagram-II} \xymatrix{
&X\ar@{-->}[d]_{\xi_{i}}\ar@{-->}[rr]^{\sigma}&&X\ar@{-->}[rr]^{\rho}&&V\ar@{->}[d]^{\nu}&\\%
&\mathbb{P}(1,1,4)\ar@{-->}[rrrr]_{\upsilon}&&&&\mathbb{P}^{2}&}
\end{equation}
for some $i\in\{1,2,3\}$, where $\phi$, $\sigma$ and $\upsilon$
are birational maps.
\end{proposition}

Let us prove Proposition~\ref{proposition:n-17}, which implies the
claim of Theorem~\ref{theorem:main} for $n=17$.

\begin{remark}
\label{remark:n-17-three-points} It easily follows from
Theorem~\ref{theorem:Kawamata} that the commutative
diagram~\ref{equation:n-17-commutative-diagram-I} exists, if the
set $\mathbb{CS}(X, \frac{1}{k}\mathcal{M})$ contains the point
$P_{1}$, $P_{2}$ and $P_{3}$.
\end{remark}

It follows from Theorem~\ref{theorem:smooth-points} and
Lemma~\ref{lemma:Ryder} that $\mathbb{CS}(X,
\frac{1}{k}\mathcal{M})\subseteq\{P_{1},P_{2},P_{3}\}$, but we may
assume that the set $\mathbb{CS}(X, \frac{1}{k}\mathcal{M})$
contains the point $P_{1}$ and does not contain the point $P_{3}$.

\begin{remark}
\label{remark:n-17-nef-and-big} The divisor $-K_{U_{1}}$ is nef
and big.
\end{remark}

Let $\mathcal{D}$ be the proper transform of the linear system
$\mathcal{M}$ on the variety $U_{1}$, $\bar{P}_{2}$ be the proper
transform of the point $P_{2}$ on the variety $U_{1}$, and
$\bar{P}_{4}$ be the singular point of the variety $U_{1}$ that is
contained in the exceptional divisor of $\alpha_{1}$. Then
$\mathcal{D}\sim_{\mathbb{Q}}-kK_{U_{1}}$ by
Theorem~\ref{theorem:Kawamata}, and it follows from
Lemma~\ref{lemma:Noether-Fano} that the set $\mathbb{CS}(U_{1},
\frac{1}{k}\mathcal{D})$ is not empty.

\begin{remark}
\label{remark:n-17-one-point} It easily follows from
Theorem~\ref{theorem:Kawamata} that the commutative
diagram~\ref{equation:n-17-commutative-diagram-II} exists in the
case when the set $\mathbb{CS}(U_{1}, \frac{1}{k}\mathcal{D})$
contains the point $\bar{P}_{4}$.
\end{remark}

Therefore, we may assume that $\mathbb{CS}(U_{1},
\frac{1}{k}\mathcal{D})$ does not contain the point $\bar{P}_{4}$.
Hence, it follows from the proof of
Lemma~\ref{lemma:n-9-one-singular-points-of-index} that
$\mathbb{CS}(U_{1}, \frac{1}{k}\mathcal{D})$ does not contains
subvarieties of $U_{1}$ that are contained in the exceptional
divisor of $\alpha_{1}$. Thus, the set $\mathbb{CS}(U_{1},
\frac{1}{k}\mathcal{D})$ contains $\bar{P}_{3}$.

Let $\gamma\colon W\to U_{1}$ be the weighted blow up of the point
$\bar{P}^{3}$ with weights $(1,1,3)$, $\mathcal{B}$ is the proper
transform of $\mathcal{M}$ on the variety $W$, $S$ be a general
surface of in $|-K_{W}|$, and $C$ be the base curve of the pencil
$|-K_{W}|$. Then $\mathcal{B}\sim_{\mathbb{Q}}-kK_{W}$  by
Theorem~\ref{theorem:Kawamata}, the curve $C$ is irreducible, the
inequality $C^{2}=-1/24$ holds on the normal surface $S$, and the
equivalence $\mathcal{B}\vert_{S}\sim_{\mathbb{Q}} kC$ holds,
which is impossible by Lemmas~\ref{lemma:normal-surface} and
\ref{lemma:Cheltsov-II}. The claim of
Proposition~\ref{proposition:n-17} is proved.

\section{Case $n=19$, hypersurface of degree $12$ in $\mathbb{P}(1,2,3,3,4)$.}
\label{section:n-19}

We use the notations and assumptions of
Section~\ref{section:start}. Let $n=19$. Then $X$ is a  general
hypersurface in $\mathbb{P}(1,2,3,3,4)$ of degree $12$,  the
singularities of the hypersurface $X$ consist of the points
$O_{1}$, $O_{2}$, $O_{3}$ that are quotient singularities of type
$\frac{1}{2}(1,1,1)$, the points $P_{1}$, $P_{2}$, $P_{3}$,
$P_{4}$ that are quotient singularities of type
$\frac{1}{3}(1,1,2)$, and the equality $-K_{X}^{3}=1/6$ holds.

It follows from Example~\ref{example:elliptic-7-11-19} that for
every $i\in\{1,2,3,4\}$ there is a commutative diagram
$$
\xymatrix{
&&U_{i}\ar@{->}[dl]_{\pi_{i}}\ar@{->}[dr]^{\eta_{i}}&&&\\%
&X\ar@{-->}[rr]_{\xi_{i}}&&\mathbb{P}(1,2,3),&}
$$
where $\xi_{i}$ is a projection, $\pi_{i}$ is the blow up of
$P_{i}$ with weights $(1,1,2)$, and $\eta_{i}$ is
an~elliptic~fibration.

\begin{proposition}
\label{proposition:n-19} There is a commutative diagram
\begin{equation}
\label{equation:n-19-commutative-diagram} \xymatrix{
&X\ar@{-->}[d]_{\xi_{i}}\ar@{-->}[rr]^{\rho}&&V\ar@{->}[d]^{\nu}&\\%
&\mathbb{P}(1,2,3)\ar@{-->}[rr]_{\sigma}&&\mathbb{P}^{2}&}
\end{equation}
for some $i\in\{1,2,3,4\}$, where $\sigma$ is a birational map.
\end{proposition}

\begin{proof}
The set $\mathbb{CS}(X, \frac{1}{k}\mathcal{M})$ contains the
point $P_{i}$ for some $i\in\{1,2,3,4\}$. Let $\mathcal{D}$ be the
proper transform of the linear system $\mathcal{M}$ on the variety
$U_{i}$. Then $\mathcal{D}\sim_{\mathbb{Q}}-kK_{U_{i}}$ by
Theorem~\ref{theorem:Kawamata}, which implies the existence of the
commutative diagram~\ref{equation:n-19-commutative-diagram}.
\end{proof}

\section{Case   $n=20$, hypersurface of degree $13$ in $\mathbb{P}(1,1,3,4,5)$.}
\label{section:n-20}

We use the notations and assumptions of
Section~\ref{section:start}. Let $n=20$. Then
$X$~is~a~general~hyper\-surface in $\mathbb{P}(1,1,3,4,5)$ of
degree $13$, the equality $-K_{X}^{3}=13/60$ holds, and the
singularities of the hypersurface  $X$ consist of the point
$P_{1}$ that is a singularity of type $\frac{1}{3}(1,1,2)$,  the
point $P_{2}$ that is a singularity of type $\frac{1}{4}(1,1,3)$,
and the point $P_{3}$ that is a singularity of type
$\frac{1}{5}(1,1,4)$. There is a commutative diagram
$$
\xymatrix{
&&&Y\ar@{->}[dl]_{\gamma_{4}}\ar@{->}[dr]^{\gamma_{2}}\ar@{->}[drrrrrr]^{\eta}&&&&&&\\
&&U_{23}\ar@{->}[dl]_{\beta_{3}}\ar@{->}[dr]^{\beta_{2}}&&U_{34}\ar@{->}[dl]^{\beta_{4}}&&&&&\mathbb{P}(1,1,3),\\
&U_{2}\ar@{->}[dr]_{\alpha_{2}}&&U_{3}\ar@{->}[dl]^{\alpha_{3}}&&&&&&\\
&&X\ar@{-->}[rrrrrrruu]_{\psi}&&&&&&&}
$$
where $\psi$ is a projection, $\alpha_{2}$ is the weighted blow up
of $P_{2}$ with weights $(1,1,3)$, $\alpha_{3}$ is the weighted
blow up of$P_{3}$ with weights $(1,1,4)$, $\beta_{3}$ is the
weighted blow up with weights $(1,1,4)$ of the proper transform of
the point $P_{3}$ on  $U_{2}$, $\beta_{2}$ is the weighted blow up
with weights $(1,1,3)$ of the proper transform of the point
$P_{2}$ on the variety $U_{3}$, $\beta_{4}$ is the weighted blow
up with weights $(1,1,3)$ of the singular point of $U_{3}$ that is
contained in the exceptional divisor of  $\alpha_{3}$,
$\gamma_{2}$ is the weighted blow up with weights $(1,1,3)$ of the
proper transform of $P_{2}$ on $U_{34}$, $\gamma_{4}$ is the
weighted blow up with weights $(1,1,3)$ of the singular point of
the variety $U_{23}$ that is contained in the exceptional divisor
of the morphism $\beta_{3}$, and $\eta$ is an elliptic fibration.

\begin{remark}
\label{remark:n-20-big-nef-I} The divisors $-K_{U_{2}}$,
$-K_{U_{3}}$, $-K_{U_{23}}$ and $-K_{U_{34}}$ are nef and big.
\end{remark}

There is a commutative diagram
$$
\xymatrix{
&&Y\ar@{->}[dl]_{\gamma_{3}}\ar@{->}[dr]^{\gamma_{1}}\ar@{->}[drrrrrrr]^{\omega}&&&&&&&\\
&U_{1}\ar@{->}[dr]_{\alpha_{1}}&&U_{3}\ar@{->}[dl]^{\alpha_{3}}&&&&&&\mathbb{P}(1,1,4),\\
&&X\ar@{-->}[rrrrrrru]_{\xi}&&&&&&&}
$$
where $\xi$ is a projection, $\alpha_{1}$ is the blow up of
$P_{1}$ with weights $(1,1,2)$, $\alpha_{3}$ is the weighted blow
up of the point $P_{3}$ с weights $(1,1,4)$, $\gamma_{3}$ is the
weighted blow up with weights $(1,1,4)$ of the proper transform of
$P_{3}$ on the variety $U_{1}$, $\gamma_{1}$ is the weighted blow
up with weights $(1,1,2)$ of the proper transform of the point
$P_{1}$ on the variety $U_{3}$, and $\omega$ is an elliptic
fibration.

\begin{remark}
\label{remark:n-20-big-nef-II} The divisor $-K_{U_{1}}$ is nef and
big.
\end{remark}

In the rest of the section we prove the following result.

\begin{proposition}
\label{proposition:n-20} Either there is a commutative diagram
\begin{equation}
\label{equation:n-20-commutative-diagram-I} \xymatrix{
&X\ar@{-->}[d]_{\psi}\ar@{-->}[rr]^{\sigma}&&X\ar@{-->}[rr]^{\rho}&&V\ar@{->}[d]^{\nu}&\\%
&\mathbb{P}(1,1,3)\ar@{-->}[rrrr]_{\upsilon}&&&&\mathbb{P}^{2},&}
\end{equation}
or there is a commutative diagram
\begin{equation}
\label{equation:n-20-commutative-diagram-II} \xymatrix{
&X\ar@{-->}[d]_{\xi}\ar@{-->}[rr]^{\sigma}&&X\ar@{-->}[rr]^{\rho}&&V\ar@{->}[d]^{\nu}&\\%
&\mathbb{P}(1,1,4)\ar@{-->}[rrrr]_{\zeta}&&&&\mathbb{P}^{2},&}
\end{equation}
where $\sigma$, $\upsilon$ and $\zeta$ are birational maps.
\end{proposition}

It follows from Theorem~\ref{theorem:smooth-points} and
Lemma~\ref{lemma:Ryder} that $\mathbb{CS}(X,
\frac{1}{k}\mathcal{M})\subseteq\{P_{1},P_{2},P_{3}\}$.

\begin{lemma}
\label{lemma:n-20-points-P1-P3} Suppose that $\{P_{1},
P_{3}\}\subseteq\mathbb{CS}(X, \frac{1}{k}\mathcal{M})$. Then the
commutative diagram~\ref{equation:n-20-commutative-diagram-II}
exists.
\end{lemma}

\begin{proof}
Let $\mathcal{B}$ be the proper transform of the linear system
$\mathcal{M}$ on the variety $Y$, and $S$ be a general surface in
$\mathcal{B}$. Then $S\sim_{\mathbb{Q}}-kK_{Y}$ by
Theorem~\ref{theorem:Kawamata}. We see that $S\cdot C=0$, where
$C$ is a general fiber of $\omega$, which implies the existence of
the commutative
diagram~\ref{equation:n-20-commutative-diagram-II}.
\end{proof}

\begin{lemma}
\label{lemma:n-20-points-P1-P2} The set $\mathbb{CS}(X,
\frac{1}{k}\mathcal{M})$ does not contain the set $\{P_{1},
P_{2}\}$.
\end{lemma}

\begin{proof}
Suppose that $\{P_{1}, P_{3}\}\subseteq\mathbb{CS}(X,
\frac{1}{k}\mathcal{M})$. Let $\pi\colon W\to U_{1}$ be the
weighted blow up with weights $(1,1,3)$ of the proper transform of
the point $P_{2}$ on the variety $U_{1}$, and $\mathcal{B}$ be the
proper transform of the linear system $\mathcal{M}$ on the variety
$W$. Then $\mathcal{B}\sim_{\mathbb{Q}}-kK_{W}$ by
Theorem~\ref{theorem:Kawamata}.

The linear system $|-K_{W}|$ is a pencil, the base locus of the
pencil $|-K_{W}|$ is the irreducible curve $\Delta$ such that the
curve $\alpha_{1}\circ\pi(C)$ is cut on the hypersurface $X$ by
the equations $z=y=0$.

Let $S$ be a sufficiently general surface of the linear system
$|-K_{W}|$, $\bar{P}_{3}$ be the proper transform of the singular
point $P_{3}$ on the variety $W$, and $\bar{P}_{5}$ and
$\bar{P}_{6}$ be other singular points of $W$ such that
$\alpha_{1}\circ\pi(\bar{P}_{5})=P_{1}$ and
$\alpha_{1}\circ\pi(\bar{P}_{6})=P_{2}$. Then $\bar{P}_{5}$ is a
quotient singularity of type $\frac{1}{2}(1,1,1)$ on the variety
$W$, the point $\bar{P}_{6}$ is a quotient singularity of type
$\frac{1}{3}(1,1,2)$ on the variety $W$, the surface $S$ is smooth
outside of the points $\bar{P}_{3}$, $\bar{P}_{5}$ and
$\bar{P}_{6}$, the singularities of the surface $S$ in the points
$\bar{P}_{3}$, $\bar{P}_{5}$ and $\bar{P}_{6}$ are Du Val
singularities of types $\mathbb{A}_{4}$, $\mathbb{A}_{1}$ and
$\mathbb{A}_{2}$ respectively.

The equality $\Delta^{2}=-1/30$ holds on the surface $S$, but the
equivalence $\mathcal{B}\vert_{S}\sim_{\mathbb{Q}} k\Delta$ holds,
which implies that $\mathcal{B}\vert_{S}=k\Delta$. Now we can
easily get a contradiction using Lemma~\ref{lemma:Cheltsov-II}.
\end{proof}

\begin{lemma}
\label{lemma:n-20-single-point} The set $\mathbb{CS}(X,
\frac{1}{k}\mathcal{M})$ is not consists of the point $P_{i}$.
\end{lemma}

\begin{proof}
Suppose that $\mathbb{CS}(X, \frac{1}{k}\mathcal{M})=\{P_{i}\}$.
Let $\mathcal{D}$ be the proper transform of $\mathcal{M}$ on
$U_{i}$. Then the equivalence
$\mathcal{D}\sim_{\mathbb{Q}}-kK_{U_{i}}$ holds by
Theorem~\ref{theorem:Kawamata}. Moreover, it follows from
Lemma~\ref{lemma:Noether-Fano} and the proof of
Lemma~\ref{lemma:n-9-one-singular-points-of-index} that the set
$\mathbb{CS}(U_{i}, \frac{1}{k}\mathcal{D})$ contains the singular
point of $U_{i}$ that is the singular point of the exceptional
divisor of the birational morphism $\alpha_{i}$.

Let $\pi\colon W\to U_{i}$ be the weighted blow up with weights
$(1,1,i)$ of the singular point of the variety $U_{i}$ that is
contained in the exceptional divisor of the morphism $\alpha_{i}$,
and $\mathcal{B}$ be the proper transform of the linear system
$\mathcal{M}$ on the variety $W$. Then
$\mathcal{B}\sim_{\mathbb{Q}}-kK_{W}$  by
Theorem~\ref{theorem:Kawamata}.

Let $S$ be a sufficiently general surface of the pencil
$|-K_{W}|$, and $\Delta$ be the unique base curve of the pencil
$|-K_{W}|$. Then the surface $S$ is normal, but the curve $\Delta$
is irreducible, rational and smooth. Moreover, simple computations
imply that the equality
$$
\Delta^{2}=\left\{\aligned
&-9/20\ {\text {при}}\ i=1,\\
&-1/30\ {\text {при}}\ i=2,\\
&0\ {\text {при}}\ i=3,\\
\endaligned
\right.
$$
holds on the surface $S$. However, we have the equivalence
$\mathcal{B}\vert_{S}\sim_{\mathbb{Q}} k\Delta$, which implies
(see the proof of Lemma~\ref{lemma:n-20-points-P1-P2}) that the
curve $\alpha_{i}\circ\pi(\Delta)$ is contained in the set
$\mathbb{CS}(X, \frac{1}{k}\mathcal{M})$ if $i\ne 3$.

Therefore, the equality $i=3$ holds.

Let $G$ be the exceptional divisor of $\alpha_{3}$, and
$\bar{P}_{4}$ be the singular point of $U_{3}$ that is contained
in divisor $G$. Then $\bar{P}^{4}$ is a quotient singularity of
type $\frac{1}{4}(1,1,3)$ on the variety $U_{3}$, and it follows
from Lemmas~\ref{lemma:Cheltsov-Kawamata} and \ref{lemma:curves}
that the set $\mathbb{CS}(U_{3}, \frac{1}{k}\mathcal{D})$ consists
of the singular point $\bar{P}_{4}$.

The variety $W$ is the variety $U_{34}$, and $\pi$ is the morphism
$\beta_{4}$. Thus, the divisor $-K_{W}$ is nef and big. Therefore,
it follows from Lemmas~\ref{lemma:Noether-Fano} and
\ref{lemma:Cheltsov-Kawamata} that the set $\mathbb{CS}(W,
\frac{1}{k}\mathcal{B})$ contains the singular point of the
variety $W$ that is contained in the exceptional divisor of the
$\pi$.

Let $\mu\colon Z\to W$ be the weighted blow up weights $(1,1,2)$
of the singular point of the variety $W$ that is contained in the
exceptional divisor of the morphism $\pi$, and $\mathcal{P}$ be
the proper transform of the linear system $\mathcal{M}$ on the
variety $Z$. Then the equivalence
$\mathcal{P}\sim_{\mathbb{Q}}-kK_{Z}$ holds by
Theorem~\ref{theorem:Kawamata}.

Let $F$ be a sufficiently general surface of the pencil
$|-K_{Z}|$, and $\Gamma$ be the unique base curve of the pencil
$|-K_{Z}|$. Then the surface $F$ is irreducible and normal, but
the curve $\Gamma$ is irreducible, rational and smooth. The
equality $\Gamma^{2}=-1/24$ holds on $F$, but
$\mathcal{P}\vert_{S}\sim_{\mathbb{Q}} k\Gamma$, which easily
leads to a contradiction using Lemmas~\ref{lemma:Cheltsov-II} and
\ref{lemma:normal-surface}.
\end{proof}

To conclude the proof of Proposition~\ref{proposition:n-20} we may
assume that $\mathbb{CS}(X,
\frac{1}{k}\mathcal{M})=\{P_{2},P_{3}\}$.

Let $\mathcal{D}$ be the proper transform of $\mathcal{M}$ on
$U_{23}$. Then $\mathcal{D}\sim_{\mathbb{Q}}-kK_{U_{23}}$ by
Theorem~\ref{theorem:Kawamata}, and it follows from
Lemmas~\ref{lemma:Noether-Fano} and \ref{lemma:Cheltsov-Kawamata}
that the set $\mathbb{CS}(U_{23}, \frac{1}{k}\mathcal{D})$
contains either the singular point of the variety $U_{23}$ that is
contained in the exceptional divisor of the morphism $\beta_{3}$,
or the singular point of the variety $U_{23}$ that is contained in
the exceptional divisor of $\beta_{2}$.

\begin{lemma}
\label{lemma:n-20-exclusion-of-P6} The set $\mathbb{CS}(U_{23},
\frac{1}{k}\mathcal{D})$ does not contain the singular point of
the variety $U_{23}$ that is contained in the exceptional divisor
of the morphism $\beta_{2}$.
\end{lemma}

\begin{proof}
Let $E$ be the exceptional  divisor of the morphism $\beta_{2}$,
and $\bar{P}_{6}$ be the singular point of the surface $E$. Then
$\bar{P}_{6}$ is a quotient singularity of type
$\frac{1}{3}(1,1,2)$ on $U_{23}$.

Suppose that the set $\mathbb{CS}(U_{23}, \frac{1}{k}\mathcal{D})$
contains $\bar{P}_{6}$. Let $\pi\colon W\to U_{23}$ be the
weighted blow up of the point $\bar{P}_{6}$ with weights
$(1,1,2)$, $\mathcal{B}$ be the proper transform of $\mathcal{M}$
on $W$, $S$ be a general surface of the pencil $|-K_{W}|$, and
$\Delta$ be the base curve of the pencil $|-K_{W}|$. Then $S$ is
normal, the curve $\Delta$ is irreducible, and
$\mathcal{B}\vert_{S}\sim_{\mathbb{Q}} k\Delta$ by
Theorem~\ref{theorem:Kawamata}, but the equality
$\Delta^{2}=-1/24$ holds on the surface $S$, which contradicts
Lemmas~\ref{lemma:normal-surface} and \ref{lemma:Cheltsov-II}.
\end{proof}

Hence, the set $\mathbb{CS}(U_{23}, \frac{1}{k}\mathcal{D})$
contains the singular point of the variety $U_{23}$ that is
contained in the exceptional divisor of the morphism $\beta_{3}$.
Thus, the existence of the commutative
diagram~\ref{equation:n-20-commutative-diagram-II} is easily
implied by Theorem~\ref{theorem:Kawamata}. The claim of
Proposition~\ref{proposition:n-20} is proved.

\section{Case $n=23$, hypersurface of degree $14$ in $\mathbb{P}(1,2,3,4,5)$.}
\label{section:n-23}

We use the notations and assumptions of
Section~\ref{section:start}. Let $n=23$. Then $X$ is a
sufficiently general hypersurface in $\mathbb{P}(1,2,3,4,5)$ of
degree $14$, the singularities of the hypersurface $X$ consist of
the points $P_{1}$, $P_{2}$, and $P_{3}$ that are quotient
singularities of type $\frac{1}{2}(1,1,1)$, the point $P_{4}$ that
is a quotient singularity of type $\frac{1}{3}(1,1,2)$, the point
$P_{5}$ that is a quotient singularity of type
$\frac{1}{4}(1,1,3)$, the point $P_{6}$ that is quotient
singularity of type $\frac{1}{5}(1,2,3)$, and the equality
$-K_{X}^{3}=7/60$ holds.

There is a commutative diagram
$$
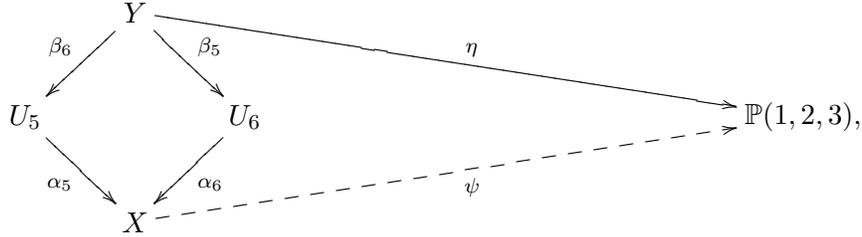
\xymatrix{
&&Y\ar@{->}[dl]_{\beta_{6}}\ar@{->}[dr]^{\beta_{5}}\ar@{->}[rrrrrrrd]^{\eta}&&&&&\\
&U_{5}\ar@{->}[dr]_{\alpha_{5}}&&U_{6}\ar@{->}[dl]^{\alpha_{6}}&&&&&&\mathbb{P}(1,2,3),\\
&&X\ar@{-->}[rrrrrrru]_{\psi}&&&&&&&}
$$
where $\psi$ is the natural projection, $\alpha_{5}$ is the
weighted blow up of the singular point $P_{5}$ with weights
$(1,1,3)$, $\alpha_{6}$ is the weighted blow up of $P_{6}$ with
weights $(1,2,3)$, $\beta_{5}$ is the weighted blow up with
weights $(1,1,3)$ of the proper transform of $P_{5}$ on the
variety $U_{6}$, $\beta_{6}$ is the weighted blow up with weights
$(1,2,3)$ of the proper transform of $P_{6}$ on $U_{5}$, and
$\eta$ is an elliptic fibration.

In the rest of the section we prove the following result.

\begin{proposition}
\label{proposition:n-23} The claim of Theorem~\ref{theorem:main}
holds for $n=23$.
\end{proposition}

It follows from Theorem~\ref{theorem:smooth-points},
Lemma~\ref{lemma:Ryder} and
Proposition~\ref{proposition:singular-points} that $\mathbb{CS}(X,
\frac{1}{k}\mathcal{M})\subseteq\{P_{5},P_{6}\}$.

\begin{lemma}
\label{lemma:n-23-P6-is-a-center} The set $\mathbb{CS}(X,
\frac{1}{k}\mathcal{M})$ contains the point $P_{6}$.
\end{lemma}

\begin{proof}
Suppose that $\mathbb{CS}(X, \frac{1}{k}\mathcal{M})$ does not
contain  $P_{6}$. Then the set $\mathbb{CS}(X,
\frac{1}{k}\mathcal{M})$ consists of the point $P_{5}$. Let
$\mathcal{D}_{5}$ be the proper transform of $\mathcal{M}$ on
$U_{5}$. Then $\mathcal{D}_{5}\sim_{\mathbb{Q}}-kK_{U_{5}}$ by
Theorem~\ref{theorem:Kawamata}, but the set $\mathbb{CS}(U_{5},
\frac{1}{k}\mathcal{D}_{5})$ is not empty by
Lemma~\ref{lemma:Noether-Fano}.

Let $\bar{P}_{7}$ be the singular point of the variety $U_{5}$
that is contained in the exceptional divisor of the morphism
$\alpha_{5}$. Then the point $\bar{P}_{7}$ is a quotient
singularity of type $\frac{1}{3}(1,1,2)$ on $U_{5}$, and it
follows from Lemma~\ref{lemma:Cheltsov-Kawamata} that
$\mathbb{CS}(U_{5}, \frac{1}{k}\mathcal{D}_{5})$ contains
$\bar{P}_{7}$. Let $\pi\colon U\to U_{5}$ be the weighted blow up
of the point $\bar{P}_{7}$ with weights $(1,1,3)$. Then the linear
system $|-2K_{U}|$ is a proper transform of the pencil
$|-2K_{X}|$, and the base locus of the pencil $|-2K_{U}|$ consists
of a single irreducible curve $Z$ such that
$\alpha_{5}\circ\pi(Z)$ is the unique base curve of the pencil
$|-2K_{X}|$.

Let $S$ be a sufficiently general surface of the pencil
$|-2K_{U}|$. Then the surface $S$ is normal, the surface $S$
contains the curve $Z$, and the inequality $Z^{2}<0$ holds on the
surface $S$, because the inequality $-K_{U}^{3}<0$ holds. However,
the equivalence $\mathcal{B}\vert_{S}\sim_{\mathbb{Q}} kZ$ holds
by Theorem~\ref{theorem:Kawamata}, where $\mathcal{B}$ is the
proper transform of $\mathcal{M}$ on $W$. It follows from
Lemma~\ref{lemma:normal-surface} that
$$
\mathrm{Supp}(S)\cap\mathrm{Supp}(D)=\mathrm{Supp}(Z),
$$
where $D$ is a sufficiently general surface of the linear system
$\mathcal{B}$, which contradicts Lemma~\ref{lemma:Cheltsov-II},
because the linear system $\mathcal{B}$ is not composed from a
pencil.
\end{proof}

It easily follows from Theorem~\ref{theorem:Kawamata} that the
claim of Theorem~\ref{theorem:main} holds for $X$ whenever the set
$\mathbb{CS}(X, \frac{1}{k}\mathcal{M})$ contains the points
$P_{5}$ and $P_{6}$. So, we may assume that
$P_{5}\not\in\mathbb{CS}(X, \frac{1}{k}\mathcal{M})$.

Let $\mathcal{D}_{6}$ be the proper transform of $\mathcal{M}$ on
$U_{6}$. Then $\mathcal{D}_{6}\sim_{\mathbb{Q}}-kK_{U_{6}}$ by
Theorem~\ref{theorem:Kawamata}, which implies that the set
$\mathbb{CS}(U_{6}, \frac{1}{k}\mathcal{D}_{6})$ is not empty by
Lemma~\ref{lemma:Noether-Fano}. Let $\bar{P}_{7}$ and
$\bar{P}_{8}$ be the singular points of the variety $U_{6}$ that
are quotient singularities of types $\frac{1}{3}(1,1,2)$ and
$\frac{1}{2}(1,1,1)$ contained in the exceptional divisor of the
morphism $\alpha_{6}$ respectively. Then the set
$\mathbb{CS}(U_{6}, \frac{1}{k}\mathcal{D}_{6})$ contains either
the point $\bar{P}_{7}$, or the point $\bar{P}_{8}$ by
Lemma~\ref{lemma:Cheltsov-Kawamata}.

\begin{lemma}
\label{lemma:n-23-P7-is-not-a-center} The set $\mathbb{CS}(U_{6},
\frac{1}{k}\mathcal{D}_{6})$ does not contain the point
$\bar{P}_{7}$.
\end{lemma}

\begin{proof}
Suppose that $\bar{P}_{7}\in\mathbb{CS}(U_{6},
\frac{1}{k}\mathcal{D}_{6})$. Let $\gamma\colon W\to U_{6}$ be the
weighted blow up of $\bar{P}_{7}$ with weights $(1,1,2)$, and $S$
be a general surface of in $|-2K_{W}|$. Then the surface $S$ is
irreducible and normal, the linear system $|-2K_{W}|$ is the
proper transform of the pencil $|-2K_{X}|$, and the base locus of
the pencil $|-2K_{W}|$ consists of the irreducible curve $\Delta$
such that the equality
$$
\Delta^{2}=-2K_{W}^{3}=-\frac{1}{6}
$$
holds on the surface $S$. Moreover, the equivalence
$\mathcal{B}\vert_{S}\sim_{\mathbb{Q}} k\Delta$ holds, where
$\mathcal{B}$ is the proper transform of $\mathcal{M}$ on $W$,
which contradicts Lemmas~\ref{lemma:normal-surface} and
\ref{lemma:Cheltsov-II}.
\end{proof}

Therefore, the set $\mathbb{CS}(U_{6},
\frac{1}{k}\mathcal{D}_{6})$ contains the point $\bar{P}_{8}$.

\begin{remark}
\label{remark:n-23-pencil-3K} The linear system $|-3K_{U_{6}}|$ is
the proper transform of the linear system $|-3K_{X}|$, the base
locus of the linear system $|-3K_{U_{6}}|$ consists of the
irreducible fiber of $\psi\circ\alpha_{6}$ that passes through the
singular point $\bar{P}_{8}$.
\end{remark}

Let $\pi\colon U\to U_{6}$ be the weighted blow up of the point
$\bar{P}_{8}$ with weights $(1,1,1)$, $F$ be the exceptional
divisor of the morphism $\pi$, $\mathcal{D}$ be the proper
transform of the linear system $\mathcal{M}$ on the variety $U$,
and $\mathcal{H}$ be the proper transform of the linear system
$|-3K_{U_{6}}|$ on $U$. Then
$$
\mathcal{D}\sim_{\mathbb{Q}}-kK_{U}\sim_{\mathbb{Q}}\pi^{*}\big(-kK_{U_{6}}\big)-\frac{k}{2}F
$$
by Theorem~\ref{theorem:Kawamata}. The simple computations imply
that
$$
\mathcal{H}\sim_{\mathbb{Q}}\pi^{*}\big(-3K_{U_{6}}\big)-\frac{1}{2}F,
$$
and the base locus of $\mathcal{H}$ consists of the irreducible
curve $Z$ such that $\alpha_{6}\circ\pi(Z)$ is the base curve of
the linear system $|-3K_{X}|$. Moreover, the equality $S\cdot
Z=1/12$ holds, where $S$ is a general surface of the linear system
$\mathcal{H}$.

Let $D_{1}$ and $D_{2}$ be general surfaces of the linear system
$\mathcal{D}$. Then
$$
-k^{2}/4=\Big(\pi^{*}\big(-3K_{U_{6}}\big)-\frac{1}{2}F\Big)\cdot
\Big(\pi^{*}\big(-kK_{U_{6}}\big)-\frac{k}{2}F\Big)^{2}
=\Big(\pi^{*}\big(-3K_{U_{6}}\big)-\frac{1}{2}F\Big)\cdot
D_{1}\cdot D_{2}\geqslant 0,
$$
which is a contradiction. Hence, the claim of
Proposition~\ref{proposition:n-23} is proved.

\section{Case   $n=25$, hypersurface of degree $15$ in $\mathbb{P}(1,1,3,4,7)$.}
\label{section:n-25}

We use the notations and assumptions of
Section~\ref{section:start}. Let $n=25$. Then $X$ is a
sufficiently general hypersurface in $\mathbb{P}(1,1,3,4,7)$ of
degree $15$, the equality $-K_{X}^{3}=5/28$ holds, and the
singularities of $X$ consist of the point $P_{1}$ that is a
quotient singularity of type $\frac{1}{4}(1,1,3)$, and the point
$P_{2}$ that is a quotient singularity of type
$\frac{1}{7}(1,3,4)$.

\begin{proposition}
\label{proposition:n-25} The claim of Theorem~\ref{theorem:main}
holds for $n=25$.
\end{proposition}

There is a commutative diagram
$$
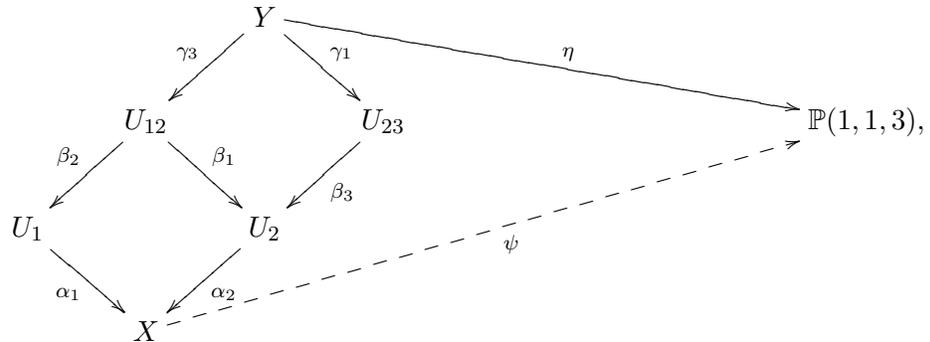
\xymatrix{
&&&Y\ar@{->}[dl]_{\gamma_{3}}\ar@{->}[dr]^{\gamma_{1}}\ar@{->}[drrrrrr]^{\eta}&&&&&&\\
&&U_{12}\ar@{->}[dl]_{\beta_{2}}\ar@{->}[dr]^{\beta_{1}}&&U_{23}\ar@{->}[dl]^{\beta_{3}}&&&&&\mathbb{P}(1,1,3),\\
&U_{1}\ar@{->}[dr]_{\alpha_{1}}&&U_{2}\ar@{->}[dl]^{\alpha_{2}}&&&&&&\\
&&X\ar@{-->}[rrrrrrruu]_{\psi}&&&&&&&}
$$
where $\psi$ is a projection, $\alpha_{1}$ is the weighted blow up
of $P_{1}$ with weights $(1,1,3)$, $\alpha_{2}$ is the weighted
blow up of $P_{2}$ with weights $(1,3,4)$, $\beta_{2}$ is the
weighted blow up with weights $(1,3,4)$ of the proper transform of
$P_{2}$ on $U_{1}$, $\beta_{1}$ is the weighted blow up with
weights $(1,1,3)$ of the proper transform of the point $P_{1}$ on
the variety $U_{2}$, $\beta_{3}$ is the weighted blow up with
weights $(1,1,3)$ of the singular point of the variety $U_{2}$
that is a quotient singularity of type $\frac{1}{4}(1,1,3)$
contained in the exceptional divisor of the morphism $\alpha_{2}$,
$\gamma_{1}$ is the weighted blow up with weights $(1,1,3)$ of the
proper transform of $P_{1}$ on  $U_{23}$, $\gamma_{3}$ is the
weighted blow up with weights $(1,1,3)$ of the point of the
variety $U_{12}$ that is a quotient singularity of type
$\frac{1}{4}(1,1,3)$ contained in the exceptional divisor of the
morphism $\beta_{2}$, and $\eta$ is an elliptic fibration.

\begin{remark}
\label{remark:n-25-big-nef} The divisors $-K_{U_{1}}$,
$-K_{U_{2}}$, $-K_{U_{12}}$ and $-K_{U_{23}}$ are nef and big.
\end{remark}

Let us prove Proposition~\ref{proposition:n-25}. It follows from
Theorem~\ref{theorem:smooth-points}, Lemma~\ref{lemma:Ryder} and
Proposition~\ref{proposition:singular-points} that $
\mathbb{CS}(X, \frac{1}{k}\mathcal{M})\subseteq\{P_{1},P_{2}\}$.
To conclude the proof of Proposition~\ref{proposition:n-25}, we
may assume that the singularities of the log pair $(X,
\frac{1}{k}\mathcal{M})$ are canonical (see
Remark~\ref{remark:canonical-singularities}).

\begin{lemma}
\label{lemma:n-25-P2-is-a-center} The set $\mathbb{CS}(X,
\frac{1}{k}\mathcal{M})$ contains the point $P_{2}$.
\end{lemma}

\begin{proof}
Suppose that the $\mathbb{CS}(X, \frac{1}{k}\mathcal{M})$ does not
contain the point $P_{2}$. Let $\mathcal{D}_{1}$ be the proper
transform of the $\mathcal{M}$ on the variety $U_{1}$. Then
$\mathbb{CS}(X, \frac{1}{k}\mathcal{M})=\{P_{1}\}$, and the set
$\mathbb{CS}(U_{1}, \frac{1}{k}\mathcal{D}_{1})$ is not empty by
Lemma~\ref{lemma:Noether-Fano}, because the equivalence
$\mathcal{D}_{1}\sim_{\mathbb{Q}}-kK_{U_{1}}$ holds by
Theorem~\ref{theorem:Kawamata}.

Let $P_{5}$ be the singular point of the variety $U_{1}$ that is
contained in the exceptional divisor of the morphism $\alpha_{1}$.
Then the point $P_{5}$ is a quotient singularity of type
$\frac{1}{3}(1,1,2)$ on $U_{1}$, and it follows from
Lemma~\ref{lemma:Cheltsov-Kawamata} that $\mathbb{CS}(U_{1},
\frac{1}{k}\mathcal{D}_{1})$ contains the point $P_{5}$.

Let $\pi\colon W\to U_{1}$ be the blow up of the point $P_{5}$
with weights $(1,1,2)$, and $S$ be a sufficiently general surface
of the pencil $|-K_{W}|$. Then the surface $S$ is irreducible and
normal, and the base locus of the pencil $|-K_{W}|$ consists of
the irreducible curve $\Delta$ such that
$$
\Delta^{2}=-K_{W}^{3}=-\frac{1}{14}
$$
on the surface $S$, but $\mathcal{B}\vert_{S}\sim_{\mathbb{Q}}
k\Delta$, where $\mathcal{B}$ is the proper transform of the
linear system $\mathcal{M}$ on the variety $W$. Therefore, we have
$\mathcal{B}\vert_{S}=k\Delta$, which implies that
$$
\mathrm{Supp}(S)\cap\mathrm{Supp}(D)=\mathrm{Supp}(\Delta),
$$
where $D$ is a general surface in $\mathcal{B}$. The latter
contradicts Lemma~\ref{lemma:Cheltsov-II}.
\end{proof}

Let $G$ be the exceptional divisor of $\alpha_{2}$,
$\mathcal{D}_{2}$ be the proper transform of the linear system
$\mathcal{M}$ on the variety $U_{2}$, $\bar{P}_{1}$ be the proper
transform of $P_{1}$ on $U_{2}$, and $\bar{P}_{3}$ and
$\bar{P}_{4}$ are the singular points of the variety $U_{2}$ that
are quotient singularities of types $\frac{1}{4}(1,1,3)$ and
$\frac{1}{3}(1,1,2)$ contained in the exceptional divisor $G$
respectively. Then $G\cong\mathbb{P}(1,3,4)$, the points
$\bar{P}_{3}$ and $\bar{P}_{4}$ are singular points of the surface
$G$, and $\mathcal{D}_{2}\sim_{\mathbb{Q}}-kK_{U_{2}}$ by
Theorem~\ref{theorem:Kawamata}. Hence, the set $\mathbb{CS}(U_{2},
\frac{1}{k}\mathcal{D}_{2})$ is not empty by
Lemma~\ref{lemma:Noether-Fano}. Moreover,  the proof of
Lemma~\ref{lemma:n-25-P2-is-a-center} implies that
$\mathbb{CS}(U_{2},
\frac{1}{k}\mathcal{D}_{2})\ne\{\bar{P}_{1}\}$.

\begin{lemma}
\label{lemma:n-25-P3-and-P4} The set $\mathbb{CS}(U_{2},
\frac{1}{k}\mathcal{D}_{2})$ does not contain both points
$\bar{P}_{3}$ and $\bar{P}_{4}$.
\end{lemma}

\begin{proof}
Suppose that $\{\bar{P}_{3},
\bar{P}_{4}\}\subseteq\mathbb{CS}(U_{2},
\frac{1}{k}\mathcal{D}_{2})$. Let $\pi\colon W\to U_{2}$ be a
composition of the weighted blow ups of the points $\bar{P}_{3}$
and $\bar{P}_{4}$ with weights $(1,1,3)$ and $(1,1,2)$
respectively, and $\mathcal{B}$ be the proper transform of
$\mathcal{M}$ on the variety $W$. Then
$\mathcal{B}\sim_{\mathbb{Q}}-kK_{W}$ by
Theorem~\ref{theorem:Kawamata}.

Let $S$ be a general surface of the pencil $|-K_{W}|$. Then the
surface $S$ is irreducible and normal, but the base locus of the
pencil $|-K_{W}|$ consists of the irreducible curves $C$ and $L$
such that the curve $\alpha_{2}\circ\pi(C)$ is the unique base
curve of the pencil $|-K_{X}|$, the curve $\pi(L)$ is contained in
the surface $G$, and $\pi(L)$ is the unique curve in
$|\mathcal{O}_{\mathbb{P}(1,\,3,\,4)}(1)|$. We have
$$
\mathcal{B}\vert_{S}\sim_{\mathbb{Q}}-kK_{W}\vert_{S}\sim_{\mathbb{Q}}kS\vert_{S}=kC+kL,
$$
but the intersection form of $L$ and $C$ on $S$ is negatively
defined, and Lemma~\ref{lemma:normal-surface} implies that
$$
\mathrm{Supp}(S)\cap\mathrm{Supp}(D)=\mathrm{Supp}(C)\cup\mathrm{Supp}(L),
$$
where $D$ is a general surface in $\mathcal{B}$, which is
impossible by Lemma~\ref{lemma:Cheltsov-II}.
\end{proof}

Thus, we have $\mathbb{CS}(U_{2},
\frac{1}{k}\mathcal{D}_{2})\subsetneq\{\bar{P}_{1}, \bar{P}_{3},
\bar{P}_{4}\}$ by Lemma~\ref{lemma:Cheltsov-Kawamata}. что .

\begin{lemma}
\label{lemma:n-25-P1-or-P4} The set $\mathbb{CS}(U_{2},
\frac{1}{k}\mathcal{D}_{2})$ contains either the point
$\bar{P}_{1}$, or the point $\bar{P}_{4}$.
\end{lemma}

\begin{proof}
Suppose that the set $\mathbb{CS}(U_{2},
\frac{1}{k}\mathcal{D}_{2})$ does not contain neither the singular
point $\bar{P}_{1}$, nor the singular point $\bar{P}_{4}$. Then
the set $\mathbb{CS}(U_{2}, \frac{1}{k}\mathcal{D}_{2})$ consists
of the point $\bar{P}_{3}$.

The linear system $|-K_{U_{2}}|$ is the proper transform of the
pencil $|-K_{X}|$, and the base locus of the pencil $|-K_{U_{2}}|$
consists of the irreducible curves $L$ and $\Delta$ such that
$\alpha_{2}(\Delta)$ is the unique base curve of the pencil
$|-K_{X}|$, the curve $L$ is contained in the  divisor $G$, the
curve $L$ is the unique curve of the linear system
$|\mathcal{O}_{\mathbb{P}(1,\,3,\,4)}(1)|$.

Let $\tilde{P}_{6}$ be the singular points of the variety $U_{23}$
that is contained in the exceptional divisor of $\beta_{3}$, and
$\mathcal{D}_{23}$ be the proper transform of $\mathcal{M}$ on
$U_{23}$. Then $\mathcal{D}_{23}\sim_{\mathbb{Q}}-kK_{U_{23}}$ by
Theorem~\ref{theorem:Kawamata}, and it follows from
Lemmas~\ref{lemma:Noether-Fano} and \ref{lemma:Cheltsov-Kawamata}
the the set $\mathbb{CS}(U_{23}, \frac{1}{k}\mathcal{D}_{23})$
contains the point $\tilde{P}_{6}$ that is a quotient singularity
of type $\frac{1}{3}(1,1,2)$ on the variety $U_{23}$.

Let $\pi\colon W\to U_{23}$ be the weighted blow up of
$\tilde{P}_{6}$ with weights $(1,1,2)$, $\mathcal{D}$ be the
proper transform of $\mathcal{M}$ on the variety $W$, and
$\bar{L}$ and $\bar{\Delta}$ be the proper transforms of $L$ and
$\Delta$ on the variety $W$ respectively. Then
$\mathcal{D}\sim_{\mathbb{Q}}-kK_{W}$ by
Theorem~\ref{theorem:Kawamata}, the pencil $|-K_{W}|$ is the
proper transform of the pencil $|-K_{U_{2}}|$, and the base locus
of $|-K_{W}|$ consists of $\bar{L}$ and $\bar{\Delta}$.

Let $S$ be a general surface of the pencil $|-K_{W}|$. Then the
surface $S$ is irreducible and normal, the equivalence
$\mathcal{D}\vert_{S}\sim_{\mathbb{Q}}k\bar{\Delta}+k\bar{L}$
holds, but the equalities
$$
\bar{\Delta}^{2}=-7/12,\ \bar{L}^{2}=-5/6,\ \bar{\Delta}\cdot\bar{L}=2/3%
$$
hold on the surface $S$. Therefore, the intersection form of the
curves $\bar{\Delta}$ and $\bar{L}$ on normal the surface $S$ is
negatively defined, which contradicts
Lemmas~\ref{lemma:normal-surface} and \ref{lemma:Cheltsov-II}.
\end{proof}

The hypersurface $X$ can be given by the equation
$$
w^{2}y+wt^{2}+wtf_{4}(x,y,z)+wf_{8}(x,y,z)+tf_{11}(x,y,z)+f_{15}(x,y,z)=0\subset\mathrm{Proj}\Big(\mathbb{C}[x,y,z,t,w]\Big),
$$
where $\mathrm{wt}(x)=1$, $\mathrm{wt}(y)=1$, $\mathrm{wt}(z)=3$,
$\mathrm{wt}(t)=4$, $\mathrm{wt}(w)=7$, and $f_{i}(x,y,t)$ is a
sufficiently general quasihomogeneous polynomial of degree $i$.

\begin{remark}
\label{remark:n-25-points-P1-P3} Suppose that the set
$\mathbb{CS}(U_{2}, \frac{1}{k}\mathcal{D}_{2})$ contains both
points $\bar{P}_{1}$ and $\bar{P}_{3}$. Then the claim of
Theorem~\ref{theorem:Kawamata} easily implies the existence of the
commutative diagram
$$
\xymatrix{
&X\ar@{-->}[d]_{\psi}\ar@{-->}[rr]^{\rho}&&V\ar@{->}[d]^{\nu}&\\%
&\mathbb{P}(1,1,3)\ar@{-->}[rr]_{\zeta}&&\mathbb{P}^{2},&}
$$
where $\zeta$ is a birational map.
\end{remark}

Therefore, we may assume that set $\mathbb{CS}(U_{2},
\frac{1}{k}\mathcal{D}_{2})$ does not contains both points
$\bar{P}_{1}$ and $\bar{P}_{3}$.

\begin{lemma}
\label{lemma:n-25-P4} The set $\mathbb{CS}(U_{2},
\frac{1}{k}\mathcal{D}_{2})$ contains the point $\bar{P}_{1}$.
\end{lemma}

\begin{proof}
Suppose that $\bar{P}_{1}\not\in\mathbb{CS}(U_{2},
\frac{1}{k}\mathcal{D}_{2})$. Then $\mathbb{CS}(U_{2},
\frac{1}{k}\mathcal{D}_{2})=\{\bar{P}_{4}\}$.

Let $\pi\colon W\to U_{2}$ be the weighted blow up of the point
$\bar{P}_{4}$ with weights $(1,1,2)$, $E$ be the exceptional
divisor of the morphism $\pi$, and $\bar{G}$ and $\mathcal{B}$ be
proper transforms of the divisor $G$ and the linear system
$\mathcal{M}$ on the variety $W$ respectively. Then it follows
from Theorem~\ref{theorem:Kawamata} that the equivalence
$\mathcal{B}\sim_{\mathbb{Q}}-kK_{W}$ holds, but the proof of
Lemma~\ref{lemma:n-25-P1-or-P4} implies that the set
$\mathbb{CS}(W, \frac{1}{k}\mathcal{B})$ does not contain the
singular point of the variety $W$ that is contained in the
exceptional divisor of the morphism $\pi$. Therefore, the
singularities of the log pair $(W, \frac{1}{k}\mathcal{B})$ are
terminal by Lemma~\ref{lemma:Cheltsov-Kawamata}.

Let $S_{x}$, $S_{y}$, $S_{z}$, $S_{t}$ and $S_{w}$ be proper
transforms on the variety $W$ of the surfaces that are cut on the
variety $X$ by the equations $x=0$, $y=0$, $z=0$, $t=0$ and $w=0$
respectively. Then
\begin{equation}
\label{equation:n-25-rational-equivalences-I} \left\{\aligned
&S_{x}\sim_{\mathbb{Q}}(\alpha_{2}\circ\pi)^{*}(-K_{X})-\frac{3}{7}E-\frac{1}{7}\bar{G},\\
&S_{y}\sim_{\mathbb{Q}}(\alpha_{2}\circ\pi)^{*}(-K_{X})-\frac{10}{7}E-\frac{8}{7}\bar{G},\\
&S_{z}\sim_{\mathbb{Q}}(\alpha_{2}\circ\pi)^{*}(-3K_{X})-\frac{2}{7}E-\frac{3}{7}\bar{G},\\
&S_{t}\sim_{\mathbb{Q}}(\alpha_{2}\circ\pi)^{*}(-4K_{X})-\frac{5}{7}E-\frac{4}{7}\bar{G},\\
&S_{w}\sim_{\mathbb{Q}}(\alpha_{2}\circ\pi)^{*}(-7K_{X}).\\
\endaligned\right.
\end{equation}

The base locus of the pencil $|-K_{W}|$ consists of the
irreducible curves $C$ and $L$ such that the curve
$\alpha_{2}\circ\pi(C)$ is cut by the equations $x=y=0$ on the
hypersurface $X$, the curve $\pi(L)$ is contained in the surface
$G$, and the curve $\pi(L)$ is contained in the linear system
$|\mathcal{O}_{\mathbb{P}(1,\,3,\,4)}(1)|$.

The equivalences~\ref{equation:n-25-rational-equivalences-I}
implies that the rational functions $y/x$, $zy/x^{4}$, $ty/x^{5}$
and $wy^{3}/x^{10}$ are contained in the linear system $|aS_{x}|$,
where $a=1$, $4$, $5$ and $10$ respectively. Therefore, the linear
system $|-20K_{W}|$ induces the birational map $\chi\colon
W\dasharrow X^{\prime}$, where $X^{\prime}$ is a hypersurface with
canonical singularities in $\mathbb{P}(1,1,4,5,10)$ of degree
$20$. In particular, the divisor $-K_{W}$ is big.

It follows from \cite{Sho93} that there is a composition of
antiflips $\zeta\colon W\dasharrow Z$ such that the rational map
$\zeta$ is regular outside of  $C\cup L$, and the divisor divisor
$-K_{Z}$ is nef. Let $\mathcal{P}$ be the proper transform of the
linear system $\mathcal{M}$ on $Z$. Then the singularities of log
pair $(Z, \frac{1}{k}\mathcal{P})$ are terminal, because $\zeta$
is a log-flop with respect to the log pair $(W,
\frac{1}{k}\mathcal{B})$, which has terminal singularities, but it
follows from Lemma~\ref{lemma:Noether-Fano} that the singularities
of $(Z, \frac{1}{k}\mathcal{P})$ are not terminal singularities.
\end{proof}

Hence, the set $\mathbb{CS}(U_{2}, \frac{1}{k}\mathcal{D}_{2})$
consists of the points $\bar{P}_{1}$ and $\bar{P}_{4}$.

Let $\pi\colon W\to U_{2}$ be a composition of the weighted blow
ups of the points $\bar{P}_{1}$ and $\bar{P}_{4}$ with weights
$(1,1,3)$ and $(1,1,2)$ respectively, $\bar{G}$ and $\mathcal{B}$
be the proper transforms of $G$ and $\mathcal{M}$ on the variety
$W$ respectively, and $F$ and $E$ be exceptional divisors of the
morphism $\pi$ that dominates the points $\bar{P}_{1}$ and
$\bar{P}_{4}$ respectively. Then the equivalence
$\mathcal{B}\sim_{\mathbb{Q}}-kK_{W}$ holds by
Theorem~\ref{theorem:Kawamata}, but it follows from the proof of
Lemma~\ref{lemma:n-25-P1-or-P4} that the singularities of $(W,
\frac{1}{k}\mathcal{B})$ are terminal.

Let $S_{x}$, $S_{y}$, $S_{z}$, $S_{t}$ and $S_{w}$ be the proper
transforms on the variety $W$ of the surfaces that are cut on $X$
by the equations $x=0$, $y=0$, $z=0$, $t=0$ and $w=0$
respectively. Then
$$
\left\{\aligned
&S_{x}\sim_{\mathbb{Q}}(\alpha_{2}\circ\pi)^{*}(-K_{X})-\frac{3}{7}E-\frac{1}{7}\bar{G}-\frac{1}{4}F\sim_{\mathbb{Q}}-K_{W},\\
&S_{y}\sim_{\mathbb{Q}}(\alpha_{2}\circ\pi)^{*}(-K_{X})-\frac{10}{7}E-\frac{8}{7}\bar{G}-\frac{1}{4}F,\\
&S_{z}\sim_{\mathbb{Q}}(\alpha_{2}\circ\pi)^{*}(-3K_{X})-\frac{2}{7}E-\frac{3}{7}\bar{G}-\frac{3}{4}F,\\
&S_{t}\sim_{\mathbb{Q}}(\alpha_{2}\circ\pi)^{*}(-4K_{X})-\frac{5}{7}E-\frac{4}{7}\bar{G},\\
&S_{w}\sim_{\mathbb{Q}}(\alpha_{2}\circ\pi)^{*}(-7K_{X})-\frac{11}{4}F,\\
\endaligned\right.
$$
which imply that the rational functions $y/x$, $zy/x^{4}$,
$twy^{4}/x^{15}$ and $wy^{3}/x^{10}$ are contained in the linear
system $|aS_{x}|$, where $a=1$, $4$, $15$ and $10$ respectively.
The linear system $|-60K_{W}|$ induces the birational map
$\chi\colon W\dasharrow X^{\prime}$ such that the variety
$X^{\prime}$ is a hypersurface in $\mathbb{P}(1,1,4,10,15)$ of
degree $30$. In particular, the divisor $-K_{W}$ is big. Now we
can obtain a contradiction in the same was as in the proof of
Lemma~\ref{lemma:n-25-P4}. The claim of
Proposition~\ref{proposition:n-25} is proved.

\section{Case  $n=26$, hypersurface of degree $15$ in $\mathbb{P}(1,1,3,5,6)$.}
\label{section:n-26}

We use the notations and assumptions of
Section~\ref{section:start}. Let $n=26$. Then $X$ is a
sufficiently general hypersurface in $\mathbb{P}(1,1,3,5,6)$ of
degree $15$, the equality $-K_{X}^{3}=1/6$ holds, and the
singularities of the hypersurface  $X$ consist of the points
$P_{1}$ and $P_{2}$ that are quotient singularities of type
$\frac{1}{3}(1,1,2)$, and the point $P_{3}$ that is a quotient
singularity of type $\frac{1}{6}(1,1,5)$.

There is a commutative diagram
$$
\xymatrix{
&U\ar@{->}[d]_{\alpha}&&W\ar@{->}[ll]_{\beta}&&Y\ar@{->}[ll]_{\gamma}\ar@{->}[d]^{\eta}&\\
&X\ar@{-->}[rrrr]_{\psi}&&&&\mathbb{P}(1,1,3)&}
$$
where $\psi$ is a projection, $\alpha$ is the weighted blow up of
$P_{3}$ with weights $(1,1,5)$, $\beta$ is the weighted blow up
with weights $(1,1,4)$ of the singular point of the variety $U$
that is contained in the exceptional divisor of the morphism
$\alpha$, $\gamma$ is the weighted blow up with weights $(1,1,3)$
of the singular point of the variety $W$ that is contained in the
exceptional divisor of the morphism $\beta$, and $\eta$ is an
elliptic fibration.

There is a commutative diagram
$$
\xymatrix{
&&U_{i}\ar@{->}[dl]_{\sigma_{i}}\ar@{->}[dr]^{\omega_{i}}&&\\
&X\ar@{-->}[rr]_{\xi_{i}}&&\mathbb{P}(1,1,6),&}
$$
where $\xi_{i}$ is a projection, $\sigma_{i}$ is the weighted blow
up of the point $P_{i}$ with weights $(1,1,2)$, and $\omega_{i}$
is an elliptic fibration, which is induced by the linear system
$|-6K_{U_{i}}|$.

It follows from \cite{CPR} that the group $\mathrm{Bir}(X)$ is
generated by biregular automorphisms of the hypersurface $X$ and a
birational involution $\tau\in\mathrm{Bir}(X)$ such that
$\psi\circ\tau=\psi$ and $\xi_{1}\circ\tau=\xi_{2}$.

In the rest of the section we prove the following result.

\begin{proposition}
\label{proposition:n-26} Either there is a commutative diagram
\begin{equation}
\label{equation:n-26-commutative-diagram-I} \xymatrix{
&X\ar@{-->}[d]_{\psi}\ar@{-->}[rr]^{\rho}&&V\ar@{->}[d]^{\nu}&\\%
&\mathbb{P}(1,1,3)\ar@{-->}[rr]_{\phi}&&\mathbb{P}^{2},&}
\end{equation}
or there is a commutative diagram
\begin{equation}
\label{equation:n-26-commutative-diagram-II} \xymatrix{
&X\ar@{-->}[d]_{\xi_{i}}\ar@{-->}[rr]^{\rho}&&V\ar@{->}[d]^{\nu}&\\%
&\mathbb{P}(1,1,6)\ar@{-->}[rr]_{\sigma}&&\mathbb{P}^{2},&}
\end{equation}
where $\phi$ and $\sigma$ are birational maps, and $i=1$ or $i=2$.
\end{proposition}

It follows from Theorem~\ref{theorem:smooth-points} and
Lemma~\ref{lemma:Ryder} that $\mathbb{CS}(X,
\frac{1}{k}\mathcal{M})\subseteq\{P_{1},P_{2},P_{3}\}$.

\begin{lemma}
\label{lemma:n-26-points-P1-and-P2} Suppose that
$P_{i}\in\mathbb{CS}(X, \frac{1}{k}\mathcal{M})$. Then the
commutative diagram~\ref{equation:n-26-commutative-diagram-II}
exists.
\end{lemma}

\begin{proof}
Let $\mathcal{B}$ be the proper transform of $\mathcal{M}$ on
$U_{i}$. Then $\mathcal{B}\sim_{\mathbb{Q}}-kK_{U_{i}}$ by
Theorem~\ref{theorem:Kawamata}, which implies the existence of the
commutative diagram~\ref{equation:n-26-commutative-diagram-II}.
\end{proof}

Therefore, we may assume that the set $\mathbb{CS}(X,
\frac{1}{k}\mathcal{M})$ consists of the point $P_{3}$.

Let $\mathcal{D}$ be the proper transform of the linear system
$\mathcal{M}$ on the variety $U$, and $\bar{P}_{4}$ be the
singular point of $U$ that is contained in the exceptional divisor
of the morphism $\alpha$. Then
$\mathcal{D}\sim_{\mathbb{Q}}-kK_{U}$ by
Theorem~\ref{theorem:Kawamata}, and the set $\mathbb{CS}(U,
\frac{1}{k}\mathcal{D})$ consists of the point $\bar{P}_{4}$ by
Lemmas~\ref{lemma:Cheltsov-Kawamata} and \ref{lemma:curves}.

Let $\mathcal{H}$ be the proper transform of $\mathcal{M}$ on $W$.
Then $\mathcal{H}\sim_{\mathbb{Q}}-kK_{W}$ by
Theorem~\ref{theorem:Kawamata}, but it follows from
Lemmas~\ref{lemma:Noether-Fano} and \ref{lemma:Cheltsov-Kawamata}
that $\mathbb{CS}(W, \frac{1}{k}\mathcal{H})$ contains the
singular point of $W$ that is contained in the exceptional divisor
of $\beta$. The existence of the
diagram~\ref{equation:n-26-commutative-diagram-I} follows from
Theorem~\ref{theorem:Kawamata}.

\section{Case $n=27$, hypersurface of degree $15$ in $\mathbb{P}(1,2,3,5,5)$.}
\label{section:n-27}

We use the notations and assumptions of
Section~\ref{section:start}. Let $n=27$. Then
$X$~is~a~general~hyper\-surface in $\mathbb{P}(1,2,3,5,5)$ of
degree $15$, whose singularities consist of the point $O$ that is
a singularity of type $\frac{1}{2}(1,1,1)$, and the points
$P_{1}$, $P_{2}$ and $P_{3}$ that are singularities of type
$\frac{1}{5}(1,2,3)$.

There is a commutative diagram
$$
\xymatrix{
&&&U_{ij}\ar@{->}[dl]_{\beta_{ij}}&&Y\ar@{->}[ll]_{\gamma_{ij}}\ar@{->}[ddr]^{\eta}&&\\%
&&U_{i}\ar@{->}[dr]_{\alpha_{i}}&&&&&\\
&&&X\ar@{-->}[rrr]_{\psi}&&&\mathbb{P}(1,2,3),&}
$$
where $\psi$ is a projection, $\alpha_{i}$ is the blow up of
$P_{i}$ with weights $(1,2,3)$, $\beta_{ij}$ is the weighted blow
up with weights $(1,2,3)$ of the proper transform of $P_{j}$ on
the variety $U_{i}$, $\gamma_{ij}$ is the weighted blow up with
weights $(1,2,3)$ of the proper transform of $P_{k}$ on $U_{ij}$,
and $\eta$ is and elliptic fibration, where $i\ne j$ and
$k\not\in\{i,j\}$.

\begin{proposition}
\label{proposition:n-27} The claim of Theorem~\ref{theorem:main}
holds for $n=27$.
\end{proposition}

Let us prove Proposition~\ref{proposition:n-27}. It follows from
Lemma~\ref{lemma:Ryder} and
Proposition~\ref{proposition:singular-points} that
$$
\varnothing\ne\mathbb{CS}\Big(X, \frac{1}{k}\mathcal{M}\Big)\subseteq\{P_{1},P_{2},P_{3}\}.%
$$

\begin{remark}
\label{remark:n-27-P1-P2-P3} In the case $\mathbb{CS}(X,
\frac{1}{k}\mathcal{M})=\{P_{1},P_{2},P_{3}\}$, the claim of
Theorem~\ref{theorem:main} holds for the hypersurface $X$ by
Theorem~\ref{theorem:Kawamata}.
\end{remark}

We may assume that the set $\mathbb{CS}(X,
\frac{1}{k}\mathcal{M})$ contains $P_{1}$ and does not contain
$P_{3}$.

\begin{lemma}
\label{lemma:n-27-P2-is-a-center} The set $\mathbb{CS}(X,
\frac{1}{k}\mathcal{M})$ contains the point $P_{2}$.
\end{lemma}

\begin{proof}
Suppose that $\mathbb{CS}(X, \frac{1}{k}\mathcal{M})$ does not
contains the point $P_{2}$. Then $\mathbb{CS}(X,
\frac{1}{k}\mathcal{M})$ consists of the point $P_{1}$. Let
$\mathcal{D}_{1}$ be the proper transform of $\mathcal{M}$ on
$U_{1}$. Then Theorem~\ref{theorem:Kawamata} implies that the
equivalence $\mathcal{D}_{1}\sim_{\mathbb{Q}}-kK_{U_{1}}$ holds.
The set $\mathbb{CS}(U_{1}, \frac{1}{k}\mathcal{D}_{1})$ is not
empty by Lemma~\ref{lemma:Noether-Fano}.

Let $G$ be is the exceptional divisor of $\alpha_{1}$, and $O$ and
$Q$ are the singular points of $G$ that are quotient singularities
of types $\frac{1}{3}(1,1,2)$ and $\frac{1}{2}(1,1,1)$ on $U_{1}$
respectively. Then it follows from the claim of
Lemma~\ref{lemma:Cheltsov-Kawamata} that $\mathbb{CS}(U_{1},
\frac{1}{k}\mathcal{D}_{1})$ contains either the point $O$, or the
point $Q$.

The linear system $|-2K_{U_{1}}|$ is a pencil, and the base locus
of the pencil $|-2K_{U_{1}}|$ consists of the irreducible curve
$C$ such that the curve $C$ passes through the point $O$, and $C$
is contracted by the rational map $\psi\circ\alpha_{1}$ to a
singular point of the surface $\mathbb{P}(1,2,3)$.

Suppose that the set $\mathbb{CS}(U_{1},
\frac{1}{k}\mathcal{D}_{1})$ contains the point $O$. Let
$\pi\colon W\to U_{1}$ be the weighted blow up of $O$ with weights
$(1,1,2)$, $\mathcal{B}$ be the proper transform of $\mathcal{M}$
on $W$, and $\bar{C}$ be the proper transform of $C$ on $W$. Then
$\mathcal{B}\sim_{\mathbb{Q}}-kK_{W}$ by
Theorem~\ref{theorem:Kawamata}, the linear system $|-2K_{W}|$ is
the proper transform of the pencil $|-2K_{U_{1}}|$, and the base
locus of $|-2K_{W}|$ consists of $\bar{C}$.

Let $S$ be a general surface in $|-2K_{W}|$. Then
$\mathcal{B}\vert_{S}\sim_{\mathbb{Q}}k\bar{C}$, but on the
surface $S$, the strict inequality $\bar{C}^{2}<0$ holds. We have
$\mathrm{Supp}(S)\cap\mathrm{Supp}(D)=\mathrm{Supp}(\bar{C})$,
where $D$ is a  general surface of the linear system
$\mathcal{B}$, which is impossible by
Lemma~\ref{lemma:Cheltsov-II}.

Hence, the set $\mathbb{CS}(U_{1}, \frac{1}{k}\mathcal{D}_{1})$
contains the point $Q$.

Let $\zeta\colon U\to U_{1}$ be the weighted blow up of the point
$Q$ with weights $(1,1,1)$, $F$ be the exceptional divisor of
$\zeta$, $\mathcal{H}$ be the proper transform of $\mathcal{M}$ on
$U$, and $\mathcal{P}$ be the proper transform of the linear
system $|-3K_{U_{1}}|$ on $U$. Then
$\mathcal{H}\sim_{\mathbb{Q}}-kK_{U}$ by
Theorem~\ref{theorem:Kawamata}, but
$$
\mathcal{P}\sim_{\mathbb{Q}}\zeta^{*}\big(-3K_{U_{1}}\big)-\frac{1}{2}F,
$$
and the base locus of $\mathcal{P}$ consists of the irreducible
curve $Z$ such that the curve $\alpha_{1}\circ\zeta(Z)$ is the
unique base curve of the linear system $|-3K_{X}|$. Therefore, we
have
$$
\Big(\zeta^{*}\big(-3K_{U_{1}}\big)-\frac{1}{2}F\Big)\cdot
Z=\frac{1}{10},
$$
which implies that
$$
-\frac{3k^{2}}{10}=\Big(\zeta^{*}\big(-3K_{U_{1}}\big)-\frac{1}{2}F\Big)\cdot
\Big(\zeta^{*}\big(-kK_{U_{1}}\big)-\frac{k}{2}F\Big)^{2}=\Big(\zeta^{*}\big(-3K_{U_{1}}\big)-\frac{1}{2}F\Big)\cdot
H_{1}\cdot H_{2}\geqslant 0,
$$
where $H_{1}$ and $H_{2}$ are general surfaces of the linear
system $\mathcal{H}$.
\end{proof}

We have $\mathbb{CS}(X, \frac{1}{k}\mathcal{M})=\{P_{1}, P_{2}\}$.
Now we can apply the proof of
Lemma~\ref{lemma:n-27-P2-is-a-center} to the proper transform of
$\mathcal{M}$ on $U_{12}$ to get a contradiction. The claim of
Proposition~\ref{proposition:n-27} is proved.

\section{Case $n=29$, hypersurface of degree $16$ in $\mathbb{P}(1,1,2,5,8)$.}
\label{section:n-29}

We use the notations and assumptions of
Section~\ref{section:start}. Let $n=29$. Then $X$ is a
sufficiently general hypersurface in $\mathbb{P}(1,1,2,5,8)$ of
degree $16$, the equality $-K_{X}^{3}=1/5$ holds, and the
singularities of the hypersurface $X$ consist of the points
$O_{1}$ and $O_{2}$ that are quotient singularities of type
$\frac{1}{2}(1,1,1)$, and the point $P$ that is a quotient
singularity of type $\frac{1}{5}(1,2,3)$.

The hypersurface $X$ is birationally superrigid, and there is a
commutative diagram
$$
\xymatrix{
&U\ar@{->}[d]_{\alpha}&&W\ar@{->}[ll]_{\beta}\ar@{->}[d]^{\eta}&\\%
&X\ar@{-->}[rr]_{\psi}&&\mathbb{P}(1,1,2),&}
$$
where $\psi$ is the natural projection, $\alpha$ is the weighted
blow up of $P$ with weights $(1,2,3)$, $\beta$ is the weighted
blow up with weights $(1,1,2)$ of the singular point of the
variety $U$ that is a quotient singularity of type
$\frac{1}{3}(1,1,2)$, and $\eta$ is an elliptic fibration.

\begin{proposition}
\label{proposition:n-29} The claim of Theorem~\ref{theorem:main}
holds for $n=29$.
\end{proposition}

\begin{proof}
Let $\mathcal{D}$ be the proper transform of $\mathcal{M}$ on the
variety $U$. Then $\mathcal{D}\sim_{\mathbb{Q}}-kK_{U}$ by
Theorem~\ref{theorem:Kawamata}, because $\mathbb{CS}(X,
\frac{1}{k}\mathcal{M})=\{P\}$ by
Theorem~\ref{theorem:smooth-points}, Lemma~\ref{lemma:Ryder} and
Proposition~\ref{proposition:singular-points}.

Let $G$ is the $\alpha$-exceptional divisor, and $Q$ and $O$ be
the singular points of $G$ that are singularities of types
$\frac{1}{3}(1,1,2)$ and $\frac{1}{2}(1,1,1)$ respectively. Then
it follows from Lemmas~\ref{lemma:Noether-Fano} and
\ref{lemma:Cheltsov-Kawamata} that either the set $\mathbb{CS}(U,
\frac{1}{k}\mathcal{D})$ consists of the point $O$, or the set
$\mathbb{CS}(U, \frac{1}{k}\mathcal{D})$ contains  $Q$.

Suppose that the set $\mathbb{CS}(U, \frac{1}{k}\mathcal{D})$
contains the point $O$. Let $\pi\colon Y\to U$ be the weighted
blow up of $O$ with weights $(1,1,1)$, $\mathcal{H}$ be the proper
transform of $\mathcal{M}$ on $Y$, $L$ be the curve on $G$ that is
contained in $|\mathcal{O}_{\mathbb{P}(1,\,2,\,3)}(1)|$, $\bar{L}$
be the proper transform of $L$ on  $Y$, and $S$ be a general
surface of the linear system $|-K_{Y}|$. Then
$\mathcal{H}\sim_{\mathbb{Q}}-kK_{Y}$ by
Theorem~\ref{theorem:Kawamata}, and the base locus of the pencil
$|-K_{Y}|$ consists of the curve $\bar{L}$ and the irreducible
curve $\Delta$ such that $\alpha\circ\pi(\Delta)$ is the base
locus of the pencil $|-K_{X}|$. Moreover, the equalities
$$
\Delta^{2}=-1,\ \bar{L}^{2}=-4/3,\ \Delta\cdot\bar{L}=1
$$
holds on the surface $S$. The intersection form of the curves
$\Delta$ and $\bar{L}$ on the surface $S$ is negatively defined.
We have $\mathcal{H}\vert_{S}\sim_{\mathbb{Q}}k\Delta+k\bar{L}$,
which contradicts Lemmas~\ref{lemma:Cheltsov-II} and
\ref{lemma:normal-surface}.

Therefore, the set $\mathbb{CS}(U, \frac{1}{k}\mathcal{D})$
contains the point $Q$. Let $\mathcal{B}$ be the proper transform
of the linear system $\mathcal{M}$ on the variety $W$. Then
$\mathcal{B}\sim_{\mathbb{Q}}-kK_{W}$ by
Theorem~\ref{theorem:Kawamata}, which implies that the linear
system $\mathcal{B}$ is contained in the fibers of the morphism
$\eta$.
\end{proof}

\section{Case $n=30$, hypersurface of degree $16$ in $\mathbb{P}(1,1,3,4,8)$.}
\label{section:n-30}

We use the notations and assumptions of
Section~\ref{section:start}. Let $n=30$. Then $X$ is a
sufficiently general hypersurface in $\mathbb{P}(1,1,3,4,8)$ of
degree $16$, the equality $-K_{X}^{3}=1/6$ holds, and the
singularities of $X$ consist of the point $O$ that is a quotient
singularity of type $\frac{1}{3}(1,1,2)$, and the points $P_{1}$
and $P_{2}$ that are singularities of type $\frac{1}{4}(1,1,3)$.
There is a commutative diagram
$$
\xymatrix{
&&Y\ar@{->}[dl]_{\gamma}\ar@{->}[dr]^{\delta}\ar@{->}[drrrrrrr]^{\eta}&&&&&&&\\
&U\ar@{->}[dr]_{\alpha}&&W\ar@{->}[dl]^{\beta}&&&&&&\mathbb{P}(1,1,3),\\
&&X\ar@{-->}[rrrrrrru]_{\psi}&&&&&&&}
$$
where $\psi$ is the natural projection, $\alpha$ is the weighted
blow up of $P_{1}$ with weights $(1,1,3)$, $\beta$ is the weighted
blow up of $P_{2}$ with weights $(1,1,3)$, $\gamma$ is the
weighted blow up with weights $(1,1,3)$ of the proper transform of
$P_{2}$ on $U$, $\delta$ is the weighted blow up with weights
$(1,1,3)$ of the proper transform of the point $P_{1}$ on the
variety $W$, and $\eta$ is and elliptic fibration.

There is a commutative diagram
$$
\xymatrix{
&&Z\ar@{->}[dl]_{\zeta}\ar@{->}[dr]^{\omega}&&\\
&X\ar@{-->}[rr]_{\xi}&&\mathbb{P}(1,1,4),&}
$$
where $\xi$ is a projection, $\zeta$ is the blow up of $O$ with
weights $(1,1,2)$, and $\omega$ is an elliptic~fibration.

\begin{proposition}
\label{proposition:n-30} Either there is a commutative diagram
\begin{equation}
\label{equation:n-30-commutative-diagram-I} \xymatrix{
&X\ar@{-->}[d]_{\psi}\ar@{-->}[rr]^{\rho}&&V\ar@{->}[d]^{\nu}&\\%
&\mathbb{P}(1,1,3)\ar@{-->}[rr]_{\phi}&&\mathbb{P}^{2},&}
\end{equation}
or there is a commutative diagram
\begin{equation}
\label{equation:n-30-commutative-diagram-II} \xymatrix{
&X\ar@{-->}[d]_{\xi}\ar@{-->}[rr]^{\theta}&&X\ar@{-->}[rr]^{\rho}&&V\ar@{->}[d]^{\nu}&\\%
&\mathbb{P}(1,1,4)\ar@{-->}[rrrr]_{\sigma}&&&&\mathbb{P}^{2},&}
\end{equation}
where $\phi$, $\theta$ and $\sigma$ are birational maps.
\end{proposition}

\begin{proof}
Suppose that $O\in\mathbb{CS}(X, \frac{1}{k}\mathcal{M})$. Then
the existence of the commutative
diagram~\ref{equation:n-30-commutative-diagram-II} follows from
Theorem~\ref{theorem:Kawamata}. Similarly, the existence of the
commutative diagram~\ref{equation:n-30-commutative-diagram-I}
follows from $\mathbb{CS}(X,
\frac{1}{k}\mathcal{M})=\{P_{1},P_{2}\}$ by
Theorem~\ref{theorem:Kawamata}. We may assume that $\mathbb{CS}(X,
\frac{1}{k}\mathcal{M})$ consists of the singular point $P_{1}$ by
Theorem~\ref{theorem:smooth-points},
Proposition~\ref{proposition:singular-points} and
Lemma~\ref{lemma:Ryder}.

Let $Q$ be the singular point of $U$ such that $\alpha(Q)=P_{1}$,
and $\mathcal{B}$ be the proper transform of the linear system
$\mathcal{M}$ on  $U$. Then $Q\in\mathbb{CS}(U,
\frac{1}{k}\mathcal{B})$ by Theorem~\ref{theorem:Kawamata} and
Lemmas~\ref{lemma:Noether-Fano} and \ref{lemma:Cheltsov-Kawamata}.

Let $\upsilon\colon \bar{U}\to U$ be the weighted blow up of the
point $Q$ with weights $(1,1,2)$, $\mathcal{D}$ be the proper
transform of the linear system $\mathcal{M}$ on the variety
$\bar{U}$, and $S$ be a sufficiently general surface of the pencil
$|-K_{\bar{U}}|$. Then $\mathcal{D}\sim_{\mathbb{Q}}-kK_{\bar{U}}$
by Theorem~\ref{theorem:Kawamata}, the surface $S$ is normal, and
the base locus of the pencil $|-K_{\bar{U}}|$ consists of the
irreducible curve $\Delta$ such that $\alpha\circ\upsilon(\Delta)$
is the unique base curve of the pencil $|-K_{X}|$. Moreover, the
inequality $\Delta^{2}<0$ holds on the surface $S$, but the
equivalence $\mathcal{D}\vert_{S}\sim_{\mathbb{Q}}k\Delta$ holds,
which contradicts Lemmas~\ref{lemma:normal-surface} and
\ref{lemma:Cheltsov-II}.
\end{proof}

\section{Case $n=31$, hypersurface of degree $16$ in $\mathbb{P}(1,1,4,5,6)$.}
\label{section:n-31}

We use the notations and assumptions of
Section~\ref{section:start}. Let $n=31$. Then $X$ is a general
hyper\-sur\-face in $\mathbb{P}(1,1,2,3,4)$ of degree $10$ and
$-K_{X}^{3}=2/15$. The singularities of $X$ consist of the points
$P_{1}$, $P_{2}$ and $P_{3}$ that are singularities of type
$\frac{1}{2}(1,1,1)$, $\frac{1}{5}(1,1,4)$ and
$\frac{1}{6}(1,1,5)$ respectively.

There is a commutative diagram
$$
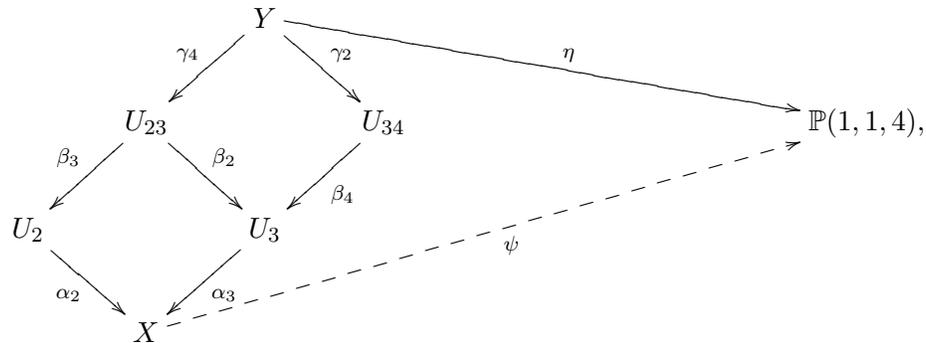
\xymatrix{
&&&Y\ar@{->}[dl]_{\gamma_{4}}\ar@{->}[dr]^{\gamma_{2}}\ar@{->}[drrrrrr]^{\eta}&&&&&&\\
&&U_{23}\ar@{->}[dl]_{\beta_{3}}\ar@{->}[dr]^{\beta_{2}}&&U_{34}\ar@{->}[dl]^{\beta_{4}}&&&&&\mathbb{P}(1,1,4),\\
&U_{2}\ar@{->}[dr]_{\alpha_{2}}&&U_{3}\ar@{->}[dl]^{\alpha_{3}}&&&&&&\\
&&X\ar@{-->}[rrrrrrruu]_{\psi}&&&&&&&}
$$
where $\psi$ is a projection, $\alpha_{2}$ is the weighted blow up
of $P_{2}$ with weights $(1,1,4)$, $\alpha_{3}$ is the weighted
blow up of $P_{3}$ with weights $(1,1,5)$, $\beta_{3}$ is the
weighted blow up with weights $(1,1,5)$ of the proper transform of
$P_{3}$ on $U_{2}$, $\beta_{2}$ is the weighted blow up  with
weights $(1,1,4)$ of the proper transform of the point $P_{2}$ on
the variety $U_{3}$, $\beta_{4}$ is the weighted blow up with
weights $(1,1,4)$ of the singular point of the variety $U_{3}$
that is contained in the exceptional divisor of the morphism
$\alpha_{3}$, $\gamma_{2}$ is the weighted blow up with weights
$(1,1,4)$ of the proper transform of $P_{2}$ on $U_{34}$,
$\gamma_{4}$ is the weighted blow up with weights $(1,1,4)$ of the
singular point of the variety $U_{23}$ that is contained in the
exceptional divisor of $\beta_{3}$, and $\eta$ is an elliptic
fibration. There is a commutative diagram
$$
\xymatrix{
&U_{2}\ar@{->}[d]_{\alpha_{2}}&&W\ar@{->}[ll]_{\beta}\ar@{->}[d]^{\omega}&\\%
&X\ar@{-->}[rr]_{\xi}&&\mathbb{P}(1,1,5),&}
$$
where $\xi$ is a projection, $\beta$ is the weighted blow up with
weights $(1,1,3)$ of the singular point of the variety $U_{2}$
that is a quotient singularity of type $\frac{1}{4}(1,1,3)$, and
$\omega$ is an elliptic fibration.

\begin{proposition}
\label{proposition:n-31} Either there is a commutative diagram
\begin{equation}
\label{equation:n-31-commutative-diagram-I} \xymatrix{
&X\ar@{-->}[d]_{\psi}\ar@{-->}[rr]^{\rho}&&V\ar@{->}[d]^{\nu}&\\%
&\mathbb{P}(1,1,4)\ar@{-->}[rr]_{\phi}&&\mathbb{P}^{2},&}
\end{equation}
or there is a commutative diagram
\begin{equation}
\label{equation:n-31-commutative-diagram-II} \xymatrix{
&X\ar@{-->}[d]_{\xi}\ar@{-->}[rr]^{\theta}&&X\ar@{-->}[rr]^{\rho}&&V\ar@{->}[d]^{\nu}&\\%
&\mathbb{P}(1,1,5)\ar@{-->}[rrrr]_{\sigma}&&&&\mathbb{P}^{2},&}
\end{equation}
where $\phi$, $\theta$ and $\sigma$ are birational maps.
\end{proposition}

Now we prove Proposition~\ref{proposition:n-31}, which implies the
claim of Theorem~\ref{theorem:main} for $n=31$. It follows from
Theorem~\ref{theorem:smooth-points}, Lemma~\ref{lemma:Ryder} and
Proposition~\ref{proposition:singular-points}~that
$$
\varnothing\ne\mathbb{CS}\Big(X,
\frac{1}{k}\mathcal{M}\Big)\subseteq\{P_{2},P_{3}\}.
$$

Let $\mathcal{D}_{2}$, $\mathcal{D}_{3}$, $\mathcal{D}_{23}$ and
$\mathcal{D}_{34}$ be the proper transforms of $\mathcal{M}$ on
$U_{2}$, $U_{3}$, $U_{23}$ and $U_{34}$ respectively, then it
follows from Lemma~\ref{lemma:Noether-Fano} that the set
$\mathbb{CS}(U_{\mu}, \frac{1}{k}\mathcal{D}_{\mu})$ is not empty,
if $\mathcal{D}_{\mu}\sim_{\mathbb{Q}} -kK_{U_{\mu}}$.

\begin{lemma}
\label{lemma:n-31-centers-on-U3} Suppose that the set
$\mathbb{CS}(X, \frac{1}{k}\mathcal{M})$ contains the point
$P_{3}$. Let $\bar{P}_{2}$ be the proper transform of the point
$P_{2}$ on the variety $U_{3}$, and $\bar{P}_{4}$ be the singular
point of $U_{3}$ that is contained in the exceptional divisor of
$\alpha_{3}$. Then $\mathcal{D}_{3}\sim_{\mathbb{Q}} -kK_{U_{3}}$
and $\mathbb{CS}(U_{3},
\frac{1}{k}\mathcal{D}_{3})\subseteq\{\bar{P}_{2},\bar{P}_{4}\}$.
\end{lemma}

\begin{proof}
The equivalence $\mathcal{D}_{3}\sim_{\mathbb{Q}}-kK_{U_{3}}$
follows from Theorem~\ref{theorem:Kawamata}. Suppose that
$$
\mathbb{CS}\Big(U_{3},
\frac{1}{k}\mathcal{D}_{3}\Big)\not\subseteq\{\bar{P}_{2},\bar{P}_{4}\},
$$
and let $G$ be the exceptional divisor of $\alpha_{3}$. Then
$G\cong\mathbb{P}(1,1,5)$, and it follows from
Lemma~\ref{lemma:Cheltsov-Kawamata} that there is a curve
$C\subset G$ of the linear system
$|\mathcal{O}_{\mathbb{P}(1,\,1,\,5)}(1)|$ that is contained in
$\mathbb{CS}(U_{3}, \frac{1}{k}\mathcal{D}_{3})$, which is
impossible by Lemma~\ref{lemma:curves}.
\end{proof}

\begin{lemma}
\label{lemma:n-31-centers-on-U2} Let $\bar{P}_{3}$ be the proper
transform of $P_{3}$ on $U_{2}$. Suppose that $\mathbb{CS}(X,
\frac{1}{k}\mathcal{M})$ contains the point $P_{2}$. Then either
$\mathbb{CS}(U_{2}, \frac{1}{k}\mathcal{D}_{2})=\{\bar{P}_{3}\}$,
or the commutative
diagram~\ref{equation:n-31-commutative-diagram-II} exists.
\end{lemma}

\begin{proof}
The equivalence $\mathcal{D}_{2}\sim_{\mathbb{Q}}-kK_{U_{2}}$ is
implied by Theorem~\ref{theorem:Kawamata}. Suppose that the set of
centers of canonical singularities $\mathbb{CS}(U_{2},
\frac{1}{k}\mathcal{D}_{2})$ does not consists of the point
$\bar{P}_{3}$. Let $\bar{P}_{5}$ be the singular point of $U_{2}$
that is contained in the exceptional divisor of $\alpha_{2}$. Then
$\bar{P}_{5}$ is a quotient singularity of type
$\frac{1}{4}(1,1,3)$ on $U_{2}$, and  $\mathbb{CS}(U_{2},
\frac{1}{k}\mathcal{D}_{2})$ contains $\bar{P}_{5}$ by
Lemma~\ref{lemma:Cheltsov-Kawamata}.

Let $\mathcal{B}$ be the proper transform of $\mathcal{M}$ on the
variety $W$. Then $\mathcal{B}\sim_{\mathbb{Q}}-kK_{W}$ by
Theorem~\ref{theorem:Kawamata}, which implies the existence of the
commutative diagram~\ref{equation:n-31-commutative-diagram-II}.
\end{proof}

It follows from Theorem~\ref{theorem:Kawamata} that we may assume
that either the equivalence
$\mathcal{D}_{23}\sim_{\mathbb{Q}}-kK_{U_{23}}$ holds, or the
equivalence $\mathcal{D}_{34}\sim_{\mathbb{Q}}-kK_{U_{34}}$ holds.
Let $\mathcal{D}$ be the proper transform of $\mathcal{M}$ on
$Y$.

\begin{lemma}
\label{lemma:n-31-centers-on-U23} Suppose that
$\mathcal{D}_{23}\sim_{\mathbb{Q}}-kK_{U_{23}}$. Then
$\mathcal{D}\sim_{\mathbb{Q}} -kK_{Y}$.
\end{lemma}

\begin{proof}
Let $F$ be the exceptional  divisor of $\beta_{2}$, $G$ be the
exceptional  divisor of $\beta_{3}$, $\check{P}_{4}$ be the
singular point of $G$, and $\check{P}_{5}$ be the singular point
of $F$. Then the proof of Lemma~\ref{lemma:n-31-centers-on-U2}
implies that $\mathbb{CS}(U_{23}, \frac{1}{k}\mathcal{D}_{23})$
does not contain the point $\check{P}_{5}$. Hence, it follows from
Lemmas~\ref{lemma:Cheltsov-Kawamata} and \ref{lemma:Noether-Fano}
that the set $\mathbb{CS}(U_{23}, \frac{1}{k}\mathcal{D}_{23})$
contains $\check{P}_{4}$, which implies
$\mathcal{D}\sim_{\mathbb{Q}}-kK_{Y}$ by
Theorem~\ref{theorem:Kawamata}.
\end{proof}

\begin{lemma}
\label{lemma:n-31-centers-on-U34} Suppose that
$\mathcal{D}_{34}\sim_{\mathbb{Q}}-kK_{U_{34}}$. Then
$\mathcal{D}\sim_{\mathbb{Q}} -kK_{Y}$.
\end{lemma}

\begin{proof}
Let $G$ be the exceptional divisor of the morphism $\beta_{4}$,
$\check{P}_{2}$ be the proper transform of the point $P_{2}$ on
the variety $U_{34}$, and $\check{P}_{6}$ be the singular point of
the surface $G$. Then $G$ is a cone over the smooth rational cubic
curve, and $\check{P}_{6}$ is a quotient singularity of type
$\frac{1}{3}(1,1,2)$ on $U_{34}$.

The set $\mathbb{CS}(U_{23}, \frac{1}{k}\mathcal{D}_{23})$ is not
empty by Lemma~\ref{lemma:Noether-Fano}, but the equivalence
$\mathcal{D}\sim_{\mathbb{Q}} -kK_{Y}$ follows from
Theorem~\ref{theorem:Kawamata}, if the set $\mathbb{CS}(U_{23},
\frac{1}{k}\mathcal{D}_{23})$ contains the point $\check{P}_{2}$.
Therefore, we may assume that the set $\mathbb{CS}(U_{23},
\frac{1}{k}\mathcal{D}_{23})$ contains the point $\check{P}_{6}$
by Lemma~\ref{lemma:Cheltsov-Kawamata}.

Let $\pi\colon W\to U_{34}$ be the weighted blow up of the point
$\check{P}_{6}$ with weights $(1,1,3)$, $\mathcal{B}$ be the
proper transform of the linear system $\mathcal{M}$ on the variety
$W$, and $S$ be a sufficiently general surface of the pencil
$|-K_{W}|$. Then the base locus of the pencil $|-K_{W}|$ consists
of the irreducible the curve $\Delta$ such that  the equivalence
$\mathcal{B}\vert_{S}\sim_{\mathbb{Q}} k\Delta$ holds, and the
inequality $\Delta^{2}<0$ holds on the surface $S$. It follows
from Lemma~\ref{lemma:normal-surface} that
$\mathcal{B}\vert_{S}=k\Delta$, which is impossible by
Lemma~\ref{lemma:Cheltsov-II}.
\end{proof}

Hence, the equivalence $\mathcal{D}\sim_{\mathbb{Q}}-kK_{Y}$
holds, which implies the existence of the
diagram~\ref{equation:n-31-commutative-diagram-I}.

\section{Case $n=32$, hypersurface of degree $16$ in $\mathbb{P}(1,2,3,4,7)$.}
\label{section:n-32}

We use the notations and assumptions of
Section~\ref{section:start}. Let $n=32$. Then $X$ is a
sufficiently general hypersurface in $\mathbb{P}(1,2,3,4,6)$ of
degree $16$, the equality $-K_{X}^{3}=2/21$ holds, and the
singularities of the hypersurface  $X$ consist of  the points
$P_{1}$, $P_{2}$, $P_{3}$ and $P_{4}$ that are quotient
singularities of type $\frac{1}{2}(1,1,1)$, the point $P_{5}$ that
is a quotient singularity of type $\frac{1}{3}(1,1,2)$, and the
point $P_{6}$ that is a singularity of type $\frac{1}{7}(1,3,4)$.
There is a commutative diagram
$$
\xymatrix{
&U\ar@{->}[d]_{\alpha}&&Y\ar@{->}[ll]_{\beta}\ar@{->}[d]^{\eta}&\\%
&X\ar@{-->}[rr]_{\psi}&&\mathbb{P}(1,2,3),&}
$$
where $\psi$ is a projection, $\alpha$ is the weighted blow up of
$P_{6}$ with weights $(1,3,4)$, $\beta$ is the weighted blow up
with weights $(1,1,3)$ of the singular point of the variety $U$
that is a quotient singularity of type $\frac{1}{4}(1,1,3)$
contained in the exceptional divisor of $\alpha$, and $\eta$ is an
elliptic fibration.

\begin{proposition}
\label{proposition:n-32} The claim of Theorem~\ref{theorem:main}
holds for $n=32$.
\end{proposition}

In the rest of the section we prove
Proposition~\ref{proposition:n-32}. It follows from
Theorem~\ref{theorem:smooth-points}, Lemma~\ref{lemma:Ryder} and
Proposition~\ref{proposition:singular-points} that $\mathbb{CS}(X,
\frac{1}{k}\mathcal{M})=\{P_{6}\}$.

Let $E$ be the exceptional  divisor of the morphism $\alpha$, and
$\mathcal{D}$ be the proper transform of the linear system
$\mathcal{M}$ on the variety $U$. Then $E\cong\mathbb{P}(1,3,4)$,
and it follows from Theorem~\ref{theorem:Kawamata} that the
equivalence $\mathcal{D}\sim_{\mathbb{Q}}-kK_{U}$ holds, but the
set $\mathbb{CS}(U, \frac{1}{k}\mathcal{D})$ is not empty by
Lemma~\ref{lemma:Noether-Fano}.

Let $P_{7}$ and $P_{8}$ be the singular points of the variety $U$
contained in the divisor $E$ that are singularities of types
$\frac{1}{3}(1,1,2)$ and $\frac{1}{4}(1,1,3)$ respecitvely.

\begin{lemma}
\label{lemma:n-32-points-P8} Suppose that $P_{8}\in\mathbb{CS}(X,
\frac{1}{k}\mathcal{M})$. Then there is a commutative diagram
\begin{equation}
\label{equation:n-32-commutative-diagramm} \xymatrix{
&X\ar@{-->}[d]_{\psi}\ar@{-->}[rr]^{\rho}&&V\ar@{->}[d]^{\nu}&\\%
&\mathbb{P}(1,2,3)\ar@{-->}[rr]_{\zeta}&&\mathbb{P}^{2},&}
\end{equation}
where $\zeta$ is a birational map.
\end{lemma}

\begin{proof}
Let $\mathcal{H}$ be the proper transform of the linear system
$\mathcal{M}$ on the variety $Y$. Then it follows from
Theorem~\ref{theorem:Kawamata} that the equivalence
$\mathcal{H}\sim_{\mathbb{Q}}-kK_{Y}$ holds. Hence, the linear
system $\mathcal{H}$ lies in the fibers of the fibration $\eta$,
which implies the existence of the commutative
diagram~\ref{equation:n-32-commutative-diagramm}.
\end{proof}

We may assume that $\mathbb{CS}(U,
\frac{1}{k}\mathcal{D})=\{P_{7}\}$ by
Lemma~\ref{lemma:Cheltsov-Kawamata}.

Let $\gamma\colon W\to U$ be the weighted blow up of the point
$P_{7}$ with weights $(1,1,2)$, $F$ be the exceptional  divisor of
the morphism $\gamma$, $\bar{E}$ be the proper transform of the
surface $E$ on the variety $W$, and $\mathcal{B}$ be the proper
transform of $\mathcal{M}$ on $W$. Then $F\cong\mathbb{P}(1,1,2)$
and $\mathcal{B}\sim_{\mathbb{Q}}-kK_{W}$.

The hypersurface $X$ can be given by the quasihomogeneous equation
$$
w^{2}y+wf_{9}(x,y,z,t)+f_{16}(x,y,z,t)=0\subset\mathbb{P}(1,2,3,4,7)\cong\mathrm{Proj}\Big(\mathbb{C}[x,y,z,t,w]\Big),
$$
where $\mathrm{wt}(x)=1$, $\mathrm{wt}(y)=2$, $\mathrm{wt}(z)=3$,
$\mathrm{wt}(t)=4$, $\mathrm{wt}(w)=7$, and $f_{9}$ and $f_{16}$
are quasihomogeneous polynomials of degree $9$ and $16$
respectively. Let $S$ be the unique surface of the linear system
$|-K_{X}|$, and $D$ be a general surface of the pencil
$|-2K_{X}|$. Then the surface $S$ is cut out by $x=0$, and $D$ is
cut out by $\lambda x^{2}+\mu y=0$, where
$(\lambda,\mu)\in\mathbb{P}^{1}$. The surface $D$ is normal, and
the base locus of the linear system $|-2K_{X}|$ consists of the
curve $C$ such that $C=D\cdot S$.

In the neighborhood of the point $P_{6}$, the monomials $x$, $z$
and $t$ can be considered as a weighted local coordinates on $X$
such that $\mathrm{wt}(x)=1$, $\mathrm{wt}(z)=3$ and
$\mathrm{wt}(z)=4$. Then in the neighborhood of the singular point
$P_{6}$, the surface $D$ can be given by equation
$$
\lambda
x^{2}+\mu\big(\epsilon_{1}x^{9}+\epsilon_{2}zx^{6}+\epsilon_{3}z^{2}x^{3}+\epsilon_{4}z^{3}+\epsilon_{5}t^{2}x+
\epsilon_{6}tx^{5}+\epsilon_{7}tzx^{2}+\cdots\big)=0,
$$
where $\epsilon_{i}\in\mathbb{C}$. In the neighborhood of $P_{7}$,
the morphism $\alpha$ can be given by the equations
$$
x=\tilde{x}\tilde{z}^{\frac{1}{7}},\ z=\tilde{z}^{\frac{3}{7}},\ t=\tilde{t}\tilde{z}^{\frac{4}{7}},%
$$
where $\tilde{x}$, $\tilde{y}$ and $\tilde{z}$ are weighted local
coordinates on the variety $U$ in the neighborhood of the singular
point $P_{7}$ such that $\mathrm{wt}(\tilde{x})=1$,
$\mathrm{wt}(\tilde{z})=2$ and $\mathrm{wt}(\tilde{t})=1$. Let
$\tilde{D}$, $\tilde{S}$ and $\tilde{C}$ be the proper transforms
on $U$ of the surface $D$, the surface $S$ and the curve $C$
respectively, and $E$ be the exceptional  divisor of the morphism
$\alpha$. Then in the neighborhood of the singular point $P_{7}$
the surface $E$ is given by the equation $\tilde{z}=0$, the
surface $\tilde{D}$ is given by the vanishing of the function
$$
\lambda
\tilde{x}^{2}+\mu\Big(\epsilon_{1}\tilde{x}^{9}\tilde{z}+\epsilon_{2}\tilde{z}\tilde{x}^{6}+\epsilon_{3}\tilde{z}\tilde{x}^{3}+\epsilon_{4}\tilde{z}+
\epsilon_{5}\tilde{t}^{2}\tilde{x}\tilde{z}+\epsilon_{6}\tilde{t}\tilde{x}^{5}\tilde{z}+\epsilon_{7}\tilde{t}\tilde{z}\tilde{x}^{2}+\cdots\Big),
$$
and the surface $S$ is given by the equation $\tilde{x}=0$.

In the neighborhood of the singular point of $F$, the morphism
$\beta$ can be given by the equations
$$
\tilde{x}=\bar{x}\bar{z}^{\frac{1}{3}},\ \tilde{z}=\bar{z}^{\frac{2}{3}},\ \tilde{t}=\bar{t}\bar{z}^{\frac{1}{3}},%
$$
where $\bar{x}$, $\bar{z}$ and $\bar{t}$ are weighted local
coordinates on the variety $W$ in the neighborhood of the singular
point of the surface $F$ such that
$\mathrm{wt}(\bar{x})=\mathrm{wt}(\bar{z})=\mathrm{wt}(\bar{t})=1$.
The surface $F$ is given by the equation $\bar{z}=0$, the proper
transform of the surface $D$ on the variety $W$ is given by the
vanishing of the analytical function
$$
\lambda
\bar{x}^{2}+\mu\Big(\epsilon_{1}\bar{x}^{9}\bar{z}^{3}+\epsilon_{2}\bar{z}^{2}\bar{x}^{6}+
\epsilon_{3}\bar{z}\bar{x}^{3}+\epsilon_{4}+
\epsilon_{5}\bar{t}^{2}\bar{x}\bar{z}+\epsilon_{6}\bar{t}\bar{x}^{5}\bar{z}^{2}+\epsilon_{7}\bar{t}\bar{z}\bar{x}^{2}+\cdots\Big),
$$
the proper transform of the surface $S$ on the variety $W$ is
given by the equation $\bar{x}=0$, and the proper transform of the
surface $E$ of the variety $W$ is given by the equation
$\bar{z}=0$.

Let $\mathcal{P}$, $\bar{D}$, $\bar{S}$ and $\bar{C}$ be the
proper transforms on the variety $W$ of the pencil $|-2K_{X}|$,
the surface $D$, the surface $S$ and the curve $C$ respectively,
and $\bar{H}$ be the proper transform on the variety $W$ of the
surface that is cut on $X$ by the equation $y=0$. Then the surface
$\bar{D}$ is a general surface of the pencil $\mathcal{P}$.
Moreover, we have
\begin{equation}
\label{equation:n-32-rational-equivalences} \left\{\aligned
&\bar{E}\sim_{\mathbb{Q}}\gamma^{*}(E)-\frac{2}{3}F,\\
&\bar{D}\sim_{\mathbb{Q}}(\alpha\circ\gamma)^{*}(-2K_{X})-\frac{2}{7}\gamma^{*}(E)-\frac{2}{3}F\sim_{\mathbb{Q}}(\alpha\circ\gamma)^{*}(-2K_{X})-\frac{2}{7}\bar{E}-\frac{6}{7}F,\\
&\bar{S}\sim_{\mathbb{Q}}(\alpha\circ\gamma)^{*}(-K_{X})-\frac{1}{7}\gamma^{*}(E)-\frac{1}{3}F\sim_{\mathbb{Q}}(\alpha\circ\gamma)^{*}(-K_{X})-\frac{1}{7}\bar{E}-\frac{3}{7}F,\\
&\bar{H}\sim_{\mathbb{Q}}\gamma^{*}\Big(\alpha^{*}(-2K_{X})-\frac{9}{7}E\Big)\sim_{\mathbb{Q}}(\alpha\circ\gamma)^{*}(-2K_{X})-\frac{9}{7}\bar{E}-\frac{6}{7}F\sim_{\mathbb{Q}}2\bar{S}-\bar{E}.
\endaligned
\right.
\end{equation}

The curve $\bar{C}$ is contained in the base locus of the pencil
$\mathcal{P}$, but the curve $\bar{C}$ is not the only curve in
the base locus of the pencil $\mathcal{P}$. Namely, let $L$ be the
curve on the surface $E$ that is contained in the linear system
$|\mathcal{O}_{\mathbb{P}(1,\,3,\,4)}(1)|$, which means that $L$
is the curve given locally by the equations
$\tilde{x}=\tilde{z}=0$, and $\bar{L}$ be the proper transform of
$L$ on $W$. Then the curve $\bar{L}$ is contained in the base
locus of the pencil $\mathcal{P}$ as well. Moreover, it follows
from the local computations that the base locus of the pencil
$\mathcal{P}$ does not contain curves outside of the union of
$\bar{C}\cup\bar{L}$.

The curve $\bar{C}$ is the intersection of the divisors $\bar{S}$
and $\bar{H}$, and the curve $\bar{L}$ is the intersections of the
divisors $\bar{S}$ and $\bar{E}$. Moreover, we have
$2\bar{C}=\bar{D}\cdot\bar{H}$,
$\bar{C}+\bar{L}=\bar{S}\cdot\bar{D}$ and
$2\bar{L}=\bar{D}\cdot\bar{E}$.

The curves $\bar{C}$ and $\bar{L}$ can be considered as divisors
on the normal surface $\bar{D}$. Then it follows from the
equivalences~\ref{equation:n-32-rational-equivalences} that
\begin{equation}
\label{equation:n-32-intersection-form} \left\{\aligned
&\bar{L}\cdot\bar{L}=\frac{\bar{E}\cdot\bar{E}\cdot\bar{D}}{4}=-\frac{5}{8},\\
&\bar{C}\cdot\bar{C}=\frac{\bar{H}\cdot\bar{H}\cdot\bar{D}}{4}=\bar{S}\cdot\bar{S}\cdot\bar{D}-\bar{S}\cdot\bar{E}\cdot\bar{D}-\frac{3}{2}\bar{E}\cdot\bar{E}\cdot\bar{D}=-\frac{7}{24},\\
&\bar{C}\cdot\bar{L}=\frac{\bar{H}\cdot\bar{E}\cdot\bar{D}}{4}=\frac{\bar{S}\cdot\bar{E}\cdot\bar{D}}{2}-\frac{\bar{E}\cdot\bar{E}\cdot\bar{D}}{4}=\frac{3}{8},\\
\endaligned
\right.
\end{equation}
which implies that the intersection forms of $\bar{C}$ and
$\bar{L}$ on $\bar{D}$ is negatively definite.

Let $G$ be a sufficiently general surface of the linear system
$\mathcal{B}$. Then
$$
G\vert_{\bar{D}}\sim_{\mathbb{Q}}-kK_{W}\vert_{\bar{D}}\sim_{\mathbb{Q}}k\bar{S}\vert_{\bar{D}}\sim_{\mathbb{Q}}k\bar{C}+k\bar{L},
$$
which is impossible by Lemmas~\ref{lemma:Cheltsov-II} and
\ref{lemma:normal-surface}.

\section{Case $n=36$, hypersurface of degree $18$ in $\mathbb{P}(1,1,4,6,7)$.}
\label{section:n-36}

We use the notations and assumptions of
Section~\ref{section:start}. Let $n=36$. Then $X$ is a
sufficiently general hypersurface in $\mathbb{P}(1,1,4,6,7)$ of
degree $18$, the equality $-K_{X}^{3}=3/28$ holds, and the
singularities of  $X$ consist of  the point $P_{1}$ that is a
quotient singularity of type $\frac{1}{2}(1,1,1)$, the point
$P_{2}$ that is a quotient singularity of type
$\frac{1}{4}(1,1,3)$, and the point $P_{3}$ that is a quotient
singularity of type $\frac{1}{7}(1,1,6)$. There is a commutative
diagram
$$
\xymatrix{
&U_{3}\ar@{->}[d]_{\alpha_{3}}&&U_{34}\ar@{->}[ll]_{\beta_{4}}&&Y\ar@{->}[ll]_{\gamma_{5}}\ar@{->}[d]^{\eta}&\\
&X\ar@{-->}[rrrr]_{\psi}&&&&\mathbb{P}(1,1,4)&}
$$
where $\psi$ is a projection, $\alpha_{3}$ is the weighted blow up
of $P_{3}$ with weights $(1,1,7)$, $\beta_{4}$ is the weighted
blow up with weights $(1,1,6)$ of the singular point of the
variety $U_{3}$ that is contained in the exceptional divisor of
the morphism $\alpha_{3}$, $\gamma_{5}$ is the weighted blow up
with weights $(1,1,4)$ of the singular point of the variety
$U_{34}$ that is contained in the exceptional divisor of the
birational morphism $\beta_{4}$, and $\eta$ is an elliptic
fibration.

\begin{remark}
\label{remark:n-36-big-nef} The divisors $-K_{U_{3}}$ and
$-K_{U_{34}}$ are nef and big.
\end{remark}

There is a commutative diagram
$$
\xymatrix{
&&W\ar@{->}[dl]_{\beta_{3}}\ar@{->}[dr]^{\beta_{2}}\ar@{->}[rrrrrrrd]^{\omega}&&&&&\\
&U_{2}\ar@{->}[dr]_{\alpha_{2}}&&U_{3}\ar@{->}[dl]^{\alpha_{3}}&&&&&&\mathbb{P}(1,1,6),\\
&&X\ar@{-->}[rrrrrrru]_{\xi}&&&&&&&}
$$
where $\xi$ is a projection, $\alpha_{2}$ is the  blow up of the
$P_{2}$ with weights $(1,1,3)$, $\alpha_{3}$ is the weighted blow
up of $P_{3}$ with weights $(1,1,6)$, $\beta_{2}$ is the weighted
blow up  with weights $(1,1,3)$ of the proper transform of $P_{2}$
on $U_{3}$, $\beta_{3}$ is the weighted blow up with weights
$(1,1,6)$ of the proper transform of the point $P_{3}$ on $U_{2}$,
and $\omega$ is an elliptic fibration.

\begin{remark}
\label{remark:n-36-big-nef-II} The divisor $-K_{U_{2}}$ is nef and
big.
\end{remark}

In the rest of the section we prove the following result.

\begin{proposition}
\label{proposition:n-36} Either there is a commutative diagram
\begin{equation}
\label{equation:n-36-commutative-diagram-I} \xymatrix{
&X\ar@{-->}[d]_{\psi}\ar@{-->}[rr]^{\zeta}&&X\ar@{-->}[rr]^{\rho}&&V\ar@{->}[d]^{\nu}&\\%
&\mathbb{P}(1,1,4)\ar@{-->}[rrrr]_{\phi}&&&&\mathbb{P}^{2},&}
\end{equation}
or there is a commutative diagram
\begin{equation}
\label{equation:n-36-commutative-diagram-II} \xymatrix{
&X\ar@{-->}[d]_{\xi}\ar@{-->}[rr]^{\rho}&&V\ar@{->}[d]^{\nu}&\\%
&\mathbb{P}(1,1,6)\ar@{-->}[rr]_{\sigma}&&\mathbb{P}^{2},&}
\end{equation}
where $\zeta$, $\phi$ and $\sigma$ are birational maps.
\end{proposition}

It follows from Theorem~\ref{theorem:smooth-points},
Proposition~\ref{proposition:singular-points} and
Lemma~\ref{lemma:Ryder} that
$$
\varnothing\ne\mathbb{CS}\Big(X,
\frac{1}{k}\mathcal{M}\Big)\subseteq\big\{P_{2},P_{3}\big\},
$$
and the existence of the commutative
diagram~\ref{equation:n-26-commutative-diagram-II} is obvious, if
$\mathbb{CS}(X, \frac{1}{k}\mathcal{M})=\{P_{2},P_{3}\}$.

\begin{lemma}
\label{lemma:n-36-points-P3} The set $\mathbb{CS}(X,
\frac{1}{k}\mathcal{M})$ contains point $P_{3}$.
\end{lemma}

\begin{proof}
Suppose that $P_{3}\not\in\mathbb{CS}(X, \frac{1}{k}\mathcal{M})$.
Then $\mathbb{CS}(X, \frac{1}{k}\mathcal{M})=\{P_{2}\}$.  Let
$\mathcal{D}_{2}$ be the proper transform of $\mathcal{M}$ on
$U_{2}$, and $P_{6}$ be the singular point of  $U_{2}$ that is
contained in the exceptional divisor of $\alpha_{2}$. Then
$\mathcal{D}_{2}\sim_{\mathbb{Q}} -kK_{U_{2}}$ by
Theorem~\ref{theorem:Kawamata}, the point $P_{6}$ is a quotient
singularity of type $\frac{1}{3}(1,1,2)$ on $U_{2}$, and it
follows from Lemmas~\ref{lemma:Noether-Fano} and
\ref{lemma:Cheltsov-Kawamata} that $P_{6}\in\mathbb{CS}(U_{2},
\frac{1}{k}\mathcal{D}_{2})$.

Let $\pi\colon Z\to U_{2}$ be the weighted blow up of  $P_{6}$
with weights $(1,1,3)$, and $\mathcal{B}$ be the proper transform
of $\mathcal{M}$ on $Z$, and $S$ be a general surface of
$|-K_{Z}|$. Then $S$ is normal, and the base locus of the
$|-K_{Z}|$ consists of the irreducible curve $\Delta$ such that
$\mathcal{B}\vert_{S}\sim_{\mathbb{Q}} k\Delta$.

The equality $\Delta^{2}=1/7$ holds on the surface $S$, which
contradicts Lemmas~\ref{lemma:Cheltsov-II} and
\ref{lemma:normal-surface}.
\end{proof}

To conclude the proof of the Proposition~\ref{proposition:n-36},
we may assume that $\mathbb{CS}(X,
\frac{1}{k}\mathcal{M})=\{P_{3}\}$.

Let $\mathcal{D}_{3}$ be the proper transform of  $\mathcal{M}$ on
$U_{3}$, and $P_{4}$ be the singular point of $U_{3}$ that is
contained in the exceptional divisor of the $\alpha_{3}$. Then
$\mathcal{D}_{3}\sim_{\mathbb{Q}} -kK_{U_{3}}$, and $P_{4}$ is a
quotient singularity of type $\frac{1}{6}(1,1,5)$ on $U_{3}$ that
is contained in $\mathbb{CS}(U_{3}, \frac{1}{k}\mathcal{D}_{3})$
by Lemmas~\ref{lemma:Noether-Fano} and
\ref{lemma:Cheltsov-Kawamata}.

\begin{lemma}
\label{lemma:n-36-points-P4} The set $\mathbb{CS}(U_{3},
\frac{1}{k}\mathcal{D}_{3})$ consists of the point $P_{4}$.
\end{lemma}

\begin{proof}
Suppose that the set  $\mathbb{CS}(U_{3},
\frac{1}{k}\mathcal{D}_{3})$ contains a subvariety $C$ of the
variety $U_{3}$ that is different from the point $P_{4}$. Let $G$
be the exceptional divisor of $\beta_{4}$. Then $C$ is a curve
that is contained in $G$ by Lemma~\ref{lemma:Cheltsov-Kawamata},
which is impossible by Lemma~\ref{lemma:curves}.
\end{proof}

Hence, the set $\mathbb{CS}(U_{3}, \frac{1}{k}\mathcal{D}_{3})$
consists of the point $P_{4}$. Let $\mathcal{D}_{34}$ be the
proper transform of the linear system $\mathcal{M}$ on the variety
$U_{34}$, and $P_{5}$ be the singular point of the variety
$U_{34}$ that is contained in exceptional divisor of the morphism
$\beta_{4}$. Then $\mathcal{D}_{34}\sim_{\mathbb{Q}} -kK_{U_{34}}$
by Theorem~\ref{theorem:Kawamata}, the point $P_{5}$ is a quotient
singularity of type $\frac{1}{5}(1,1,4)$ on $U_{34}$, but
$\mathbb{CS}(U_{34}, \frac{1}{k}\mathcal{D}_{34})$ contains the
point $P_{5}$ by Lemmas~\ref{lemma:Noether-Fano} and
\ref{lemma:Cheltsov-Kawamata}. It follows from
Theorem~\ref{theorem:Kawamata} that the proper transform of the
linear system $\mathcal{M}$ on $Y$ is contained in fibers of the
elliptic fibration $\eta$, which implies the existence of the
commutative diagram~\ref{equation:n-36-commutative-diagram-I}. The
claim of Proposition~\ref{proposition:n-36} is proved.

\section{Case $n=38$, hypersurface of degree $18$ in $\mathbb{P}(1,2,3,5,8)$.}
\label{section:n-38}

We use the notations and assumptions of
Section~\ref{section:start}. Let $n=38$. Then $X$ is a
hypersurface of degree $18$ in $\mathbb{P}(1,2,3,5,8)$, the
singularities of $X$ consist of the points $P_{1}$ and $P_{2}$
that are quotient singularities of type $\frac{1}{2}(1,1,1)$, the
point $P_{3}$ that is a quotient singularity of type
$\frac{1}{5}(1,2,3)$, and the  point $P_{4}$ that is a quotient
singularity of type $\frac{1}{8}(1,3,5)$, and the equality
$-K_{X}^{3}=3/40$ holds.

There is a commutative diagram
$$
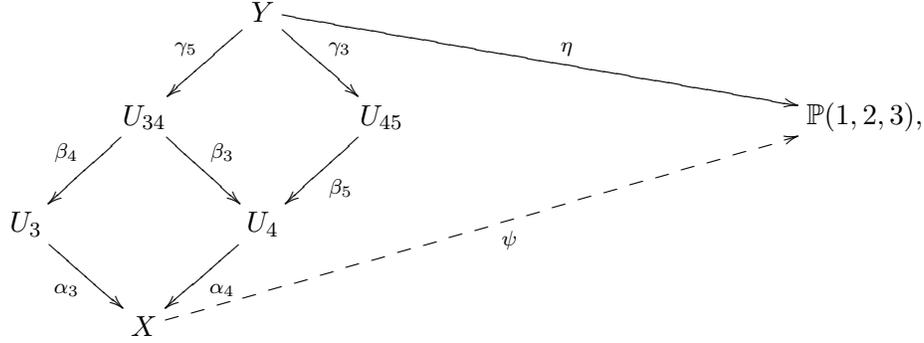
\xymatrix{
&&&Y\ar@{->}[dl]_{\gamma_{5}}\ar@{->}[dr]^{\gamma_{3}}\ar@{->}[drrrrrr]^{\eta}&&&&&&\\
&&U_{34}\ar@{->}[dl]_{\beta_{4}}\ar@{->}[dr]^{\beta_{3}}&&U_{45}\ar@{->}[dl]^{\beta_{5}}&&&&&\mathbb{P}(1,2,3),\\
&U_{3}\ar@{->}[dr]_{\alpha_{3}}&&U_{4}\ar@{->}[dl]^{\alpha_{4}}&&&&&&\\
&&X\ar@{-->}[rrrrrrruu]_{\psi}&&&&&&&}
$$
where $\psi$ is a projection, $\alpha_{3}$ is the weighted blow up
of $P_{3}$ with weights $(1,2,3)$, $\alpha_{4}$ is the weighted
blow up of $P_{4}$ with weights $(1,3,5)$, $\beta_{4}$ is the
weighted blow up with weights $(1,3,5)$ of the proper transform of
$P_{4}$ on  $U_{3}$, $\beta_{3}$ is the weighted blow up with
weights $(1,2,3)$ of the proper transform of the point $P_{3}$ on
$U_{4}$, $\beta_{5}$ is the weighted blow up with weights
$(1,2,3)$ of the singular point of the variety $U_{4}$ that is a
quotient singularity of type $\frac{1}{5}(1,2,3)$ contained in the
exceptional divisor of the morphism $\alpha_{4}$, $\gamma_{3}$ is
the weighted blow up with weights $(1,2,3)$ of the proper
transform of the point $P_{3}$ on the $U_{45}$, $\gamma_{5}$ is
the weighted blow up with weights $(1,2,3)$ of the singular point
of the variety $U_{34}$ that is a quotient singularity of type
$\frac{1}{5}(1,2,3)$ contained in the exceptional divisor of the
morphism $\beta_{4}$, and $\eta$ is an elliptic fibration.

\begin{proposition}
\label{proposition:n-38} The claim of Theorem~\ref{theorem:main}
holds for $n=38$.
\end{proposition}

\begin{proof}
It follows from Theorem~\ref{theorem:smooth-points},
Lemma~\ref{lemma:Ryder} and
Proposition~\ref{proposition:singular-points} that
$$
\varnothing\ne\mathbb{CS}\Big(X,
\frac{1}{k}\mathcal{M}\Big)\subseteq\big\{P_{3},P_{4}\big\},
$$
but the proof of Proposition~\ref{proposition:n-27} implies that
$P_{4}\in\mathbb{CS}(X, \frac{1}{k}\mathcal{M})$.

Let $\mathcal{D}_{4}$ be the proper transform of the linear system
$\mathcal{M}$ on the variety $U_{4}$, $\bar{P}_{3}$ be  the proper
transform of the point $P_{3}$ on $U_{4}$, and $P_{5}$ and $P_{6}$
be the singular points of the variety $U_{4}$ that are quotient
singularities of types $\frac{1}{5}(1,2,3)$ and
$\frac{1}{3}(1,1,2)$  contained in exceptional divisor of the
morphism $\alpha_{4}$ respectively. Then the arguments of the
proof of Proposition~\ref{proposition:n-32} imply that
$$
\mathbb{CS}\Big(U_{4}, \frac{1}{k}\mathcal{D}_{4}\Big)\cap\big\{\bar{P}_{3}, P_{5}\big\}\ne\varnothing.%
$$

Suppose that $\bar{P}_{3}\in\mathbb{CS}(U_{4},
\frac{1}{k}\mathcal{D}_{5})$. Then the proofs of
Propositions~\ref{proposition:n-27} and \ref{proposition:n-32}
implies that the set $\mathbb{CS}(U_{4},
\frac{1}{k}\mathcal{D}_{5})$ contains the singular point $P_{5}$.
Therefore, the claim of Theorem~\ref{theorem:Kawamata} implies
that the claim of Theorem~\ref{theorem:main} holds for the
hypersurface $X$.

We may assume that the set $\mathbb{CS}(U_{4},
\frac{1}{k}\mathcal{D}_{5})$ contains the point $P_{5}$.

Let $\mathcal{D}_{45}$ be the proper transform of $\mathcal{M}$ on
$U_{45}$, and $P_{7}$ and $P_{8}$ be the singular points of the
variety $U_{45}$ that are quotient singularities of types
$\frac{1}{2}(1,1,1)$ and $\frac{1}{3}(1,1,2)$ contained in the
exceptional divisor of $\beta_{5}$ respectively. Then if follows
from Lemma~\ref{lemma:Noether-Fano} that $\mathbb{CS}(U_{45},
\frac{1}{k}\mathcal{D}_{45})\ne\varnothing$, and
Theorem~\ref{theorem:Kawamata} implies that the claim of
Theorem~\ref{theorem:main} holds for the hypersurface $X$ in the
case when the set $\mathbb{CS}(U_{45},
\frac{1}{k}\mathcal{D}_{45})$ contains the proper transform of the
point $P_{3}$ on $U_{45}$.

Therefore, it follows from Lemma~\ref{lemma:Cheltsov-Kawamata}
that to conclude the proof of Proposition~\ref{proposition:n-38}
we may assume that the set $\mathbb{CS}(U_{45},
\frac{1}{k}\mathcal{D}_{45})$ contains either the point $P_{7}$,
or the point $P_{8}$.

Suppose that $P_{7}\in\mathbb{CS}(U_{45},
\frac{1}{k}\mathcal{D}_{45})$. Then considering the proper
transform of the complete linear system $|-3K_{X}|$ on the
weighted blow up of the point $P_{7}$ with weights $(1,1,1)$, we
easily obtain a contradiction as in the proof of
Lemma~\ref{lemma:n-27-P2-is-a-center}.

Thus, the set $\mathbb{CS}(U_{45}, \frac{1}{k}\mathcal{D}_{45})$
contains the point $P_{8}$.

Let $\pi\colon W\to U_{45}$ be the weighted blow up of $P_{8}$,
$\mathcal{B}$ be the proper transform of $\mathcal{M}$ on the
variety $W$, and $D$ be a general surface in $|-2K_{W}|$. Then
$\mathcal{B}\sim_{\mathbb{Q}}-kK_{W}$ by
Theorem~\ref{theorem:Kawamata}, the surface $D$ is normal, the
pencil $|-2K_{W}|$ is the proper transform of the pencil
$|-2K_{X}|$, and the base locus of the pencil $|-2K_{W}|$ consists
of the curves $C$ and $L$ such that
$\alpha_{4}\circ\beta_{5}\circ\pi(C)$ is the unique base curve of
the pencil $|-2K_{X}|$, and the curve $\beta_{5}\circ\pi(L)$ is
contained in the exceptional divisor of the morphism $\alpha_{4}$.

The the intersection form of the curves $C$ and $L$ on the surface
$D$ is negatively definite, but the equivalence
$\mathcal{B}\vert_{D}\sim_{\mathbb{Q}}kC+kL$ holds, which is
impossible by Lemmas~\ref{lemma:normal-surface} and
\ref{lemma:Cheltsov-II}.
\end{proof}

\section{Case $n=40$, hypersurface of degree $19$ in $\mathbb{P}(1,3,4,5,7)$.}
\label{section:n-40}

We use the notations and assumptions of
Section~\ref{section:start}. Let $n=40$. Then $X$ is a
sufficiently general hypersurface in $\mathbb{P}(1,3,4,5,7)$ of
degree $19$, the singularities of the hypersurface $X$ consist of
the point $P_{1}$ that is a quotient singularity of type
$\frac{1}{3}(1,1,2)$, the point $P_{2}$ that is a quotient
singularity of type $\frac{1}{4}(1,1,3)$, the point $P_{3}$ that
is a quotient singularity of type $\frac{1}{5}(1,2,3)$, and the
point $P_{4}$ that is a quotient singularity of type
$\frac{1}{7}(1,3,4)$, and $-K_{X}^{3}=19/420$.

There is a commutative diagram
$$
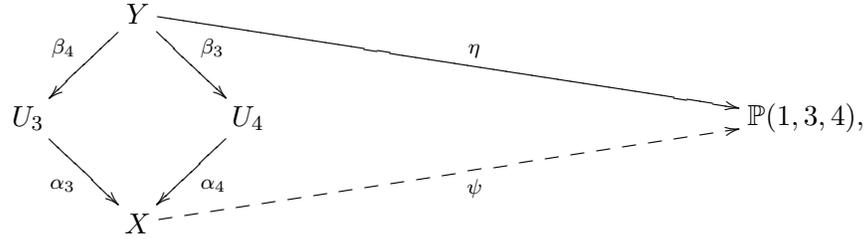
\xymatrix{
&&Y\ar@{->}[dl]_{\beta_{4}}\ar@{->}[dr]^{\beta_{3}}\ar@{->}[rrrrrrrd]^{\eta}&&&&&\\
&U_{3}\ar@{->}[dr]_{\alpha_{3}}&&U_{4}\ar@{->}[dl]^{\alpha_{4}}&&&&&&\mathbb{P}(1,3,4),\\
&&X\ar@{-->}[rrrrrrru]_{\psi}&&&&&&&}
$$
where $\psi$ is the natural projection, $\alpha_{3}$ is the
weighted blow up of the singular point $P_{3}$ with weights
$(1,2,3)$, $\alpha_{4}$ is the weighted blow up of $P_{4}$ with
weights $(1,3,4)$, $\beta_{3}$ is the weighted blow up with
weights $(1,2,3)$ of the proper transform of $P_{3}$ on $U_{4}$,
$\beta_{4}$ is the weighted blow up with weights $(1,3,4)$ of the
proper transform of $P_{4}$ on  $U_{3}$, and $\eta$ is an elliptic
fibration.

In the rest of the section we prove the following result.

\begin{proposition}
\label{proposition:n-40} The claim of Theorem~\ref{theorem:main}
holds for $n=40$.
\end{proposition}

It follows from Theorem~\ref{theorem:smooth-points},
Lemma~\ref{lemma:Ryder} and
Proposition~\ref{proposition:singular-points} that
$$
\varnothing\ne\mathbb{CS}\Big(X,
\frac{1}{k}\mathcal{M}\Big)\subseteq\big\{P_{3}, P_{4}\big\}.
$$

Let $\mathcal{D}_{3}$ and $\mathcal{D}_{4}$ be the proper
transforms of the linear system $\mathcal{M}$ on $U_{3}$ and
$U_{4}$ respectively.

\begin{lemma}
\label{lemma:n-40-points-P3-P4} Suppose that $\mathbb{CS}(X,
\frac{1}{k}\mathcal{M})=\{P_{3}, P_{4}\}$. Then there is a
commutative diagram
\begin{equation}
\label{equation:n-40-commutative-diagramm} \xymatrix{
&X\ar@{-->}[d]_{\psi}\ar@{-->}[rr]^{\rho}&&V\ar@{->}[d]^{\nu}&\\%
&\mathbb{P}(1,3,4)\ar@{-->}[rr]_{\zeta}&&\mathbb{P}^{2},&}
\end{equation}
where $\zeta$ is a birational map.
\end{lemma}

\begin{proof}
Let $\mathcal{H}$ be the proper transform of the linear system
$\mathcal{M}$ on the variety $Y$. Then it follows from
Theorem~\ref{theorem:Kawamata} that the equivalence
$\mathcal{H}\sim_{\mathbb{Q}}-kK_{Y}$ holds, which implies that
$\mathcal{H}$ lies in the fibers of the elliptic fibration $\eta$,
which implies the existence of the commutative
diagram~\ref{equation:n-40-commutative-diagramm}.
\end{proof}

Let $P_{5}$ and $P_{6}$ be the singular points of the variety
$U_{3}$  that are contained in the exceptional divisor of
$\alpha_{3}$ such that $P_{5}$ and $P_{6}$ are quotient
singularities of types $\frac{1}{2}(1,1,1)$ and
$\frac{1}{3}(1,1,2)$ respectively, and $P_{7}$ and $P_{8}$ are the
singular points of the variety $U_{4}$ that are quotient
singularities of types $\frac{1}{3}(1,1,2)$ and
$\frac{1}{4}(1,1,3)$ contained in the exceptional divisor of the
morphism $\alpha_{4}$ respectively. Then it follows from
Theorem~\ref{theorem:Kawamata} and Lemmas~\ref{lemma:Noether-Fano}
and Lemma~\ref{lemma:Cheltsov-Kawamata} that
$$
\left\{\aligned
&P_{4}\not\in\mathbb{CS}\Big(X, \frac{1}{k}\mathcal{M}\Big)\Rightarrow\mathbb{CS}\Big(U_{3}, \frac{1}{k}\mathcal{D}_{3}\Big)\cap\big\{P_{5}, P_{6}\big\}\ne\varnothing,\\%
&P_{3}\not\in\mathbb{CS}\Big(X, \frac{1}{k}\mathcal{M}\Big)\Rightarrow\mathbb{CS}\Big(U_{4}, \frac{1}{k}\mathcal{D}_{4}\Big)\cap\big\{P_{7}, P_{8}\big\}\ne\varnothing.\\%
\endaligned
\right.
$$

\begin{lemma}
\label{lemma:n-40-point-P5} Suppose that
$P_{4}\not\in\mathbb{CS}(X, \frac{1}{k}\mathcal{M})$. Then
$P_{5}\not\in\mathbb{CS}(U_{3}, \frac{1}{k}\mathcal{D}_{3})$.
\end{lemma}

\begin{proof}
Suppose that the set $\mathbb{CS}(U_{3},
\frac{1}{k}\mathcal{D}_{3})$ contains the point $P_{5}$. Let
$\gamma\colon W\to U_{3}$ be the weighted blow up of $P_{5}$ with
weights $(1,1,1)$, $F$ be the exceptional  divisor of $\gamma$,
$\mathcal{D}$ be the proper transform of the linear system
$\mathcal{M}$ on the variety $W$, $\mathcal{H}$ be the proper
transform of the pencil $|-3K_{X}|$ on the variety $W$, $\bar{D}$
be a sufficiently general surface of the linear system
$\mathcal{D}$, and $\bar{H}$ be a sufficiently general surface of
the pencil $\mathcal{H}$. Then we have the equivalence
$$
\bar{H}\sim_{\mathbb{Q}}(\alpha_{3}\circ\gamma)^{*}(-3K_{X})-\frac{3}{5}\gamma^{*}(E)-\frac{1}{2}F,
$$
where $E$ is the exceptional divisor of the morphism $\alpha_{3}$.
The base locus of the pencil $\mathcal{H}$ consists of the curve
$\bar{C}$ such that $\alpha_{3}\circ\gamma(C)$ is the unique base
curve of the pencil $|-3K_{X}|$. On the other hand, it follows
from Theorem~\ref{theorem:Kawamata} that
$\bar{D}\sim_{\mathbb{Q}}-kK_{W}$.

The equivalence $\bar{D}\vert_{\bar{H}}\sim_{\mathbb{Q}}k\bar{C}$
and the inequality $\bar{C}^{2}<0$ hold on the surface $\bar{H}$.
It follows from Lemma~\ref{lemma:normal-surface} that the support
of $\bar{H}\cdot\bar{D}$ consists of $\bar{C}$, which contradicts
Lemma~\ref{lemma:Cheltsov-II}.
\end{proof}

\begin{lemma}
\label{lemma:n-40-point-P6} Suppose that
$P_{4}\not\in\mathbb{CS}(X, \frac{1}{k}\mathcal{M})$. Then
$P_{6}\not\in\mathbb{CS}(U_{3}, \frac{1}{k}\mathcal{D}_{3})$.
\end{lemma}

\begin{proof}
Suppose that $P_{6}\in\mathbb{CS}(U_{3},
\frac{1}{k}\mathcal{D}_{3})$. Let $\gamma\colon W\to U_{3}$ be the
weighted blow up of $P_{6}$ with weights $(1,1,2)$, $F$ and $G$ be
the exceptional divisors of $\alpha_{3}$ and $\gamma$
respectively, $\mathcal{B}$ and $\mathcal{D}$ be the proper
transforms of  $\mathcal{M}$ and $|-7K_{X}|$ on the variety $W$
respectively, and $D$ be a general surface of the linear system
$\mathcal{D}$. Then it follows from Theorem~\ref{theorem:Kawamata}
that $\mathcal{B}\sim_{\mathbb{Q}}-kK_{W}$, but the base locus of
the linear system $\mathcal{D}$ does not contain curves. Moreover,
we have
$$
D\sim_{\mathbb{Q}}(\alpha_{3}\circ\gamma)^{*}\big(-7K_{X}\big)-\frac{2}{5}\gamma^{*}(F)-\frac{2}{3}G,
$$
the divisor $D$ is nef, but the explicit calculations imply that
$$
D\cdot B_{1}\cdot
B_{2}=\Big((\alpha_{3}\circ\gamma)^{*}\big(-7K_{X}\big)-\frac{2}{5}F-G\Big)
\Big((\alpha_{3}\circ\gamma)^{*}\big(-kK_{X}\big)-\frac{k}{5}F-\frac{k}{2}G\Big)^{2}=-\frac{1}{12}k^{2},%
$$
where $B_{1}$ and $B_{2}$ are general surfaces in $\mathcal{B}$,
which is a contradiction.
\end{proof}

\begin{lemma}
\label{lemma:n-40-point-P7} Suppose that
$P_{3}\not\in\mathbb{CS}(X, \frac{1}{k}\mathcal{M})$. Then
$P_{7}\not\in\mathbb{CS}(U_{4}, \frac{1}{k}\mathcal{D}_{4})$.
\end{lemma}

\begin{proof}
Suppose that the set $\mathbb{CS}(U_{4},
\frac{1}{k}\mathcal{D}_{4})$ contains the point $P_{7}$. Let
$\gamma\colon W\to U_{4}$ be the weighted blow up of  $P_{7}$ with
weights $(1,1,2)$, $F$ be the exceptional  divisor of  $\gamma$,
$\mathcal{D}$ be the proper transform of the linear system
$\mathcal{M}$ on $W$, $\mathcal{H}$ be the proper transform of
$|-4K_{X}|$ on $W$, $\bar{D}$ be a general surface of the linear
system $\mathcal{D}$, and $\bar{H}$ be a general surface of the
linear system $\mathcal{H}$. Then
$$
\bar{H}\sim_{\mathbb{Q}}(\alpha_{3}\circ\gamma)^{*}(-4K_{X})-\frac{4}{7}\gamma^{*}(E)-\frac{1}{3}F,
$$
and the base locus of $\mathcal{H}$ consists of the curve
$\bar{C}$ such that $\alpha_{3}\circ\gamma(C)$ is the  base curve
of the linear system $|-4K_{X}|$. It follows from
Theorem~\ref{theorem:Kawamata} that the equivalence
$\bar{D}\sim_{\mathbb{Q}}-kK_{W}$ holds.

The equality $\bar{C}^{2}=-1/30$ holds on the normal surface
$\bar{H}$, which implies that the support of the cycle
$\bar{H}\cdot\bar{D}$ consists of $\bar{C}$, because
$\bar{D}\vert_{\bar{H}}\sim_{\mathbb{Q}}k\bar{C}$, which
contradicts Lemma~\ref{lemma:Cheltsov-II}.
\end{proof}

\begin{lemma}
\label{lemma:n-40-point-P8} Suppose that
$P_{3}\not\in\mathbb{CS}(X, \frac{1}{k}\mathcal{M})$. Then
$P_{8}\not\in\mathbb{CS}(U_{4}, \frac{1}{k}\mathcal{D}_{4})$.
\end{lemma}

\begin{proof}
Suppose that the set $\mathbb{CS}(U_{4},
\frac{1}{k}\mathcal{D}_{4})$ contains the point $P_{8}$. Let
$\gamma\colon W\to U_{4}$ be the weighted blow up of $P_{8}$ with
weights $(1,1,3)$, $\mathcal{D}$ be the proper transform of
$\mathcal{M}$ on $W$, $\bar{H}$ be a general surface of the pencil
$|-3K_{W}|$, and $\bar{D}$ be a general surface in $\mathcal{D}$.
Then $\bar{D}\sim_{\mathbb{Q}}-kK_{W}$, but the base locus of the
pencil $|-3K_{W}|$ consists of  the irreducible curve $\bar{C}$
such that $\alpha_{4}\circ\gamma(C)$ is the base curve of
$|-3K_{X}|$. Moreover, the equality $\bar{C}^{2}=-1/20$ holds on
$\bar{H}$, but $\bar{D}\vert_{\bar{H}}\sim_{\mathbb{Q}}k\bar{C}$,
which implies that the support of $\bar{H}\cdot\bar{D}$ consists
of $\bar{C}$, which is impossible by
Lemma~\ref{lemma:Cheltsov-II}.
\end{proof}

The claim Proposition~\ref{proposition:n-40} is proved.

\section{Case $n=43$, hypersurface of degree $20$ in $\mathbb{P}(1,2,4,5,9)$.}
\label{section:n-43}

We use the notations and assumptions of
Section~\ref{section:start}. Let $n=43$. Then $X$ is a
sufficiently general hypersurface in $\mathbb{P}(1,2,4,5,9)$ of
degree $20$, the equality $-K_{X}^{3}=1/18$ holds, and the
singularities of the hypersurface $X$ consist of the points
$P_{1}$, $P_{2}$, $P_{3}$, $P_{4}$ and $P_{5}$ that are quotient
singularities of type $\frac{1}{2}(1,1,1)$, and the point $P_{6}$
that is a quotient singularity of type $\frac{1}{9}(1,4,5)$.

\begin{proposition}
\label{proposition:n-43} The claim of Theorem~\ref{theorem:main}
holds for $n=43$.
\end{proposition}

In the rest of the section we prove
Proposition~\ref{proposition:n-43}. It follows from
Theorem~\ref{theorem:smooth-points}, Lemma~\ref{lemma:Ryder} and
Proposition~\ref{proposition:singular-points} that $\mathbb{CS}(X,
\frac{1}{k}\mathcal{M})=\{P_{6}\}$. There is a commutative diagram
$$
\xymatrix{
&U\ar@{->}[d]_{\alpha}&&Y\ar@{->}[ll]_{\beta}\ar@{->}[d]^{\eta}&\\%
&X\ar@{-->}[rr]_{\psi}&&\mathbb{P}(1,2,4),&}
$$
where $\psi$ is a projection, $\alpha$ is the weighted blow up of
$P_{6}$ with weights $(1,4,5)$, $\beta$ is the weighted blow up
with weights $(1,1,4)$ of the singular point of the variety $U$
that is a quotient singularity of type $\frac{1}{5}(1,1,4)$
contained in the exceptional divisor of $\alpha$, and $\eta$ is an
elliptic fibration.

Let $\mathcal{D}$ be the proper transform of the linear system
$\mathcal{M}$ on the variety $U$, and $P_{7}$ and $P_{8}$ be the
singular points of $U$ that are quotient singularities of types
$\frac{1}{4}(1,1,3)$ and $\frac{1}{5}(1,1,4)$ contained in the
exceptional divisor of $\alpha$ respectively. Then
$\mathcal{D}\sim_{\mathbb{Q}}-kK_{U}$ by
Theorem~\ref{theorem:Kawamata}.

We must show that the proper transform of the linear system
$\mathcal{M}$ on the variety $Y$ is contained in the fibers of the
$\eta$, which is implied by Theorem~\ref{theorem:Kawamata}, if
$P_{8}\in\mathbb{CS}(U, \frac{1}{k}\mathcal{D})$. Thus, to
conclude the proof of Proposition~\ref{proposition:n-43} we may
assume that $P_{8}\not\in\mathbb{CS}(U, \frac{1}{k}\mathcal{D})$.

\begin{remark}
\label{remark:n-43-P7} The set $\mathbb{CS}(U,
\frac{1}{k}\mathcal{D})$ contains $P_{7}$ by
Lemma~\ref{lemma:Cheltsov-Kawamata}, because $-K_{U}$ is nef and
big.
\end{remark}

Let $\gamma\colon W\to U$ be the weighted blow up of the singular
point $P_{7}$ with weights $(1,1,3)$, $\mathcal{B}$ be the proper
transform of the linear system $\mathcal{M}$ on the variety $W$,
and $P_{9}$ be the singular point of the variety $W$ that is a
quotient singularity of type $\frac{1}{3}(1,1,2)$ contained in the
exceptional divisor of the morphism $\gamma$. Then the equivalence
$\mathcal{B}\sim_{\mathbb{Q}}-kK_{W}$ holds by
Theorem~\ref{theorem:Kawamata}.

\begin{lemma}
\label{lemma:n-43-P9} The set $\mathbb{CS}(W,
\frac{1}{k}\mathcal{B})$ does not contains the point $P_{9}$.
\end{lemma}

\begin{proof}
Suppose that the set $\mathbb{CS}(W, \frac{1}{k}\mathcal{B})$
contains the point $P_{9}$. Let $\pi\colon Z\to W$ be the weighted
blow up of the singular point $P_{9}$ with weights $(1,1,2)$,
$\mathcal{H}$ be the proper transform of the linear system
$\mathcal{M}$ on the variety $Z$, and $\mathcal{P}$ be the proper
transform of the linear system $|-5K_{X}|$ on the variety $Z$.
Then $\mathcal{H}\sim_{\mathbb{Q}}-kK_{Z}$ by
Theorem~\ref{theorem:Kawamata}, but the base locus of
$\mathcal{P}$ consists of the irreducible curve $\Gamma$ such that
$\alpha\circ\gamma\circ\pi(\Gamma)$ is the base curve in of
$|-5K_{X}|$.

Let $H_{1}$ and $H_{2}$ be general surfaces of the linear system
$\mathcal{H}$, and $D$ be general surface of the linear system
$\mathcal{P}$. Then $D\cdot\Gamma=1$ and $D^{3}=6$. Therefore, the
divisor $D$ is nef and big, but the elementary computations imply
that $D\cdot H_{1}\cdot H_{2}=0$, which is impossible by
Corollary~\ref{corollary:Cheltsov}.
\end{proof}

Therefore, the claim of Lemma~\ref{lemma:Cheltsov-Kawamata}
implies the following corollary.

\begin{corollary}
\label{corollary:n-43-terminal-singularities} The singularities of
the log pair $(W, \frac{1}{k}\mathcal{B})$ are terminal.
\end{corollary}

The hypersurface $X$ can be given by the quasihomogeneous equation
of degree $20$
$$
w^{2}y+wf_{11}(x,y,z,t)+f_{20}(x,y,z,t)=0\subset\mathbb{P}(1,2,4,5,9)\cong\mathrm{Proj}\Big(\mathbb{C}[x,y,z,t,w]\Big),
$$
where $\mathrm{wt}(x)=1$, $\mathrm{wt}(y)=2$, $\mathrm{wt}(z)=4$,
$\mathrm{wt}(t)=5$, $\mathrm{wt}(w)=9$, and $f_{i}(x,y,z,t)$ is a
quasihomogeneous polynomial of degree $i$. Let $D$ be a general
surface in $|-2K_{X}|$, and $S$ be a surface that is cut on the
hypersurface $X$ by the equation $x=0$. Then $D$ is cut on $X$ by
the quasihomogeneous equation $\lambda x^{2}+\mu y=0$, where
$(\lambda,\mu)\in\mathbb{P}^{1}$, and the base locus of
$|-2K_{X}|$ consists of the irreducible curve $C$ that is cut on
the hypersurface  $X$ by the equations $x=y=0$.

In the neighborhood of the point $P_{6}$, the monomials $x$, $z$
and $t$ can be considered as weighted local coordinates on $X$
such that $\mathrm{wt}(x)=1$, $\mathrm{wt}(z)=4$ and
$\mathrm{wt}(t)=5$. Then in the neighborhood of the singular point
$P_{7}$, the weighted blow up $\alpha$ is given by the equations
$$
x=\tilde{x}\tilde{z}^{\frac{1}{9}},\ z=\tilde{z}^{\frac{4}{9}},\ t=\tilde{t}\tilde{z}^{\frac{5}{9}},%
$$
where $\tilde{x}$, $\tilde{z}$ and $\tilde{t}$ are weighted local
coordinated on the variety $U$ in the neighborhood of the singular
point $P_{7}$ such that $\mathrm{wt}(\tilde{x})=1$,
$\mathrm{wt}(\tilde{z})=3$ and $\mathrm{wt}(\tilde{t})=1$.

Let $E$ be the exceptional  divisor of the morphism  $\alpha$, and
$\tilde{D}$, $\tilde{S}$ and $\tilde{C}$ be the proper transforms
on the variety $U$ of the surface $D$, the surface $S$ and the
curve $C$ respectively. Then $E$ is given by the equation
$\tilde{z}=0$, and the surface $\tilde{S}$ is given by the
equation $\tilde{x}=0$. Moreover, it follows from the local
equation of the surface $\tilde{D}$ that
$\tilde{D}\cdot\tilde{S}=\tilde{C}+2\tilde{L}_{1}$, where
$\tilde{L}_{1}$ is the curve that is locally given by the
equations $\tilde{z}=\tilde{x}=0$. Moreover, the surface
$\tilde{D}$ is not normal in a general point of the curve
$\tilde{L}_{1}$. Nevertheless, we have the equivalences
$$
\tilde{D}\sim_{\mathbb{Q}}2\tilde{S}\sim_{\mathbb{Q}}\alpha^{*}\big(-2K_{X}\big)-\frac{2}{9}E.
$$

In the neighborhood of the point $P_{9}$ the morphism $\gamma$ is
given by the equations
$$
\tilde{x}=\bar{x}\bar{z}^{\frac{1}{4}},\ \tilde{z}=\bar{z}^{\frac{3}{4}},\ \tilde{t}=\bar{t}\bar{z}^{\frac{1}{4}},%
$$
where $\bar{x}$, $\bar{z}$ and $\bar{t}$ are weighted local
coordinates on the variety $W$ in the neighborhood of the point
$P_{9}$ such that $\mathrm{wt}(\bar{x})=1$,
$\mathrm{wt}(\bar{z})=2$ and $\mathrm{wt}(\bar{t})=1$. In
particular, the exceptional divisor of the morphism $\gamma$ is
given by the equation $\bar{z}=0$, and the proper transform of the
surface $S$ on the variety $W$ is given by the equation
$\bar{x}=0$.

Let $F$ be the exceptional divisor of the morphism $\gamma$, and
$\bar{D}$, $\bar{S}$, $\bar{E}$, $\bar{C}$ and $\bar{L}_{1}$ be
the proper transforms on the variety $W$ of the surface $D$, the
surface $S$, the surface $E$, the curve $C$ and the curve
$\tilde{L}_{1}$ respectively. Then we the equivalence
\begin{equation}
\label{equation:n-43-rational-equivalences-I} \left\{\aligned
&\bar{S}\sim_{\mathbb{Q}}-K_{W}\sim_{\mathbb{Q}}(\alpha\circ\gamma)^{*}\big(-K_{X}\big)-\frac{1}{9}\gamma^{*}(E)-\frac{1}{4}F,\\
&\bar{D}\sim_{\mathbb{Q}}-2K_{W}\sim_{\mathbb{Q}}(\alpha\circ\gamma)^{*}\big(-2K_{X}\big)-\frac{2}{9}\gamma^{*}(E)-\frac{1}{2}F,\\
&\bar{E}\sim_{\mathbb{Q}}\gamma^{*}(E)-\frac{3}{4}F.\\
\endaligned
\right.
\end{equation}

Let $\bar{L}_{2}$ be the curve on the variety $W$ that is given by
the equation $\bar{z}=\bar{x}=0$. Then
$$
\bar{D}\cdot\bar{S}=\bar{C}+2\bar{L}_{1}+\bar{L}_{2},\ \bar{D}\cdot\bar{E}=2\bar{L}_{1},\ \bar{D}\cdot F=2\bar{L}_{2},%
$$
but the base locus of $|-2K_{W}|$ consists of the curves
$\bar{C}$, $\bar{L}_{1}$ and $\bar{L}_{2}$. The
equivalences~\ref{equation:n-43-rational-equivalences-I} imply
$$
\bar{D}\cdot\bar{C}=0,\ \bar{D}\cdot\bar{L}_{1}=-\frac{2}{5},\ \bar{D}\cdot\bar{L}_{2}=\frac{2}{3}.%
$$

Let $\bar{H}$ and $\bar{T}$ be the proper transforms on the
variety $W$ of the surfaces that are cut on the hypersurface $X$
by the equations $y=0$ and $t=0$ respectively. Then
\begin{equation}
\label{equation:n-43-rational-equivalences-II} \left\{\aligned
&\bar{T}\sim_{\mathbb{Q}}(\alpha\circ\gamma)^{*}\big(-5K_{X}\big)-\frac{5}{9}\gamma^{*}(E)-\frac{1}{4}F\sim_{\mathbb{Q}}(\alpha\circ\gamma)^{*}\big(-5K_{X}\big)-\frac{5}{9}\bar{E}-\frac{2}{3}F,\\%
&\bar{H}\sim_{\mathbb{Q}}(\alpha\circ\gamma)^{*}\big(-2K_{X}\big)-\frac{11}{9}\gamma^{*}(E)-\frac{3}{2}F\sim_{\mathbb{Q}}(\alpha\circ\gamma)^{*}\big(-2K_{X}\big)-\frac{11}{9}\bar{E}-\frac{5}{3}F,\\%
\endaligned
\right.
\end{equation}
which implies that
$$
-14K_{W}\sim_{\mathbb{Q}}14\bar{D}\sim_{\mathbb{Q}}2\bar{T}+2\bar{H}+2\bar{E},
$$
and the support of the cycle $\bar{T}\cdot\bar{H}$ does not
contain the curves $\bar{L}_{2}$ and $\bar{C}$. Therefore, the
base locus of the linear system $|-14K_{W}|$ does not contain
curves except the curve $\bar{L}_{1}$.

The singularities of the mobile log pair $(W, \lambda |-14K_{W}|)$
are log-terminal for some rational number $\lambda>1/14$, but the
divisor $K_{W}+\lambda |-14K_{W}|$ has non-negative intersection
with all curves on the variety $W$ except the curve $\bar{L}_{1}$.
It follows from \cite{Sho93} that the log-flip $\zeta\colon
W\dasharrow Z$ in the curve $\bar{L}_{1}$ with respect to the log
pair $(W, \lambda |-14K_{W}|)$ exists.

Let $\mathcal{P}$ be the proper transform of the linear system
$\mathcal{M}$ on the variety $Z$. Then the singularities of the
log pair $(Z, \frac{1}{k}\mathcal{P})$ are terminal, because the
singularities of the log pair $(W, \frac{1}{k}\mathcal{B})$ are
terminal, but the rational map $\zeta$ is a log flop with respect
to the log pair $(W, \frac{1}{k}\mathcal{B})$, but $-K_{Z}$ is
numerically effective, because the base locus of the linear system
$|-14K_{W}|$ does not contain curves outside the curve
$\bar{L}_{1}$, and the inequality $-K_{W}\cdot\bar{L}_{1}<0$
holds.

In the rest of the section we show that $-K_{Z}$ is big, which
contradicts Lemma~\ref{lemma:Noether-Fano}.

The rational functions $y/x^{2}$ and $ty/x^{7}$ are contained in
$|2S|$ and $|7S|$ respectively, but the
equivalences~\ref{equation:n-43-rational-equivalences-II} implies
that $y/x^{2}$ and $ty/x^{7}$ are contained in $|2\bar{S}|$ and
$|7\bar{S}|$ respectively.

Let $\bar{Z}$ be the proper transform on the variety $W$ of the
irreducible surface that is cut on the hypersurface $X$ by the
equation $z=0$. Then the equivalences
$$
\bar{Z}\sim_{\mathbb{Q}}(\alpha\circ\gamma)^{*}\big(-4K_{X}\big)-\frac{4}{9}\gamma^{*}(E)\sim_{\mathbb{Q}}(\alpha\circ\gamma)^{*}\big(-4K_{X}\big)-\frac{4}{9}\bar{E}-\frac{1}{3}F
$$
hold, which imply that
$-6K_{W}\sim_{\mathbb{Q}}\bar{Z}+\bar{H}+\bar{E}$. Thus, the
rational function $zy/x^{6}$ is contained in the linear system
$|6\bar{S}|$. Thus, the linear system $|-42K_{W}|$ maps the
variety $W$ dominantly on some three-dimensional variety, which
implies that the divisor $-K_{Z}$ is big.

\section{Case  $n=44$, hypersurface of degree $20$ in $\mathbb{P}(1,2,5,6,7)$.}
\label{section:n-44}

We use the notations and assumptions of
Section~\ref{section:start}. Let $n=44$. Then $X$ is a
sufficiently general hypersurface in $\mathbb{P}(1,2,5,6,7)$ of
degree $20$, the equality $-K_{X}^{3}=1/21$ holds, and the
singularities of the hypersurface $X$ consist of the points
$P_{1}$, $P_{2}$ and $P_{3}$ that are quotient singularities of
type $\frac{1}{2}(1,1,1)$, the point $P_{4}$ that is a quotient
singularity of type $\frac{1}{6}(1,1,5)$, and the point $P_{5}$
that is a quotient singularity of type $\frac{1}{7}(1,2,5)$. There
is a commutative diagram
$$
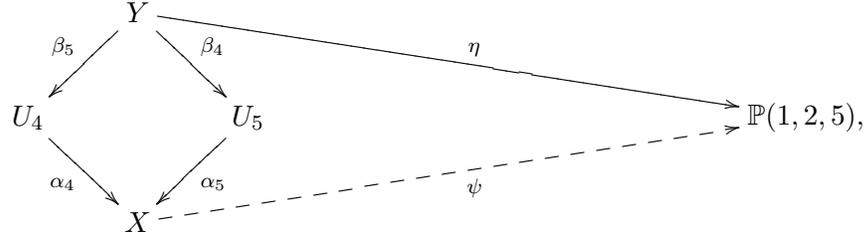
\xymatrix{
&&Y\ar@{->}[dl]_{\beta_{5}}\ar@{->}[dr]^{\beta_{4}}\ar@{->}[rrrrrrrd]^{\eta}&&&&&\\
&U_{4}\ar@{->}[dr]_{\alpha_{4}}&&U_{5}\ar@{->}[dl]^{\alpha_{5}}&&&&&&\mathbb{P}(1,2,5),\\
&&X\ar@{-->}[rrrrrrru]_{\psi}&&&&&&&}
$$
where $\psi$ is a projection, $\alpha_{4}$ is the weighted blow up
of $P_{4}$ with weights $(1,1,5)$, $\alpha_{5}$ is the weighted
blow up of $P_{5}$ with weights $(1,2,5)$, $\beta_{4}$ is the
weighted blow up with weights $(1,1,5)$ of the proper transform of
$P_{4}$ on $U_{5}$, $\beta_{5}$ is the weighted blow up with
weights $(1,2,5)$ of the proper transform of the point $P_{5}$ on
the variety $U_{4}$, and $\eta$ is an elliptic fibration. There is
a commutative diagram
$$
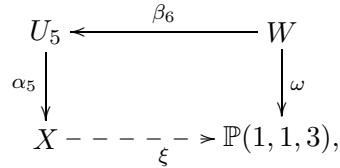
\xymatrix{
&U_{5}\ar@{->}[d]_{\alpha_{5}}&&W\ar@{->}[ll]_{\beta_{6}}\ar@{->}[d]^{\omega}&\\%
&X\ar@{-->}[rr]_{\xi}&&\mathbb{P}(1,1,3),&}
$$
where $\xi$ is a projection, $\beta_{6}$ is the weighted blow up
with weights $(1,2,3)$ of the singular point of $U_{5}$ that is a
singularity of type $\frac{1}{5}(1,2,3)$  contained in the
$\alpha_{5}$-exceptional divisor, and $\omega$ is an elliptic
fibration.

\begin{proposition}
\label{proposition:n-44} Either there is a commutative diagram
\begin{equation}
\label{equation:n-44-commutative-diagram-I} \xymatrix{
&X\ar@{-->}[d]_{\psi}\ar@{-->}[rr]^{\rho}&&V\ar@{->}[d]^{\nu}&\\%
&\mathbb{P}(1,2,5)\ar@{-->}[rr]_{\phi}&&\mathbb{P}^{2},&}
\end{equation}
or there is a commutative diagram
\begin{equation}
\label{equation:n-44-commutative-diagram-II} \xymatrix{
&X\ar@{-->}[d]_{\xi}\ar@{-->}[rr]^{\zeta}&&X\ar@{-->}[rr]^{\rho}&&V\ar@{->}[d]^{\nu}&\\%
&\mathbb{P}(1,1,3)\ar@{-->}[rrrr]_{\sigma}&&&&\mathbb{P}^{2}&}
\end{equation}
where $\phi$, $\zeta$ and $\sigma$ are birational maps.
\end{proposition}

Let us prove Proposition~\ref{proposition:n-44}. It follows from
Lemma~\ref{lemma:Ryder} and
Proposition~\ref{proposition:singular-points}~that
$$
\varnothing\ne\mathbb{CS}\Big(X,
\frac{1}{k}\mathcal{M}\Big)\subseteq\big\{P_{4},P_{5}\big\},
$$
but the proof of Lemma~\ref{lemma:n-23-P6-is-a-center} implies
that the set $\mathbb{CS}(X, \frac{1}{k}\mathcal{M})$ contains the
point $P_{5}$.

The existence of the commutative
diagram~\ref{equation:n-44-commutative-diagram-I} easily follows
from Theorem~\ref{theorem:Kawamata} in the case when
$\mathbb{CS}(X, \frac{1}{k}\mathcal{M})=\{P_{4},P_{5}\}$. Thus, to
conclude the proof of Proposition~\ref{proposition:n-44}, we may
assume that the set $\mathbb{CS}(X, \frac{1}{k}\mathcal{M})$
consists of the point $P_{5}$.

In the rest of the section we prove the existence of the
commutative diagram~\ref{equation:n-44-commutative-diagram-II}.

Let $\mathcal{D}_{5}$ be the proper transform of $\mathcal{M}$ on
$U_{5}$. Then $\mathcal{D}_{5}\sim_{\mathbb{Q}}-kK_{U_{5}}$ by
Theorem~\ref{theorem:Kawamata}, which implies that the set
$\mathbb{CS}(U_{5}, \frac{1}{k}\mathcal{D}_{5})$ is not empty by
Lemma~\ref{lemma:Noether-Fano}. Let $G$ be the exceptional divisor
of the morphism $\alpha_{5}$, and $\bar{P}_{6}$ and $\bar{P}_{7}$
are the singular points of $G$ that are quotient singularities of
types $\frac{1}{5}(1,2,3)$ and $\frac{1}{2}(1,1,1)$ on the variety
$U_{5}$ respectively.

In the case when the set $\mathbb{CS}(U_{5},
\frac{1}{k}\mathcal{D}_{5})$ contains the point $\bar{P}_{6}$, the
existence of the commutative
diagram~\ref{equation:n-44-commutative-diagram-II} follows from
Theorem~\ref{theorem:Kawamata}. Therefore, to conclude the proof
of Proposition~\ref{proposition:n-44}, we may assume that the set
$\mathbb{CS}(U_{5}, \frac{1}{k}\mathcal{D}_{5})$ contains the
point $\bar{P}_{7}$ by Lemma~\ref{lemma:Cheltsov-Kawamata}.

\begin{remark}
\label{remark:n-44-pencil-5K} The linear system $|-5K_{U_{5}}|$ is
a proper transform of $|-5K_{X}|$, and the base locus of the
linear system $|-5K_{U_{5}}|$ consists of the irreducible curve
that is the fiber of the rational map $\psi\circ\alpha_{5}$
passing through the point $\bar{P}_{7}$.
\end{remark}

Let $\pi\colon U\to U_{5}$ be the weighted blow up of
$\bar{P}_{7}$ with weights $(1,1,1)$, $F$ be the exceptional
divisor of $\pi$, $\mathcal{D}$ be the proper transform of
$\mathcal{M}$ on the variety $U$, and $\mathcal{H}$ be the proper
transform of the linear system $|-5K_{U_{5}}|$ on the variety $U$.
Then $\mathcal{D}\sim_{\mathbb{Q}}-kK_{U}$ by
Theorem~\ref{theorem:Kawamata}, but
$$
\mathcal{H}\sim_{\mathbb{Q}}\pi^{*}\big(-5K_{U_{5}}\big)-\frac{1}{2}F,
$$
and the base locus of $\mathcal{H}$ consists of the irreducible
curve $Z$ such that $\alpha_{5}\circ\pi(Z)$ is the unique curve in
the base locus of the linear system $|-5K_{X}|$.

Let $S$ be a general surface of the linear system $\mathcal{H}$.
Then the equality $S\cdot Z=1/3$ holds, which implies that the
divisor $\pi^{*}(-10K_{U_{5}})-F$ is nef. Let $D_{1}$ and $D_{2}$
be general surfaces in $\mathcal{D}$. Then
$$
-\frac{2k^{2}}{3}=\Big(\pi^{*}\big(-10K_{U_{5}}\big)-F\Big)\cdot
\Big(\pi^{*}\big(-kK_{U_{5}}\big)-\frac{k}{2}F\Big)^{2}
=\Big(\pi^{*}\big(-10K_{U_{5}}\big)-F\Big)\cdot D_{1}\cdot
D_{2}\geqslant 0,
$$
which is a contradiction. The claim of
Proposition~\ref{proposition:n-44} is proved.

\section{Case $n=47$, hypersurface of degree $21$ in $\mathbb{P}(1,1,5,7,8)$.}
\label{section:n-47}

We use the notations and assumptions of
Section~\ref{section:start}. Let $n=47$. Then
$X$~is~a~general~hyper\-surface in $\mathbb{P}(1,1,5,7,8)$ of
degree $21$, whose singularities consist of the point $P_{1}$ that
is a singularity of type $\frac{1}{5}(1,2,3)$, and the point
$P_{2}$ that is a singularity of type $\frac{1}{8}(1,1,7)$.

There is a commutative diagram
$$
\xymatrix{
&U\ar@{->}[d]_{\alpha}&&W\ar@{->}[ll]_{\beta}&&Y\ar@{->}[ll]_{\gamma}\ar@{->}[d]^{\eta}&\\%
&X\ar@{-->}[rrrr]_{\psi}&&&&\mathbb{P}(1,1,5),&}
$$
where $\alpha$ is the weighted blow up of the point $P_{2}$ with
weights $(1,1,7)$, $\beta$ is the weighted blow up with weights
$(1,1,6)$ of the singular point of the variety $U$ that is a
quotient singularity of type $\frac{1}{7}(1,1,6)$, $\gamma$ is the
weighted blow up with weights $(1,1,5)$ of the singular point $W$
that is a quotient singularity of type $\frac{1}{6}(1,1,5)$, and
$\eta$ is and elliptic fibration.

\begin{proposition}
\label{proposition:n-47} The claim of Theorem~\ref{theorem:main}
holds for $n=47$.
\end{proposition}

In the rest of the section we prove
Proposition~\ref{proposition:n-47}. It follows from
Theorem~\ref{theorem:smooth-points}, Lemma~\ref{lemma:Ryder} and
Proposition~\ref{proposition:singular-points} that $\mathbb{CS}(X,
\frac{1}{k}\mathcal{M})\subseteq\{P_{1}, P_{2}\}$.

The hypersurface $X$ can be given by the equation
$$
w^{2}z+\sum_{i=0}^{2}wz^{i}g_{13-5i}(x,y,t)+\sum_{i=0}^{3}z^{i}g_{21-5i}(x,y,t)=0\subset\mathrm{Proj}\Big(\mathbb{C}[x,y,z,t,w]\Big),
$$
where $\mathrm{wt}(x)=1$, $\mathrm{wt}(y)=1$, $\mathrm{wt}(z)=5$,
$\mathrm{wt}(t)=7$, $\mathrm{wt}(w)=8$, and $g_{i}(x,y,t)$ is a
quasihomogeneous polynomial of degree $i$.

\begin{lemma}
\label{lemma:n-47-P2} The set $\mathbb{CS}(X,
\frac{1}{k}\mathcal{M})$ contains the point $P_{2}$.
\end{lemma}

\begin{proof}
Suppose that $\mathbb{CS}(X, \frac{1}{k}\mathcal{M})$ does not
contains $P_{2}$. Then $\mathbb{CS}(X,
\frac{1}{k}\mathcal{M})=\{P_{1}\}$.

Let $\pi\colon Z\to X$ be the weighted blow up of $P_{1}$ with
weights $(1,2,3)$, $E$ be the exceptional divisor of  $\pi$, and
$\mathcal{B}$ be the proper transform of $\mathcal{M}$ on $Z$.
Then $E\cong\mathbb{P}(1,2,35)$ and
$\mathcal{B}\sim_{\mathbb{Q}}-kK_{Z}$.

Let $\bar{P}_{3}$ and $\bar{P}_{4}$ be the singular points of the
variety $Z$ contained in $E$ that are singularities of types
$\frac{1}{2}(1,1,1)$ and $\frac{1}{3}(1,1,2)$ respectively. Then
the proof of Proposition~\ref{proposition:n-29} implies that the
singularities of the log pair $(Z, \frac{1}{k}\mathcal{B})$ are
terminal. However, the divisor $-K_{Z}$ is not nef.

The base locus of the pencil $|-K_{Z}|$ consists of the
irreducible curves $C$ and $L$ such that the curve $\pi(C)$ is cut
out by $x=y=0$, the curve $L$ is contained in the divisor $E$, the
curve $L$ is contained in
$|\mathcal{O}_{\mathbb{P}(1,\,2,\,3)}(1)|$, the inequalities
$-K_{Z}\cdot C<0$ and $-K_{Z}\cdot L>0$ hold.

It follows from \cite{Sho93} that the antiflip $\zeta\colon
Z\dasharrow \bar{Z}$ in the curve $C$ exists, and  $-K_{\bar{Z}}$
is nef.

Let $\mathcal{P}$ be the proper transform of the linear system
$\mathcal{M}$ on the variety $\bar{Z}$. Then the singularities of
the log pair $(\bar{Z}, \frac{1}{k}\mathcal{P})$ are terminal,
because the singularities of the log pair $(Z,
\frac{1}{k}\mathcal{B})$ are terminal, and the antiflip $\zeta$ is
a log flop with respect to the log pair $(Z,
\frac{1}{k}\mathcal{B})$.

One can easily check that the rational functions $y/x$,
$zy/x^{6}$, $ty/x^{8}$ and $yw/x^{9}$ are contained in the linear
system $|-aK_{Z}|$, where $a=1$, $6$, $8$ and $9$ respectively.
Therefore, the complete linear system $|-72K_{Z}|$ induces the
birational map $\chi\colon Z\dasharrow \bar{X}$ such that
$\bar{X}$ is a hypersurface of degree $24$ in
$\mathbb{P}(1,1,6,8,9)$. Hence, the divisor $-K_{\bar{Z}}$ is big,
which contradicts Lemma~\ref{lemma:Noether-Fano}.
\end{proof}

Let $G$ be the exceptional divisor of the morphism $\alpha$,
$\mathcal{D}$ be the proper transform of the linear system
$\mathcal{M}$ on the variety $U$, $\bar{P}_{1}$ be the proper
transform of $P_{1}$ on $U$, and $\bar{P}_{5}$ be the singular
point of the variety $U$ that is contained in $G$. Then
$\mathcal{D}\sim_{\mathbb{Q}}-kK_{U}$ by
Theorem~\ref{theorem:Kawamata}.

\begin{lemma}
\label{lemma:n-47-P5} The set $\mathbb{CS}(U,
\frac{1}{k}\mathcal{D})$ contains the point $\bar{P}_{5}$.
\end{lemma}

\begin{proof}
Suppose that $\bar{P}_{5}\not\in\mathbb{CS}(U,
\frac{1}{k}\mathcal{D})$. Then $\mathbb{CS}(U,
\frac{1}{k}\mathcal{D})=\{\bar{P}_{1}\}$ by
Lemmas~\ref{lemma:Noether-Fano} and \ref{lemma:Cheltsov-Kawamata}.

Let $\pi\colon Z\to U$ be the weighted blow up of the point
$\bar{P}_{1}$ with weights $(1,2,3)$, $E$ be the
ex\-cep\-ti\-o\-nal divisor of the morphism $\pi$, and $\bar{G}$
and $\mathcal{B}$ be the proper transforms of $G$ and
$\mathcal{M}$ on the variety $Z$ respectively. Then
$\mathcal{B}\sim_{\mathbb{Q}}-kK_{Z}$ by
Theorem~\ref{theorem:Kawamata}, but the proof of
Lemma~\ref{lemma:n-47-P2} implies that the singularities of the
log pair $(Z, \frac{1}{k}\mathcal{B})$ are terminal.

Let $S_{x}$, $S_{y}$, $S_{z}$, $S_{t}$ and $S_{w}$ be the proper
transforms on the variety $Z$ of the surfaces that are cut on $X$
by the equations $x=0$, $y=0$, $z=0$, $t=0$ and $w=0$
respectively. Then
\begin{equation}
\label{equation:n-47-rational-equivalences-II} \left\{\aligned
&S_{x}\sim_{\mathbb{Q}}(\alpha\circ\pi)^{*}\big(-K_{X}\big)-\frac{1}{5}E-\frac{1}{8}\bar{G},\\
&S_{y}\sim_{\mathbb{Q}}(\alpha\circ\pi)^{*}\big(-K_{X}\big)-\frac{6}{5}E-\frac{1}{8}\bar{G},\\
&S_{z}\sim_{\mathbb{Q}}(\alpha\circ\pi)^{*}\big(-5K_{X}\big)-\frac{1}{5}E-\frac{13}{8}\bar{G},\\
&S_{t}\sim_{\mathbb{Q}}(\alpha\circ\pi)^{*}\big(-7K_{X}\big)-\frac{2}{5}E-\frac{7}{8}\bar{G},\\
&S_{w}\sim_{\mathbb{Q}}(\alpha\circ\pi)^{*}\big(-8K_{X}\big)-\frac{3}{5}E.\\
\endaligned\right.
\end{equation}

The base locus of the pencil $|-K_{Z}|$ consists of the
irreducible curves $C$ and $L$ such that the curve
$\alpha\circ\pi(C)$ is cut by the equation $x=y=0$ on the
hypersurface $X$, the curve $L$ is contained in the divisor $E$,
and the curve $L$ is the unique curve of the linear system
$|\mathcal{O}_{\mathbb{P}(1,\,2,\,3)}(1)|$.

It follows from \cite{Sho93} that there is the antiflip
$\zeta\colon Z\dasharrow \bar{Z}$ in $C$ such that $-K_{\bar{Z}}$
is nef.

Let $\mathcal{P}$ be the proper transform of the linear system
$\mathcal{M}$ on the variety $\bar{Z}$. Then the singularities of
the log pair $(\bar{Z}, \frac{1}{k}\mathcal{P})$ are terminal.

The equivalences~\ref{equation:n-47-rational-equivalences-II}
imply that the functions $y/x$, $zy/x^{6}$, $ty/x^{8}$ and
$wzy^{2}/x^{15}$ are contained in the complete linear system
$|aS_{x}|$, where $a=1$, $6$, $8$ and $15$ respectively. Hence,
the linear system $|-120K_{Z}|$ induces the birational map
$\chi\colon Z\dasharrow \bar{X}$ such that $\bar{X}$ is a
hypersurface of degree $30$ in $\mathbb{P}(1,1,6,8,15)$, which
implies that $-K_{\bar{Z}}$ is big, which contradicts
Lemma~\ref{lemma:Noether-Fano}.
\end{proof}

\begin{remark}
\label{remark:n-47-P5} It follows from Lemma~\ref{lemma:curves}
that $\mathbb{CS}(U, \frac{1}{k}\mathcal{D})=\{\bar{P}_{5}\}$.
\end{remark}

Let $\mathcal{H}$ be the proper transform of the linear system
$\mathcal{M}$ on the variety $W$, $F$ be the exceptional divisor
of $\beta$, $\tilde{P}_{1}$ be the proper transform of $P_{1}$ on
$W$, and $\tilde{P}_{6}$ be the singular point of  $W$ that is
contained in $F$. Then $\mathcal{H}\sim_{\mathbb{Q}}-kK_{W}$ by
Theorem~\ref{theorem:Kawamata}, but $\mathbb{CS}(W,
\frac{1}{k}\mathcal{H})\ne\varnothing$.

\begin{lemma}
\label{lemma:n-47-points-P6} Suppose that
$\tilde{P}_{6}\in\mathbb{CS}(W, \frac{1}{k}\mathcal{H})$. Then
there is a commutative diagram
\begin{equation}
\label{equation:n-47-commutative-diagramm} \xymatrix{
&X\ar@{-->}[d]_{\psi}\ar@{-->}[rr]^{\rho}&&V\ar@{->}[d]^{\nu}&\\%
&\mathbb{P}(1,1,5)\ar@{-->}[rr]_{\zeta}&&\mathbb{P}^{2},&}
\end{equation}
where $\zeta$ is a birational map.
\end{lemma}

\begin{proof}
The existence of the
diagram~\ref{equation:n-47-commutative-diagramm} follows from
Theorem~\ref{theorem:Kawamata}.
\end{proof}

We may assume that the set $\mathbb{CS}(W,
\frac{1}{k}\mathcal{H})$ consists of the point $\tilde{P}_{1}$ by
Lemma~\ref{lemma:Cheltsov-Kawamata}.

Let $\pi\colon Z\to W$ be the weighted blow up of $\tilde{P}_{1}$
with weights $(1,2,3)$, $E$ be the exceptional divisor of $\pi$,
$\mathcal{B}$ be the proper transform of $\mathcal{M}$ on $Z$, and
$\tilde{G}$ and $\tilde{F}$ be the proper transforms on the
variety $Z$ of the surfaces $G$ and $F$ respectively. Then it
follows from the proof of Lemma~\ref{lemma:n-47-P2} that the
singularities of the log pair $(Z, \frac{1}{k}\mathcal{B})$ are
terminal, but $\mathcal{B}\sim_{\mathbb{Q}}-kK_{Z}$ by
Theorem~\ref{theorem:Kawamata}.

Let $S_{x}$, $S_{y}$, $S_{z}$, $S_{t}$ and $S_{w}$ be proper
transforms on the variety $Z$ of the surfaces that are cut on the
variety $X$ by the equations $x=0$, $y=0$, $z=0$, $t=0$ and $w=0$
respectively. Then
\begin{equation}
\label{equation:n-47-rational-equivalences-III} \left\{\aligned
&S_{x}\sim_{\mathbb{Q}}(\alpha\circ\beta\circ\pi)^{*}\big(-K_{X}\big)-\frac{1}{5}E-\frac{1}{8}\tilde{G}-\frac{1}{4}\tilde{F},\\
&S_{y}\sim_{\mathbb{Q}}(\alpha\circ\beta\circ\pi)^{*}\big(-K_{X}\big)-\frac{6}{5}E-\frac{1}{8}\tilde{G}-\frac{1}{4}\tilde{F},\\
&S_{z}\sim_{\mathbb{Q}}(\alpha\circ\beta\circ\pi)^{*}\big(-5K_{X}\big)-\frac{1}{5}E-\frac{13}{8}\tilde{G}-\frac{9}{4}\tilde{F},\\
&S_{t}\sim_{\mathbb{Q}}(\alpha\circ\beta\circ\pi)^{*}\big(-7K_{X}\big)-\frac{2}{5}E-\frac{7}{8}\tilde{G}-\frac{3}{4}\tilde{F},\\
&S_{w}\sim_{\mathbb{Q}}(\alpha\circ\beta\circ\pi)^{*}\big(-8K_{X}\big)-\frac{3}{5}E.\\
\endaligned\right.
\end{equation}

The equivalences~\ref{equation:n-47-rational-equivalences-III}
imply that the functions $y/x$, $zy/x^{6}$, $tzy^{2}/x^{14}$ and
$wz^{2}y^{3}/x^{21}$ are contained in the linear systems
$|S_{x}|$, $|6S_{x}|$, $|14S_{x}|$ and $|21S_{x}|$ respectively.
Therefore, the complete linear system $|-42K_{Z}|$ induces the
birational map $\chi\colon Z\dasharrow \bar{X}$ such that
$\bar{X}$ is a threefold\footnote{The proofs of
Lemmas~\ref{lemma:n-47-P2} and \ref{lemma:n-47-P5} give the
birational transformations of the hypersurface $X$ into
hypersurfaces in $\mathbb{P}(1,1,6,8,9)$ and
$\mathbb{P}(1,1,6,8,15)$ of degrees $24$ and $30$ respectively.
The anticanonical models of the varieties $U$ and $W$ are
hypersurfaces in $\mathbb{P}(1,1,5,7,13)$ and
$\mathbb{P}(1,1,5,12,18)$ of degrees $26$ and $36$ respectively,
and the threefold $\bar{X}$ is a hypersurface in
$\mathbb{P}(1,1,6,14,21)$ of degree $42$. Up to to the action of
$\mathrm{Bir}(X)$, there are no other non-trivial birational
transformations of $X$ into Fano threefolds with canonical
singularities.}.

The base locus of the linear system $|-42K_{Z}|$ consists of the
irreducible curve $C$ such that the curve
$\alpha\circ\beta\circ\pi(C)$ is cut on $X$ by the equations
$x=y=0$. Therefore, the existence of the antiflip $\zeta\colon
Z\dasharrow \bar{Z}$ in the curve $C$ follows from \cite{Sho93},
which implies that $-K_{\bar{Z}}$ is nef and big.

The rational map $\zeta$ is a log flip with respect to the log
pair $(Z, \frac{1}{k}\mathcal{B})$. Therefore, we see that the
singularities of the mobile log pair $(\bar{Z},
\frac{1}{k}\mathcal{P})$ are terminal, where $\mathcal{P}$ is the
proper transform of the linear system $\mathcal{M}$ on the variety
$Z$, which is impossible by Lemma~\ref{lemma:Noether-Fano}.

The claim of Proposition~\ref{proposition:n-47} is proved.

\section{Case  $n=48$, hypersurface of degree $21$ in $\mathbb{P}(1,2,3,7,9)$.}
\label{section:n-48}

We use the notations and assumptions of
Section~\ref{section:start}. Let $n=48$. Then $X$ is
a~general~hyper\-surface in $\mathbb{P}(1,2,3,7,9)$ of degree
$21$, whose singularities consist of the point $P_{1}$ that is a
singularity of type $\frac{1}{2}(1,1,1)$, the points $P_{2}$ and
$P_{3}$ that are singularity of type $\frac{1}{3}(1,1,2)$, and the
point $P_{4}$ that is a singularity of type $\frac{1}{9}(1,2,7)$.
There is a commutative diagram
$$
\xymatrix{
&U\ar@{->}[d]_{\alpha}&&W\ar@{->}[ll]_{\beta}&&Y\ar@{->}[ll]_{\gamma}\ar@{->}[d]^{\eta}&\\%
&X\ar@{-->}[rrrr]_{\psi}&&&&\mathbb{P}(1,2,3),&}
$$
where $\psi$ is a projection, $\alpha$ is the weighted blow up of
$P_{4}$ with weights $(1,2,7)$, $\beta$ is the weighted blow up
with weights $(1,2,5)$ of the singular point of the variety $U$
that is a quotient singularity of type $\frac{1}{7}(1,2,5)$
contained in the exceptional divisor of $\alpha$, $\gamma$ is the
weighted blow up with weights $(1,2,3)$ of the singular point of
$W$ that is a singularity of type $\frac{1}{5}(1,2,3)$  contained
in the exceptional divisor of $\beta$, and $\eta$ is an elliptic
fibration.

\begin{proposition}
\label{proposition:n-48} The claim of Theorem~\ref{theorem:main}
holds for $n=48$.
\end{proposition}

In the rest of the section we prove
Proposition~\ref{proposition:n-48}. It follows from
Theorem~\ref{theorem:smooth-points}, Lemma~\ref{lemma:Ryder} and
Proposition~\ref{proposition:singular-points} that $\mathbb{CS}(X,
\frac{1}{k}\mathcal{M})=\{P_{4}\}$.

Let $E$ be the $\alpha$-exceptional divisor, $\mathcal{D}$ be the
proper transform of $\mathcal{M}$ on $U$, and $P_{5}$ and
$P_{6}$~be the singular points of $U$ contained in $E$ that are
singularities of types $\frac{1}{2}(1,1,1)$ and
$\frac{1}{7}(1,2,5)$ res\-pec\-tive\-ly. Then
$E\cong\mathbb{P}(1,2,7)$, but
$\mathcal{D}\sim_{\mathbb{Q}}-kK_{U}$ by
Theorem~\ref{theorem:Kawamata}.

It follows from Lemmas~\ref{lemma:Noether-Fano},
\ref{lemma:Cheltsov-Kawamata} and \ref{lemma:curves} that
$\mathbb{CS}(U, \frac{1}{k}\mathcal{D})\subseteq\{P_{5}, P_{6}\}$.

\begin{lemma}
\label{lemma:n-48-points-P5} The set $\mathbb{CS}(U,
\frac{1}{k}\mathcal{D})$ does not contain the point $P_{5}$.
\end{lemma}

\begin{proof}
Suppose that the set $\mathbb{CS}(U, \frac{1}{k}\mathcal{D})$
contains the point $P_{5}$. Let $\pi\colon Z\to U$ be the weighted
blow up of $P_{5}$ with weights $(1,1,1)$, $G$ be the exceptional
divisor of  $\pi$, and $\mathcal{B}$ and $\mathcal{P}$ be the
proper transforms of the linear systems $\mathcal{M}$ and
$|-7K_{X}|$ on the variety $Z$ respectively. Then
$$
\mathcal{B}\sim_{\mathbb{Q}}-kK_{Z}\sim_{\mathbb{Q}}(\alpha\circ\pi)^{*}\big(-kK_{X}\big)-\frac{k}{9}\pi^{*}(E)-\frac{k}{2}G
$$
by Theorem~\ref{theorem:Kawamata}, but the base locus of the
linear system $\mathcal{P}$ does not contain curves. Let $H$ be a
general divisor of the linear system $\mathcal{P}$. Then the
divisor $H$ is numerically effective, but
$$
H\cdot B_{1}\cdot
B_{2}=\Big((\alpha\circ\pi)^{*}\big(-kK_{X}\big)-\frac{k}{9}\pi^{*}(E)-\frac{k}{2}G\Big)^{2}
\Big((\alpha\circ\pi)^{*}\big(-7K_{X}\big)-\frac{7}{9}\pi^{*}(E)-\frac{1}{2}G\Big)=-\frac{1}{6}k^{2},
$$
where $B_{1}$ and $B_{2}$ are general surfaces of the linear
system $\mathcal{B}$, which is a contradiction.
\end{proof}

Hence, the set $\mathbb{CS}(U, \frac{1}{k}\mathcal{D})$ consists
of the point $P_{6}$.

Let $F$ be the $\beta$-exceptional divisor, $\mathcal{H}$ be the
proper transform of $\mathcal{M}$ on $W$, and $P_{7}$ and
$P_{8}$~be the singular points of $W$ contained in $F$ that are
singularities of types $\frac{1}{2}(1,1,1)$ and
$\frac{1}{5}(1,2,3)$  respectively. Then
$F\cong\mathbb{P}(1,2,5)$, but
$\mathcal{H}\sim_{\mathbb{Q}}-kK_{W}$ by
Theorem~\ref{theorem:Kawamata}.

It follows from Lemmas~\ref{lemma:Noether-Fano},
\ref{lemma:Cheltsov-Kawamata} and \ref{lemma:curves} that
$\mathbb{CS}(W, \frac{1}{k}\mathcal{H})\subseteq\{P_{7}, P_{8}\}$.

\begin{lemma}
\label{lemma:n-48-points-P8} Suppose that $P_{8}\in\mathbb{CS}(W,
\frac{1}{k}\mathcal{H})$. Then there is a commutative diagram
\begin{equation}
\label{equation:n-48-commutative-diagramm} \xymatrix{
&X\ar@{-->}[d]_{\psi}\ar@{-->}[rr]^{\rho}&&V\ar@{->}[d]^{\nu}&\\%
&\mathbb{P}(1,2,3)\ar@{-->}[rr]_{\zeta}&&\mathbb{P}^{2},&}
\end{equation}
where $\zeta$ is a birational map.
\end{lemma}

\begin{proof}
Let $S$ be the proper transform on the variety $Y$ of a
sufficiently general surface of the linear system $\mathcal{M}$,
and $\Gamma$ be a  general fiber of  $\eta$. Then
$S\sim_{\mathbb{Q}}-kK_{Y}$ by Theorem~\ref{theorem:Kawamata},
which implies that $S\cdot \Gamma=0$. Therefore, the surface $S$
lies in the fibers of the fibration $\eta$, which implies the
existence of the commutative
diagram~\ref{equation:n-48-commutative-diagramm}.
\end{proof}

To conclude the proof of Proposition~\ref{proposition:n-48}, we
may assume that $P_{7}\in\mathbb{CS}(W, \frac{1}{k}\mathcal{H})$.

Let $\pi\colon Z\to W$ be the weighted blow up of $P_{7}$ with
weights $(1,1,1)$, $G$ be the exceptional  divisor of $\pi$, and
$\mathcal{B}$ be the proper transform of $\mathcal{M}$ on $Z$.
Then $\mathcal{B}\sim_{\mathbb{Q}}-kK_{Z}$ by
Theorem~\ref{theorem:Kawamata}.

The hypersurface $X$ can be given by the equation
$$
w^{2}z+wf_{12}(x,y,z,t)+f_{21}(x,y,z,t)=0\subset\mathbb{P}(1,2,3,7,9)\cong\mathrm{Proj}\Big(\mathbb{C}[x,y,z,t,w]\Big),
$$
where $\mathrm{wt}(x)=1$, $\mathrm{wt}(y)=2$, $\mathrm{wt}(z)=3$,
$\mathrm{wt}(t)=7$, $\mathrm{wt}(w)=9$, and $f_{i}(x,y,z,t)$ is a
quasihomogeneous polynomial of degree $i$.

Let $\mathcal{P}$ be the proper transform on the variety $Z$ of
the pencil of surfaces that are cut on the hypersurface $X$ by the
equations $\lambda x^{3}+\mu z=0$, where
$(\lambda,\mu)\in\mathbb{P}^{1}$. Then the base locus of the
pencil $\mathcal{P}$ consists of the irreducible curves $C$,
$L_{1}$ and $L_{2}$ such that $\alpha\circ\beta\circ\pi(C)$ is the
curve that is cut on the hypersurface $X$ by the equations
$x=z=0$, the curve $\beta\circ\pi(L_{1})$ is contained in the
exceptional divisor $E$, the curve $\beta\circ\pi(L_{1})$ is the
unique curve in the base locus of the linear system
$|\mathcal{O}_{\mathbb{P}(1,\,2,\,7)}(1)|$, the curve $\pi(L_{2})$
is contained in $F$, and the curve $\pi(L_{2})$ is the unique
curve of the linear system
$|\mathcal{O}_{\mathbb{P}(1,\,2,\,5)}(1)|$.

Let $D$ be a general surface of the pencil $\mathcal{P}$,
$\bar{E}$ and $\bar{F}$ be the proper transforms of the
exceptional divisors $E$ and $F$ on the variety $Z$ respectively,
and $S$ be the proper transform on $Z$ of the surface that is cut
on the hypersurface $X$ by the equation $x=0$. Then
$$
S\cdot D=C+L_{1}+L_{2},\ \bar{E}\cdot D=3L_{1},\ \bar{F}\cdot D=3L_{2}, %
$$
the surface $D$ is normal, and
\begin{equation}
\label{equation:n-48-rational-equivalences} \left\{\aligned
&\bar{F}\sim_{\mathbb{Q}}\pi^{*}(F)-\frac{1}{2}G,\\
&\bar{E}\sim_{\mathbb{Q}}(\beta\circ\pi)^{*}(E)-\frac{5}{7}\pi^{*}(F)-\frac{1}{2}G,\\
&D\sim_{\mathbb{Q}}(\alpha\circ\beta\circ\pi)^{*}\big(-3K_{X}\big)-\frac{3}{9}(\beta\circ\pi)^{*}(E)-\frac{3}{7}\pi^{*}(F)-\frac{3}{2}G,\\
&S\sim_{\mathbb{Q}}(\alpha\circ\beta\circ\pi)^{*}\big(-K_{X}\big)-\frac{1}{9}(\beta\circ\pi)^{*}(E)-\frac{1}{7}\pi^{*}(F)-\frac{1}{2}G.\\
\endaligned
\right.
\end{equation}

Consider the curves $C$, $L_{1}$ and $L_{2}$ as divisors on
$\bar{D}$. The
equivalences~\ref{equation:n-48-rational-equivalences} imply that
$$
C\cdot C=L_{1}\cdot L_{1}=-\frac{1}{2},\ L_{2}\cdot L_{2}=-\frac{2}{5},\ C\cdot L_{1}=C\cdot L_{2}=L_{1}\cdot L_{2}=0,%
$$
which implies that the intersection form of the curves $C$,
$L_{1}$ and $L_{2}$ on the normal surface $D$ is negatively
defined. On the other hand, we have
$$
B\vert_{D}\sim_{\mathbb{Q}}-kK_{Z}\vert_{D}\sim_{\mathbb{Q}}kS\vert_{D}\sim_{\mathbb{Q}}kC+kL_{1}+kL_{2},
$$
where $B$ is a general surface of the linear system $\mathcal{B}$,
which contradicts Lemmas~\ref{lemma:normal-surface} and
\ref{lemma:Cheltsov-II}.

The claim of Proposition~\ref{proposition:n-48} is proved.

\section{Case $n=49$, hypersurface of degree $21$ in $\mathbb{P}(1,3,5,6,7)$.}
\label{section:n-49}

We use the notations and assumptions of
Section~\ref{section:start}. Let $n=49$. Then $X$ is
a~general~hyper\-surface in $\mathbb{P}(1,3,5,6,7)$ of degree
$21$, whose singularities consist of the points $P_{1}$, $P_{2}$
and $P_{3}$ that are singularities of type $\frac{1}{3}(1,1,2)$,
the point $P_{4}$ that is a singularity of type
$\frac{1}{5}(1,2,3)$, and the point $P_{5}$ that is a singularity
of type $\frac{1}{6}(1,1,5)$.

There is a commutative diagram
$$
\xymatrix{
&&Y\ar@{->}[dl]_{\alpha}\ar@{->}[dr]^{\eta}&&\\%
&X\ar@{-->}[rr]_{\psi}&&\mathbb{P}(1,3,5),&}
$$
where $\psi$ is a projection, $\alpha$ is the blow up of $P_{5}$
with weights $(1,1,5)$, and $\eta$ is an elliptic~fibration.

There is a commutative diagram
$$
\xymatrix{
&&U\ar@{->}[dl]_{\beta}\ar@{->}[dr]^{\omega}&&\\
&X\ar@{-->}[rr]_{\xi}&&\mathbb{P}(1,3,6),&}
$$
where $\xi$ is a projection, $\beta$ is the blow up of $P_{4}$
with weights $(1,2,3)$, and $\omega$ is an elliptic~fibration.

\begin{proposition}
\label{proposition:n-49} Either there is a commutative diagram
\begin{equation}
\label{equation:n-49-commutative-diagram-I} \xymatrix{
&X\ar@{-->}[d]_{\psi}\ar@{-->}[rr]^{\rho}&&V\ar@{->}[d]^{\nu}&\\%
&\mathbb{P}(1,3,5)\ar@{-->}[rr]_{\phi}&&\mathbb{P}^{2},&}
\end{equation}
or there is a commutative diagram
\begin{equation}
\label{equation:n-49-commutative-diagram-II} \xymatrix{
&X\ar@{-->}[d]_{\xi}\ar@{-->}[rr]^{\rho}&&V\ar@{->}[d]^{\nu}&\\%
&\mathbb{P}(1,3,6)\ar@{-->}[rr]_{\sigma}&&\mathbb{P}^{2},&}
\end{equation}
where $\phi$ and $\sigma$ are birational maps.
\end{proposition}

\begin{proof}
It follows from Proposition~\ref{proposition:singular-points} and
Lemma~\ref{lemma:Ryder} that $\mathbb{CS}(X,
\frac{1}{k}\mathcal{M})\subseteq\{P_{4}, P_{5}\}$.

Suppose that the set $\mathbb{CS}(X, \frac{1}{k}\mathcal{M})$
contains point $P_{4}$. Let $\mathcal{D}$ be the proper transform
of the linear system $\mathcal{M}$ on the variety $U$. Then
$\mathcal{D}\sim_{\mathbb{Q}}-kK_{U}$ by
Theorem~\ref{theorem:Kawamata}. Intersecting a general surface of
the linear system $\mathcal{D}$ with a general fiber of $\omega$,
we see that $\mathcal{D}$ lies in the fibers of the elliptic
$\omega$, which implies the existence of the commutative
diagram~\ref{equation:n-49-commutative-diagram-II}.

Similarly, the commutative
diagram~\ref{equation:n-49-commutative-diagram-I} exists when
$P_{4}\in\mathbb{CS}(X, \frac{1}{k}\mathcal{M})$.
\end{proof}

\section{Case $n=51$, hypersurface of degree $22$ in $\mathbb{P}(1,1,4,6,11)$.}
\label{section:n-51}

We use the notations and assumptions of
Section~\ref{section:start}. Let $n=51$. Then $X$ is a
hypersurface of degree $22$ in $\mathbb{P}(1,1,4,6,11)$, the
equality $-K_{X}^{3}=1/12$ holds, and the singularities of $X$
consist of the point $P_{1}$ that is a quotient singularity of
type $\frac{1}{2}(1,1,1)$, the point $P_{2}$ that is a quotient
singularity of type $\frac{1}{3}(1,1,4)$, and the point $P_{3}$
that is a quotient singularity of type $\frac{1}{6}(1,1,5)$.

The hypersurface $X$ is birationally superrigid. There is a
commutative diagram
$$
\xymatrix{
&U\ar@{->}[d]_{\alpha}&&Y\ar@{->}[ll]_{\beta}\ar@{->}[d]^{\eta}&\\%
&X\ar@{-->}[rr]_{\psi}&&\mathbb{P}(1,1,4),&}
$$
where $\psi$ is a projection, $\alpha$ is the weighted blow up of
$P_{3}$ with weights $(1,1,5)$, $\beta$ is the weighted blow up
with weights $(1,1,4)$ of the singular point of the variety $U$
that is contained in the exceptional divisor of $\alpha$, and
$\eta$ is an elliptic fibration. There is a commutative diagram
$$
\xymatrix{
&&W\ar@{->}[dl]_{\gamma}\ar@{->}[dr]^{\omega}&&\\
&X\ar@{-->}[rr]_{\xi}&&\mathbb{P}(1,1,6),&}
$$
where $\xi$ is a projection, $\gamma$ is the blow up of $P_{2}$
with weights $(1,1,3)$, and $\omega$ is an elliptic~fibration.

\begin{proposition}
\label{proposition:n-51} Either there is a commutative diagram
\begin{equation}
\label{equation:n-51-commutative-diagram-I} \xymatrix{
&X\ar@{-->}[d]_{\psi}\ar@{-->}[rr]^{\rho}&&V\ar@{->}[d]^{\nu}&\\%
&\mathbb{P}(1,1,4)\ar@{-->}[rr]_{\phi}&&\mathbb{P}^{2},&}
\end{equation}
or there is a commutative diagram
\begin{equation}
\label{equation:n-51-commutative-diagram-II} \xymatrix{
&X\ar@{-->}[d]_{\xi}\ar@{-->}[rr]^{\rho}&&V\ar@{->}[d]^{\nu}&\\%
&\mathbb{P}(1,1,6)\ar@{-->}[rr]_{\sigma}&&\mathbb{P}^{2},&}
\end{equation}
where $\phi$ and $\sigma$ are birational maps.
\end{proposition}

\begin{proof}
Suppose that $P_{2}$ is contained in $\mathbb{CS}(X,
\frac{1}{k}\mathcal{M})$. Let $\mathcal{D}$ be the proper
transform of $\mathcal{M}$ on the variety $W$. Then the
equi\-va\-lence $\mathcal{D}\sim_{\mathbb{Q}}-kK_{W}$ holds by
Theorem~\ref{theorem:Kawamata}, which implies  the existence of
the commutative
diagram~\ref{equation:n-51-commutative-diagram-II}.

We may assume that $\mathbb{CS}(X,
\frac{1}{k}\mathcal{M})=\{P_{3}\}$ by
Theorem~\ref{theorem:smooth-points},
Proposition~\ref{proposition:singular-points} and
Lemma~\ref{lemma:Ryder}.

Let $\mathcal{B}$ be the proper transform of $\mathcal{M}$ on $U$.
Then $\mathcal{B}\sim_{\mathbb{Q}}-kK_{U}$ by
Theorem~\ref{theorem:Kawamata}, but the anticanonical divisor
$-K_{U}$ is nef and big. It follows from
Lemmas~\ref{lemma:Noether-Fano} and \ref{lemma:Cheltsov-Kawamata}
that the set of centers of canonical singularities $\mathbb{CS}(U,
\frac{1}{k}\mathcal{B})$ contains the singular point of the
variety $U$ that is contained in the exceptional divisor of the
birational morphism $\alpha$.

Let $\mathcal{H}$ be the proper transform of $\mathcal{M}$ on $Y$.
Then $\mathcal{H}\sim_{\mathbb{Q}}-kK_{Y}$ by
Theorem~\ref{theorem:Kawamata}.

Intersecting a general surface of the linear system $\mathcal{H}$
with a general fiber of the fibration $\eta$, we see that the
linear system $\mathcal{H}$ lies in the fibers of the elliptic
fibration $\eta$, which implies the existence of the commutative
diagram~\ref{equation:n-51-commutative-diagram-I}.
\end{proof}

\section{Case $n=56$, hypersurface of degree $24$ in $\mathbb{P}(1,2,3,8,11)$.}
\label{section:n-56}

We use the notations and assumptions of
Section~\ref{section:start}. Let $n=56$. Then $X$
is~a~general~hyper\-surface in $\mathbb{P}(1,2,3,8,11)$ of degree
$24$, whose singularities consist of the points $P_{1}$,
$P_{2}$~and~$P_{3}$~that are singularities of type
$\frac{1}{2}(1,1,1)$, and the point $P_{4}$ that is a singularity
of type $\frac{1}{11}(1,3,8)$.

There is a commutative diagram
$$
\xymatrix{
&U\ar@{->}[d]_{\alpha}&&W\ar@{->}[ll]_{\beta}&&Y\ar@{->}[ll]_{\gamma}\ar@{->}[d]^{\eta}&\\%
&X\ar@{-->}[rrrr]_{\psi}&&&&\mathbb{P}(1,2,3),&}
$$
where $\psi$ is a projection, $\alpha$ is the weighted blow up of
$P_{4}$ with weights $(1,3,8)$, $\beta$ is the weighted blow up
with weights $(1,3,5)$ of the point of $U$ that is a quotient
singularity of type $\frac{1}{8}(1,3,5)$ contained in the
exceptional divisor of $\alpha$, $\gamma$ is the weighted blow up
with weights $(1,2,3)$ of the point of $W$ that is a quotient
singularity of type $\frac{1}{5}(1,2,3)$ contained in the
exceptional divisor of the morphism $\beta$, and $\eta$ is an
elliptic fibration.

\begin{proposition}
\label{proposition:n-56} The claim of Theorem~\ref{theorem:main}
holds for $n=56$.
\end{proposition}

Let us prove Proposition~\ref{proposition:n-56}. It follows from
Proposition~\ref{proposition:singular-points} that $\mathbb{CS}(X,
\frac{1}{k}\mathcal{M})=\{P_{4}\}$.

Let $E$ be the exceptional divisor of the morphism $\alpha$,
$\mathcal{D}$ be the proper transform of the linear system
$\mathcal{M}$ on the variety $U$, and $P_{5}$ and $P_{6}$ be the
singular point of $U$ contained in $E$ that are singularities of
types $\frac{1}{3}(1,1,2)$ and $\frac{1}{8}(1,3,5)$ respectively.
Then $\mathcal{D}\sim_{\mathbb{Q}}-kK_{U}$ by
Theorem~\ref{theorem:Kawamata}.

\begin{lemma}
\label{lemma:n-56-points-P5} The set $\mathbb{CS}(U,
\frac{1}{k}\mathcal{D})$ does not contain the point $P_{5}$.
\end{lemma}

\begin{proof}
Suppose that the set $\mathbb{CS}(U, \frac{1}{k}\mathcal{D})$
contains the point $P_{5}$. Let $\pi\colon Z\to U$ be the weighted
blow up of $P_{5}$ with weights $(1,1,2)$, $G$ be the exceptional
divisor of  $\pi$, and $\mathcal{B}$ and $\mathcal{P}$ be the
proper transforms of $\mathcal{M}$ and $|-8K_{X}|$ on the variety
$Z$ respectively. Then $\mathcal{B}\sim_{\mathbb{Q}}-kK_{Z}$, but
the base locus of the linear system $\mathcal{P}$ does not contain
curves.

Let $H$ be a general surface in $\mathcal{P}$. Then the divisor
$H$ is nef and big. In particular, we have
$$
H\cdot B_{1}\cdot
B_{2}=\Big((\alpha\circ\pi)^{*}\big(-kK_{X}\big)-\frac{k}{11}\pi^{*}(E)-\frac{k}{3}G\Big)^{2}
\Big((\alpha\circ\pi)^{*}\big(-8K_{X}\big)-\frac{8}{11}\pi^{*}(E)-\frac{2}{3}G\Big)=0,
$$
where $B_{1}$ and $B_{2}$ are general surfaces in $\mathcal{B}$,
which contradicts Corollary~\ref{corollary:Cheltsov}.
\end{proof}

Hence, the set $\mathbb{CS}(U, \frac{1}{k}\mathcal{D})$ consists
of the point $P_{6}$ by Lemma~\ref{lemma:Cheltsov-Kawamata}.

Let $F$ be the exceptional divisor of the morphism $\beta$,
$\mathcal{H}$ be the proper transform of the linear system
$\mathcal{M}$ on the variety $W$, and $P_{7}$ and $P_{8}$ be the
singular points of $W$ that are singularities of types
$\frac{1}{3}(1,1,2)$ and $\frac{1}{5}(1,2,3)$ contained in $F$
respectively. Then $\mathcal{H}\sim_{\mathbb{Q}}-kK_{W}$ by
Theorem~\ref{theorem:Kawamata}.

The proof of Lemma~\ref{lemma:n-48-points-P8} implies the
existence of the commutative diagram
$$
\xymatrix{
&X\ar@{-->}[d]_{\psi}\ar@{-->}[rr]^{\rho}&&V\ar@{->}[d]^{\nu}&\\%
&\mathbb{P}(1,2,3)\ar@{-->}[rr]_{\zeta}&&\mathbb{P}^{2}&}
$$
in the case when $P_{8}\in\mathbb{CS}(W, \frac{1}{k}\mathcal{H})$,
where $\zeta$ is a birational map.

We may assume that $\mathbb{CS}(W,
\frac{1}{k}\mathcal{H})=\{P_{7}\}$ due to
Lemmas~\ref{lemma:Noether-Fano} and \ref{lemma:Cheltsov-Kawamata}.

Let $\pi\colon Z\to W$ be the weighted blow up of $P_{7}$ with
weights $(1,1,2)$, $G$ be the exceptional  divisor of $\pi$, and
$\mathcal{B}$ be the proper transform of $\mathcal{M}$ on  $Z$.
Then $\mathcal{B}\sim_{\mathbb{Q}}-kK_{Z}$ by
Theorem~\ref{theorem:Kawamata}.

\begin{lemma}
\label{lemma:n-56-terminal-singularities} The singularities of the
log pair $(Z, \frac{1}{k}\mathcal{B})$ are terminal.
\end{lemma}

\begin{proof}
Suppose that $\mathbb{CS}(Z,
\frac{1}{k}\mathcal{B})\ne\varnothing$. Let $P_{9}$ be the
singular point of $G$. Then $P_{9}$ is a singularity of type
$\frac{1}{2}(1,1,1)$ on $Z$, the set $\mathbb{CS}(Z,
\frac{1}{k}\mathcal{B})$ contains $P_{9}$ by
Lemma~\ref{lemma:Cheltsov-Kawamata}, and
$G\cong\mathbb{P}(1,1,2)$.

Let $\bar{\pi}:\bar{Z}\to Z$ be the weighted blow up of the point
$P_{9}$ with weights $(1,1,1)$, and $\bar{G}$ be the exceptional
divisor of $\bar{\pi}$. Take any divisor $D$ on the variety
$\bar{Z}$ such that the equivalence
$$
D\sim_{\mathbb{Q}}
-2K_{\bar{Z}}-(\beta\circ\pi\circ\bar{\pi})^{*}\big(16K_{U}\big)-(\pi\circ\bar{\pi})^{*}\big(18K_{W}\big)
$$
holds. Analyzing the base locus of the pencil $|-2K_{\bar{Z}}|$,
we see that $D$ is nef, but $D^{3}>0$.

The divisor $D$ is nef and big, but
$D\cdot\bar{H}_{1}\cdot\bar{H}_{2}=0$, where $\bar{H}_{1}$ and
$\bar{H}_{2}$ are the proper transforms on $\bar{Z}$ of general
surfaces in $\mathcal{M}$, which is impossible by
Corollary~\ref{corollary:Cheltsov}.
\end{proof}

The hypersurface $X$ can be given by the equation
$$
w^{2}y+wf_{13}(x,y,z,t)+f_{24}(x,y,z,t)=0\subset\mathbb{P}(1,2,3,8,11)\cong\mathrm{Proj}\Big(\mathbb{C}[x,y,z,t,w]\Big),
$$
where $\mathrm{wt}(x)=1$, $\mathrm{wt}(y)=2$, $\mathrm{wt}(z)=3$,
$\mathrm{wt}(t)=8$, $\mathrm{wt}(w)=11$, and $f_{i}(x,y,z,t)$ is a
sufficiently general quasihomogeneous polynomial of degree $i$.

Let $\bar{E}$ and $\bar{F}$ be the proper transforms on the
variety $Z$ of the divisors $E$ and $F$ respectively, and $S_{1}$,
$S_{2}$, $S_{3}$ and $S_{8}$ be the proper transforms on the
variety $Z$ of the surfaces that are cut on the hypersurface $X$
by the equations $x=0$, $y=0$, $z=0$ and $t=0$ respectively. Then
\begin{equation}
\label{equation:n-56-rational-equivalences-I} \left\{\aligned
&\bar{F}\sim_{\mathbb{Q}}\pi^{*}(F)-\frac{1}{3}G,\\
&\bar{E}\sim_{\mathbb{Q}}(\beta\circ\pi)^{*}(E)-\frac{5}{8}\pi^{*}(F)-\frac{2}{3}G,\\
&S_{1}\sim_{\mathbb{Q}}(\alpha\circ\beta\circ\pi)^{*}\big(-K_{X}\big)-\frac{1}{11}(\beta\circ\pi)^{*}(E)-\frac{1}{8}\pi^{*}(F)-\frac{1}{3}G,\\
&S_{2}\sim_{\mathbb{Q}}(\alpha\circ\beta\circ\pi)^{*}\big(-2K_{X}\big)-\frac{13}{11}(\beta\circ\pi)^{*}(E)-\frac{5}{8}\pi^{*}(F)-\frac{2}{3}G,\\
&S_{3}\sim_{\mathbb{Q}}(\alpha\circ\beta\circ\pi)^{*}\big(-3K_{X}\big)-\frac{3}{11}(\beta\circ\pi)^{*}(E)-\frac{3}{8}\pi^{*}(F)-\frac{1}{3}G,\\
&S_{8}\sim_{\mathbb{Q}}(\alpha\circ\beta\circ\pi)^{*}\big(-8K_{X}\big)-\frac{8}{11}(\beta\circ\pi)^{*}(E).\\
\endaligned\right.
\end{equation}

The base locus of $|-2K_{Z}|$ consists of the irreducible curves
$C$, $L_{1}$, $L_{2}$ and $L_{3}$ such that the curve
$\alpha\circ\beta\circ\pi(C)$ is cut on $X$ by the equations
$x=y=0$, the curve $\beta\circ\pi(L_{1})$ is contained in the
divisor $E$, the curve $\beta\circ\pi(L_{1})$ is contained in
$|\mathcal{O}_{\mathbb{P}(1,\,3,\,8)}(1)|$, the curve $\pi(L_{2})$
is contained in the divisor $F$, the curve $\pi(L_{2})$ is
contained in $|\mathcal{O}_{\mathbb{P}(1,\,3,\,5)}(1)|$, the curve
$L_{3}$ is contained in the divisor $G$, the curve $L_{3}$ is
contained in $|\mathcal{O}_{\mathbb{P}(1,\,1,\,2)}(1)|$. Moreover,
we have
$$
S_{1}\cdot D=C+2L_{1}+2L_{2}+L_{3},\ \bar{E}\cdot D=2L_{1},\
\bar{F}\cdot D=2L_{2},\ G\cdot D=2L_{3},
$$
where $D$ is a general surface in $|-2K_{Z}|$. It follows from the
equivalences~\ref{equation:n-56-rational-equivalences-I} that
$$
-K_{Z}\cdot C=\frac{1}{10},\ -K_{Z}\cdot L_{1}=-\frac{1}{3},\
-K_{Z}\cdot L_{2}=-\frac{1}{10},\ -K_{Z}\cdot L_{3}=\frac{1}{2},
$$
which implies that the curves $L_{1}$ and $L_{2}$ are the only
curves on the variety $Z$ that have negative intersection with the
divisor $-K_{Z}$.

The singularities of the log pair $(Z, \lambda |-2K_{Z}|)$ are
log-terminal for some rational числа $\lambda>1/2$, but the
divisor $K_{Z}+\lambda |-2K_{Z}|$ has nonnegative intersection
with all curves on the variety $Z$ except the curves $L_{1}$ and
$L_{2}$. It follows from \cite{Sho93} that there is a composition
of antiflips $\zeta\colon Z\dasharrow\bar{Z}$ such that the
divisor $-K_{\bar{Z}}$ is numerically effective.

Let $\mathcal{P}$ be the proper transform of the linear system
$\mathcal{M}$ on the variety $\bar{Z}$. Then the singularities of
the log pair $(\bar{Z}, \frac{1}{k}\mathcal{P})$ are terminal,
because the singularities of the log pair $(Z,
\frac{1}{k}\mathcal{B})$ are terminal, and the rational map
$\zeta$ is a log flop with respect to the log pair $(Z,
\frac{1}{k}\mathcal{B})$.

It follows from the
equivalences~\ref{equation:n-56-rational-equivalences-I} that
$-K_{Z}\sim_{\mathbb{Q}} S_{1}$ and
\begin{equation}
\label{equation:n-56-rational-equivalences-II} \left\{\aligned
&S_{1}\sim_{\mathbb{Q}}(\alpha\circ\beta\circ\pi)^{*}\big(-K_{X}\big)-\frac{1}{11}\bar{E}-\frac{2}{11}\bar{F}-\frac{5}{11}G,\\
&S_{2}\sim_{\mathbb{Q}}(\alpha\circ\beta\circ\pi)^{*}\big(-2K_{X}\big)-\frac{13}{11}\bar{E}-\frac{15}{11}\bar{F}-\frac{21}{11}G,\\
&S_{3}\sim_{\mathbb{Q}}(\alpha\circ\beta\circ\pi)^{*}\big(-3K_{X}\big)-\frac{3}{11}\bar{E}-\frac{6}{11}\bar{F}-\frac{4}{11}G,\\
&S_{8}\sim_{\mathbb{Q}}(\alpha\circ\beta\circ\pi)^{*}\big(-8K_{X}\big)-\frac{8}{11}\bar{E}-\frac{5}{11}\bar{F}-\frac{7}{11}G.\\
\endaligned
\right.
\end{equation}

The equivalences~\ref{equation:n-56-rational-equivalences-II}
implies that the rational functions $y/x^{2}$, $zy/x^{5}$ and
$ty^{3}/x^{14}$ are contained in the linear systems $|2S_{1}|$,
$|5S_{1}|$ and $|14S_{1}|$ respectively. In particular, the
complete linear system $|-70K_{Z}|$ induces the dominant rational
map $Z\dasharrow\mathbb{P}(1,2,5,14)$, which implies that the
anticanonical divisor $-K_{\bar{Z}}$ is nef and big, which
contradicts Lemma~\ref{lemma:Noether-Fano}.

\section{Case   $n=58$, hypersurface of degree $24$ in $\mathbb{P}(1,3,4,7,10)$.}
\label{section:n-58}

We use the notations and assumptions of
Section~\ref{section:start}. Let $n=58$. Then $X$ is a
hypersurface of degree $24$ in $\mathbb{P}(1,3,4,7,10)$, the
equality $-K_{X}^{3}=1/35$ holds, and the singularities of $X$
consist of the points $P_{1}$, $P_{2}$ and $P_{3}$ of type
$\frac{1}{2}(1,1,1)$, $\frac{1}{7}(1,3,4)$ and
$\frac{1}{10}(1,3,7)$ respectively.

There is a commutative diagram
$$
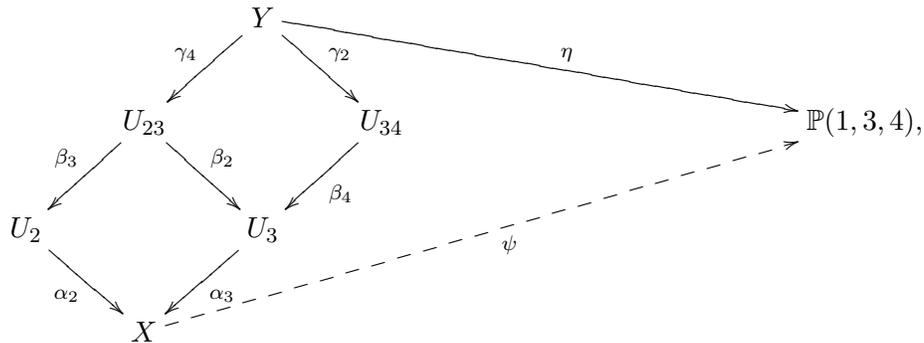
\xymatrix{
&&&Y\ar@{->}[dl]_{\gamma_{4}}\ar@{->}[dr]^{\gamma_{2}}\ar@{->}[drrrrrr]^{\eta}&&&&&&\\
&&U_{23}\ar@{->}[dl]_{\beta_{3}}\ar@{->}[dr]^{\beta_{2}}&&U_{34}\ar@{->}[dl]^{\beta_{4}}&&&&&\mathbb{P}(1,3,4),\\
&U_{2}\ar@{->}[dr]_{\alpha_{2}}&&U_{3}\ar@{->}[dl]^{\alpha_{3}}&&&&&&\\
&&X\ar@{-->}[rrrrrrruu]_{\psi}&&&&&&&}
$$
where $\psi$ is a projection, $\alpha_{2}$ is the weighted blow up
of  $P_{2}$ with weights $(1,3,4)$, $\alpha_{3}$ is the weighted
blow up of  $P_{3}$ with weights $(1,3,7)$, $\beta_{3}$ is the
weighted blow up with weights $(1,3,7)$ of the proper transform of
$P_{3}$ on  $U_{2}$, $\beta_{2}$ is the weighted blow up with
weights $(1,3,4)$ of the proper transform of the point $P_{2}$ on
$U_{2}$, $\beta_{4}$ is the weighted blow up with weights
$(1,3,4)$ of the singular point of the variety $U_{3}$ that is a
quotient singularity of type $\frac{1}{7}(1,3,4)$ contained in the
exceptional divisor of the morphism $\alpha_{3}$, $\gamma_{2}$ is
the weighted blow up with weights $(1,3,4)$ of the proper
transform of the point $P_{2}$ on $U_{34}$, $\gamma_{4}$ is the
weighted blow up with weights $(1,3,4)$ of the singular point of
the variety $U_{23}$ that is a quotient singularity of type
$\frac{1}{7}(1,3,4)$ contained in the exceptional divisor of the
morphism $\beta_{3}$, and $\eta$ is an elliptic fibration.

In the rest of the section we prove the following result.

\begin{proposition}
\label{proposition:n-58} The claim of Theorem~\ref{theorem:main}
holds for $n=58$.
\end{proposition}

It follows from Theorem~\ref{theorem:smooth-points},
Lemma~\ref{lemma:Ryder} and
Proposition~\ref{proposition:singular-points} that $\mathbb{CS}(X,
\frac{1}{k}\mathcal{M})\subseteq\{P_{2},P_{3}\}$.

\begin{lemma}
\label{lemma:n-58-P3-is-a-center} The set $\mathbb{CS}(X,
\frac{1}{k}\mathcal{M})$ contains the point $P_{3}$.
\end{lemma}

\begin{proof}
Suppose that the set $\mathbb{CS}(X, \frac{1}{k}\mathcal{M})$ does
not contain the point $P_{3}$. Let $\mathcal{D}_{2}$ be the proper
transform of the linear system $\mathcal{M}$ on the variety
$U_{2}$, and $O$ and $Q$ be the singular points of the variety
$U_{2}$ contained in the exceptional divisor of $\alpha_{2}$ such
that the points $O$ and $Q$ are sin\-gu\-la\-ri\-ties of types
$\frac{1}{3}(1,1,2)$ and $\frac{1}{4}(1,1,3)$ respectively. Then
$\mathcal{D}_{2}\sim_{\mathbb{Q}}-kK_{U_{2}}$ by
Theorem~\ref{theorem:Kawamata}.

Suppose that $O\in \mathbb{CS}(U_{2},
\frac{1}{k}\mathcal{D}_{2})$. Let $\pi\colon W\to U_{2}$ be the
weighted blow up of the singular point $O$ with weights $(1,1,2)$,
$\mathcal{B}$ and $\mathcal{P}$ be the proper transforms of
$\mathcal{M}$ and $|-4K_{X}|$ on the variety $W$ respectively, and
$S$ be a general surface in $\mathcal{P}$. Then the base locus of
$\mathcal{P}$ consists of an irreducible curve $C$ such that
$\alpha_{2}(C)$ is the base curve of $|-4K_{X}|$, but $C^{2}<0$ on
the normal surface $S$, which contradicts
$\mathcal{B}\vert_{S}\sim_{\mathbb{Q}}kC$ by
Lemmas~\ref{lemma:Cheltsov-II} and \ref{lemma:normal-surface}.

We see that $Q\in \mathbb{CS}(U_{2}, \frac{1}{k}\mathcal{D}_{2})$
by Lemma~\ref{lemma:Cheltsov-Kawamata}, because
$\mathbb{CS}(U_{2}, \frac{1}{k}\mathcal{D}_{2})\ne\varnothing$ by
Lemma~\ref{lemma:Noether-Fano}.

Let $\zeta\colon U\to U_{2}$ be the weighted blow up of the point
$Q$ with weights $(1,1,3)$, $\mathcal{H}$ be the proper transform
of $\mathcal{M}$ on the variety $U$, $H$ be a general surface of
the linear system $\mathcal{H}$, and $D$ be a general surface in
$|-3K_{U}|$. Then $D$ is normal, and the base locus of the pencil
пучка $|-3K_{U}|$ consists of the irreducible curve $Z$ such that
$\alpha_{2}(Z)$ is the unique base curve of the pencil
$|-3K_{X}|$.

The equivalence $\mathcal{H}\vert_{D}\sim_{\mathbb{Q}}kZ$ holds by
Theorem~\ref{theorem:Kawamata}, and the inequality $Z^{2}<0$ holds
on the surface $D$, which is impossible by
Lemmas~\ref{lemma:Cheltsov-II} and \ref{lemma:normal-surface}.
\end{proof}

Let $\mathcal{D}_{2}$ and $\mathcal{D}_{23}$ be the proper
transforms of $\mathcal{M}$ on $U_{2}$ and $U_{23}$ respectively.
The arguments of the proof of
Lemma~\ref{lemma:n-58-P3-is-a-center} implies the following two
corollaries.

\begin{corollary}
\label{corollary:n-58-U2-centers} Suppose that
$\mathcal{D}_{2}\sim_{\mathbb{Q}}-kK_{U_{2}}$. Then the set
$\mathbb{CS}(U_{2}, \frac{1}{k}\mathcal{D}_{2})$ does not
contain~sub\-varieties of the variety $U_{2}$ that are contained
in the exceptional divisor of the morphism $\alpha_{2}$.
\end{corollary}

\begin{corollary}
\label{corollary:n-58-U23-centers} Suppose that
$\mathcal{D}_{23}\sim_{\mathbb{Q}}-kK_{U_{23}}$. Then the set
$\mathbb{CS}(U_{23}, \frac{1}{k}\mathcal{D}_{23})$ does not
contain sub\-va\-rieties of the variety $U_{23}$ that are
contained in the exceptional divisor of the morphism $\beta_{2}$.
\end{corollary}

Let $\mathcal{D}_{3}$ and $\mathcal{D}_{23}$ be the proper
transforms of $\mathcal{M}$ on $U_{3}$ and $U_{34}$ respectively.
Then we have the equivalence
$\mathcal{D}_{3}\sim_{\mathbb{Q}}-kK_{U_{3}}$ by
Theorem~\ref{theorem:Kawamata}.

\begin{lemma}
\label{lemma:n-58-P2-is-a-center} The set $\mathbb{CS}(X,
\frac{1}{k}\mathcal{M})$ contains the point $P_{2}$.
\end{lemma}

\begin{proof}
Suppose that $P_{2}\not\in\mathbb{CS}(X, \frac{1}{k}\mathcal{M})$.
Then $\mathcal{D}_{3}\sim_{\mathbb{Q}}-kK_{U_{2}}$ by
Theorem~\ref{theorem:Kawamata}, which implies~that the set
$\mathbb{CS}(U_{3}, \frac{1}{k}\mathcal{D}_{3})$ is not empty by
Lemma~\ref{lemma:Noether-Fano}.

Let $P_{4}$ and $P_{5}$ be the singular points of the variety
$U_{3}$ that are contained in the exceptional divisor of the
birational morphism $\alpha_{3}$ such that the points $P_{4}$ and
$P_{5}$ are quotient singularities of types $\frac{1}{7}(1,3,4)$
and $\frac{1}{3}(1,1,2)$ respectively. Then the set
$\mathbb{CS}(U_{3}, \frac{1}{k}\mathcal{D}_{3})$ contains either
the singular point $P_{4}$, or the singular point $P_{5}$ by
Lemma~\ref{lemma:Cheltsov-Kawamata}.

It follows from Lemma~\ref{lemma:n-58-P3-is-a-center} that the set
$\mathbb{CS}(U_{3}, \frac{1}{k}\mathcal{D}_{3})$ does not contain
the point $P_{5}$, which implies that $\mathbb{CS}(U_{3},
\frac{1}{k}\mathcal{D}_{3})$ contains $P_{4}$. Hence, the
equivalence $\mathcal{D}_{34}\sim_{\mathbb{Q}}-kK_{U_{34}}$ holds
by Theorem~\ref{theorem:Kawamata}, and the set
$\mathbb{CS}(U_{34}, \frac{1}{k}\mathcal{D}_{34})$ is not empty by
Lemma~\ref{lemma:Noether-Fano}.

Let $P_{6}$ and $P_{7}$ be the singular points of the variety
$U_{34}$ that are contained in the exceptional divisor of the
birational morphism $\beta_{4}$ such that the points $P_{6}$ and
$P_{7}$ are quotient singularities of types $\frac{1}{3}(1,1,2)$
and $\frac{1}{4}(1,1,3)$ respectively. Then the set
$\mathbb{CS}(U_{34}, \frac{1}{k}\mathcal{D}_{34})$ contains either
the singular point $P_{6}$, or the singular point $P_{7}$ by
Lemma~\ref{lemma:Cheltsov-Kawamata}.

Suppose that the set $\mathbb{CS}(U_{34},
\frac{1}{k}\mathcal{D}_{34})$ contains the point $P_{7}$. Let
$\zeta\colon U\to U_{34}$ be the weighted blow up of $P_{7}$ with
weights $(1,1,3)$, $\mathcal{H}$ be the proper transform of
$\mathcal{M}$ on $U$, $H$ be a general surface of $\mathcal{H}$,
and $D$ be a general surface of the pencil $|-3K_{U}|$. Then $D$
is normal, and the base locus of the pencil $|-3K_{U}|$ consists
of the irreducible curve $Z$ such that
$\alpha_{3}\circ\beta_{4}(Z)$ is the unique base curve of the
pencil $|-3K_{X}|$. Moreover, the equivalence
$\mathcal{H}\vert_{D}\sim_{\mathbb{Q}}kZ$ holds, and the
inequality  $Z^{2}<0$ holds on the surface $D$, which is
impossible by Lemmas~\ref{lemma:Cheltsov-II} and
\ref{lemma:normal-surface}.

Therefore, the set $\mathbb{CS}(U_{34},
\frac{1}{k}\mathcal{D}_{34})$ contains the point $P_{6}$.

The hypersurface $X$ can be given by the quasihomogeneous equation
$$
w^{2}z+wf_{14}(x,y,z,t)+f_{24}(x,y,z,t)=0\subset\mathrm{Proj}\Big(\mathbb{C}[x,y,z,t,w]\Big),
$$
where $\mathrm{wt}(x)=1$, $\mathrm{wt}(y)=3$, $\mathrm{wt}(z)=4$,
$\mathrm{wt}(t)=7$, $\mathrm{wt}(w)=10$, and $f_{i}(x,y,z,w)$ is a
quasihomogeneous polynomial of degree $i$. Let $\mathcal{P}$ be a
pencil consisting of the surfaces that are cut on the hypersurface
$X$ by the equations $\lambda x^{4}+\mu z=0$, where
$(\lambda,\mu)\in\mathbb{P}^{1}$. Then the base locus of the
pencil $\mathcal{P}$ consists of the irreducible curve that are
cut on $X$ by the equations $x=z=0$.

Let $\pi\colon W\to U_{34}$ be the weighted blow up of $P_{6}$
with weights $(1,1,2)$, $\mathcal{B}$ be the proper transform of
$\mathcal{M}$ on $W$, $\mathcal{H}$ be the proper transform of
$\mathcal{P}$ on $W$, $H$ be a general surface of the pencil
$\mathcal{H}$, and $E$, $F$ and $G$ be the exceptional divisors of
$\alpha_{3}$, $\beta_{4}$ and $\pi$ respectively. Then
$$
H\sim_{\mathbb{Q}}(\alpha_{3}\circ\beta_{4}\circ\pi)^{*}\big(-4K_{X}\big)-\frac{1}{5}(\beta_{4}\circ\pi)^{*}(E)-\frac{4}{7}(\pi)^{*}(F)-\frac{1}{3}G,\\
$$
the surface $H$ normal, and the base locus of the pencil
$\mathcal{H}$ consists of the curves $C$ and $L$ such that
$\alpha_{3}\circ\beta_{4}\circ\pi(C)$ is the unique base curve of
the pencil $\mathcal{P}$, and the curve $\beta_{4}\circ\pi(L)$ is
the unique curve on the surface $E\cong\mathbb{P}(1,3,7)$ that is
contained in the linear system
$|\mathcal{O}_{\mathbb{P}(1,\,3,\,7)}(1)|$.

Let $S$ be the surface of the linear system $|-K_{W}|$, and
$\bar{E}$ be the proper transform of the surface $E$ on the
variety $W$. Then $S\cdot H=C+L$ and $\bar{E}\cdot H=4L$, but
$$
\bar{E}\sim_{\mathbb{Q}}(\beta_{4}\circ\pi)^{*}(E)-\frac{4}{7}(\pi)^{*}(F)-\frac{1}{3}G,\\
$$
which implies that the intersection form of the curves $L$ and $C$
on the surface $H$ is negatively defined. The equivalence
$\mathcal{B}\vert_{H}\sim_{\mathbb{Q}}kC+kL$ holds, which
contradicts Lemmas~\ref{lemma:normal-surface} and
\ref{lemma:Cheltsov-II}.
\end{proof}

Hence, we have $\mathbb{CS}(X,
\frac{1}{k}\mathcal{M})=\{P_{2},P_{3}\}$, and
$\mathcal{D}_{23}\sim_{\mathbb{Q}}-kK_{U_{23}}$ by
Theorem~\ref{theorem:Kawamata}.

It easily follows from Lemmas~\ref{lemma:Noether-Fano} and
\ref{lemma:Cheltsov-Kawamata}, the proof of
Lemma~\ref{lemma:n-58-P2-is-a-center} and
Corollary~\ref{corollary:n-58-U23-centers} that the set
$\mathbb{CS}(U_{23}, \frac{1}{k}\mathcal{D}_{23})$ contains the
singular point of the variety $U_{23}$ that is a quotient
singularity of type $\frac{1}{7}(1,3,4)$ contained in the
exceptional divisor of $\beta_{3}$, which implies that the proper
transform of the linear system $\mathcal{M}$ on the variety $Y$
lies in the fibers of the fibration $\eta$.

\section{Case $n=64$, hypersurface of degree $26$ in $\mathbb{P}(1,2,5,6,13)$.}
\label{section:n-64}

We use the notations and assumptions of
Section~\ref{section:start}. Let $n=64$. Then $X$ is
a~general~hyper\-surface in $\mathbb{P}(1,2,5,6,13)$ of degree
$26$. The singularities of the hypersurface $X$ consist of the
points $P_{1}$, $P_{2}$, $P_{3}$~and~$P_{4}$~that are
singularities of type $\frac{1}{2}(1,1,1)$, the point $P_{5}$ that
is a singularity of type $\frac{1}{5}(1,2,3)$, and the point
$P_{6}$ that is a singularity of type $\frac{1}{6}(1,1,5)$.

There is a commutative diagram
$$
\xymatrix{
&&Y\ar@{->}[dl]_{\alpha}\ar@{->}[dr]^{\eta}&&\\%
&X\ar@{-->}[rr]_{\psi}&&\mathbb{P}(1,2,5),&}
$$
where $\psi$ is a projection, $\alpha$ is the blow up of $P_{6}$
with weights $(1,1,5)$, and $\eta$ is an elliptic~fibration.

There is a commutative diagram
$$
\xymatrix{
&&U\ar@{->}[dl]_{\beta}\ar@{->}[dr]^{\omega}&&\\
&X\ar@{-->}[rr]_{\xi}&&\mathbb{P}(1,1,3),&}
$$
where $\xi$ is a projection, $\beta$ is the blow up of $P_{5}$
with weights $(1,2,3)$, and $\omega$ is an elliptic~fibration.

\begin{proposition}
\label{proposition:n-64} Either there is a commutative diagram
$$\xymatrix{
&X\ar@{-->}[d]_{\psi}\ar@{-->}[rr]^{\rho}&&V\ar@{->}[d]^{\nu}&\\%
&\mathbb{P}(1,2,5)\ar@{-->}[rr]_{\phi}&&\mathbb{P}^{2},&}
$$
or there is a commutative diagram
$$
\xymatrix{
&X\ar@{-->}[d]_{\xi}\ar@{-->}[rr]^{\rho}&&V\ar@{->}[d]^{\nu}&\\%
&\mathbb{P}(1,1,3)\ar@{-->}[rr]_{\sigma}&&\mathbb{P}^{2},&}
$$
where $\phi$ and $\sigma$ are birational maps.
\end{proposition}

\begin{proof}
See the proof of Proposition~\ref{proposition:n-49}.
\end{proof}

\section{Case $n=65$, hypersurface of degree $27$ in $\mathbb{P}(1,2,5,9,11)$.}
\label{section:n-65}

We use the notations and assumptions of
Section~\ref{section:start}. Let $n=65$. Then $X$ is a
hypersurface  of degree $27$ in $\mathbb{P}(1,2,5,9,11)$, the
equality $-K_{X}^{3}=3/110$ holds, and the singularities of $X$
consist of the point $P_{1}$ that is a quotient singularity of
type $\frac{1}{2}(1,1,1)$, the point $P_{2}$ that is a quotient
singularity of type $\frac{1}{5}(1,1,4)$, and the point $P_{3}$
that is a quotient singularity of type $\frac{1}{11}(1,2,9)$.

There is a commutative diagram
$$
\xymatrix{
&U\ar@{->}[d]_{\alpha}&&W\ar@{->}[ll]_{\beta}&&Y\ar@{->}[ll]_{\gamma}\ar@{->}[d]^{\eta}&\\%
&X\ar@{-->}[rrrr]_{\psi}&&&&\mathbb{P}(1,2,5),&}
$$
where $\psi$ is a projection, $\alpha$ is the weighted blow up of
$P_{3}$ with weights $(1,2,9)$, $\beta$ is the weighted blow up
with weights $(1,2,7)$ of the point of $U$ that is a singularity
of type $\frac{1}{9}(1,2,7)$ contained in the $\alpha$-exceptional
divisor, $\gamma$ is the weighted blow up with weights $(1,2,5)$
of the singular point of type $\frac{1}{7}(1,2,5)$ contained in
the $\beta$-excep\-ti\-onal divisor, and $\eta$ is an elliptic
fibra\-ti\-on.

\begin{proposition}
\label{proposition:n-65} The claim of Theorem~\ref{theorem:main}
holds for $n=65$.
\end{proposition}

Let us prove Proposition~\ref{proposition:n-65}. It follows from
Theorem~\ref{theorem:smooth-points}, Lemma~\ref{lemma:Ryder} and
Proposition~\ref{proposition:singular-points} that $\mathbb{CS}(X,
\frac{1}{k}\mathcal{M})=\{P_{3}\}$. Let $E$ be the exceptional
divisor of  $\alpha$, $\mathcal{D}$ be the proper transform of the
linear system $\mathcal{M}$ on $U$, and $P_{4}$ and $P_{5}$ be the
singular points of $U$ that are singularities of types
$\frac{1}{2}(1,1,1)$ and $\frac{1}{9}(1,2,7)$ contained in $E$
respectively.
Then~$\mathcal{D}\sim_{\mathbb{Q}}-kK_{U}$~by~Theorem~\ref{theorem:Kawamata}.

It follows from Lemmas~\ref{lemma:Noether-Fano} and
\ref{lemma:Cheltsov-Kawamata} and the proof of
Lemma~\ref{lemma:n-48-points-P5} that $\mathbb{CS}(U,
\frac{1}{k}\mathcal{D})=\{P_{5}\}$.

Let $F$ be the exceptional  divisor of the morphism $\beta$,
$\mathcal{H}$ be the proper transform of the linear system
$\mathcal{M}$ on the variety $W$, and $P_{6}$ and $P_{7}$ be the
singular points of  $W$ that are singularities of types
$\frac{1}{2}(1,1,1)$ and $\frac{1}{7}(1,2,5)$ contained in  $F$
respectively. Then $\mathcal{H}\sim_{\mathbb{Q}}-kK_{W}$ by
Theorem~\ref{theorem:Kawamata}.

It follows from Lemmas~\ref{lemma:Noether-Fano} and
\ref{lemma:Cheltsov-Kawamata} that $\mathbb{CS}(U,
\frac{1}{k}\mathcal{D})\cap\{P_{6}, P_{7}\}\ne\varnothing$.

\begin{remark}
\label{remark:n-65-points-P7} In the case when the set
$\mathbb{CS}(W, \frac{1}{k}\mathcal{H})$ contains the point
$P_{7}$, it follows from Theorem~\ref{theorem:Kawamata} that there
is a commutative diagram
$$
\xymatrix{
&X\ar@{-->}[d]_{\psi}\ar@{-->}[rr]^{\rho}&&V\ar@{->}[d]^{\nu}&\\%
&\mathbb{P}(1,2,5)\ar@{-->}[rr]_{\zeta}&&\mathbb{P}^{2},&}
$$
where $\zeta$ is a birational map.
\end{remark}

We may assume that $P_{6}\in\mathbb{CS}(W,
\frac{1}{k}\mathcal{H})$. Note, that $E\cong\mathbb{P}(1,2,9)$ and
$F\cong\mathbb{P}(1,2,7)$.

Let $\pi\colon Z\to W$ be the weighted blow up of $P_{6}$ with
weights $(1,1,1)$, $G$ be the exceptional divisor of $\pi$, and
$\mathcal{B}$ be the proper transform of $\mathcal{M}$ on the
variety $Z$. Then $G\cong\mathbb{P}^{2}$, and the equivalence
$\mathcal{B}\sim_{\mathbb{Q}}-kK_{Z}$ holds by
Theorem~\ref{theorem:Kawamata}.

The hypersurface $X$ can be given by equation
$$
w^{2}z+wf_{16}(x,y,z,t)+f_{27}(x,y,z,t)=0\subset\mathbb{P}(1,2,5,9,11)\cong\mathrm{Proj}\Big(\mathbb{C}[x,y,z,t,w]\Big),
$$
where $\mathrm{wt}(x)=1$, $\mathrm{wt}(y)=2$, $\mathrm{wt}(z)=5$,
$\mathrm{wt}(t)=9$, $\mathrm{wt}(w)=11$, and $f_{i}(x,y,z,t)$ is a
quasihomogeneous polynomial of degree $i$. Let $\bar{E}$ and
$\bar{F}$ be the proper transforms of the exceptional divisors $E$
and $F$ on the variety $Z$ respectively, and $\mathcal{P}$ be the
proper transform on $Z$ of the pencil of surfaces that are cut out
on the hypersurface $X$ by $\lambda x^{5}+\mu z=0$, where
$(\lambda,\mu)\in\mathbb{P}^{1}$.

The base locus of $\mathcal{P}$ consists of the irreducible curves
$C$, $L_{1}$, $L_{2}$, $\Delta_{1}$, $\Delta_{2}$ and $\Delta$
such that the curve $\alpha\circ\beta\circ\pi(C)$ is cut out on
$X$ by $x=z=0$, the curve $\beta\circ\pi(L_{1})$ is contained in
the excep\-ti\-onal divisor $E$, the curve $\beta\circ\pi(L_{1})$
is the unique curve in $|\mathcal{O}_{\mathbb{P}(1,\,2,\,9)}(1)|$,
the curve $\pi(L_{1})$~is contained in $F$, the curve
$\pi(L_{1})$~is contained in
$|\mathcal{O}_{\mathbb{P}(1,\,2,\,7)}(1)|$, the curves
$\Delta_{1}$ and $\Delta_{2}$ are the lines on $G$ that are cut
out by $\bar{E}$ and $\bar{F}$ respectively, and the curve
$\Delta$ is a line on $G$, which is different from the lines
$\Delta_{1}$ and $\Delta_{2}$.

Let $D$ be a general surface of the pencil $\mathcal{P}$, and $S$
be the proper transform on $Z$ of the surface that is cut on the
hypersurface $X$ by the equation $x=0$. Then
$$
S\cdot D=C+L_{1}+L_{2},\ \bar{E}\cdot D=5L_{1}+\Delta_{1},\ \bar{F}\cdot D=5L_{2}+\Delta_{2}, %
$$
the surface $D$ is normal, and $D$ is smooth in the neighborhood
of $G$. In particular, it follows from the local computations and
the adjunction formula that the equalities
\begin{equation}
\label{equation:n-65-equalities}
\Delta_{1}\cdot\Delta_{2}=\Delta_{1}\cdot L_{2}=\Delta_{2}\cdot L_{1}=1\ \Delta_{1}\cdot C=\Delta_{2}\cdot C=0,\ \Delta_{1}^{2}=\Delta_{2}^{2}=-4 %
\end{equation}
hold on the surface $D$. However, we have
\begin{equation}
\label{equation:n-65-rational-equivalences} \left\{\aligned
&\bar{F}\sim_{\mathbb{Q}}\pi^{*}(F)-\frac{1}{2}G,\\
&\bar{E}\sim_{\mathbb{Q}}(\beta\circ\pi)^{*}(E)-\frac{7}{9}\pi^{*}(F)-\frac{1}{2}G,\\
&D\sim_{\mathbb{Q}}(\alpha\circ\beta\circ\pi)^{*}(-5K_{X})-\frac{5}{11}(\beta\circ\pi)^{*}(E)-\frac{5}{9}\pi^{*}(F)-\frac{3}{2}G,\\
&S\sim_{\mathbb{Q}}(\alpha\circ\beta\circ\pi)^{*}(-K_{X})-\frac{1}{11}(\beta\circ\pi)^{*}(E)-\frac{1}{9}\pi^{*}(F)-\frac{1}{2}G.\\
\endaligned
\right.
\end{equation}

It follows from the equalities~\ref{equation:n-65-equalities} and
equivalences~\ref{equation:n-65-rational-equivalences} that the
equalities
$$
C\cdot C=L_{1}\cdot L_{1}=-\frac{1}{2},\ L_{2}\cdot L_{2}=-\frac{3}{7},\ C\cdot L_{1}=C\cdot L_{2}=L_{1}\cdot L_{2}=0%
$$
hold on $D$. The intersection form of the curves $C$, $L_{1}$ and
$L_{2}$ on the surface $D$ is negatively defined, but
$\mathcal{B}\vert_{D}\sim_{\mathbb{Q}}kC+kL_{1}+kL_{2}$, which
contradicts Lemmas~\ref{lemma:normal-surface} and
\ref{lemma:Cheltsov-II}.

\section{Case $n=68$, hypersurface of degree $28$ in $\mathbb{P}(1,3,4,7,14)$.}
\label{section:n-68}

We use the notations and assumptions of
Section~\ref{section:start}. Let $n=68$. Then $X$ is a
hypersurface of degree $28$ in $\mathbb{P}(1,3,4,7,14)$, the
singularities of the $X$ consist of the point $P_{1}$ that is a
singularity of type $\frac{1}{2}(1,1,1)$, the point $P_{2}$ that
is a singularity of type $\frac{1}{3}(1,1,2)$, the points $P_{3}$
and $P_{4}$ that are singularities of type $\frac{1}{7}(1,3,4)$,
and $-K_{X}^{3}=1/42$.

\begin{proposition}
\label{proposition:n-68} The claim of Theorem~\ref{theorem:main}
holds for $n=68$.
\end{proposition}

\begin{proof}
It follows from Lemma~\ref{lemma:Ryder} and
Proposition~\ref{proposition:singular-points} that $\mathbb{CS}(X,
\frac{1}{k}\mathcal{M})\subseteq\{P_{3},P_{4}\}$, but the proof is
trivial if $\mathbb{CS}(X, \frac{1}{k}\mathcal{M})=\{P_{3},
P_{4}\}$. Hence, we may assume that $P_{4}\not\in\mathbb{CS}(X,
\frac{1}{k}\mathcal{M})$.

There is a commutative diagram
$$
\xymatrix{
&&Y\ar@{->}[dl]_{\gamma}\ar@{->}[dr]^{\delta}\ar@{->}[drrrrrrr]^{\eta}&&&&&&&\\
&U\ar@{->}[dr]_{\alpha}&&W\ar@{->}[dl]^{\beta}&&&&&&\mathbb{P}(1,2,3),\\
&&X\ar@{-->}[rrrrrrru]_{\psi}&&&&&&&}
$$
where $\psi$ is a projection, $\alpha$ is the weighted blow up of
$P_{3}$ with weights $(1,3,4)$, $\beta$ is the weighted blow up of
the point $P_{4}$ with weights $(1,3,4)$, $\gamma$ is the weighted
blow up with weights $(1,3,4)$ of the proper transform of the
singular point $P_{4}$ on the variety $U$, $\delta$ is the
weighted blow up with weights $(1,3,4)$ of the proper transform of
$P_{3}$ on $W$, and $\eta$ is an elliptic fibration.

Let $\mathcal{B}$ be the proper transform of $\mathcal{M}$ on $U$,
and $P_{5}$ and $P_{6}$ be the singular points of $U$ that are
singularities of types $\frac{1}{3}(1,1,2)$ and
$\frac{1}{4}(1,1,3)$ contained in the exceptional divisor of the
birational morphism $\alpha$ respectively. Then $\mathbb{CS}(U,
\frac{1}{k}\mathcal{B})\cap\{P_{5}, P_{6}\}\ne\varnothing$ by
Lemma~\ref{lemma:Cheltsov-Kawamata}.

Suppose that $P_{6}\in\mathbb{CS}(U, \frac{1}{k}\mathcal{B})$. Let
$\pi\colon W\to U$ be the weighted blow up of the point $P_{6}$
with weights $(1,1,3)$, and $\mathcal{B}$ and $\mathcal{P}$ be the
proper transforms of $\mathcal{M}$ and the pencil $|-3K_{X}|$ on
the variety $W$ respectively. Then
$\mathcal{B}\sim_{\mathbb{Q}}-kK_{W}$ by
Theorem~\ref{theorem:Kawamata}, and the base locus of
$\mathcal{P}$ consists of  the irreducible curve $Z$ such that
$\alpha\circ\pi(Z)$ is the base curve of the pencil $|-3K_{X}|$.

Let $S$ be a general surface in $\mathcal{P}$, and $B$ be a
general surface in $\mathcal{B}$. Then the surface $S$ is normal,
the inequality $Z^{2}<0$ holds on the surface $S$, but
$\mathcal{B}\vert_{S}\sim_{\mathbb{Q}} kZ$. Therefore, the support
of the cycle $B\cdot S$ is contained in $Z$ by
Lemma~\ref{lemma:normal-surface}, which is impossible by
Lemma~\ref{lemma:Cheltsov-II}.

Hence, the set $\mathbb{CS}(U, \frac{1}{k}\mathcal{B})$ contains
the point $P_{5}$.

Let $\zeta\colon Z\to U$ be the weighted blow up of $P_{5}$ with
weights $(1,1,2)$, and $\mathcal{D}$ and $\mathcal{H}$ be the
proper transforms of  $\mathcal{M}$ and $|-4K_{X}|$ on $Z$
respectively. Then $\mathcal{D}\sim_{\mathbb{Q}}-kK_{Z}$ by
Theorem~\ref{theorem:Kawamata}, and the base locus of
$\mathcal{H}$ consists of a curve $C$ such that
$\alpha\circ\zeta(C)$ is the base curve of $|-4K_{X}|$.

Let $H$ be a  general surface in $\mathcal{H}$. Then $H\cdot C=0$
and $H^{3}>0$. Thus, the divisor $H$ is nef and big. On the other
hand, the equality $H\cdot D_{1}\cdot D_{2}=0$ holds, where
$D_{1}$ and $D_{2}$ are general surfaces of the linear system
$\mathcal{D}$, which is impossible by
Corollary~\ref{corollary:Cheltsov}.
\end{proof}

\section{Case  $n=74$, hypersurface of degree $30$ in $\mathbb{P}(1,3,4,10,13)$.}
\label{section:n-74}

We use the notations and assumptions of
Section~\ref{section:start}. Let $n=74$. Then $X$ is a
hypersurface of degree $30$ in $\mathbb{P}(1,3,4,10,13)$ , the
equality $-K_{X}^{3}=1/52$ holds, and the singularities of $X$
consist of the point $P_{1}$ that is a quotient singularity of
type $\frac{1}{2}(1,1,1)$, the point $P_{2}$ that is a quotient
singularity of type $\frac{1}{4}(1,1,3)$, and the point $P_{3}$
that is a quotient singularity of type $\frac{1}{13}(1,3,10)$.

There is a commutative diagram
$$
\xymatrix{
&U\ar@{->}[d]_{\alpha}&&W\ar@{->}[ll]_{\beta}&&Y\ar@{->}[ll]_{\gamma}\ar@{->}[d]^{\eta}&\\%
&X\ar@{-->}[rrrr]_{\psi}&&&&\mathbb{P}(1,3,4),&}
$$
where $\psi$ is a projection, $\alpha$ is the weighted blow up of
$P_{3}$ with weights $(1,3,10)$, $\beta$ is the weighted blow up
with weights $(1,3,7)$ of the singular point of the variety $U$
that is a quotient singularity of type $\frac{1}{10}(1,3,7)$
contained in the exceptional divisor of the morphism $\alpha$,
$\gamma$ is the weighted blow up with weights $(1,3,4)$ of the
singular point of the variety $W$ that is a quotient singularity
of type $\frac{1}{7}(1,3,4)$ contained in the exceptional divisor
of $\beta$, and $\eta$ is an elliptic fibration.

In the rest of the section we prove the following result.

\begin{proposition}
\label{proposition:n-74} The claim of Theorem~\ref{theorem:main}
holds for $n=74$.
\end{proposition}

It follows from Theorem~\ref{theorem:smooth-points},
Lemma~\ref{lemma:Ryder} and
Proposition~\ref{proposition:singular-points} $\mathbb{CS}(X,
\frac{1}{k}\mathcal{M})=\{P_{3}\}$.

Let $E$ be the exceptional  divisor of the morphism $\alpha$,
$\mathcal{D}$ be the proper transform of the linear system
$\mathcal{M}$ on the variety $U$, and $P_{4}$ and $P_{5}$ be the
singular points of  $U$ that are quotient singularities of types
$\frac{1}{3}(1,1,2)$ and $\frac{1}{10}(1,3,7)$ contained in the
divisor $E$ respectively. Then the equivalence
$\mathcal{D}\sim_{\mathbb{Q}}-kK_{U}$ holds by
Theorem~\ref{theorem:Kawamata}.

\begin{lemma}
\label{lemma:n-74-points-P4} The set $\mathbb{CS}(U,
\frac{1}{k}\mathcal{D})$ does not contain the point $P_{4}$.
\end{lemma}

\begin{proof}
Suppose that the set $\mathbb{CS}(U, \frac{1}{k}\mathcal{D})$
contains the point $P_{4}$. Let $\pi\colon Z\to U$ be the weighted
blow up of $P_{4}$ with weights $(1,1,2)$, $G$ be the exceptional
divisor of $\pi$, and $\mathcal{B}$ be the proper transforms of
the linear system $\mathcal{M}$ on the variety $Z$. Then
$\mathcal{B}\sim_{\mathbb{Q}}-kK_{Z}$ by
Theorem~\ref{theorem:Kawamata}.

Let $D$ be a divisor on $Z$ such that the equivalence
$$
D\sim_{\mathbb{Q}} -4K_{Z}-\pi^{*}\big(36K_{U}\big)
$$
holds. Analyzing the base locus of the pencil $|-4K_{Z}|$, we see
that $D$ is nef and big, but
$$
D\cdot B_{1}\cdot
B_{2}=\Big((\alpha\circ\pi)^{*}\big(-kK_{X}\big)-\frac{k}{13}\pi^{*}(E)-\frac{k}{3}G\Big)^{2}
\Big(-4K_{Z}-\pi^{*}\big(36K_{U}\big)\Big)=0,
$$
where $B_{1}$ and $B_{2}$ are general surfaces in $\mathcal{B}$.
The latter is impossible by Corollary~\ref{corollary:Cheltsov}.
\end{proof}

\begin{corollary}
\label{corollary:n-74-points-P5} The set $\mathbb{CS}(U,
\frac{1}{k}\mathcal{D})$ consists of the point $P_{5}$ by
Lemmas~\ref{lemma:Noether-Fano} and \ref{lemma:Cheltsov-Kawamata}.
\end{corollary}

Let $F$ be the exceptional divisor of the morphism $\beta$,
$\mathcal{H}$ be the proper transform of the linear system
$\mathcal{M}$ on the variety $W$, and $P_{6}$ and $P_{7}$ be the
singular points of  $W$ that are quotient singularities of types
$\frac{1}{3}(1,1,2)$ and $\frac{1}{7}(1,3,4)$ contained in $F$
respectively. Then $\mathcal{H}\sim_{\mathbb{Q}}-kK_{W}$.

Suppose that the set $\mathbb{CS}(W, \frac{1}{k}\mathcal{H})$
contains the point $P_{7}$. Then it easily follows from the proof
of Lemma~\ref{lemma:n-48-points-P8} that the claim of
Theorem~\ref{theorem:main} holds for the hypersurface $X$. Hence,
we may assume that the set $\mathbb{CS}(W,
\frac{1}{k}\mathcal{H})$ does not contain the point $P_{7}$. Thus,
the set $\mathbb{CS}(W, \frac{1}{k}\mathcal{H})$ consists of the
singular point $P_{6}$ by Lemmas~\ref{lemma:Noether-Fano} and
\ref{lemma:Cheltsov-Kawamata}

Let $\pi\colon Z\to W$ be the weighted blow up of the point
$P_{6}$ with weights $(1,1,2)$, $G$ be the exceptional  divisor of
the birational morphism $\pi$, and $\mathcal{B}$ be the proper
transform of the linear system $\mathcal{M}$ on the variety $Z$.
Then $\mathcal{B}\sim_{\mathbb{Q}}-kK_{Z}$ by
Theorem~\ref{theorem:Kawamata}.

Let $D$ be a divisor on $Z$ such that the equivalence
$$
D\sim_{\mathbb{Q}}
-4K_{Z}-(\beta\circ\pi)^{*}\big(20K_{U}\big)-\pi^{*}\big(24K_{W}\big)
$$
holds. The divisor $D$ is nef and big, but $D\cdot B_{1}\cdot
B_{2}=0$, where $B_{1}$ and $B_{2}$ are general surfaces of the
linear system $\mathcal{B}$, which is impossible by
Corollary~\ref{corollary:Cheltsov}.

\section{Case $n=79$, hypersurface of degree $33$ in $\mathbb{P}(1,3,5,11,14)$.}
\label{section:n-79}

We use the notations and assumptions of
Section~\ref{section:start}. Let $n=79$. Then $X$ is a
hypersurface of degree $33$ in $\mathbb{P}(1,3,5,11,14)$, the
equality $-K_{X}^{3}=1/70$ holds, and the singularities of $X$
consist of the points $P_{1}$ and $P_{2}$ that are singularities
of type $\frac{1}{5}(1,1,4)$ and $\frac{1}{14}(1,3,11)$
respectively.

\begin{proposition}
\label{proposition:n-79} The claim of Theorem~\ref{theorem:main}
holds for $n=79$.
\end{proposition}

\begin{proof}
There is a commutative diagram
$$
\xymatrix{
&U\ar@{->}[d]_{\alpha}&&W\ar@{->}[ll]_{\beta}&&Y\ar@{->}[ll]_{\gamma}\ar@{->}[d]^{\eta}&\\%
&X\ar@{-->}[rrrr]_{\psi}&&&&\mathbb{P}(1,3,5),&}
$$
where $\psi$ is a projection, $\alpha$ is the weighted blow up of
$P_{2}$ with weights $(1,3,11)$, $\beta$ is the weighted blow up
with weights $(1,3,8)$ of the singular point of type
$\frac{1}{11}(1,3,8)$ contained in the exceptional divisor of
$\alpha$, $\gamma$ is the weighted blow up with weights $(1,3,5)$
of the singular point of type $\frac{1}{8}(1,3,5)$ contained in
the exceptional divisor of  $\beta$, and $\eta$ is an elliptic
fibration.

It follows from Theorem~\ref{theorem:smooth-points},
Lemma~\ref{lemma:Ryder} and
Proposition~\ref{proposition:singular-points} that $\mathbb{CS}(X,
\frac{1}{k}\mathcal{M})=\{P_{2}\}$.

Let $E$ be the exceptional divisor of the morphism $\alpha$,
$\mathcal{D}$ be the proper transform of the linear system
$\mathcal{M}$ on $U$, and $P_{3}$ and $P_{4}$ be the singular
points of $U$ that are quotient singularities of types
$\frac{1}{3}(1,1,2)$ and $\frac{1}{11}(1,3,8)$ contained in $E$
respectively. Then $\mathcal{D}\sim_{\mathbb{Q}}-kK_{U}$, and it
follows from Lemmas~\ref{lemma:Noether-Fano} and
\ref{lemma:Cheltsov-Kawamata} and the proof of
Lemma~\ref{lemma:n-74-points-P4} that $\mathbb{CS}(U,
\frac{1}{k}\mathcal{D})=\{P_{4}\}$.

Let $F$ be the exceptional  divisor of the morphism $\beta$,
$\mathcal{H}$ be the proper transform of the linear system
$\mathcal{M}$ on the variety $W$, and $P_{5}$ and $P_{6}$ be the
singular points of $W$ that are quotient singularities of types
$\frac{1}{3}(1,1,2)$ and $\frac{1}{8}(1,3,4)$ contained in $F$ .
Then $\mathcal{H}\sim_{\mathbb{Q}}-kK_{W}$.

In the case when the set $\mathbb{CS}(W, \frac{1}{k}\mathcal{H})$
contains the point $P_{6}$, it  easily follows from
Theorem~\ref{theorem:Kawamata} that the claim of
Theorem~\ref{theorem:main} holds for the hypersurface $X$.
Therefore, we may assume that the set $\mathbb{CS}(W,
\frac{1}{k}\mathcal{H})$ consists of the point $P_{5}$ by
Lemmas~\ref{lemma:Noether-Fano} and \ref{lemma:Cheltsov-Kawamata}.

Let $\pi\colon Z\to W$ be the weighted blow up of t $P_{5}$ with
weights $(1,1,2)$, $G$ be the exceptional divisor of $\pi$, and
$\mathcal{B}$ be the proper transform of $\mathcal{M}$ on  $Z$.
Then $\mathcal{B}\sim_{\mathbb{Q}}-kK_{Z}$ by
Theorem~\ref{theorem:Kawamata}.

The hypersurface $X$ can be given by the equation
$$
w^{2}z+wf_{19}(x,y,z,t)+f_{33}(x,y,z,t)=0\subset\mathbb{P}(1,2,3,8,11)\cong\mathrm{Proj}\Big(\mathbb{C}[x,y,z,t,w]\Big),
$$
where $\mathrm{wt}(x)=1$, $\mathrm{wt}(y)=3$, $\mathrm{wt}(z)=5$,
$\mathrm{wt}(t)=11$, $\mathrm{wt}(w)=14$, and $f_{i}(x,y,z,t)$ is
a  quasihomogeneous polynomial  of degree $i$. Let $\mathcal{P}$
be the linear system on the hypersurface  $X$ that is generated by
the monomials $x^{30}$, $y^{10}$, $z^{6}$, $t^{2}x^{8}$,
$t^{2}y^{2}x^{2}$, $ty^{6}x$ and $wtz$, $\mathcal{R}$ be the
proper transform of $\mathcal{P}$ on the variety $Z$, and $R$ be a
general surface in $\mathcal{R}$. Then $R$ is nef and big, but
$$
R\sim_{\mathbb{Q}}(\alpha\circ\beta\circ\pi)^{*}\big(-30K_{X}\big)-\frac{30}{11}(\beta\circ\pi)^{*}(E)-\frac{8}{11}\pi^{*}(F)-\frac{2}{3}G,
$$
which implies that the equality $R\cdot B_{1}\cdot B_{2}=0$ holds,
where $B_{1}$ and $B_{2}$ are general surfaces of the linear
system $\mathcal{B}$. The latter contradicts
Corollary~\ref{corollary:Cheltsov}.
\end{proof}

\section{Case $n=80$, hypersurface of degree $34$ in $\mathbb{P}(1,3,4,10,17)$.}
\label{section:n-80}

We use the notations and assumptions of
Section~\ref{section:start}. Let $n=80$. Then $X$ is a
hypersurface of degree $34$ in $\mathbb{P}(1,3,4,10,17)$, whose
singularities consist of the point $P_{1}$ that is a singularity
of type $\frac{1}{2}(1,1,1)$, the point $P_{2}$ that is a
singularity of type $\frac{1}{3}(1,1,2)$, the point $P_{3}$ that
is a singularity of type $\frac{1}{4}(1,1,3)$, and the point
$P_{4}$ that is a singularity of type $\frac{1}{10}(1,3,7)$.

\begin{proposition}
\label{proposition:n-80} The claim of Theorem~\ref{theorem:main}
holds for $n=80$.
\end{proposition}

\begin{proof}
It follows from Proposition~\ref{proposition:singular-points} that
$\mathbb{CS}(X, \frac{1}{k}\mathcal{M})=\{P_{4}\}$. There is a
commutative diagram
$$
\xymatrix{
&U\ar@{->}[d]_{\alpha}&&Y\ar@{->}[ll]_{\beta}\ar@{->}[d]^{\eta}&\\%
&X\ar@{-->}[rr]_{\psi}&&\mathbb{P}(1,3,4),&}
$$
where $\psi$ is a projection, $\alpha$ is the weighted blow up of
$P_{4}$ with weights $(1,3,7)$, $\beta$ is the weighted blow up
with weights $(1,3,4)$ of the singular point the variety $U$ that
is a quotient singularity of type $\frac{1}{7}(1,3,4)$ contained
in the exceptional divisor of $\alpha$, and $\eta$ is an elliptic
fibration.

Let $\mathcal{D}$ be the proper transform of the linear system
$\mathcal{M}$ on the variety $U$, and $P_{5}$ and $P_{6}$ be the
singular points of the variety $U$ that are quotient singularities
of types $\frac{1}{3}(1,1,2)$ and $\frac{1}{4}(1,1,3)$ contained
in  the exceptional divisor of  $\alpha$ respectively. Then the
equivalence $\mathcal{D}\sim_{\mathbb{Q}}-kK_{U}$.

It follows from Lemmas~\ref{lemma:Noether-Fano} and
\ref{lemma:Cheltsov-Kawamata} that either the set $\mathbb{CS}(U,
\frac{1}{k}\mathcal{D})$ contains the point $P_{5}$, or the set
$\mathbb{CS}(U, \frac{1}{k}\mathcal{D})$ contains the point
$P_{6}$. In the latter case it follows from
Theorem~\ref{theorem:Kawamata} that the claim of
Theorem~\ref{theorem:main} holds for $X$. Therefore, we may assume
$P_{5}\in\mathbb{CS}(U, \frac{1}{k}\mathcal{D})$.

The hypersurface $X$ can be given by the quasihomogeneous equation
$$
t^{3}z+t^{2}f_{14}(x,y,z,w)+tf_{24}(x,y,z,w)+f_{34}(x,y,z,w)=0\subset\mathrm{Proj}\Big(\mathbb{C}[x,y,z,t,w]\Big),
$$
where $\mathrm{wt}(x)=1$, $\mathrm{wt}(y)=3$, $\mathrm{wt}(z)=4$,
$\mathrm{wt}(t)=10$, $\mathrm{wt}(w)=17$, and $f_{i}(x,y,z,w)$ is
a general quasihomogeneous polynomial  of degree $i$. Let
$\mathcal{P}$ be a pencil consisting of the surfaces that are cut
on  $X$ by the equations $\lambda x^{4}+\mu z=0$, where
$(\lambda,\mu)\in\mathbb{P}^{1}$. Then the base locus of
$\mathcal{P}$ consists of the irreducible curve that is cut on $X$
by the equations $x=z=0$.

Let $\gamma\colon W\to U$ be the weighted blow up of $P_{5}$ with
weights $(1,1,2)$, $\mathcal{B}$ be the proper transform of the
linear system $\mathcal{M}$ on the variety $W$, $\mathcal{H}$ be
the proper transform of the pencil $\mathcal{P}$ on the variety
$W$, $D$ be a sufficiently general surface of the pencil
$\mathcal{H}$, and $E$ and $F$ be the exceptional divisors of the
morphisms $\alpha$ and $\gamma$ respectively. Then the surface $D$
is normal, the equivalences
$$
D\sim_{\mathbb{Q}}-4K_{W}\sim_{\mathbb{Q}}(\alpha\circ\gamma)^{*}\big(-4K_{X}\big)-\frac{2}{5}\gamma^{*}(E)-\frac{4}{3}F\\
$$
hold, and the base locus of the pencil $\mathcal{H}$ consists of
the curves $C$ and $L$ such that $\alpha\circ\gamma(C)$ is the
base curve of the pencil $\mathcal{P}$, and $\gamma(L)$ is the
curve on the surface $E\cong\mathbb{P}(1,3,7)$ that is contained
in the linear system $|\mathcal{O}_{\mathbb{P}(1,\,3,\,7)}(1)|$.

Let $S$ be the surface of the linear system $|-K_{W}|$, and
$\bar{E}$ be the proper transform of the surface $E$ on the
variety $W$. Then $S\cdot D=C+L$ and $\bar{E}\cdot D=4L$, which
implies that
$$
C\cdot C=-\frac{1}{3},\ C\cdot L=0,\ L\cdot L=-\frac{2}{7}
$$
on the surface $D$. Hence, the intersection form of the curves $C$
and $L$ on the surface~$D$~is~nega\-tively defined, but
$\mathcal{B}\vert_{D}\sim_{\mathbb{Q}} kC+kL$, which contradicts
Lemmas~\ref{lemma:normal-surface} and \ref{lemma:Cheltsov-II}.
\end{proof}

\section{Case $n=82$, hypersurface of degree $36$ in $\mathbb{P}(1,1,5,12,18)$.}
\label{section:n-82}

We use the notations and assumptions of
Section~\ref{section:start}. Let $n=82$. Then $X$ is a
hypersurface of degree $36$ in $\mathbb{P}(1,1,5,12,18)$, the
equality $-K_{X}^{3}=1/30$ holds, and the singularities of $X$
consist of the point $P_{1}$ that is a quotient singularity of
type $\frac{1}{5}(1,2,3)$, and the point $P_{2}$ that is a
quotient singularity of type $\frac{1}{6}(1,1,5)$. There is a
commutative diagram
$$
\xymatrix{
&&U\ar@{->}[dl]_{\alpha}\ar@{->}[dr]^{\eta}&&\\
&X\ar@{-->}[rr]_{\psi}&&\mathbb{P}(1,1,5),&}
$$
where $\psi$ is a projection, $\alpha$ is the weighted blow up of
$P_{2}$ with weights $(1,1,6)$, and $\eta$ is an elliptic
fibration. The hypersurface $X$ is birationally superrigid.

\begin{proposition}
\label{proposition:n-82} The claim of Theorem~\ref{theorem:main}
holds for $n=82$.
\end{proposition}

\begin{proof}
Suppose that $P_{2}\in\mathbb{CS}(X, \frac{1}{k}\mathcal{M})$.
Then Theorem~\ref{theorem:Kawamata} implies the existence of the
commutative diagram
$$
\xymatrix{
&X\ar@{-->}[d]_{\psi}\ar@{-->}[rr]^{\rho}&&V\ar@{->}[d]^{\nu}&\\%
&\mathbb{P}(1,1,5)\ar@{-->}[rr]_{\sigma}&&\mathbb{P}^{2},&}
$$
where $\sigma$ is a birational map. Thus, it follows from
Theorem~\ref{theorem:smooth-points}, Lemma~\ref{lemma:Ryder} and
Proposition~\ref{proposition:singular-points} that we may assume
that the set $\mathbb{CS}(X, \frac{1}{k}\mathcal{M})$ consists of
the point $P_{1}$.

Let $\pi\colon W\to X$ be the weighted blow up of $P_{1}$ with
weights $(1,2,3)$, and $\mathcal{B}$ be the proper transform of
$\mathcal{M}$ on $W$. Then $\mathcal{B}\sim_{\mathbb{Q}}-kK_{W}$
by Theorem~\ref{theorem:Kawamata}, and the singularities of the
exceptional divisor of the morphism $\pi$ consist of the points
$Q$ and $O$ that are quotient singularities of types
$\frac{1}{2}(1,1,1)$ and $\frac{1}{3}(1,1,2)$ on the variety $W$
respectively.

Suppose that the set $\mathbb{CS}(W, \frac{1}{k}\mathcal{B})$ is
not empty. Then the set $\mathbb{CS}(W, \frac{1}{k}\mathcal{B})$
contains either the point $O$, or the point $Q$. On the other
hand, it easily follows from the proof of
Proposition~\ref{proposition:n-29} that the set $\mathbb{CS}(W,
\frac{1}{k}\mathcal{B})$ does not contain neither the point $O$,
nor the point $Q$, which implies that the singularities of the log
pair $(W, \frac{1}{k}\mathcal{B})$ are terminal.

It follows from the proof of Proposition~\ref{proposition:n-43}
that there is a birational map $\gamma\colon W\dasharrow Y$ such
that $\gamma$ is an antiflip, the divisor $-K_{Y}$ is nef and big,
and the linear system $|-rK_{Y}|$ induced a birational map
$Y\dasharrow X^{\prime}$ such that $X^{\prime}$ is a hypersurface
in $\mathbb{P}(1,1,6,14,21)$ of degree $42$ with canonical
singularities (see the proof of Theorem~5.5.1 in \cite{CPR}),
where $r\gg 0$.

Let $\mathcal{H}$ be the proper transform of $\mathcal{M}$ on the
variety $Y$. Then $\mathcal{H}\sim_{\mathbb{Q}} -kK_{Y}$, because
$\gamma$ is an isomorphism in codimension one. The map $\gamma$ is
a log flip with respect to $(W, \lambda\mathcal{B})$ for some
rational number $\lambda>1/k$. Thus, the singularities of the
mobile log pair $(Y, \frac{1}{k}\mathcal{H})$ are terminal, which
contradicts Lemma~\ref{lemma:Noether-Fano}.
\end{proof}

\section{Case $n\in\{21, 24, 33, 35, 41, 42, 46, 50, 54, 55, 61, 62, 63, 67, 69, 71, 76, 77, 83, 85, 91\}$.} %
\label{section:n-21-24-85-91}

We use the notations and assumptions of
Section~\ref{section:start}.

\begin{proposition}
\label{proposition:n-21-24-85-91} Suppose that  $n\in\{21, 24, 33,
35, 41, 42, 46, 50, 54, 55, 61, 62, 63, 67, 69, 71, 76,\allowbreak
77, 83, 85, 91\}$. Then there is a commutative diagram
\begin{equation}
\label{equation:n-21-24-85-91-commutative-diagram} \xymatrix{
&X\ar@{-->}[d]_{\psi}\ar@{-->}[rr]^{\rho}&&V\ar@{->}[d]^{\nu}&\\%
&\mathbb{P}(1,a_{1},a_{2})\ar@{-->}[rr]_{\sigma}&&\mathbb{P}^{2},&}
\end{equation}
where $\psi$ is the natural projection, and $\sigma$ is a
birational map.
\end{proposition}

\begin{proof}
It follows from Theorem~\ref{theorem:smooth-points} and
Lemma~\ref{lemma:Ryder} that the set $\mathbb{CS}(X,
\frac{1}{k}\mathcal{M})$ consists of singular points of the
hypersurface $X$. We prove the existence of the
diagram~\ref{equation:n-21-24-85-91-commutative-diagram} case by
case.

\medskip

\texttt{Case $n=21$.}

The variety $X$ is a  hypersurface in $\mathbb{P}(1,1,2,4,7)$ of
degree $14$, the equality $-K_{X}^{3}=1/4$ holds, and the
singularities of the hypersurface $X$ consist of  the points
$O_{1}$ and $O_{2}$ that are quotient singularities of type
$\frac{1}{2}(1,1,1)$, and the point $P$ that is a quotient
singularity of type $\frac{1}{4}(1,1,3)$.

There is a commutative diagram
$$
\xymatrix{
&U\ar@{->}[d]_{\alpha}&&W\ar@{->}[ll]_{\beta}\ar@{->}[d]^{\eta}&\\%
&X\ar@{-->}[rr]_{\psi}&&\mathbb{P}(1,1,2),&}
$$
where $\psi$ is a projection, $\alpha$ is the weighted blow up of
$P$ with weights $(1,1,3)$, $\beta$ is the weighted blow up with
weights $(1,1,2)$ of the singular point of $U$ that is a
singularity of type $\frac{1}{3}(1,1,2)$, and $\eta$ is an
elliptic fibration. It follows from
Proposition~\ref{proposition:singular-points} that $\mathbb{CS}(X,
\frac{1}{k}\mathcal{M})=\{P\}$.

Let $\mathcal{D}$ be the proper transform of $\mathcal{M}$ on $U$.
Then $\mathcal{D}\sim_{\mathbb{Q}}-kK_{U}$  by
Theorem~\ref{theorem:Kawamata}, and it follows from
Lemmas~\ref{lemma:Noether-Fano} and \ref{lemma:Cheltsov-Kawamata}
that the set $\mathbb{CS}(U, \frac{1}{k}\mathcal{D})$ contains the
singular point of $U$ that is contained in the exceptional divisor
of the morphism $\alpha$.

Let $\mathcal{B}$ be the proper transform of $\mathcal{M}$ on $Y$.
Then $\mathcal{B}\sim_{\mathbb{Q}}-kK_{Y}$ by
Theorem~\ref{theorem:Kawamata}, which implies the existence of the
the commutative
diagram~\ref{equation:n-21-24-85-91-commutative-diagram}.

\medskip

\texttt{Case $n=24$.}

The variety $X$ is a hypersurface in $\mathbb{P}(1,1,2,5,7)$ of
degree $15$, the equality $-K_{X}^{3}=3/14$ holds, and the
singularities of the hypersurface $X$ consist of the point $P_{1}$
that is a quotient singularity of type $\frac{1}{2}(1,1,1)$, and
the point $P_{2}$ that is a quotient singularity of type
$\frac{1}{7}(1,2,5)$.

It follows from Proposition~\ref{proposition:singular-points} that
$\mathbb{CS}(X, \frac{1}{k}\mathcal{M})=\{P_{2}\}$. There is a
commutative diagram
$$
\xymatrix{
&U\ar@{->}[d]_{\alpha}&&W\ar@{->}[ll]_{\beta}&&Y\ar@{->}[ll]_{\gamma}\ar@{->}[d]^{\eta}&\\%
&X\ar@{-->}[rrrr]_{\psi}&&&&\mathbb{P}(1,1,2),&}
$$
where $\alpha$ is the weighted blow up of the point $P$ with
weights $(1,2,5)$,  $\beta$ is the weighted blow up with weights
$(1,2,3)$ of the singular point of $U$ that is a singularity of
type $\frac{1}{5}(1,2,3)$, $\gamma$ is the weighted blow up with
weights $(1,1,2)$ of the singular point of the variety $W$ that is
a quotient singularity of type $\frac{1}{3}(1,1,2)$, and $\eta$ is
an elliptic fibration.

The proof of Proposition~\ref{proposition:n-29} implies that
$\mathcal{D}\sim_{\mathbb{Q}}-kK_{W}$, where $\mathcal{D}$ is the
proper transform of the linear system $\mathcal{M}$ on the variety
$W$, and the set $\mathbb{CS}(W, \frac{1}{k}\mathcal{D})$ does not
contain subvarieties of the variety $W$ that are not contained in
the exceptional divisor of the morphism $\beta$.

We can apply arguments of the proof of
Proposition~\ref{proposition:n-29} to the log pair $(W,
\frac{1}{k}\mathcal{D})$ to prove that the set $\mathbb{CS}(W,
\frac{1}{k}\mathcal{D})$ contains the singular point of $W$ that
is a singularity of type $\frac{1}{3}(1,1,2)$, which implies the
existence of the commutative
diagram~\ref{equation:n-21-24-85-91-commutative-diagram} by
Theorem~\ref{theorem:Kawamata}.

\medskip

\texttt{Case $n=33$.}

The variety $X$ is a hypersurface in $\mathbb{P}(1,2,3,5,7)$ of
degree $17$, the singularities of $X$ consist of the point $P_{1}$
that is a quotient singularity of type $\frac{1}{2}(1,1,1)$, the
point $P_{2}$ that is a quotient singularity of type
$\frac{1}{3}(1,1,2)$, the point $P_{3}$ that is a quotient
singularity of type $\frac{1}{5}(1,2,3)$, and the point $P_{4}$
that is a  singularity of type $\frac{1}{7}(1,2,5)$. There is a
commutative diagram
$$
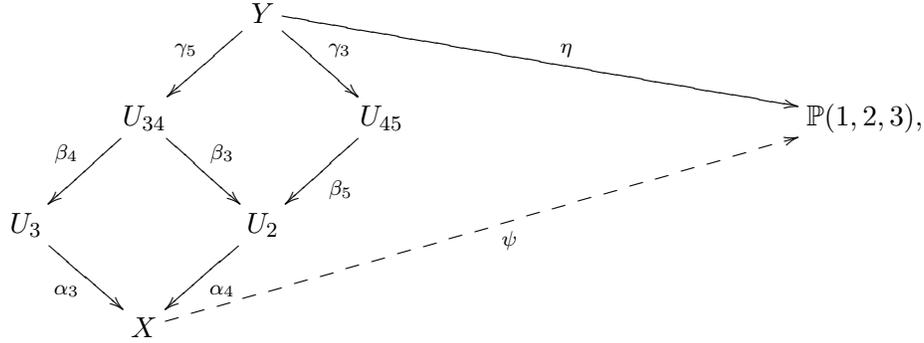
\xymatrix{
&&&Y\ar@{->}[dl]_{\gamma_{5}}\ar@{->}[dr]^{\gamma_{3}}\ar@{->}[drrrrrr]^{\eta}&&&&&&\\
&&U_{34}\ar@{->}[dl]_{\beta_{4}}\ar@{->}[dr]^{\beta_{3}}&&U_{45}\ar@{->}[dl]^{\beta_{5}}&&&&&\mathbb{P}(1,2,3),\\
&U_{3}\ar@{->}[dr]_{\alpha_{3}}&&U_{2}\ar@{->}[dl]^{\alpha_{4}}&&&&&&\\
&&X\ar@{-->}[rrrrrrruu]_{\psi}&&&&&&&}
$$
where $\alpha_{3}$ is the weighted blow up of the point $P_{3}$
with weights $(1,2,3)$, $\alpha_{4}$ is the weighted blow up of
the point $P_{4}$ with weights $(1,2,5)$, $\beta_{4}$ is the
weighted blow up with weights $(1,2,5)$ of the proper transform of
$P_{4}$ on $U_{3}$, $\beta_{3}$ is the weighted blow up with
weights $(1,2,3)$ of the proper transform of $P_{3}$ on $U_{4}$,
$\beta_{5}$ is the weighted blow up with weights $(1,2,3)$ of the
singular point of the variety $U_{4}$ that is a quotient
singularity of type $\frac{1}{5}(1,2,3)$ contained  in the
exceptional divisor of the morphism $\alpha_{4}$, $\gamma_{3}$ is
the weighted blow up with weights $(1,2,3)$ of the proper
transform of the point $P_{3}$ on the variety $U_{45}$,
$\gamma_{5}$ is the weighted blow up with weights $(1,2,3)$ of the
singular point of $U_{34}$ that is a quotient singularity of type
$\frac{1}{5}(1,2,3)$ contained in the exceptional divisor of the
morphism $\beta_{4}$, and $\eta$ is an elliptic fibration.

It follows from Proposition~\ref{proposition:singular-points} and
the proof of Proposition~\ref{proposition:n-23} that
$\mathbb{CS}(X, \frac{1}{k}\mathcal{M})=\{P_{3}, P_{4}\}$.

Let $\mathcal{D}_{34}$ be the proper transform of $\mathcal{M}$ on
$U_{34}$, and $\bar{P}_{5}$ and $\bar{P}_{6}$ be the singular
points of the variety $U_{34}$ that are quotient singularities of
types $\frac{1}{5}(1,2,3)$ and $\frac{1}{2}(1,1,1)$ contained in
the exceptional divisor of the morphism $\beta_{4}$ respectively.
Then $\mathcal{D}_{34}\sim_{\mathbb{Q}}-kK_{U_{34}}$ by
Theorem~\ref{theorem:Kawamata}.

It follows from Lemma~\ref{lemma:Cheltsov-Kawamata} and the proof
of Proposition~\ref{proposition:n-23} that
$\bar{P}_{5}\in\mathbb{CS}(U_{34}, \frac{1}{k}\mathcal{D}_{34})$,
and the existence of the commutative
diagram~\ref{equation:n-21-24-85-91-commutative-diagram} follows
from Theorem~\ref{theorem:Kawamata}.

\medskip

\texttt{Case $n=35$.}

The variety $X$ is a  hypersurface in $\mathbb{P}(1,1,3,5,9)$ of
degree $18$, whose singularities consist of the points $P_{1}$ and
$P_{2}$ that are singularities of type $\frac{1}{3}(1,1,2)$, and
the point $P_{3}$ that is a singularity of type
$\frac{1}{5}(1,1,4)$. There is a commutative diagram
$$
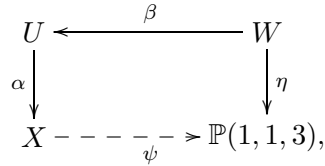
\xymatrix{
&U\ar@{->}[d]_{\alpha}&&W\ar@{->}[ll]_{\beta}\ar@{->}[d]^{\eta}&\\%
&X\ar@{-->}[rr]_{\psi}&&\mathbb{P}(1,1,3),&}
$$
where $\alpha$ is the blow up of $P_{3}$ with weights $(1,1,4)$,
$\beta$ is the  blow up with weights $(1,1,3)$ of the point of $U$
that is  singularity of type $\frac{1}{4}(1,1,3)$, and $\eta$ is
an elliptic fibration.

It follows from Proposition~\ref{proposition:singular-points} that
$\mathbb{CS}(X, \frac{1}{k}\mathcal{M})=\{P_{3}\}$.

Let $\mathcal{D}$ is the proper transform $\mathcal{M}$ on $U$.
Then $\mathcal{D}\sim_{\mathbb{Q}}-kK_{U}$ by
Theorem~\ref{theorem:Kawamata}, and it follows from
Lemma~\ref{lemma:Cheltsov-Kawamata} that $\mathbb{CS}(U,
\frac{1}{k}\mathcal{D})$ contains the singular point of $U$ that
is contained in the exceptional divisor of $\alpha$. The existence
of the diagram~\ref{equation:n-21-24-85-91-commutative-diagram}
follows from Theorem~\ref{theorem:Kawamata}.

\medskip

\texttt{Case $n=41$.}

The variety $X$ is a hypersurface in $\mathbb{P}(1,1,4,5,10)$ of
degree $20$, whose singularitiesconsist of  the point $O$ that is
a singularity of type $\frac{1}{2}(1,1,1)$, the points $P_{1}$ and
$P_{2}$ of type $\frac{1}{5}(1,1,4)$.

The indeterminacies of  $\psi$ are resolved by blowing up $P_{1}$
and $P_{2}$ with weights $(1,1,4)$, it follows from
Proposition~\ref{proposition:singular-points} that $\mathbb{CS}(X,
\frac{1}{k}\mathcal{M})\subseteq\{P_{1},P_{2}\}$, and the
existence of the
diagram~\ref{equation:n-21-24-85-91-commutative-diagram} follows
from the proof of Proposition~\ref{proposition:n-15}.

\medskip

\texttt{Case $n=42$.}

The variety $X$ is a hypersurface in $\mathbb{P}(1,2,3,5,10)$ of
degree $20$, the singularities of $X$ consist of the points
$P_{1}$, $P_{2}$ and $P_{3}$ that are quotient singularities of
type $\frac{1}{2}(1,1,1)$,  the point $P_{4}$ that is a
singularity of type $\frac{1}{3}(1,1,2)$, the points $P_{5}$ and
$P_{6}$, that are singularities of types $\frac{1}{5}(1,2,3)$.

The indeterminacies of  $\psi$ are resolved by blowing up $P_{5}$
and $P_{6}$ with weights $(1,2,3)$,

It follows from Proposition~\ref{proposition:singular-points} that
$\mathbb{CS}(X, \frac{1}{k}\mathcal{M})\subseteq\{P_{5},P_{6}\}$,
but the existence of the
diagram~\ref{equation:n-21-24-85-91-commutative-diagram} is
obvious if $\mathbb{CS}(X,
\frac{1}{k}\mathcal{M})=\{P_{5},P_{6}\}$. We may assume that
$\mathbb{CS}(X, \frac{1}{k}\mathcal{M})=\{P_{5}\}$.

Let $\alpha\colon U\to X$ be the weighted blow up of the point
$P_{5}$ with weights $(1,2,3)$, $\mathcal{B}$ be the proper
transform of the linear system $\mathcal{M}$ on $U$, and $O$ and
$Q$ be the singular points of $U$ that are quotient singularities
of types $\frac{1}{3}(1,1,2)$ and $\frac{1}{2}(1,1,1)$ contained
in the exceptional divisor of the morphism $\alpha$ respectively.
Then it follows from the proof of
Lemma~\ref{lemma:n-27-P2-is-a-center} that $\mathbb{CS}(U,
\frac{1}{k}\mathcal{B})=\{Q\}$.

Let $\zeta\colon Z\to U$ be the weighted blow up of the point $Q$
with weights $(1,1,1)$, $\mathcal{D}$ be the proper transform of
the linear system $\mathcal{M}$ on the variety $Z$, and
$\mathcal{H}$ be the proper transform of the linear system
$|-3K_{X}|$ on the variety $Z$. Then
$\mathcal{D}\sim_{\mathbb{Q}}-kK_{Z}$ by
Theorem~\ref{theorem:Kawamata}, and the base locus of the linear
system $\mathcal{H}$ consists of a curve $C$ such that
$\alpha\circ\zeta(C)$ is the base curve  of $|-3K_{X}|$.

Let $S$ be a general surface in $\mathcal{H}$. Then $S$ is normal,
and the inequality  $C^{2}<0$ holds on the surface $S$, but the
equivalence $\mathcal{D}\vert_{S}\sim_{\mathbb{Q}} kC$ holds,
which contradicts Lemmas~\ref{lemma:Cheltsov-II} and
\ref{lemma:normal-surface}.

\medskip

\texttt{Case $n=45$.}

The variety $X$ is a hypersurface in $\mathbb{P}(1,3,4,5,8)$ of
degree $20$, whose singularities consist of the point $P_{1}$ of
type $\frac{1}{3}(1,1,2)$, the points $P_{2}$ and $P_{3}$ of type
$\frac{1}{4}(1,1,3)$,~the~point~$P_{4}$~of~type~$\frac{1}{8}(1,3,5)$.

There is a commutative diagram
$$
\xymatrix{
&U\ar@{->}[d]_{\alpha}&&Y\ar@{->}[ll]_{\beta}\ar@{->}[d]^{\eta}&\\%
&X\ar@{-->}[rr]_{\psi}&&\mathbb{P}(1,3,4),&}
$$
where $\alpha$ is the weighted blow up of the point $P_{4}$ with
weights $(1,3,5)$, $\beta$ is the weighted blow up with weights
$(1,2,3)$ of the singular point of $U$ that is a quotient
singularity of type $\frac{1}{5}(1,2,3)$ contained in the
exceptional divisor of the morphism $\alpha$, and $\eta$ is an
elliptic fibration.

It follows from Proposition~\ref{proposition:singular-points} that
$\mathbb{CS}(X, \frac{1}{k}\mathcal{M})=\{P_{4}\}$.

Let $\mathcal{B}$ be the proper transform of the linear system
$\mathcal{M}$ on the variety $U$, $\bar{P}_{5}$ be the singular
point of the variety $U$ that is a quotient singularity of type
$\frac{1}{5}(1,2,3)$ contained in the exceptional divisor of the
morphism $\alpha$, and $\bar{P}_{6}$ be the singular point of the
the variety $U$ that is a quotient singularity of type
$\frac{1}{2}(1,1,1)$ contained in the exceptional divisor of the
morphism $\alpha$.

It follows from Lemma~\ref{lemma:Cheltsov-Kawamata} that
$\mathbb{CS}(U, \frac{1}{k}\mathcal{B})\cap\{\bar{P}_{5},
\bar{P}_{6}\}\ne\varnothing$, but the proof of
Proposition~\ref{proposition:n-44} easily implies that
$\bar{P}_{5}\in\mathbb{CS}(U, \frac{1}{k}\mathcal{B})$, which
implies the existence of the commutative
diagram~\ref{equation:n-21-24-85-91-commutative-diagram}.

\medskip

\texttt{Case $n=46$.}

The variety $X$ is a  hypersurface in $\mathbb{P}(1,1,3,7,10)$ of
degree $21$, the singularities of $X$ consist of the point $P$
that is a singularity of type $\frac{1}{10}(1,3,7)$. There is a
commutative diagram
$$
\xymatrix{
&U\ar@{->}[d]_{\alpha}&&W\ar@{->}[ll]_{\beta}&&Y\ar@{->}[ll]_{\gamma}\ar@{->}[d]^{\eta}&\\%
&X\ar@{-->}[rrrr]_{\psi}&&&&\mathbb{P}(1,1,3),&}
$$
where $\alpha$ is the blow up of $P$ with weights $(1,3,7)$,
$\beta$ is the blow up with weights $(1,3,4)$ of the singular
point of $U$ that is a singularity of type $\frac{1}{7}(1,3,4)$,
$\gamma$ is the blow up with weights $(1,1,3)$ of the singular
point of $W$ that is a singularity of type $\frac{1}{4}(1,1,3)$,
and $\eta$ is an elliptic fibration.

It follows from Proposition~\ref{proposition:singular-points} that
$\mathbb{CS}(X, \frac{1}{k}\mathcal{M})=\{P\}$.

Let $\mathcal{D}$ be the proper transform of the linear system
$\mathcal{M}$ on the variety $W$. Then it follows from the proof
of Proposition~\ref{proposition:n-29} that
$\mathcal{D}\sim_{\mathbb{Q}}-kK_{W}$, and the set $\mathbb{CS}(W,
\frac{1}{k}\mathcal{D})$ does not contain subvarieties of the
variety $W$ that are not contained in the exceptional divisor of
$\beta$.

The set $\mathbb{CS}(W, \frac{1}{k}\mathcal{D})$ is not empty,
because the divisor $-K_{W}$ is nef and big. Therefore, we can
apply the arguments of the proof of
Proposition~\ref{proposition:n-29} the the log pair $(W,
\frac{1}{k}\mathcal{D})$, which implies that the set
$\mathbb{CS}(W, \frac{1}{k}\mathcal{D})$ contains the singular
point of $W$ that is a singularity of type $\frac{1}{4}(1,1,3)$.

The existence of the commutative
diagram~\ref{equation:n-21-24-85-91-commutative-diagram} is
implied by Theorem~\ref{theorem:Kawamata}.

\medskip

\texttt{Case $n=50$.}

The variety $X$ is a hypersurface in $\mathbb{P}(1,1,3,7,11)$ of
degree $22$, the singularities of $X$ consist of the point $P_{1}$
that is a quotient singularity of type $\frac{1}{3}(1,1,2)$, and
the point $P_{2}$ that is a quotient singularity of type
$\frac{1}{7}(1,3,4)$. The equality $-K_{X}^{3}=2/21$ holds. There
is a commutative diagram
$$
\xymatrix{
&U\ar@{->}[d]_{\alpha}&&W\ar@{->}[ll]_{\beta}\ar@{->}[d]^{\eta}&\\%
&X\ar@{-->}[rr]_{\psi}&&\mathbb{P}(1,1,3),&}
$$
where $\alpha$ is the blow up of  $P_{2}$ with weights $(1,3,7)$,
$\beta$ is the blow up with weights $(1,1,3)$ of the singular
point of $U$ that is a singularity of type $\frac{1}{4}(1,1,3)$,
and $\eta$ is an elliptic fibration.

The proof of Proposition~\ref{proposition:n-29} implies the
existence of the
diagram~\ref{equation:n-21-24-85-91-commutative-diagram}.

\medskip

\texttt{Case $n=54$.}

The variety $X$ is a  hypersurface in $\mathbb{P}(1,1,6,8,9)$ of
degree $24$, whose singularities consist of the point $P_{1}$ type
$\frac{1}{2}(1,1,1)$, the point $P_{2}$ of type
$\frac{1}{3}(1,1,2)$, the point $P_{3}$ of type
$\frac{1}{9}(1,1,8)$.

There is a commutative diagram
$$
\xymatrix{
&U\ar@{->}[d]_{\alpha}&&W\ar@{->}[ll]_{\beta}&&Y\ar@{->}[ll]_{\gamma}\ar@{->}[d]^{\eta}&\\%
&X\ar@{-->}[rrrr]_{\psi}&&&&\mathbb{P}(1,1,6),&}
$$
where $\alpha$ is the blow up of $P_{3}$ with weights $(1,1,8)$,
$\beta$ is the blow up with weights $(1,1,7)$ of the singular
point of $U$ that is a singularity of type $\frac{1}{8}(1,1,7)$,
$\gamma$ is the blow up with weights $(1,1,3)$ of the singular
point of $W$ that is a  singularity of type $\frac{1}{7}(1,1,6)$,
and $\eta$ is an elliptic fibration.

The existence of the
diagram~\ref{equation:n-21-24-85-91-commutative-diagram} follows
from the proof of Proposition~\ref{proposition:n-16}.

\medskip

\texttt{Case $n=55$.}

The variety $X$ is a hypersurface in $\mathbb{P}(1,2,3,7,12)$ of
degree $24$, whose singularities consist of the points $P_{1}$ and
$P_{2}$ of type $\frac{1}{2}(1,1,1)$, the points $P_{3}$ and
$P_{4}$ of type $\frac{1}{3}(1,1,2)$, and the singular point
$P_{5}$~of type $\frac{1}{7}(1,2,5)$.  There is a commutative
diagram
$$
\xymatrix{
&U\ar@{->}[d]_{\alpha}&&W\ar@{->}[ll]_{\beta}\ar@{->}[d]^{\eta}&\\%
&X\ar@{-->}[rr]_{\psi}&&\mathbb{P}(1,2,3),&}
$$
where $\alpha$ is the blow up of $P_{5}$ with weights $(1,2,5)$,
$\beta$ is the blow up with weights $(1,2,3)$ of the singular
point of $U$ that is a singularity of type
$\frac{1}{5}(1,2,3)$,~and~$\eta$~is an elliptic fibration.

It follows from Proposition~\ref{proposition:singular-points} that
$\mathbb{CS}(X, \frac{1}{k}\mathcal{M})=\{P_{5}\}$.

Let $\mathcal{B}$ be the proper transform of the linear system
$\mathcal{M}$ on the variety $U$, and $P_{6}$ and $P_{7}$ be the
singular points of $U$ that are singularities of types
$\frac{1}{2}(1,1,1)$ and $\frac{1}{5}(1,2,3)$ contained in the
exceptional divisor of $\alpha$ respectively. Then $\mathbb{CS}(U,
\frac{1}{k}\mathcal{B})\subseteq\{P_{6}, P_{7}\}$ by
Lemmas~\ref{lemma:Noether-Fano}, \ref{lemma:Cheltsov-Kawamata} and
\ref{lemma:curves}.

The existence of the
diagram~\ref{equation:n-21-24-85-91-commutative-diagram} follows
from Theorem~\ref{theorem:Kawamata} in the case when
$P_{7}\in\mathbb{CS}(U, \frac{1}{k}\mathcal{B})$, which implies
that we may assume that the set $\mathbb{CS}(U,
\frac{1}{k}\mathcal{B})$ consists of the point $P_{6}$.

The hypersurface $X$ can be given by the equation
$$
t^{3}z+t^{2}f_{10}(x,y,z,w)+tf_{17}(x,y,z,w)+f_{24}(x,y,z,w)=0\subset\mathrm{Proj}\Big(\mathbb{C}[x,y,z,t,w]\Big),
$$
where $\mathrm{wt}(x)=1$, $\mathrm{wt}(y)=2$, $\mathrm{wt}(z)=3$,
$\mathrm{wt}(t)=7$, $\mathrm{wt}(w)=12$, and $f_{i}(x,y,z,w)$ is a
sufficiently general quasihomogeneous polynomial of degree $i$.
Let $\gamma\colon W\to U$ be the weighted blow up of the point
$P_{6}$ with weights $(1,1,1)$, $\mathcal{H}$ be the proper
transform of the linear system $\mathcal{M}$ on the variety $W$,
$\mathcal{P}$ be the proper transform on the variety $W$ of the
pencil of surfaces that are cut on the hypersurface $X$ by the
equations $\lambda x^{3}+\mu z=0$, where
$(\lambda,\mu)\in\mathbb{P}^{1}$, and $D$ be a sufficiently
general surface of the pencil $\mathcal{P}$. Then the base locus
of the pencil $\mathcal{P}$ consists of the irreducible curves
$C$, $L$ and $\Delta$ such that $\alpha\circ\gamma(C)$ is the base
curve of $|-3K_{X}|$, the curve $\gamma(L)$ is contained in the
exceptional divisor of  $\alpha$, and the curve $\Delta$ is
contained in the exceptional divisor of $\gamma$.

The surface $D$ is normal, the intersection form of the curves $L$
and $C$ on the surface $D$ is negatively defined, but
$\mathcal{H}\vert_{D}\sim_{\mathbb{Q}}kC+kL$, which is impossible
by Lemmas~\ref{lemma:normal-surface} and \ref{lemma:Cheltsov-II}.

\medskip

\texttt{Case $n=61$.}

The variety $X$ is a hypersurface in $\mathbb{P}(1,4,5,7,9)$ of
degree $25$, the singularities of $X$ consist of the point $P_{1}$
that is a quotient singularity of type $\frac{1}{4}(1,1,3)$, the
point $P_{2}$ that is a quotient singularity of type
$\frac{1}{7}(1,2,5)$, and the point $P_{3}$ that is a quotient
singularity of type $\frac{1}{9}(1,4,5)$.

The indeterminacies the map rational $\psi$ is resolved by the
weighted blow ups of the singular points $P_{2}$ and $P_{3}$ with
weights $(1,2,5)$, and $(1,4,5)$ respectively.

It follows from Proposition~\ref{proposition:singular-points} that
$\mathbb{CS}(X, \frac{1}{k}\mathcal{M})\subseteq\{P_{2},P_{3}\}$,
the proof of Lemma~\ref{lemma:n-58-P3-is-a-center} implies that
$\mathbb{CS}(X, \frac{1}{k}\mathcal{M})=\{P_{2}, P_{3}\}$, and the
existence of the
diagram~\ref{equation:n-21-24-85-91-commutative-diagram} follows
from Theorem~\ref{theorem:Kawamata}.

\medskip

\texttt{Case $n=62$.}

The variety $X$ is a hypersurface in $\mathbb{P}(1,1,5,7,13)$ of
degree $26$, the singularities of  $X$ consist of the point
$P_{1}$ that is a quotient singularity of type
$\frac{1}{5}(1,2,3)$, and the point $P_{2}$ that is a singularity
of type $\frac{1}{7}(1,1,6)$. The equality $-K_{X}^{3}=2/35$
holds. There is a commutative diagram
$$
\xymatrix{
&U\ar@{->}[d]_{\alpha}&&Y\ar@{->}[ll]_{\beta}\ar@{->}[d]^{\eta}&\\%
&X\ar@{-->}[rr]_{\psi}&&\mathbb{P}(1,1,5),&}
$$
where $\alpha$ is the weighted blow up of the point $P_{2}$ with
weights $(1,1,6)$, $\beta$ is the weighted blow up with weights
$(1,1,5)$ of the singular point of the variety $U$ that is a
quotient singularity of type $\frac{1}{6}(1,1,5)$, and $\eta$ is
an elliptic fibration.

It follows from Proposition~\ref{proposition:singular-points} that
$\mathbb{CS}(X, \frac{1}{k}\mathcal{M})\subseteq\{P_{1}, P_{2}\}$.
Therefore, the existence of the commutative
diagram~\ref{equation:n-21-24-85-91-commutative-diagram} follows
from the proof of Proposition~\ref{proposition:n-47}.

\medskip

\texttt{Case $n=63$.}

The variety $X$ is a hypersurface in $\mathbb{P}(1,2,3,8,13)$ of
degree $26$, the singularities of the hypersurface $X$ consist of
the points $P_{1}$, $P_{2}$ and $P_{3}$ that are quotient
singularities of type $\frac{1}{2}(1,1,1)$, the point $P_{4}$ that
is a quotient singularity of type $\frac{1}{3}(1,1,2)$, and the
point $P_{5}$ that is a singularity of type $\frac{1}{8}(1,3,5)$.
The equality $-K_{X}^{3}=1/24$ holds.

There is a commutative diagram
$$
\xymatrix{
&U\ar@{->}[d]_{\alpha}&&Y\ar@{->}[ll]_{\beta}\ar@{->}[d]^{\eta}&\\%
&X\ar@{-->}[rr]_{\psi}&&\mathbb{P}(1,2,3),&}
$$
where $\alpha$ is the weighted blow up of the point $P_{5}$ with
weights $(1,3,5)$, and $\beta$ is the weighted blow up with
weights $(1,2,3)$ of the singular point of the variety $U$ that is
a quotient singularity of type $\frac{1}{5}(1,2,3)$ contained in
the exceptional divisor of $\alpha$, and $\eta$ is an elliptic
fibration.

It follows from Proposition~\ref{proposition:singular-points} that
$\mathbb{CS}(X, \frac{1}{k}\mathcal{M})=\{P_{5}\}$.

Let $\mathcal{B}$ be the proper transform of $\mathcal{M}$ on $U$,
and $P_{6}$ and $P_{7}$ be the singular points of $U$ that are
quotient singularities of types $\frac{1}{5}(1,2,3)$ and
$\frac{1}{2}(1,1,1)$ contained in the exceptional divisor of the
morphism $\alpha$ respectively. Then it follows from
Lemma~\ref{lemma:Cheltsov-Kawamata} that the set  $\mathbb{CS}(U,
\frac{1}{k}\mathcal{B})$ contains either the point $P_{6}$, or the
point $P_{7}$.

Suppose that the set $\mathbb{CS}(U, \frac{1}{k}\mathcal{B})$
contains the point $P_{7}$. Let $\gamma\colon W\to U$ be the
weighted blow up of the point $P_{7}$ with weights $(1,1,1)$, and
$\mathcal{H}$ and $\mathcal{P}$ be the proper transforms of the
linear system $\mathcal{M}$ and the pencil $|-2K_{X}|$ on the
variety $W$ respectively. Then the base locus of the pencil
$\mathcal{P}$ consists of the irreducible curves $C$ and $L$ such
that the curve $\alpha\circ\gamma(C)$ is the unique curve in the
base locus of the pencil $|-2K_{X}|$, and the curve $\gamma(L)$ is
contained in the exceptional divisor of the birational morphism
$\alpha$. Moreover, the surface $D$ is normal, and the equalities
$$
C\cdot C=-\frac{5}{12},\ L\cdot L=-\frac{7}{20},\ C\cdot L=\frac{1}{4}%
$$
hold on $D$. Hence, the intersection form of the curves $L$ and
$C$ on the surface $D$ is negatively defined. On the other hand,
the equivalence $\mathcal{H}\vert_{D}\sim_{\mathbb{Q}}kC+kL$ holds
on the surface $D$, which contradicts
Lemmas~\ref{lemma:normal-surface} and \ref{lemma:Cheltsov-II}.

Therefore, the set $\mathbb{CS}(U, \frac{1}{k}\mathcal{B})$
contains the point $P_{6}$. Now the claim of
Theorem~\ref{theorem:Kawamata} easily implies the existence of the
commutative
diagram~\ref{equation:n-21-24-85-91-commutative-diagram}.

\medskip

\texttt{Case $n=67$.}

The variety $X$ is a hypersurface in $\mathbb{P}(1,1,4,9,14)$ of
degree $28$, and the singularities of the hypersurface $X$ consist
of the point $P_{1}$ that is a quotient singularity of type
$\frac{1}{2}(1,1,1)$, and the point $P_{2}$ that is a quotient
singularity of type $\frac{1}{9}(1,4,5)$. The equality
$-K_{X}^{3}=1/18$ holds.

It follows from Proposition~\ref{proposition:singular-points} that
$\mathbb{CS}(X, \frac{1}{k}\mathcal{M})=\{P_{2}\}$.

There is a commutative diagram
$$
\xymatrix{
&U\ar@{->}[d]_{\alpha}&&W\ar@{->}[ll]_{\beta}\ar@{->}[d]^{\eta}&\\%
&X\ar@{-->}[rr]_{\psi}&&\mathbb{P}(1,1,4),&}
$$
where $\alpha$ is the weighted blow up of the point $P_{2}$ with
weights $(1,4,5)$, $\beta$ is the weighted blow up with weights
$(1,1,4)$ of the singular point of $U$ that is quotient
singularity of type $\frac{1}{5}(1,1,4)$, and $\eta$ is an
elliptic fibration.

The existence of the commutative
diagram~\ref{equation:n-21-24-85-91-commutative-diagram} follows
from the proof of Proposition~\ref{proposition:n-29}.

\medskip

\texttt{Case $n=69$.}

The variety $X$ is a hypersurface in $\mathbb{P}(1,4,6,7,11)$ of
degree $28$, and the singularities of the hypersurface $X$ consist
of the points $P_{1}$ and $P_{2}$ that are quotient singularities
of type $\frac{1}{2}(1,1,1)$, and the points $P_{3}$ and $P_{4}$
that are quotient singularities of types $\frac{1}{6}(1,1,5)$ и
$\frac{1}{11}(1,4,7)$ respectively.

The equality $-K_{X}^{3}=1/66$ holds, and there is a commutative
diagram
$$
\xymatrix{
&U\ar@{->}[d]_{\alpha}&&W\ar@{->}[ll]_{\beta}\ar@{->}[d]^{\eta}&\\%
&X\ar@{-->}[rr]_{\psi}&&\mathbb{P}(1,3,4),&}
$$
where $\alpha$ is the weighted blow up of the point $P_{4}$ with
weights $(1,4,7)$, $\beta$ is the weighted blow up with weights
$(1,3,4)$ of the singular point of type $\frac{1}{7}(1,3,4)$, and
$\eta$ is an elliptic fibration.

It follows from Proposition~\ref{proposition:singular-points} that
$\mathbb{CS}(X, \frac{1}{k}\mathcal{M})=\{P_{4}\}$.

Let $\mathcal{B}$ be the proper transform of $\mathcal{M}$ on $U$,
and $P_{5}$ and $P_{6}$ be the singular points of  $U$ that are
quotient singularities of types $\frac{1}{7}(1,3,4)$ and
$\frac{1}{4}(1,1,3)$ contained in the exceptional divisor of the
morphism $\alpha$ respectively. Then it follows from
Lemmas~\ref{lemma:Noether-Fano} and \ref{lemma:Cheltsov-Kawamata}
that $\mathbb{CS}(U, \frac{1}{k}\mathcal{B})$ contains either  the
point $P_{5}$, or the point $P_{6}$.

It follows from the proof of Lemma~\ref{lemma:n-40-point-P7} that
the set $\mathbb{CS}(U, \frac{1}{k}\mathcal{B})$ does not contain
the singular point $P_{6}$. Thus, the set $\mathbb{CS}(U,
\frac{1}{k}\mathcal{B})$ contains the point $P_{5}$.

The existence of the commu\-ta\-tive
diagram~\ref{equation:n-21-24-85-91-commutative-diagram} follows
from Theorem~\ref{theorem:Kawamata}.

\medskip

\texttt{Case $n=71$.}

The variety $X$ is a hypersurface in $\mathbb{P}(1,1,6,8,15)$ of
degree $30$, whose singularities consist of the point $P_{1}$,
$P_{2}$ and $P_{3}$ that are singularities of type
$\frac{1}{2}(1,1,1)$, $\frac{1}{3}(1,1,2)$ and
$\frac{1}{8}(1,1,7)$ respectively. The equality $-K_{X}^{3}=1/24$
holds. There is a commutative diagram
$$
\xymatrix{
&U\ar@{->}[d]_{\alpha}&&W\ar@{->}[ll]_{\beta}\ar@{->}[d]^{\eta}&\\%
&X\ar@{-->}[rr]_{\psi}&&\mathbb{P}(1,1,6),&}
$$
where $\alpha$ is the weighted blow up of the point $P_{3}$ with
weights $(1,1,7)$, $\beta$ is the weighted blow up with weights
$(1,1,6)$ of the singular point of $U$ that is contained in the
exceptional divisor of  $\alpha$, and $\eta$ is an elliptic
fibration. It follows from
Proposition~\ref{proposition:singular-points} that $\mathbb{CS}(X,
\frac{1}{k}\mathcal{M})=\{P_{3}\}$.

Let $\mathcal{D}$ be the proper transform of $\mathcal{M}$ on $U$.
Then $\mathcal{D}\sim_{\mathbb{Q}}-kK_{U}$  by
Theorem~\ref{theorem:Kawamata}, and it follows from
Lemmas~\ref{lemma:Noether-Fano} and \ref{lemma:Cheltsov-Kawamata}
that the set $\mathbb{CS}(U, \frac{1}{k}\mathcal{D})$ contains the
singular point of  $U$ that is contained in the exceptional
divisor of the morphism $\alpha$.

The existence of the commutative
diagram~\ref{equation:n-21-24-85-91-commutative-diagram} is
implied by Theorem~\ref{theorem:Kawamata}.

\medskip

\texttt{Case $n=76$.}

The variety $X$ is a hypersurface in $\mathbb{P}(1,5,6,8,11)$ of
degree $30$, whose singularities  consist~of the point $P_{1}$,
$P_{2}$ and $P_{3}$ that are  singularities of type
$\frac{1}{2}(1,1,1)$, $\frac{1}{8}(1,3,5)$ and
$\frac{1}{11}(1,5,6)$~respectively. The indeterminacies of the
projection $\psi$ are resolved by the weighted blow up of the
sin\-gu\-lar points  $P_{2}$ and $P_{3}$ with weights $(1,3,5)$
and $(1,5,6)$ respectively.

The proofs of Lemmas~\ref{lemma:n-43-P9} and
\ref{lemma:n-58-P3-is-a-center} implies that $\mathbb{CS}(X,
\frac{1}{k}\mathcal{M})=\{P_{2}, P_{3}\}$, which implies the
existence of the commutative
diagram~\ref{equation:n-21-24-85-91-commutative-diagram} due to
Theorem~\ref{theorem:Kawamata}.

\medskip

\texttt{Case $n=77$.}

The variety $X$ is a hypersurface in $\mathbb{P}(1,2,5,9,16)$ of
degree $32$, the singularities of $X$ consist of the points
$P_{1}$ and $P_{2}$ that are quotient singularities of type
$\frac{1}{2}(1,1,1)$, the point $P_{3}$ that is a quotient
singularity of type $\frac{1}{5}(1,1,4)$, the point $P_{4}$ that
is a quotient singularity of type $\frac{1}{9}(1,2,7)$.

The equality $-K_{X}^{3}=1/45$ holds, and there is a commutative
diagram
$$
\xymatrix{
&U\ar@{->}[d]_{\alpha}&&Y\ar@{->}[ll]_{\beta}\ar@{->}[d]^{\eta}&\\%
&X\ar@{-->}[rr]_{\psi}&&\mathbb{P}(1,2,5),&}
$$
where $\alpha$ is the blow up of $P_{4}$ with weights $(1,2,7)$,
$\beta$ is the blow up with weights $(1,2,5)$ of the singular
point of  $U$ that is a  singularity of type $\frac{1}{7}(1,2,5)$,
and $\eta$ is an elliptic fibration.

It follows from Proposition~\ref{proposition:singular-points} that
$\mathbb{CS}(X, \frac{1}{k}\mathcal{M})=\{P_{4}\}$.

Let $\mathcal{B}$ be the proper transform the linear system
$\mathcal{M}$ on the variety $U$, and $P_{5}$ and $P_{6}$ be the
singular points of $U$ that are singularities of types
$\frac{1}{7}(1,2,5)$ and $\frac{1}{2}(1,1,1)$ contained in the
exceptional divisor of $\alpha$ respectively. Then $\mathbb{CS}(U,
\frac{1}{k}\mathcal{B})\cap\{P_{5}, P_{6}\}\ne\varnothing$ by
Lemmas~\ref{lemma:Noether-Fano} and \ref{lemma:Cheltsov-Kawamata}.

Suppose that the set $\mathbb{CS}(U, \frac{1}{k}\mathcal{B})$
contains the point $P_{6}$. Let $\gamma\colon W\to U$ be the
weighted blow up of the singular point $P_{5}$ with weights
$(1,1,1)$, $\mathcal{H}$ and $\mathcal{D}$ be the proper
transforms of the linear systems $\mathcal{M}$ and $|-16K_{X}|$ on
the variety $W$ respectively, $D$ be a general surface of the
linear system $\mathcal{D}$, and $H_{1}$ and $H_{2}$ be general
surfaces of the linear system $\mathcal{H}$. Then the base locus
of the linear system $\mathcal{D}$ does not contain curves. In
particular, the divisor $D$ is nef, but $D\cdot H_{1}\cdot
H_{1}<0$.

Therefore, the set $\mathbb{CS}(U, \frac{1}{k}\mathcal{B})$
contains the point $P_{5}$. Now the claim of
Theorem~\ref{theorem:Kawamata} implies the existence of the
commutative
diagram~\ref{equation:n-21-24-85-91-commutative-diagram}.

\medskip

\texttt{Case $n=83$.}

The variety $X$ is a general hypersurface in
$\mathbb{P}(1,3,4,11,18)$ of degree $36$, and the singularities of
the hypersurface $X$ consist of the point $P_{1}$ that is a
quotient singularity of type $\frac{1}{2}(1,1,1)$, the points
$P_{2}$ and $P_{3}$ that are quotient singularities of type
$\frac{1}{3}(1,1,2)$, the point $P_{4}$ that is a quotient
singularity of type $\frac{1}{11}(1,4,7)$. There is a commutative
diagram
$$
\xymatrix{
&U\ar@{->}[d]_{\alpha}&&W\ar@{->}[ll]_{\beta}\ar@{->}[d]^{\eta}&\\%
&X\ar@{-->}[rr]_{\psi}&&\mathbb{P}(1,3,4),&}
$$
where $\alpha$ is the weighted blow up of the point $P_{4}$ with
weights $(1,4,7)$, $\beta$ is the weighted blow up with weights
$(1,3,4)$ of the singular point of $U$ that is a quotient
singularity of type $\frac{1}{7}(1,3,4)$, and $\eta$ is an
elliptic fibration. It follows from
Proposition~\ref{proposition:singular-points} that $\mathbb{CS}(X,
\frac{1}{k}\mathcal{M})=\{P_{4}\}$.

Let $\mathcal{B}$ be the proper transform of $\mathcal{M}$ on $U$,
and $P_{5}$ and $P_{6}$ be the singular points of  $U$~that are
singularities of types $\frac{1}{7}(1,3,4)$ and
$\frac{1}{4}(1,1,3)$ contained in the exceptional divisor of the
morphism $\alpha$ respectively. Then the proof of
Proposition~\ref{proposition:n-32} implies that $\mathbb{CS}(U,
\frac{1}{k}\mathcal{B})$~does not contain the singular point
$P_{6}$. Thus, the set $\mathbb{CS}(U,
\frac{1}{k}\mathcal{B})$~contains the  point $P_{5}$ by
Lemmas~\ref{lemma:Noether-Fano} and \ref{lemma:Cheltsov-Kawamata},
and the existence of the
diagram~\ref{equation:n-21-24-85-91-commutative-diagram} follows
from Theorem~\ref{theorem:Kawamata}.

\medskip

\texttt{Case $n=85$.}

The variety $X$ is a  hypersurface in $\mathbb{P}(1,3,5,11,19)$ of
degree $38$, the singularities of $X$ consist of the point $P_{1}$
that is a quotient singularity of type $\frac{1}{3}(1,1,2)$, the
point $P_{2}$ that is a quotient singularity of type
$\frac{1}{5}(1,1,4)$, and the point $P_{3}$ that is a quotient
singularity of type $\frac{1}{11}(1,3,8)$.

The equality $-K_{X}^{3}=2/165$ holds, and there is a commutative
diagram
$$
\xymatrix{
&U\ar@{->}[d]_{\alpha}&&Y\ar@{->}[ll]_{\beta}\ar@{->}[d]^{\eta}&\\%
&X\ar@{-->}[rr]_{\psi}&&\mathbb{P}(1,3,5),&}
$$
where $\alpha$ is the weighted blow up of the point $P_{3}$ with
weights $(1,3,8)$, $\beta$ is the weighted blow up with weights
$(1,3,5)$ of the singular point of $U$ that is a singularity of
type $\frac{1}{8}(1,3,5)$, and $\eta$ is an elliptic fibration. It
follows from Proposition~\ref{proposition:singular-points} that
$\mathbb{CS}(X, \frac{1}{k}\mathcal{M})=\{P_{3}\}$.

Let $\mathcal{B}$ be the proper transform of $\mathcal{M}$ on $U$,
and $P_{4}$ and $P_{5}$ be the singular points of $U$ that are
quotient singularities of types $\frac{1}{8}(1,3,5)$ and
$\frac{1}{3}(1,1,2)$ contained in the exceptional divisor of the
morphism $\alpha$ respectively. Then
 $\mathbb{CS}(U, \frac{1}{k}\mathcal{B})\cap\{P_{4}, P_{5}\}\ne\varnothing$
 by Lemmas~\ref{lemma:Noether-Fano} and
\ref{lemma:Cheltsov-Kawamata}.

Suppose that the set $\mathbb{CS}(U, \frac{1}{k}\mathcal{B})$
contains the point $P_{5}$. Let $\gamma\colon W\to U$ be the
weighted blow up of the singular point $P_{5}$ with weights
$(1,1,2)$, $\mathcal{H}$ and $\mathcal{D}$ be the proper
transforms of the linear systems $\mathcal{M}$ and $|-19K_{X}|$ on
the variety $W$ respectively, $D$ be a general surface of the
linear system $\mathcal{D}$, and $H_{1}$ and $H_{2}$ be general
surfaces of the linear system $\mathcal{H}$. Then the base locus
of the linear system $\mathcal{D}$ does not contain curves.

The divisor $D$ is nef. In particular, the inequality $D\cdot
H_{1}\cdot H_{1}\geqslant 0$ holds, but it follows from the simple
explicit computations that $D\cdot H_{1}\cdot H_{1}=-2k^{2}/15$,
which is a contradiction.

Hence, the set $\mathbb{CS}(U, \frac{1}{k}\mathcal{B})$ contains
the singular point $P_{4}$. Thus, the existence of the commutative
diagram~\ref{equation:n-21-24-85-91-commutative-diagram} follows
from Theorem~\ref{theorem:Kawamata}.

\medskip

\texttt{Case $n=91$.}

The variety $X$ is a hypersurface in $\mathbb{P}(1,4,5,13,22)$ of
degree $44$, the singularities of $X$ consist of the point $P_{1}$
that is a quotient singularity of type $\frac{1}{2}(1,1,1)$, the
point $P_{2}$ that is a quotient singularity of type
$\frac{1}{5}(1,2,3)$, the point $P_{3}$ that is a quotient
singularity of type $\frac{1}{13}(1,4,9)$.

The equality $-K_{X}^{3}=1/130$ holds, and there is a commutative
diagram
$$
\xymatrix{
&U\ar@{->}[d]_{\alpha}&&W\ar@{->}[ll]_{\beta}\ar@{->}[d]^{\eta}&\\%
&X\ar@{-->}[rr]_{\psi}&&\mathbb{P}(1,4,5),&}
$$
where $\alpha$ is the weighted blow up of the point $P_{3}$ with
weights $(1,4,9)$, $\beta$ is the weighted blow up  with weights
$(1,4,5)$ of the singular point of the variety $U$ that is a
quotient singularity of type $\frac{1}{9}(1,4,5)$, and $\eta$ is
an elliptic fibration.

It follows from Proposition~\ref{proposition:singular-points} that
$\mathbb{CS}(X, \frac{1}{k}\mathcal{M})=\{P_{3}\}$.

Let $\mathcal{B}$ be the proper transform of $\mathcal{M}$ on $U$,
and $P_{4}$ and $P_{5}$ be the singular points of $U$ that are
quotient singularities of types $\frac{1}{9}(1,4,5)$ and
$\frac{1}{4}(1,1,3)$ contained in the exceptional divisor of the
morphism $\alpha$ respectively. Then it follows from
Lemmas~\ref{lemma:Noether-Fano} and \ref{lemma:Cheltsov-Kawamata}
that $\mathbb{CS}(U, \frac{1}{k}\mathcal{B})$ contains either the
point $P_{4}$, or the point $P_{5}$.

It follows from the proof of Proposition~\ref{proposition:n-32}
that the set $\mathbb{CS}(U, \frac{1}{k}\mathcal{B})$ does not
contain the singular point $P_{5}$. Thus, the set $\mathbb{CS}(U,
\frac{1}{k}\mathcal{B})$ contains the singular point $P_{4}$, and
the existence of the commutative
diagram~\ref{equation:n-21-24-85-91-commutative-diagram} follows
from Theorem~\ref{theorem:Kawamata}.
\end{proof}

The claim of Theorem~\ref{theorem:main} is proved.

\medskip

\address{\begin{tabbing}
\hspace*{28 em}\=\kill
Steklov Institute of Mathematics \>School of Mathematics\\%
8 Gubkin street, Moscow 117966   \>The University of Edinburgh\\%
Russia                           \>Kings Buildings,  Mayfield Road\\%
                                 \> Edinburgh EH9 3JZ, UK\\%
\texttt{cheltsov@yahoo.com}      \>                      \\%
                                 \>\texttt{I.Cheltsov@ed.ac.uk}%
\end{tabbing}}


\begin{thebibliography}{99}

\bibitem{Ch00a}
I.\,Cheltsov, \emph{Log pairs on birationally rigid varieties}\\
Journal of Mathematical Sciences \textbf{102} (2000), 3843--3875

\bibitem{Ch03a}
I.\,Cheltsov, \emph{Anticanonical models of Fano $3$-folds of degree four}\\
Sbornik: Mathematics \textbf{194} (2003), 617--640%

\bibitem{Ch05umn}
I.\,Cheltsov, \emph{Birationally rigid Fano varieties}\\
Russian Mathematical Surveys \textbf{60} (2005), 71--160%

\bibitem{ChPa05}
I.\,Cheltsov, J.\,Park, \emph{Weighted Fano threefold hypersurfaces}\\
Journal fur die Reine und Angewandte Mathematik, принято в печать%

\bibitem{Co00}
A.\,Corti, \emph{Singularities of linear systems and 3-fold birational geometry}\\
L.M.S. Lecture Note Series \textbf{281} (2000), 259--312%

\bibitem{CPR}
A.\,Corti, A.\,Pukhlikov, M.\,Reid, \emph{Fano 3-fold hypersurfaces}\\
L.M.S. Lecture Note Series \textbf{281} (2000), 175--258%

\bibitem{IF00}
A.\,R.\,Iano-Fletcher, \emph{Working with weighted complete intersections}\\
L.M.S. Lecture Note Series \textbf{281} (2000), 101--173%

\bibitem{Ka96}
Y.\,Kawamata,  \emph{Divisorial contractions to $3$-dimensional terminal quotient singularities}\\
Higher-dimensional complex varieties (Trento, 1994), de Gruyter, Berlin (1996), 241--246%

\bibitem{KMM}
Y.\,Kawamata, K.\,Matsuda, K.\,Matsuki, \emph{Introduction to the minimal model problem}\\
Advanced Studies in Pure Mathematics \textbf{10} (1987), 283--360%

\bibitem{Ko97}
J.\,Koll\'ar, \emph{Singularities of pairs}\\
Proceedings of Symposia in Pure Mathematics \textbf{62} (1997), 221--287%

\bibitem{Ry05}
D.\,Ryder, \emph{Classification of elliptic and K3 fibrations birational to some $\mathbb{Q}$-Fano $3$-folds}\\
Journal of Mathematical Sciences, The University of Tokyo, \textbf{13} (2006),  13--42%

\bibitem{Ry02}
D.\,Ryder, \emph{The Curve Exclusion Theorem for elliptic and $K3$ fibrations birational to Fano $3$-fold hypersurfaces}\\
arXiv:math.AG/0606177 (2006)%

\bibitem{Sho93}
V.\,Shokurov, \emph{$3$-fold log flips}\\
Russian Academy of Sciences. Izvestiya. Mathematics \textbf{40} (1993), 95--202%
\end{thebibliography}
\end{document}